\documentclass[12pt,twoside,a4paper]{amsart}

\usepackage[outline]{contour}                                
\newcommand*{\fancy}[1]{{\color{white}\contour{black}{#1}}}  %
\usepackage[dvipsnames]{xcolor}                     
\usepackage{hyperref} 
\hypersetup{breaklinks, colorlinks,                 
    linkcolor = {Blue},                             
    citecolor = {Blue},                             
    urlcolor  = {OliveGreen}                        
}                                                   %
\newcommand{\urlro}[1]{\hypersetup{urlcolor = {YellowOrange
}}\url{#1}\hypersetup{urlcolor = {OliveGreen}}}%
\newcommand{\urlp}[1]{\hypersetup{urlcolor = {BrickRed%
}}\url{#1}\hypersetup{urlcolor = {OliveGreen}}}%

\usepackage{tikz}
\usepackage{pgfplots} 
\pgfplotsset{compat=1.11}
\usetikzlibrary{intersections, pgfplots.fillbetween}
 
\usepackage[OT2, T1]{fontenc}\usepackage[russian, english]{babel}
\usepackage[T1]{fontenc}
\usepackage[latin9]{inputenc}
\usepackage[margin=20mm,top=30mm,bottom=30mm]{geometry}
\usepackage{amsmath,amssymb,amsthm,amscd,stmaryrd,graphicx,bm,latexsym}

\makeatletter
\def\@tvsp{\mathchoice{{}\mkern-4.5mu}{{}\mkern-4.5mu}{{}\mkern-2.5mu}{}}
\def\ltrivert{\left|\@tvsp\left|\@tvsp\left|}
\def\rtrivert{\right|\@tvsp\right|\@tvsp\right|}
\makeatother
\usepackage{mathrsfs}\usepackage{mathtools}\usepackage{lmodern}\usepackage{natbib}\usepackage{enumitem}\usepackage{soul}
\newlist{parts}{enumerate}{5}       
\setlist[parts]{wide,label=\textbf{\upshape(\alph*)},ref=(\alph*),align=left,labelindent=0pt,labelsep=.5em,itemsep=.3ex,topsep=0ex,%
}

\DeclareFontFamily{U}{mathb}{\hyphenchar\font45}
\DeclareFontShape{U}{mathb}{m}{n}{
<-6> mathb5 <6-7> mathb6 <7-8> mathb7
<8-9> mathb8 <9-10> mathb9
<10-12> mathb10 <12-> mathb12
}{}
\DeclareSymbolFont{mathb}{U}{mathb}{m}{n}
\DeclareMathSymbol{\llcurly}{\mathrel}{mathb}{"CE}
\DeclareMathSymbol{\ggcurly}{\mathrel}{mathb}{"CF}

\usepackage[
final,
scrtime
]{prelim2e}

\renewcommand{\PrelimText}{\footnotesize[\,
\texttt{\jobname.tex}\hfill \today\ at \thistime\,]}
  \theoremstyle{plain}
    \newtheorem{Thm}{Theorem}[section]  
      \newtheorem{Lem}[Thm]{Lemma}\newtheorem{Cor}[Thm]{Corollary}\newtheorem{Rem}[Thm]{Remark}
      \newtheorem{Example}[Thm]{Example}
      \newtheorem{Examples}[Thm]{Examples}
   
    \theoremstyle{definition} 
      \newtheorem{Def}[Thm]{Definition} \newtheorem{Question}[Thm]{Question}
      \newtheorem{Problem}[Thm]{Problem}
      
\renewcommand{\iff}{\Leftrightarrow}

\renewcommand{\cite}{\citet}
\newcommand*{\nameadjunct}{\relax}
\makeatletter
\renewcommand*{\NAT@nmfmt}[1]{\NAT@up #1\nameadjunct}
\makeatother

\newcommand*{\citeposs}[2][]{%
    \begingroup
    \renewcommand*{\nameadjunct}{'s}%
    \citet[#1]{#2}%
    \endgroup
}

\newcounter{aequiv}

 {\end{list}}  
\let\oldsqrt\sqrt
\def\sqrt{\mathpalette\DHLhksqrt}
\def\DHLhksqrt#1#2{%
\setbox0=\hbox{$#1\oldsqrt{#2\,}$}\dimen0=\ht0
\advance\dimen0-0.2\ht0
\setbox2=\hbox{\vrule height\ht0 depth -\dimen0}%
{\box0\lower0.4pt\box2}}
\renewcommand{\limsup}{\varlimsup}\renewcommand{\liminf}{\varliminf}

\newcommand{\C}{{\mathbb C}}%
\newcommand{\N}{{\mathbb N}}\newcommand{\R}{{\mathbb R}}%
\newcommand{\Z}{{\mathbb Z}}
\newcommand{\cA}{{\mathcal A}}\newcommand{\cB}{{\mathcal B}}%
\newcommand{\cE}{{\mathcal E}}\newcommand{\cF}{{\mathcal F}}%
\newcommand{\cG}{{\mathcal G}}%
\newcommand{\cK}{{\mathcal K}}\newcommand{\cL}{{\mathcal L}}%
\newcommand{\cM}{{\mathcal M}}%
\newcommand{\cP}{{\mathcal P}}%
\newcommand{\cS}{{\mathcal S}}%
\newcommand{\cX}{{\mathcal X}}%
\renewcommand{\epsilon}{\varepsilon}\renewcommand{\phi}{\varphi} %
\renewcommand{\theta}{\vartheta}     %
\newcommand{\Alpha}{\mathrm{A}}    %
\newcommand{\Eta}{\mathrm{H}}      %
\DeclareMathOperator{\sgn}{sgn}  
\DeclareMathOperator*{\bigconv}{\mbox{\LARGE$\ast$}}
\DeclareFontEncoding{LS2}{}{\noaccents@}                                                 %
\DeclareFontSubstitution{LS2}{stix}{m}{n}                                                %
\DeclareSymbolFont{largesymbols_stix}{LS2}{stixex}{m}{n}                                 %
\DeclareMathDelimiter{\lpp}{\mathopen}{largesymbols_stix}{"DE}{largesymbols_stix}{"02}   %
\DeclareMathDelimiter{\rpp}{\mathclose}{largesymbols_stix}{"DF}{largesymbols_stix}{"03}  %
\renewcommand{\[}{\begin{eqnarray*}}\renewcommand{\]}{\end{eqnarray*}}
\newcommand{\la}{\begin{eqnarray}}\newcommand{\al}{\end{eqnarray}}
  
\renewcommand{\Prob}{\mbox{\rm Prob}} 
\newcommand{\dd}{{\mathrm{d}}}
\newcommand{\leb}{\fancy{$\lambda$}}
\DeclareMathOperator*{\esssup}{ess\,sup}
\renewcommand{\P}{{\mathbb P}}
\newcommand{\toitself}{\,\mbox{\rotatebox[x=0.5\width,y=0.34\height]{270}
{\scalebox{0.4}[2]{$\curvearrowleft$}}}\,\,}
\newcommand{\sdot}{\bullet}  
\newcommand{\bdot}[1]{\overset{\,{}_{{}_{\bullet}}}{#1}} 

\newcommand{\im}{\mathbin{\text{\raisebox{.2ex}{\scalebox{.6}{$\!\square$\!}}}}} 

\newcommand{\bed}[2]{\pmb{\left(\right.}#1\,\pmb{|}\, #2\pmb{\left.\right)}}

\newcommand{\stle}{\le_{\mathrm{st}}}\newcommand{\stge}{\ge_{\mathrm{st}}}

\newcommand{\1}{\fancy{$1$}} 

\begin{document}
\author[Mattner]{Lutz Mattner}
\address{Universit\"at Trier, Fachbereich IV -- Mathematik, 54286~Trier, Germany}
\email{mattner@uni-trier.de}
\title
[A convolution inequality, and Berry-Esseen for summands 
 close to normal]
{A convolution inequality,\\ yielding a sharper Berry-Esseen theorem\\ for summands Zolotarev-close to normal}
\begin{abstract} 
The classical Berry-Esseen error bound,
for the normal approximation to the law of a sum of independent 
and identically distributed random variables, 
is here improved by replacing 
the standardised third absolute moment by a weak norm distance to normality.
We thus sharpen and simplify two results of \citet{Ulyanov1976} and of \citet{Senatov1998}, 
each of them previously optimal, 
in the line of research initiated by \citet{Zolotarev1965} and \citet{Paulauskas1969}.

Our proof is based on a seemingly incomparable normal approximation theorem of \citet{Zolotarev1986},
combined with our main technical result:

The Kolmogorov distance (supremum norm of difference of distribution functions) 
between a convolution of two laws and a convolution of two Lipschitz laws 
is bounded homogeneously of degree~1 in the pair of the Kantorovich distances
(often called Wasserstein distances, the
L$^1$ norms of differences of distribution functions)
of the corresponding factors, 
and also in the pair of the Lipschitz constants.

Side results include 
a short introduction to $\zeta$ norms on the real line,
simpler inequalities for various probability distances,
slight improvements of the theorem of \citet{Zolotarev1986} and of a lower bound theorem of~\citet{BCG2012},
an application to sampling from finite populations, auxiliary results on rounding
and on winsorisation,
and computations of a few examples.

The introductory section in particular 
is aimed at analysts in general rather than specialists in probability approximations.
\end{abstract} 
\dedicatory{Dedicated to Ukraine}
\maketitle
\tableofcontents

\section{Introduction, from Berry-Esseen to its sharpening
Theorem~\ref{Thm:Berry-Esseen_for_Z-close_to_normal}}         \label{sec:Introduction}
\subsection{Aim}             \label{subsec:Aim}
The main purpose of this paper is to prove Theorem~\ref{Thm:BE_K_Z}, stated on page~\pageref{Thm:BE_K_Z}
below, which is a Berry-Esseen type central limit theorem, 
for sums of $n$ independent and identically distributed random variables, 
taking  also a closeness of the summands to normality into account, namely by
bounding the normal approximation error 
in the usual Kolmogorov norm $\left\|\,\cdot\,\right\|_\mathrm{K}$
from~\eqref{Eq:Def_Kolmogorov-norm} by $\frac{1}{\sqrt{n}}$ times a weak norm distance of the law of one
standardised summand to the standard normal law. 
This strictly improves four of five similar and apparently mutually incomparable results, 
each as far as known to the present author previously optimal of its kind,
of, in a certain logical rather than historical order,
\citet{Shiganov1987}, \citet{Ulyanov1976},
\citeauthor{Zolotarev1973} (implicit in his papers from \citeyear{Zolotarev1973} 
and \citeyear{Zolotarev1976}), and \citet{Senatov1998}, namely
by having on the right hand side  weaker norms with all exponents equal to~$1$ for $n\ge2$.
As its precursors from \citet{Paulauskas1969} onwards, 
Theorem~\ref{Thm:BE_K_Z} contains
the classical \citet{Berry1941}--\citet{Esseen1942} theorem~\eqref{Eq:Berry-Esseen_inequality} as a corollary, 
albeit in its present version with some rather large constant, namely with $7.2$,  obtained from
combining Theorem~\ref{Thm:BE_K_Z} or 
inequality~\eqref{Eq:Berry-Esseen_for_Z-close_to_normal_also_for_n=1}  
with inequality~\eqref{Eq:zeta_1_vee_zeta_3_distance_to_normal_bounded_by_nu_3},
rather than with the up to now best value $0.469$ from~\eqref{Eq:Def_cSH}
announced by \citet{Shevtsova2013}. 

We prove Theorem~\ref{Thm:BE_K_Z}, perhaps somewhat surprisingly, 
by reducing it to the Berry-Esseen type theorem of \citet{Zolotarev1986,Zolotarev1997}, recalled and slightly refined 
as Theorem~\ref{Fact:Zolotarev1986} below, which has a norm incomparable to Kolmogorov's on the {\em left} hand side.
The reduction is possible by using the present Corollary~\ref{Cor:Main_for_n=2} 
to our main technical result, Theorem~\ref{Thm:F_star_G_vs_H_star_H} on page~\pageref{Thm:F_star_G_vs_H_star_H}, 
which bounds the Kolmogorov distance by Zolotarev's $\zeta_1$ (or Kantorovich, or
Wasserstein)
distance for certain convolution products.   
Theorem~\ref{Thm:F_star_G_vs_H_star_H} in turn is proved by what seems to us to be, 
in the field of probability approximation theorems, a not quite standard use of the Krein-Milman theorem,  
or  more precisely of the closely related \citet{Bauer1958} maximum principle. 

Ignoring constant factors, abbreviated as
\la                                    \label{Eq:i.c.f.}
 \text{{\em i.c.f.}}
\al
in this paper, Zolotarev's Berry-Esseen type theorem
just mentioned is actually stronger than the i.i.d.~case of the more recent result of 
\citet{Goldstein2010} and \citet{Tyurin2010} recalled 
as Theorem~\ref{Thm:Goldstein-Tyurin} below, but, as a side result of this paper, 
we use Goldstein-Tyurin in the obvious way to improve a bit the constant in Zolotarev's theorem.
Again i.c.f., as already indicated above, the present Theorem~\ref{Thm:Berry-Esseen_for_Z-close_to_normal}
improves the Berry-Esseen theorem~\eqref{Eq:Berry-Esseen_inequality}, but we analogously use Shevtsova's constant 
from~\eqref{Eq:Def_cSH} for the latter  to get a smaller constant 
than otherwise obtainable here in the former. We also emphasise the asymptotically optimal,
for $\zeta \llcurly \varkappa \rightarrow0$,
error bound~(\ref{Eq:Zolotarev's_zeta_1-B-E-Thm_sharper},\ref{Eq:xi_near_zero})
in Zolotarev's Theorem~\ref{Fact:Zolotarev1986}, obtained here by essentially
his proof, but not pointed out by him.

In subsections~\ref{subsec:Classical_B-E}--\ref{subsec:Known_solutions_with_zeta}
below we explain in more detail the 
development leading to the present Theorem~\ref{Thm:BE_K_Z}.
The length of these subsections may be excused by our aim of writing  
there for analysts in general, rather than for experts in Berry-Esseen refinements,
with indeed the hope of attracting some of the former to this 
fascinating area of probability theory. 
Readers already knowing~$\zeta_3$ may jump to Theorem~\ref{Thm:BE_K_Z} on page~\pageref{page:Thm_BE_K_Z}
immediately, and readers interested just in norm inequalities for convolutions
with two factors
may jump to Theorem~\ref{Thm:F_star_G_vs_H_star_H} on page~\pageref{page:Thm_convolution_inequality}.

\subsection{Some notation and conventions}        \label{subsec:Some_notation} 
For stating our results and comparisons more precisely, let us introduce here some notation. 
For just reading the convolution inequality
Theorem~\ref{Thm:F_star_G_vs_H_star_H}, however, it suffices to recall the standard notation 
from the paragraph around~\eqref{Eq:Def_Lip-norm}, and to accept
the perhaps not so standard notation~\eqref{Eq:Def_conv_distr_fcts}. 
Throughout this paper, we have tried to ``recall'' any unfamiliar notation,
usually by pointing at appropriate places in the present subsection, 
which hence might perhaps be skipped for now and 
consulted only when needed.

In addition to~\eqref{Eq:i.c.f.}, we use some more common abbreviations:
{\em i.i.d.}~for independent and identically distributed,
{\em a.e.}~for almost everywhere,
{\em w.l.o.g.} for without loss of generality,
{\em w.r.t.} for with respect to,
{\em iff} for if an only if,
L.H.S.~for left hand side,
and R.H.S.~for right hand side.

We use the indicator notation of \citet[p.~11]{Iverson1962}\,-\,\citet[pp.~xx--xxi in the English
translation~\citeyear{deFinetti1972}]{deFinetti1967} 
for propositions,  and also a more common one for sets, 
\la                             \label{Eq:indicator_notation}
  (\text{statement}) &\coloneqq& \begin{Bmatrix}1 \\ 0 \end{Bmatrix}
  \text{ if statement is } \begin{Bmatrix}\text{true} \\ \text{false} \end{Bmatrix},
  \qquad \1_A(x) \,\ \coloneqq \,\ (x\in A).
\al

We use the standard lattice theoretical notation
$x\wedge y\coloneqq \inf\{x,y\}$,
$x\vee y \coloneqq\sup\{x,y\}$, $x_+\coloneqq x\vee0$, $x_-\coloneqq (-x)\vee 0$,
and $|x|\coloneqq x_++x_-$\,, where $x$ and $y$ may be
real numbers, functions, or signed measures.

If $f$ is a 
$\C$-valued Borel function  defined on $\R$ almost everywhere with respect to Lebesgue measure $\leb$, 
we write as usual $\left\|f\right\|_1\coloneqq \int|f|\,\dd\leb$ and 
$\left\|f\right\|_\infty
\coloneqq \esssup_{x\in\R}|f(x)|$ $\coloneqq$ the $\leb$-essential  
supremum of $|f|$. 
For an everywhere defined function $\in\C^\R$, the ordinary supremum of $|f|$ 
will be written just as $\sup_{x\in\R}|f(x)|$, and its 
Lipschitz constant as 
\la                                       \label{Eq:Def_Lip-norm}
 \|f\|^{}_{\mathrm{L}} &\coloneqq& \sup\left\{\,\left|\tfrac{f(y)-f(x)}{y-x}\right| 
   : x,y\in\R, x\neq y\right\}\,,
\al
and $f$ is called {\em Lipschitz} if $\|f\|^{}_{\mathrm{L}} <\infty$.
For example we have 
\la                     \label{Eq:Lip=sup_|f'|}
 \|f\|^{}_{\mathrm{L}} &=& \left\|f'\right\|_\infty \quad\text{ for }f\in\C^\R\text{ Lipschitz},
\al
by the fundamental theorem of calculus for absolutely continuous functions, 
see for example~\citet[Theorem 7.20]{Rudin1987}.

We denote the vector space 
of all bounded signed measures on the Borel sets of $\R$ simply by
\la                                                \label{Eq:Def_cM}
 \cM\,, 
\al 
but use the standard notation
$\Prob(\R)$ for the subset of all probability measures, or {\em laws} for brevity.
Writing $\mathrm{N}_{\mu,\sigma^2}$ for the normal law with mean $\mu\in\R$ 
and standard deviation $\sigma\in[0,\infty[\,$, we abbreviate our notation  
in the centred or even standard case to
\la                           \label{Eq:Def_N_sigma_N_1}
  \mathrm{N}_\sigma &\coloneqq& \mathrm{N}_{0,\sigma^2}\quad\text{ for }\sigma\in[0,\infty[\,,  
     \qquad \mathrm{N} \,\ \coloneqq\,\ \mathrm{N}_1 \,.
\al
Further special laws occurring below include the Dirac measures $\delta_a$ for $a\in\R$, 
and the Bernoulli laws $\mathrm{B}_p\coloneqq(1-p)\delta_0+p\delta_1$ for $p\in[0,1]$.

For $M\in\cM$ 
we define its ordinary and complementary distribution functions 
$F_M$ and $\overline{F}_{\!M}$ by
\la                      \label{Eq:Def_F_M_and_overline{F}_M}
  && F^{}_M(x) \, \coloneqq\, M(\,\mathopen]-\infty,x\mathclose]\,),
  \quad \overline{F}^{}_{\!M}(x) \, \coloneqq\, M(\,\mathopen[x,\infty\mathclose[\,) 
     \, =\, M(\R)- F^{}_M(x-)
  \quad\text{ for }x\in \R,
\al
write as usual ${\Phi} \coloneqq F^{}_{\mathrm{N}}$ for the standard normal distribution function,
so $\Phi(x)=\int_{-\infty}^x\phi(y)\,\dd y$ with $\phi(y)\coloneqq\frac{1}{\sqrt{2\pi}}\exp(-\frac{y^2}{2})$,
call 
\la                                                                \label{Eq:Def_Kolmogorov-norm}
 {\big\|} M {\big\|_{\mathrm{K}}} &\coloneqq& \big\| F^{}_M\big\|^{}_\infty \vee \big\| \overline{F}^{}_{\!M}\big\|^{}_\infty \,,
 \quad\text{ hence }\  {\big\|} M {\big\|_{\mathrm{K}}}
   \,\ =\,\ \big\| F^{}_M\big\|^{}_\infty \ \text{ if } M(\R)=0\,,
\al
the (reflection invariant version of the) \textit{Kolmogorov norm} of $M$,
write ${|}M{|}$ for its variation measure,  and consider
\la                                       \label{Eq:Def_nu_r}
  \nu_r(M) \coloneqq \int |x|^r \,\dd|M|(x)     \quad\text{ for }r\in[0,\infty[\,, 
\al  
the $r$th absolute moment of $M$, put 
\la                              \label{Eq:Def_M_Lipschitz}
  \|M\|^{}_{\mathrm{L}} &\coloneqq& \|F_M\|^{}_{\mathrm{L}}
    \,\ =\,\ \|\overline{F}_{\!M}\|^{}_{\mathrm{L}} \,,
\al
and call $M$ {\em Lipschitz} if $\|M\|^{}_{\mathrm{L}}<\infty$.
                                                                              
We use the conventions of measure theory about $\pm\infty$, in particular $0\cdot\infty\coloneqq0$.  

Let us call any subadditive and absolutely homogeneous $[0,\infty]$-valued
function on a vector space over~$\R$ or~$\C$  
an {\em eqnorm}. An eqnorm $\|\cdot\|$ 
with $\|x\|=0$ implying $x=0$ is here called an {\em enorm}, and a
$[0,\infty[$-valued eqnorm is called a {\em qnorm}. 
Thus an enorm is a norm except that it may assume the enormous value $\infty$,
and qnorm is short for {\em quasinorm} (with the more common name ``seminorm''
avoided by us, since usually no complementing other half justifying the ``semi'' is in sight).
On $\cM$ for example 
$\left\|\,\cdot\,\right\|_{\mathrm{K}}$ and  the usual (unweighted) total variation norm $\nu_0$
included in~\eqref{Eq:Def_nu_r} are indeed norms,
$\|\cdot \|^{}_{\mathrm{L}}$ is an enorm,
while $\nu_r$ with $r>0$ is merely an eqnorm, 
with then in particular $\nu_r(M)=0$ iff $M$ is a multiple of $\delta_0$.

We write $T\im M$ for the image measure of a signed measure $M$, on any measurable space,
under a measurable function $T$, that is
\la                             \label{Eq:Def_image_measure}
 \big(T\im M\big)(B) &\coloneqq& M(T^{-1}[B])\,,
\al
where $T^{-1}[B]$ denotes a preimage.

The convolution of signed measures $M_1,M_2\in \cM$,
namely $\big(\R^2\ni (x,y)\mapsto x+y\big)\im (M_1\otimes M_2)$ with $\otimes$
indicating a product measure, 
is denoted by $M_1 {\ast} M_2$ as usual, while  the convolution 
of the distribution functions $F^{}_{M_1} , F^{}_{M_2}$ is defined to be the distribution function 
of~$M_1 {\ast} M_2$, and is written with a star $\star$ instead of an asterisk $\ast$, that is, 
\la                                                 \label{Eq:Def_conv_distr_fcts}
 F^{}_{M_1} {\star} F^{}_{M_2} &\coloneqq& F^{}_{M_1\ast M_2},
\al 
to avoid confusion with the more usual convolution of $\leb$-integrable functions $f_1,f_2$, 
which yields, up to equality a.e.,  a $\leb$-density, commonly denoted by $f_1\ast f_2$, 
of the convolution of the signed measures $f_1\leb, f_2\leb$,
with the latter defined by $\big(f_j\leb\big)(B)\coloneqq\int_Bf_j\,\dd\leb$ for $B\subseteq\R$ Borel.
Convolutions with $n\in\N_0\coloneqq\{0\}\cup\N$ factors are written
like $\bigconv_{j=1}^n M_j$, which in case of $n=0$ means $\delta_0$,
and convolution powers with exponent $n$ as $M^{\ast n} \coloneqq \bigconv_{j=1}^n M$.

Image measures of an $M\in\cM$ under translations or scalings are 
often written like 
$M(\frac{\cdot-a}{\lambda})\coloneqq (x\mapsto \lambda x +a)\im M$ 
for $a\in\R$ and $\lambda\in\R\!\setminus\!\{0\}$, and we also write 
$\check{M}\coloneqq M(\frac{\cdot}{-1})$  \label{page:Def_reflection_of_M} 
for the reflection of $M$.
If $M=\check{M}$, then $M$ is called {\em symmetric}.

For  $r\in [0,\infty[$ we put 
\la
 \cM_r & \coloneqq& \{ M\in\cM : \nu_r(M)<\infty\},    \label{Eq:Def_cM_r}
 \qquad {\Prob_r(\R)} \,\ \coloneqq\,\  \cM_r\cap \Prob(\R), \\
   \cP_r &\coloneqq& \big\{P\in\Prob_r(\R): P\text{ is no Dirac measure}\big\}.  \label{Eq:Def_cP_r}
\al 
We put ${\mu_k}(M) \coloneqq \int x^k\,\dd M(x)$ for $k\in\N_0$ and $M\in \cM_k$,
and 
\la                                                      \label{Eq:Def_cM_r,k}
 \cM_{r,k} &\coloneqq&  \big\{ M \in \cM_r : \mu_j(M)=0\text{ for }j\in\{0,\ldots,k\}\big\}
\al
for $r\in [0,\infty[$ and $k\in\{0,\ldots,\lfloor r \rfloor\}$. We further write 
${\mu} \coloneqq \mu_1$, ${\sigma}(P)\coloneqq\sqrt{\mu_2(P)-(\mu_1(P))^2}$ 
for the standard deviation of $P\in\Prob_2(\R)$, 
$\widetilde{P}$ for the standardisation of $P\in\cP_2$, that is, 
the law of the standardisation $\widetilde{X}\coloneqq\frac{X-\mu(P)}{\sigma(P)}$ 
of any random variable~$X$ with law~$P$, 
equivalently $\widetilde{P}(B) = P\big(\sigma(P)B+\mu(P)\big)$ for $B\subseteq\R$ Borel,
and correspondingly ${\widetilde{F^{}_P}} \coloneqq F^{}_{\widetilde{P}}$.
We accordingly put
\[
 \widetilde{\cP_r}
  &\coloneqq& \{ \widetilde{P} : P\in \cP_r \}
  \,\ =\,\ \{ P \in\cP_r : \mu(P)=0,\sigma(P)=1\}
  \quad\text{ for }r\in [2,\infty[\,.
\]
The centring at the mean of a law $P\in\Prob_1(\R)$ is 
\[
  \bdot{P}  &\coloneqq& \big(x\mapsto x-\mu(P)\big)\im P
       \,\ =\,\  P\big({\cdot} + \mu(P)\big)  \,\ =\,\ \delta_{-\mu(P)}\ast P \,.
\]

Let further
\la                           \label{Eq:Def_lattice_span_h(P)}
 h(P)&\coloneqq&\sup\bigcup_{a\in\R}\{\eta\in\mathopen]0, \infty\mathclose[:P(a+\eta\mathbb Z)=1\}
 \quad \text{ for } P\in\Prob(\R)\,,
\al
the lattice span of $P$;  this is of course to be read as
$h(P)\coloneqq0$ if $P$ is non-lattice, that is, if
$P(a+\eta\mathbb Z)<1$ for each choice of~$a$ and~$\eta$.

Zolotarev's enorms $\zeta_r$, occurring for $r\in\{1,3\}$ in Theorem~\ref{Thm:BE_K_Z},
are defined by~(\ref{Eq:Def_cF_r},\ref{Eq:Def_cF_r^infty},\ref{Eq:Def_zeta_r}),
with alternative representations provided by~\eqref{Eq:Four_zeta_expressions} 
and (\ref{Eq:zeta_1=kappa_1},\ref{Eq:varkappa_1=L^1-enorm}).

The standard asymptotic comparison notation $\preccurlyeq, \asymp,\llcurly,\sim$ 
is recalled in section~\ref{sec:Asymptotic_comparison_terminology}, 
where also our use of ``i.c.f.'' is explained.

\subsection{The classical Berry-Esseen theorem}        \label{subsec:Classical_B-E}
With the above notation, 
a classical form of the central limit theorem, namely for the standardised partial sums of each fixed 
sequence of independent and identically distributed real-valued random variables with finite and nonzero variances, 
as first proved implicitly by \citet[p.~219, Satz III]{Lindeberg1922}, 
presented more explicitly perhaps first by \citet[num\'ero~45 on p.~233, combined with pp.~192--193 which provide
the \citet{BuchananHildebrandt1908} type theorem  often named after {\citet[p.~173, Satz I]{Polya1920}}]{Levy1925},
and presumably hence sometimes attributed to L\'evy as for example 
in the standard monograph~\citet[p.~126, Theorem 4.8]{Petrov1995},  can be succinctly stated as 
\la                                                   \label{Eq:Lindeberg-Levy-Thm}
  \lim_{n\rightarrow\infty}\left\|\widetilde{P^{\ast n}}-\mathrm{N}\right\|_{\mathrm{K}} 
  &=& 0\quad\text{ for }P\in\cP_2\,.
\al
The typical asymptotics of the approximation error $\left\|\widetilde{P^{\ast n}}-\mathrm{N}\right\|_{\mathrm{K}}$ was, 
under the sole additional condition $\nu_3(P)<\infty$, provided by \citet[p.~162]{Esseen1956} as
\la                                                                                   \label{Eq:Esseen1956_asymptotics}
   \lim_{n\rightarrow\infty} \sqrt{n}\left\|\widetilde{P^{\ast n}}-\mathrm{N}\right\|_{\mathrm{K}} 
   &=& \frac{1}{\sqrt{2\pi}}\left(\frac{h(\widetilde{P})}2 +  \frac{|\mu_3(\widetilde{P})|}{6} \right)  
   \quad\text{ for }P\in\cP_3\,.
\al
It therefore seems natural to ask for finite sample error bounds of the form
\la                                \label{Eq:c(P)/sqrt(n)_bound}
  \left\|\widetilde{P^{\ast n}}-\mathrm{N}\right\|_{\mathrm{K}} &\le& \frac{c(P)}{\sqrt{n}}  
  \quad \text{ for }P\in\cP_3\text{ and } n\in\N 
\al
for some  appropriate choice of $c(P)$, 
which should in particular be not too difficult to compute or to bound from above,
and such a bound is provided by the celebrated \citet{Berry1941}-\citet{Esseen1942} theorem
\la                                                                     \label{Eq:Berry-Esseen_inequality}
  \left\|\widetilde{P^{\ast n}}-\mathrm{N}\right\|_{\mathrm{K}} &\le& \frac{c}{\sqrt{n}}\nu_3(\widetilde{P}) 
  \quad \text{ for } P\in\cP_3\text{ and } n\in\N 
\al
with some universal constant $c<\infty$. 
According to \citet{Shevtsova2013}, we can choose here $c$ as
\newcommand{\cSh}{ c^{}_{\text{\foreignlanguage{russian}{Sh}}}} 
\la                                  \label{Eq:Def_cSH}
  \cSh &\coloneqq& 0.469\,,
\al  
with some improvements to be expected in the future, but certainly not beyond
\la                                                                      \label{Eq:Def_c_E}
  c^{}_{\mathrm{E}}&\coloneqq&  \frac{3+\sqrt{10}}{6\sqrt{2\pi}}\,\ =\,\ 0.4097\ldots 
\al  
as follows easily from using~\eqref{Eq:Esseen1956_asymptotics} just for Bernoulli laws, 
$P=\mathrm{B}_p$ with $p\in\mathopen]0,1\mathclose[$\,:

For $p\in[0,1]$, 
we easily compute
$\mu(\mathrm{B}_p)=p$, $\sigma(\mathrm{B}_p)=\sqrt{p(1-p)}$, 
$\mu_3(\bdot{\mathrm{B}}_p)=p(1-p)(1-2p)$, 
$\nu_3(\bdot{\mathrm{B}}_p)=p(1-p)(\frac{1}{2}+2(p-\frac{1}{2})^2)$,
and if $p\in\mathopen]0,1\mathclose[$ hence
$h(\mathrm{B}_p)=1$, $h(\widetilde{\mathrm{B}_p})=\frac{1}{\sqrt{p(1-p)}}$,
$\mu_3(\widetilde{\mathrm{B}_p}) = \frac{1-2p}{\sqrt{p(1-p)}}$,
$\nu_3(\widetilde{\mathrm{B}_p}) = \frac{\frac12+2(p-\frac12)^2 }{\sqrt{p(1-p)}}$,
so that~\eqref{Eq:Esseen1956_asymptotics} and~\eqref{Eq:Berry-Esseen_inequality} specialise to 
\la                                                                   \label{Eq:Esseen1956_asymptotics_Bernoulli}
 \lim_{n\rightarrow\infty} \sqrt{n}\left\|\widetilde{\mathrm{B}_p^{\ast n}}-\mathrm{N}\right\|_{\mathrm{K}} 
   &=& \frac{3+|1-2p|}{6\sqrt{2\pi p(1-p)}}  
     \quad\text{ for }p\in\mathopen]0,1\mathclose[\,, \\
 \left\|\widetilde{\mathrm{B}_p^{\ast n}}-\mathrm{N}\right\|_{\mathrm{K}} 
   &\le& c\, \frac{\frac12+2(p-\frac12)^2  }{\sqrt{np(1-p)}}               \label{Eq:Berry-Esseen_Bernoulli}
     \quad \text{ for }p\in\mathopen]0,1\mathclose[  \text{ and } n\in\N\,,  
\al
implying by elementary calculations 
\[
 c &\ge& \sup_{p\in\mathopen]0,1\mathclose[} \frac{3+|1-2p|}{6\sqrt{2\pi}(\frac12+2(p-\frac12)^2)}
   \,\ =\,\ c^{}_{\mathrm{E}}
\]
with the supremum attained at $p=p^{}_{\mathrm{E}}\coloneqq \frac12(4-\sqrt{10})=0.418861\ldots$
and at $p=1-p^{}_{\mathrm{E}}$.

\citet[p.~161, Theorem]{Esseen1956} actually proved the more interesting and less trivial result
that $c^{}_{\mathrm{E}}$ is the best lower bound for $c$ in~\eqref{Eq:Berry-Esseen_inequality}
obtainable from~\eqref{Eq:Esseen1956_asymptotics} even without restricting attention to $P$ 
being Bernoulli, that is, 
\la                                                                        \label{Eq:c_E_really_is_...}
 c^{}_{\mathrm{E}} &=&\sup\left\{ 
        \frac{\text{R.H.S.\eqref{Eq:Esseen1956_asymptotics}}}
              {\phantom{\int\limits^{.}}\nu_3(\widetilde{P})\phantom{\int\limits^{.}}} 
     : P\in\cP_3 \right\}, 
\al 
with the supremum 
attained exactly at $P=\mathrm{B}_{p^{}_\mathrm{E}}$ and at its nondegenerate affine-linear images.
Esseen's result (\ref{Eq:Esseen1956_asymptotics},\ref{Eq:c_E_really_is_...})
was sharpened significantly, and generalised to convolutions of not necessary identical laws, by
\citet[in particular part I, p.~231 in the English version, inequality~(2.14) in 
Theorem~2.2]{Chistyakov2001} and by \citet[in particular p.~303, Corollary 4.18]{Shevtsova2012}, 
lending further support to the still open conjecture 
that~\eqref{Eq:Berry-Esseen_inequality} might hold with $c=c^{}_\mathrm{E}$. 
\citet[p.~1, Theorem~1]{Schulz2016} proved this conjecture to be true at least in the Bernoulli case
(that is, \eqref{Eq:Berry-Esseen_Bernoulli} holds with $c=c^{}_{\mathrm{E}}$), containing in particular 
the case of $\widetilde{P}=\widetilde{\mathrm{B}_{p^{}_\mathrm{E}}}$
asymptotically worst
according to~(\ref{Eq:Esseen1956_asymptotics},\ref{Eq:c_E_really_is_...}).

Let us mention here two asides. First,  \citet[pp.~15--16, Theorem 1 and Remark 4.3]{Schulz2016}
also showed that, for $p\in[\frac13,\frac23]$ at least, but not for every $p\in\mathopen]0,1\mathclose[$,
\eqref{Eq:c(P)/sqrt(n)_bound} holds for $P=\mathrm{B}_p$ with the then obviously optimal 
$c(\mathrm{B}_p)=\text{R.H.S.\eqref{Eq:Esseen1956_asymptotics_Bernoulli}}
 =\text{R.H.S.\eqref{Eq:Esseen1956_asymptotics}}$.
This result is for $p\in[\frac13,\frac{2}{3}]\setminus\{ p^{}_\mathrm{E} , 1- p^{}_\mathrm{E}\}$  
strictly sharper than \eqref{Eq:Berry-Esseen_Bernoulli} with $c=c^{}_{\mathrm{E}}$, and 
much more difficult to prove than the special case of $p=\frac12$ 
obtained earlier by \citet{HippMattner2007}.
Second, the analogue of~\eqref{Eq:c_E_really_is_...} 
for the approximation error taken as
$ \sup\big\{ \left| \widetilde{P^{\ast n}}(I) -\mathrm{N}(I)\right|   
    : I\subseteq \R \text{ an 
                              interval}  \big\}$
instead of $\text{L.H.S.\eqref{Eq:Berry-Esseen_inequality}}
  = \sup\big\{ \left| \widetilde{P^{\ast n}}(I) -\mathrm{N}(I)\right|   
     : I\subseteq \R      \text{ an unbounded interval}  \big\}$ 
looks a bit more elegant, with then the symmetric Bernoulli law $\mathrm{B}_\frac{1}{2}$ being 
extremal and with $\frac{2}{\sqrt{2\pi}}<2\,c^{}_{\mathrm{E}}$ playing the role of $c^{}_{\mathrm{E}}$, 
and is easier to prove as \citet{DinevMattner} showed.
   
To prepare for a return to our discussion of~\eqref{Eq:Berry-Esseen_inequality}, 
and for a later use in the proof of Theorem~\ref{Thm:BoChiGoe}, let us recall 
the equivalence
\la                         \label{Eq:Elementary_Cramer_Levy}
 \widetilde{P^{\ast n}} = \mathrm{N}  &\iff& \widetilde{P}=\mathrm{N}\quad\text{ for }n\in\N\text{ and }P\in\cP_2\,,
\al
where the elementary converse is due to $\mathrm{N}^{\ast n}=\mathrm{N}_{\sqrt{n}}$ and hence 
$\widetilde{P^{\ast n}} =\widetilde{\widetilde{P}^{\ast n}} = \widetilde{\mathrm{N}_{\sqrt{n}}}=\mathrm{N}$,
and the not completely trivial direct half is a simple special case of the 
\citet{Cramer1936}-L\'evy theorem, and is obtainable by assuming w.l.o.g.~$P=\widetilde{P}$
and observing that the corresponding Fourier transforms 
$\widehat{\mathrm{N}},\widehat{P}:\R\rightarrow \C$
are continuous with $\widehat{P}(0)=1>0$ and $(\widehat{P}(\frac{t}{\sqrt{n}}))^n=\widehat{\mathrm{N}}(t)=\exp(-\frac{t^2}{2})>0$
and hence $\widehat{P}(t)=(\widehat{\mathrm{N}}(\sqrt{n}t))^\frac{1}{n} =\widehat{\mathrm{N}}(t)$ for, 
respectively, $t\in\R$. Here ``not completely trivial'' refers to nonuniqueness in general of 
convolution roots in $\Prob(\R)$, as for example in \citet[p.~506, Curiosity (iii)]{Feller1971}.

\subsection{Zolotarev's Problem~\ref{Zolotarev-Problem}}  \label{subsec:Zolotarev-Problem}
We now observe that, compared to the asymptotic central limit theorem error
in~R.H.S.\eqref{Eq:Esseen1956_asymptotics}, 
the quantity $c\,\nu_3(\widetilde{P})$ on the right in
the Berry-Esseen theorem~\eqref{Eq:Berry-Esseen_inequality} 
has the defect of never being small, since we have 
\la                                                       \label{Eq:nu_3(widetilde(P))_ge_1}
 \nu_3(\widetilde{P})&\ge& (\nu_2(\widetilde{P}))^{3/2}\,\ =\,\ 1 \quad\text{ for }P\in\cP_3\,,
\al 
by, say,  Jensen's inequality applied to the convex function $\mathopen[0,\infty\mathclose[\ni t\mapsto t^{3/2}$,
and~\eqref{Eq:nu_3(widetilde(P))_ge_1} in particular holds for $\widetilde{P}=\mathrm{N}$, 
with
\la            \label{Eq:nu_3_von_N}  \label{Eq:nu_1_and_nu_3_of_N}
 && \nu_1(\mathrm{N}) \,=\,  \frac{2}{\sqrt{2\pi}} \ =\,  0.79788\ldots  \,, \quad
 \nu_2(\mathrm{N}) \, =\, 1\,,\quad 
 \nu_3(\mathrm{N}) \, =\, \frac{4}{\sqrt{2\pi}} \,\ =\,\ 1.595769\ldots,
\al
in which case 
L.H.S.\eqref{Eq:Berry-Esseen_inequality} actually vanishes by~\eqref{Eq:Elementary_Cramer_Levy}.
On the other hand, we have by a simple argument,
as mentioned by \citet{Zolotarev1972}
and associated somewhat imprecisely to \citet{Levy1937}
by \citet[p.~531]{Zolotarev1973}, 
\la &&                                          \label{Eq:simple_Kolomogorov-telescoping}
 \left\|\widetilde{P^{\ast n}}-\mathrm{N}\right\|_{\mathrm{K}}
 \ =\ \left\|\widetilde{P}^{\ast n}-\mathrm{N}_{\sqrt{n}}\right\|_{\mathrm{K}}
 \ =\ \left\|\widetilde{P}^{\ast n}-\mathrm{N}^{\ast n}\right\|_{\mathrm{K}}
 \ \le\  n \left\|\widetilde{P}-\mathrm{N}\right\|_{\mathrm{K}} 
 \quad\text{ for } P\in\cP_2, n\in\N,
\al
namely by using scale invariance of the Kolmogorov distance in the first step,
and in the final step a simple telescoping argument given more generally 
as~\eqref{Eq:semiadditivity} below.
The inequality in~\eqref{Eq:simple_Kolomogorov-telescoping}
makes precise in some way the idea that 
$\left\|\widetilde{P^{\ast n}}-\mathrm{N}\right\|_{\mathrm{K}}$ 
should be small if $P$ is close to normal, but of course, in contrast to~\eqref{Eq:Berry-Esseen_inequality},
the dependence on~$n$ of $\text{R.H.S\eqref{Eq:simple_Kolomogorov-telescoping}}$     
is rather unhelpful.
Aiming then at combining the virtues of~\eqref{Eq:Berry-Esseen_inequality} and~\eqref{Eq:simple_Kolomogorov-telescoping},
it appears natural to pose a problem like~\ref{Zolotarev-Problem}  
below. We suggest to name it after its apparent originator, 
Vladimir Mikhailovich Zolotarev (1931--2019, see \citet{ZolotarevObituary}),
who in any case was a main contributor to its successively better solutions, 
and who in particular provided basic ingredients (the definition and basic properties of $\zeta$ metrics,
and essentially Theorem~\ref{Thm:Zolotarev's_zeta_1-B-E-Thm} below)
for the proof of this paper's purpose, Theorem~\ref{Thm:Berry-Esseen_for_Z-close_to_normal}
on page~\pageref{Thm:Berry-Esseen_for_Z-close_to_normal}. 
The actual wording of the problem is here chosen as to fit the solutions we can report below:

\begin{Problem}[Zolotarev, \citeyear{Zolotarev1965} and several further works cited below]
                                                                \label{Zolotarev-Problem}
Find a nice sequence of metrics $d_n$ on~$\widetilde{\cP_3}$, perhaps decreasing in $n$,
and perhaps simply $d_n=d$ constant in~$n$, such that there exists a constant $c<\infty$ with 
\la                                                            \label{Eq:Problem_with_d_n}
 \left\|\widetilde{P^{\ast n}}-\mathrm{N}\right\|_{\mathrm{K}}
  &\le& \frac{c}{\sqrt{n}}d_n(\widetilde{P},\mathrm{N}) \quad \text{ for }P\in\cP_3\text{ and }n\in\N\,. 
\al
Or, more ambitiously, find a nice norm $\|\cdot\|$ on 
$                                                        
 \cM_{3,2} = \big\{M\in\cM_3: \mu_0(M)=\mu_1(M)=\mu_2(M)$ 
 $=$ $0\big\}
$
and a corresponding starting point $n_0\in\N$ such that there exists a constant $c<\infty$ with 
\la                                                                      \label{Eq:Problem_with_norm}
 \left\|\widetilde{P^{\ast n}}-\mathrm{N}\right\|_{\mathrm{K}}
  &\le& \frac{c}{\sqrt{n}}\left\|\widetilde{P}-\mathrm{N}\right\| 
  \quad \text{ for }P\in\cP_3\text{ and }n\ge n_0\,.
\al  
\end{Problem}

Here the somewhat vague adjective ``nice'' could be made a bit more precise as 
``rather easy to compute or bound, and then as weak as possible'' 
for the arguments actually occurring, that is, on $\widetilde{\cP_3}\times\{\mathrm{N}\}$
for $d_n$, and on 
\[
  \widetilde{\cP_3}-\mathrm{N} &\coloneqq& \{ P-\mathrm{N} : P\in \widetilde{\cP_3}\}
\]
for $\|\cdot\|$.
For example we will see already at~\eqref{Eq:Paulauskas1969-distance}
that ``as weak as possible'' excludes the case of 
\la                                                  \label{Eq:Trivial_distance}
 d_n(\widetilde{P}, \mathrm{N}) &\coloneqq& 
   \frac{4}{\sqrt{2\pi}}(\widetilde{P}\neq \mathrm{N})+ \nu_3(\widetilde{P}-\mathrm{N}), 
\al  
where we have written the usual discrete metric using~\eqref{Eq:indicator_notation};
this indeed yields a solution to~\eqref{Eq:Problem_with_d_n},
and trivially so given the Berry-Esseen theorem~\eqref{Eq:Berry-Esseen_inequality}, 
using~\eqref{Eq:nu_3_von_N} and hence 
$ \nu_3(\widetilde{P}) = \frac{4}{\sqrt{2\pi}} + \nu_3(\widetilde{P}) -\nu_3(\mathrm{N})
 \le d_n(\widetilde{P},\mathrm{N})$ for $P$ not normal.
By contrast, 
\la
 d_n(\widetilde{P}, \mathrm{N}) &\coloneqq& \text{R.H.S.\eqref{Eq:Esseen1956_asymptotics}}
 \,\ = \,\ \frac{1}{\sqrt{2\pi}}\left(\frac{|h(\widetilde{P})-h(\mathrm{N})| }2 
  +  \frac{|\mu_3(\widetilde{P})-\mu_3(\mathrm{N})  |}{6} \right) 
\al
may initially look like a perhaps suitable (quasi-)metric distance, but is obviously too weak to
make~\eqref{Eq:Problem_with_d_n} generally true;
in fact for  $P\in\cP_3$ nonnormal, $\text{L.H.S.\eqref{Eq:Problem_with_d_n}}>0$ 
by~\eqref{Eq:Elementary_Cramer_Levy},
but $\text{R.H.S.\eqref{Eq:Problem_with_d_n}}=0$ whenever $P$ is nonlattice and, for example, symmetric.
Concerning~\eqref{Eq:Problem_with_norm}, we will see below 
in~(\ref{Eq:Shiganov1987-distance_1},\ref{Eq:Shiganov1987-distance_2},%
\ref{Eq:Shiganov1987-distance_3},\ref{Eq:Berry-Esseen_for_Z-close_to_normal},%
\ref{Eq:Berry-Esseen_for_varkappa-close_to_normal})
that allowing here an $n_0>1$ admits 
weaker norms $\|\cdot\|$ than it would otherwise be the case.

We proceed to review known nontrivial solutions to Problem~\ref{Zolotarev-Problem}
in subsections~\ref{subsec:Known_sulutions_with_nu},  \ref{subsec:Known_sulutions_with_kappa},  
and~\ref{subsec:Known_solutions_with_zeta}, 
but for simplicity mention there papers treating more general or related questions 
only as far as their specialisations contribute to the present setting.
So we do not explicitly review related results for 
higher dimensions, distributions possibly nonidentical or 
without third moments, the Kolmogorov norm $\|\cdot\|_{\mathrm{K}}$ on the left hand side
replaced by \citeposs[the second theorem]{Nagaev1965}
weighted version R.H.S.\eqref{Eq:BE_Nag_Z?} often called ``nonuniform'',
or by the total variation norm $\nu_0$ as for example in 
\citet[p.~1254, Theorem~4]{Boutsikas2011},
error bounds for  
short Edgeworth expansions as provided by
\citet[p.~54, Corollary~3.1]{Yaroslavtseva2008b},
or for gamma approximations as in \citet[p.~594, Theorem~3.2]{Boutsikas2015},
or for stable rather than normal approximations as in \citet{ChristophWolf1992}.
Clearly we thus can provide at best a partial picture of the relevant literature. 

\subsection{Known solutions with $\nu$ distances (weighted total variation norms)
                                to normality}     \label{subsec:Known_sulutions_with_nu}
After a pioneering result of \citet{Zolotarev1965}, who obtained the bound
$\left\|\widetilde{P^{\ast n}}-\mathrm{N}\right\|_{\mathrm{K}}\le c\,(\nu_3(\widetilde{P}-\mathrm{N})/\sqrt{n})^\frac{1}{4}$, 
yielding in $n$ the rate $n^{-\frac{1}{8}}$ rather than $n^{-\frac12}$,  
and after related seminar talks of Zolotarev in Vilnius as recalled
in \citet[p.~430]{BloznelisRackauskas2019},
the apparently first nontrivial solution to Problem~\ref{Zolotarev-Problem} as stated here was given 
by \citet{Paulauskas1969}: \eqref{Eq:Problem_with_d_n} holds with 
\la                                                                         \label{Eq:Paulauskas1969-distance}
 d_n(\widetilde{P},\mathrm{N}) &\coloneqq& \left(\nu_3^{\,\frac14} \vee \nu_3 \right)\!(\widetilde{P}-\mathrm{N})
\al
with $c$ unspecified as, unless the contrary is stated, in all further results reviewed here, 
and where $\vee$ indicates the usual supremum of functions (that is, pointwise maximum).
The distance in~\eqref{Eq:Paulauskas1969-distance}, unlike R.H.S.\eqref{Eq:Trivial_distance}, 
can be arbitrarily close to zero also for $P\in\cP_3$ not normal, and so with it~\eqref{Eq:Problem_with_d_n}
strictly improves i.c.f.~the Berry-Esseen theorem~\eqref{Eq:Berry-Esseen_inequality},
since we have,
recalling first~\eqref{Eq:nu_3_von_N} and then~\eqref{Eq:nu_3(widetilde(P))_ge_1}  
for deducing~\eqref{Eq:pseudomoment_le_c_times_moment} from
\eqref{Eq:nu_r_subadditive_special},
\la                                                        
 \nu_r(\widetilde{P}-\mathrm{N})                       \label{Eq:nu_r_subadditive_special}
   &\le&  \nu_r(\widetilde{P}) + \nu_r(\mathrm{N})    \quad\text{ for }r\in[0,\infty[\,, \\
 \nu_3(\widetilde{P}-\mathrm{N})                  \label{Eq:pseudomoment_le_c_times_moment}
  &\le&    \left(1+\frac{4}{\sqrt{2\pi}}\right)\nu_3(\widetilde{P})\,,
\al 
and hence 
$\text{R.H.S.\eqref{Eq:Paulauskas1969-distance}}\le 1\!\vee\! \nu_3(\widetilde{P}-\mathrm{N})
 \le \text{R.H.S.\eqref{Eq:pseudomoment_le_c_times_moment}}$.
On the other hand, \eqref{Eq:Paulauskas1969-distance} has two obvious defects, 
namely the bad exponent $\frac{1}{4}$,  
which is however not simply omittable if $n\le3$, by 
\eqref{Eq:Zolotarev_example_norm_error_n_le_4_K_nu_0} 
and~\eqref{Eq:Zolotarev_example_r-norm_errors_n=1} with $r=3$
in Example~\ref{Example:Zolotarev_1973}, 
and the strength of the norm $\nu_3$: We have equality in~\eqref{Eq:nu_r_subadditive_special}
for example whenever $P\in\cP_3$ is discrete, and in this case~\eqref{Eq:Problem_with_d_n}
with~\eqref{Eq:Paulauskas1969-distance} 
is i.c.f.~just equivalent to~\eqref{Eq:Berry-Esseen_inequality}.

The problem of the bad exponent in \eqref{Eq:Paulauskas1969-distance} was solved by 
\citet{Sazonov1972}, \citet{Zolotarev1973}, \citet{Salakhutdinov1978}, \citet{Ulyanov1978},
and \citet{Shiganov1987}: \eqref{Eq:Problem_with_norm} holds with 
\la                                                   \label{Eq:Sazonov1972-Zolotaref1973-distance}
  \|\cdot\| &\coloneqq& \nu_0\vee\nu_3\,, \quad n_0\,\ \coloneqq\,\ 1\,,  \quad c\,\ \coloneqq\,\ 1.8\,,
\al 
and \eqref{Eq:Problem_with_d_n} holds with each of the following three choices
\la
   d_n(\widetilde{P},\mathrm{N}) &\coloneqq&  \left(\nu_1^{1\wedge\frac{n}{2}}\vee\nu_3\right)\!(\widetilde{P}-\mathrm{N})\,,
    \quad c \,\ \coloneqq\,\ 4.2\,,  \label{Eq:Shiganov1987-distance_1} \\
   d_n(\widetilde{P},\mathrm{N}) &\coloneqq&  \left(\nu_2^{1\wedge\frac{n}{3}}\vee\nu_3\right)\!(\widetilde{P}-\mathrm{N})\,,
    \quad c \,\ \coloneqq\,\ 13.5\,, \label{Eq:Shiganov1987-distance_2} \\ 
   d_n(\widetilde{P},\mathrm{N}) &\coloneqq&  \left(\nu_3^{1\wedge\frac{n}{4}}\vee\nu_3\right)\!(\widetilde{P}-\mathrm{N})\,,
    \quad c \,\ \coloneqq\,\ 35\,;  \label{Eq:Shiganov1987-distance_3}
\al
so in particular, for each $r\in\{0,1,2,3\}$,  \eqref{Eq:Problem_with_norm} holds with 
$\|\cdot\| \coloneqq \nu_r\vee \nu_3$ and \mbox{$n_0\coloneqq r+1$}.
Here~\eqref{Eq:Sazonov1972-Zolotaref1973-distance} but with $c$ unspecified was proved,  apparently independently
and at any rate differently, 
by \citet[p.~570, Theorem 3.1 for dimension $k=1$]{Sazonov1972} and by~\citet[p.~533, Theorem, inequality~(7)]{Zolotarev1973},
\eqref{Eq:Shiganov1987-distance_3} but with $c$ unspecified was published for $n\ge 9$ without proof  by
\citet[special case of Corollary~1]{Salakhutdinov1978},
and obtained for general $n$, with a sketch of a proof, by~\citet[Theorem~3 for dimension $k=1$, Lemma~2(a)]{Ulyanov1978},
and the rest 
of (\ref{Eq:Sazonov1972-Zolotaref1973-distance}--\ref{Eq:Shiganov1987-distance_3}) is due to \citet{Shiganov1987}.

I.c.f., \eqref{Eq:Problem_with_d_n} 
with (\ref{Eq:Sazonov1972-Zolotaref1973-distance}--\ref{Eq:Shiganov1987-distance_3})
combined by taking a minimum is equivalent to~\eqref{Eq:Problem_with_d_n} with  
\la                              \label{Eq:Shiganov_combined}
 d_n(\widetilde{P},\mathrm{N}) &\coloneqq& \min_{r=3\wedge(n-1)}^3
  \left(\nu_r^{1\wedge\frac{n}{r+1}}\vee\nu_3\right)\!(\widetilde{P}-\mathrm{N})\,,
\al   
namely obviously so if $n=1$ or $n\ge 4$, and if  $n\in\{2,3\}$
by applying Lyapunov's inequality,
\la                                    \label{Eq:Lyapunov_three_moments_inequality}
  \nu_s &\le& \nu_r^\frac{t-s}{t-r} \nu_t^\frac{s-r}{t-r}
     \,\ \le \,\ \nu_r\vee\nu_t
  \quad \text{ on $\,\cM\,$ for $\,0\le r \le s\le t<\infty\,, \tfrac{0}{0}\coloneqq 1$\,},
\al
with $s\coloneqq3\wedge(n-1)=n-1$ and $t\coloneqq 3$ to show that for $r\in\{0,\ldots,s-1\}$
we have $\nu_s \le \nu_r\vee\nu_3$, and hence, using $\frac{n}{r+1}>\frac{n}{s+1}\ge 1$, 
we get
$  
 \nu_s^{1\wedge\frac{n}{s+1}}\vee\nu_3 
 = \nu_s\vee\nu_3\le \nu_r\vee\nu_3 = \nu_r^{1\wedge\frac{n}{r+1}}\vee\nu_3
$
and see that indeed $\left(\nu_r^{1\wedge\frac{n}{r+1}}\vee\nu_3\right)\!(\widetilde{P}-\mathrm{N})$
with the present $r$ is irrelevant for the minimum.

Similarly to~\eqref{Eq:nu_3(widetilde(P))_ge_1}, we get $\nu_r(\widetilde{P})\le 1$ for
$r\in[0,2]$ , and with~\eqref{Eq:nu_r_subadditive_special} hence
\la                                                   \label{Eq:nu_r_bounded_for_r_le_2}
 \nu_r(\widetilde{P}-\mathrm{N}) &\le& 1+\nu_r(\mathrm{N}) 
    \,\ \le \,\ 2\quad\text{ for }r\in[0,2]\text{ and }P\in\cP_2\,.
\al
This shows in particular that
in an inequality~\eqref{Eq:Problem_with_d_n}
with $d_n(\widetilde{P},\mathrm{N})
= \left(\nu_r^{\alpha_n }\vee\nu_3\right)\!(\widetilde{P}-\mathrm{N})$
with some $r\in\{0,1,2,3\}$, it is i.c.f.~preferable to have $\alpha_n$ as large as
possible, but $\le$~$1$ in case of $r=3$.
Analogous remarks apply below to bounds involving $\varkappa_r$ or $\zeta_r$
in place of $\nu_r$.

As essentially known from \citet[upper bounds for $r(n)$ in (5)]{Zolotarev1972} 
and proved more explicitly by \citet[Examples 1.2 and 1.3]{Yaroslavtseva2008b},
each of the exponents $1\wedge\frac{n}{r+1}$ in~\eqref{Eq:Shiganov_combined}
is optimal: This is trivially so if $n\ge4$.
If $n\le 3$, we observe first that decreasing $1 \wedge\frac{n}{r+1}$ would i.c.f.\ 
worsen~\eqref{Eq:Problem_with_d_n}  with~\eqref{Eq:Shiganov_combined},
by~\eqref{Eq:nu_r_bounded_for_r_le_2} if $r\le2$ and trivially if $r=3$,
and second that increasing $1 \wedge\frac{n}{r+1}$ is inadmissible 
due to $n\le 3$ 
and~(\ref{Eq:Zolotarev_example_norm_error_n_le_4_K_nu_0},\ref{Eq:Zolotarev_example_r-norm_errors_n=1})
in Example~\ref{Example:Zolotarev_1973}.
Hence in particular the starting points $n_0=r+1$ for 
\eqref{Eq:Problem_with_norm} with $\|\cdot\| \coloneqq \nu_r\vee \nu_3$,
given above after~\eqref{Eq:Shiganov1987-distance_3}, are optimal.

\subsection{Known solutions with $\varkappa$ distances (weighted $\mathrm{L}^1$ 
norms of distribution functions)}                             \label{subsec:Known_sulutions_with_kappa}  
The problem of the strength of the eq\-norms~$\nu_r$ in 
(\ref{Eq:Paulauskas1969-distance}--
\ref{Eq:Shiganov1987-distance_3})
was
attacked by \citeposs{Zolotarev1970,Zolotarev1971,Zolotarev1972,Zolotarev1973}
introduction, to the area of normal approximation error bounds for convolution powers on~$\R$
(but see \citet[p.~31]{ChristophWolf1992}
for a few earlier references concerning assumptions for asymptotic expansions), 
of weaker eqnorms $\varkappa_r$.
Since the strong eqnorms~$\nu_r$ in (\ref{Eq:Sazonov1972-Zolotaref1973-distance}--\ref{Eq:Shiganov1987-distance_3})
may be thought of arising through
\la                               \label{Eq:|mu_r|_le_nu_r}
 &&\left|\mu_r(M)\right| \,\ =\,\ \left|\int x^r\,\dd M(x)\right|
  \,\ \le\,\ \int |x|^r\,\dd|M|(x)  \,\ =\,\ \nu_r(M) 
  \quad\text{ for }r\in\N\text{ and }M\in\cM_r\,,
\al
thus bounding in particular the too weak eqnorm $|\mu_3|$ occurring 
in~\eqref{Eq:Esseen1956_asymptotics}, 
the idea is to get weaker but hopefully still strong enough eqnorms $\varkappa_r$
by preparing a triangle inequality as in~\eqref{Eq:|mu_r|_le_nu_r}
by an integration by parts:

Recalling the notation (\ref{Eq:indicator_notation},\ref{Eq:Def_F_M_and_overline{F}_M}), let 
\la                    \label{Eq:Def_h_M}
  h_M &\coloneqq& \1_{]0,\infty[}\overline{F}_{\!M}- \1_{]-\infty,0[}F_M \quad\text{ for }M\in\cM\,.
\al
We then have and put,
referring to the proof of Lemma~\ref{Lem:Generalised_signed_moments}
for the easy justifications of~(\ref{Eq:Def_varkappa},\ref{Eq:|mu|_le_varkappa_le_nu}), 
\la                        \label{Eq:Def_varkappa}        
 \nu_r(M)  &=& \int r|x|^{r-1}\big|h_{|M|}(x)\big|\,\dd x    \\ \nonumber        
   &\ge &  \int r|x|^{r-1}\left| h_M(x)\right|\dd x     
     \,\ \eqqcolon \,\ \varkappa_r(M)\qquad \text{ for }
      r\in\mathopen]0,\infty\mathclose[\text{ and }
     M\in \cM\,,
\al
with equality throughout iff $M$ is of the same sign on each of $]-\infty,0[$ 
and $]0,\infty[\,$, and hence
\la                                                \label{Eq:|mu|_le_varkappa_le_nu}
 \quad && \mu_r(M) \,=\,  \int rx^{r-1}h_{M}(x)\,\dd x\, ,  \quad 
 \left|\mu_r(M)\right| 
     \, \le \, \varkappa_r(M) \, \le \, \nu_r(M) 
   \quad\text{ for }r\in\N,\,M\in\cM_r
\al
as desired. Introducing now the assumption $M(\R)=0$, 
and recalling the notation~\eqref{Eq:Def_cM_r,k} and hence 
$ \cM_{0,0}=\{M\in\cM : M(\R)=0\}$, we get 
\la                                          \label{Eq:varkappa_if_M(IR)=0}
  \varkappa_r(M) &=& \int_\R r|x|^{r-1}\big|F_M(x)\big| \,\dd x
  \quad \text{ for }r\in\mathopen]0,\infty\mathclose[\,\text{ and }\, M\in \cM_{0,0} \,,  
\al
and in particular 
\la                 \label{Eq:varkappa_1=L^1-enorm}
 \varkappa_1(M) &=& \left\| F_M \right\|_1  \quad \text{ for } M\in \cM_{0,0}\,, 
\al
where on the right in~\eqref{Eq:varkappa_1=L^1-enorm} 
we have the usual $\mathrm{L}^1$ enorm, with respect to Lebesgue measure~$\leb$ 
on~$\R$, of the distribution function $F_M$. 
For $M=P-Q$ with $P,Q\in\Prob(\R)$, $\varkappa_1(M)$ is also known as
the Kantorovich or, historically less appropriately but still more commonly,
the Wasserstein distance between $P$ and $Q$.

As an aside, let us mention four early theorems,
             \label{page:aside_on_early_theorems_with_varkappa_1}
in the probabilistic literature, where the not immediately
probabilistically interpretable quantity
$\varkappa_1(P-Q)=\int|F_P(x)-F_Q(x)|\,\dd x$ occurs for $P,Q\in\Prob_1(\R)$.
First, \citet[p.~30, Theorem 1]{Esseen1945} bounds  $\varkappa_1(P-Q)$ by $\frac{\pi}{T}$ 
if the Fourier transforms
$\hat{P},\hat{Q}$ coincide on the interval $[-T,T]$. 
Second,  \cite[p.~277, (4.7)]{FortetMourier1953} prove that   
$\varkappa_1(P-Q)$ is $\zeta_1(P-Q)$ from~\eqref{Eq:Def_zeta_r} below,
but rather defined by the second expression in~\eqref{Eq:Four_zeta_expressions},
apparently motivated by a desire to elementarise their uniform functional strong law of large 
numbers, compare \citet[p.~277, ``Ainsi les th\'eor\`emes g\'eneraux \ldots'']{FortetMourier1953}.
Third, \cite[p.~801, (1.8) with $r=1$]{Agnew1954} gives a 
central limit theorem, namely the present~\eqref{Eq:Lindeberg-Levy-Thm}
with $\|\cdot\|_{\mathrm{K}}$ replaced by $\varkappa_1$, 
which actually is an obvious corollary to \citet[p.~70, Theorem 1]{Esseen1945}
combined with~\eqref{Eq:Lindeberg-Levy-Thm}.
For these first three theorems, apparently    
no probabilistic interpretation was either obvious or supplied, 
as remarked for the third by \citet{Morgenstern1955}. 
Fourth, \citet[p.~42, Teorema I]{Dall'Aglio1956} proves that $\varkappa_1(P-Q)$
is the minimal transport 
distance
$W(P,Q)\coloneqq \inf\{ \int |x-y| \,\dd R(x,y): R \in\Prob(\R\times\R)\text{ with marginals }P,Q \}$.
The second and fourth theorems combined yield for the real line 
the theorem of \citet{KantorovichRubinstein1958} 
as presented by \citet[p.~421, Theorem 11.8.2]{Dudley2003}, 
in the present notation $W(P,Q)=\zeta_1(P-Q)$.

For some further references and more detailed historical remarks about
the minimal transport distance $W(P,Q)$
one might start with \citet{Rueschendorf2000}, \citet[p.~435]{Dudley2003}, 
\citet[pp.~106--107]{Villani2009},
and \citet[Introduction]{BogachevKolesnikov2012}.
In particular Villani and also \citet[pp.~788--789]{BogachevKolesnikov2012}
state that \citet{Kantorovich1942} introduced
$W(P,Q)$, for any compact metric space  rather than the real line,
and their statement is accurate if one
takes into account an equivalence
between bimeasures and product measures,  as given by
\citet[p.~166, section 1, claim (i)]{MarczewskiRyll-Nardzewski1953}.
Later,
\citet{Dobrushin1970} named $W(P,Q)$ after \citet{Vasershtein1969},
presumably in ignorance of earlier occurrences,
and most followers of this naming, when using the Latin alphabet,
seem to prefer ``Wasserstein'' over ``Vasershtein'',
presumably since they imagine the latter as the result of two not exactly inverse
transliterations, from Latin to Cyrillic and back,
of the former German name.

Back to our main theme,
each $\varkappa_r$ is not only bounded from above by $\nu_r$
as noted in~\eqref{Eq:Def_varkappa}, but quite obviously  strictly weaker,
and this even as far as just convergence in $\widetilde{\cP_3}$ to normality is concerned,
as shown by, say, the binomial central limit example
\la                                \label{Eq:Binomial_CLT_varkappa_r_but_not_nu_r}
 && \lim_{n\rightarrow\infty} \varkappa_r\big(\widetilde{\mathrm{B}_p^{\ast n}} - \mathrm{N} \big)
  \ =\ 0 \ \neq \ 2\,\nu_r(\mathrm{N}) \ =\
   \lim_{n\rightarrow\infty} \nu_r\big(\widetilde{\mathrm{B}_p^{\ast n}} - \mathrm{N} \big)
   \quad\text{ for }p\in\mathopen]0,1\mathclose[\,,\ r\in \mathopen]0,\infty\mathclose[\,,
\al
where the first equality holds, for example, 
by Remark~\ref{Rem:CLT_rate_varkappa_r}\ref{part:Osipov_applied_to_varkappa_r},
and it implies 
$\nu_r(\widetilde{\mathrm{B}_p^{\ast n}}) =  \varkappa_r(\widetilde{\mathrm{B}_p^{\ast n}}) 
\rightarrow  \varkappa_r(\mathrm{N})=\nu_r(\mathrm{N})$
by equality in~\eqref{Eq:Def_varkappa} for $M\ge0$
and by the qnorm property of~$\varkappa_r$ on~$\cM_r$, 
and the final equality hence follows using
$\nu_r\big(\widetilde{\mathrm{B}_p^{\ast n}} - \mathrm{N} \big)
 = \nu_r\big(\widetilde{\mathrm{B}_p^{\ast n}}\big) +\nu_r(\mathrm{N} )$.

Obtaining then a theorem like~\eqref{Eq:Problem_with_d_n} 
with~\eqref{Eq:Shiganov1987-distance_1}, but with $\nu_1$ and $\nu_3$ improved,
albeit at the cost of reobtaining somewhat less than ideal exponents,
\citet{Ulyanov1976} sharpened another result of \citet{Zolotarev1973} to the before 
Corollary~\ref{Cor:BE_K_varkappa} below
best solution to Problem~\ref{Zolotarev-Problem} in terms of $\varkappa_1$ and $\varkappa_3$ 
known to us: \eqref{Eq:Problem_with_d_n} holds with 
\la                                                                  \label{Eq:Ulyanov1976-distance}
 d_n(\widetilde{P},\mathrm{N}) &\coloneqq& 
  \left( (\varkappa_1\!\vee\!\varkappa_3)^{1-2^{-n}}\vee\varkappa_3\right)(\widetilde{P}-\mathrm{N}).
\al 
More precisely, \citet[p.~533, the definition of $\varkappa_0$]{Zolotarev1973} 
and \citet[p.~270, the definition of $\kappa$ in case of $g=|\cdot|$]{Ulyanov1976}
considered $\int\max(1,cx^2)|\widetilde{F}(x)-\Phi(x)|\dd x$
with $c\in\{3,1\}$, which can  be replaced i.c.f.~equivalently by
$\big(\varkappa_1\!\vee\!\varkappa_3\big)(\widetilde{P}-\mathrm{N})$,
and with this replacement \citet[p.~533, Theorem, (6)]{Zolotarev1973} 
yields~\eqref{Eq:Problem_with_d_n} with 
$d_n(\widetilde{P},\mathrm{N}) \coloneqq 
\left( (\varkappa_1\!\vee\!\varkappa_3)^{\frac{n}{n+1}}\vee\varkappa_1\vee\varkappa_3\right)(\widetilde{P}-\mathrm{N})$, 
\citet[pp.~271 and 282, Theorem 1 with $g\coloneqq|\cdot|$]{Ulyanov1976}
improves the exponent $\frac{n}{n+1}$ to $1-2^{-n}$, and finally 
the $\varkappa_1$ with exponent $1$ can be omitted due to boundedness of $\varkappa_1$ on 
$\widetilde{\cP_2}-\widetilde{\cP_2}$, or more precisely  by
\la                                         \label{Eq:zeta_1_le}
 \varkappa_1(\widetilde{P}-\mathrm{N}) &\le& 1+\nu_1(\mathrm{N})
  \,\ =\,\ 1+\frac{2}{\sqrt{2\pi}}  
  \,\ =\,\  1.79788\ldots \quad\text{ for }P\in\cP_2\,,   
\al
which holds by (\ref{Eq:Def_varkappa},\ref{Eq:nu_r_bounded_for_r_le_2},\ref{Eq:nu_1_and_nu_3_of_N}).

Further, \citet[p.~661, Corollary]{Ulyanov1978} yields, by further specialisation, 
that~\eqref{Eq:Problem_with_d_n} 
also holds with 
\la                                                       \label{Eq:Ulyanov1978-distance}
 d_n(\widetilde{P},\mathrm{N}) &\coloneqq& 
  \left(\varkappa_3^{\frac{1}{2}(1-2^{-n})} \vee \varkappa_3\right)(\widetilde{P}-\mathrm{N}),
\al 
improving an exponent in \citet[p.~533, Theorem, (5)]{Zolotarev1973}.
But at each of the here considered arguments $\widetilde{P}-\mathrm{N}$,
we have $\varkappa_1 \le \sqrt{\frac{4}{\sqrt{2\pi}}\varkappa_3}$
by~\eqref{Eq:varkappa_1_vs_varkappa_r_and_Lip_centred} with $r\coloneqq 3$ in 
the quite simple Lemma~\ref{Lem:varkappa_1_vs_varkappa_r_and_Lip}  below,
whereas $\sqrt{\varkappa_3}/\varkappa_1$ can be arbitrarily large
even if $\varkappa_3$, and hence also $\varkappa_1$, is small,
for example by using
$\varkappa_3(\widetilde{P}-\mathrm{N})\ge 6\,\zeta_3(\widetilde{P}-\mathrm{N})$
from~\eqref{Eq:zeta_vs_varkappa_etc}, and either one
of~\eqref{Eq:zeta_3_1_truncated_normal_standardised} or~\eqref{Eq:zeta_3_1_truncated_normal_winsorised}
from Example~\ref{Examples:zeta_3=mu_3}(a,b).
Hence, i.c.f., \eqref{Eq:Problem_with_d_n} with~\eqref{Eq:Ulyanov1978-distance}
not only follows easily from~\eqref{Eq:Problem_with_d_n} with~\eqref{Eq:Ulyanov1976-distance}, 
but is also strictly worse.  

\subsection{Introducing Zolotarev's $\zeta$ distances (dual to smooth function norms)}                          \label{subsec:Intro_zeta}
Leaving aside for a moment the problem of nonideal exponents in bounds 
like~\eqref{Eq:Problem_with_d_n} with~\eqref{Eq:Ulyanov1976-distance},
it turns out that bounds with 
$\varkappa_3$ on~$\cM_{3,2}$ replaced by an even weaker norm, 
namely \citeposs{Zolotarev1976} $\zeta_3$, can be obtained easily, given 
\citeposs{Zolotarev1973} paper. We may, as \citet{Christoph1979} essentially did, 
apparently independently of \citet{Zolotarev1976}, introduce~$\zeta_3$,
and more generally~$\zeta_r$ with in this paper for simplicity $r\in\N_0$,
similarly to~$\varkappa_r$ in~\eqref{Eq:Def_varkappa} above, 
roughly speaking by performing $r$~integrations by parts on $\int x^r\,\dd M(x)$,
rather than just zero as in~\eqref{Eq:|mu_r|_le_nu_r}
or one as in~\eqref{Eq:Def_varkappa}.
This leads to the expression $\int\!\left| h_{M,r} \right|\dd\leb$
in~\eqref{Eq:Four_zeta_underbar_expressions}, with $h_{M,r}$ defined 
by~(\ref{Eq:Def_F_M,k},\ref{Eq:Def_h_M,k}), 
and to an alternative representation in~\eqref{Eq:Def_zeta_r} via~\eqref{Eq:int_g_dM_formal}.  

To be more precise, let us put 
\la                       \label{Eq:Def_cG_k,alpha}
 &&\,\cG_{k,\alpha} \,\coloneqq\, 
   \big\{ 
    g \in\C^\R : g^{(k-1)}\text{ absolutely continuous},\, 
    \left\| \tfrac{g^{(k)}}{1+|\,\cdot\,|^\alpha}\right\|_\infty < \infty
   \big\}  \text{ for }k\in\N,\alpha\in[0,\infty[\,,
\al
where of course the derivative $g^{(k-1)}$ of order $k\!-\!1$ is assumed to exist everywhere,
and the down-weighted $\mathrm{L}_{}^\infty$ norm 
is taken of the in general only $\leb$-a.e.~defined $k$th derivative $g^{(k)}$.
%
%
We also recall the definition of the sets $\cM_{r,k}$ from~\eqref{Eq:Def_cM_r,k},
and the notation~{\rm(\ref{Eq:Def_F_M_and_overline{F}_M},\ref{Eq:Def_h_M})}.

\begin{Lem}[On successive integration by parts with signed measures]  
                        \label{page:Lemma_F_M,k_and_int_g_dM}        \label{Lem:F_M,k_and_int_g_dM}
Let  $k\in\N$, $M\in \cM_{k-1}$, and 
\la    \quad                                  \label{Eq:Def_F_M,k}
 F^{}_{M,k}(x) &\!\coloneqq \!& \int\limits_{]-\infty,x]}\!\!\frac{(y-x)^{k-1}}{(k-1)!}\dd M(y),\quad
 \overline{F}^{}_{\!M,k}(x) \!\,\ \coloneqq\!\,\ 
   \int\limits_{[x,\infty[}\!\!\frac{(y-x)^{k-1}}{(k-1)!}\dd M(y)
   \quad\text{for }x\in\R, \\
 h^{}_{M,k} &\!\coloneqq\!&       \label{Eq:Def_h_M,k}
  \1^{}_{]0,\infty[}\overline{F}^{}_{\!M,k} -\1^{}_{]-\infty,0[} F^{}_{M,k}\,.
\al
\begin{parts}
\item  In case of $k=1$ we have    \label{part:F_M,k_with_k=1}
 $F_{M,1} = F_M$, $\overline{F}_{\!M,1}= \overline{F}_{\!M} = M(\R)- F_M(\cdot-)$, 
 and $h_{M,1}=h_M$.
\item                           \label{part:F_M_overline_complementary}
 Let $M\in\cM_{k-1,k-1}$. Then
 \la                     \label{Eq:h_M,k=-F_M,k}
    h^{}_{M,k} & = &  \overline{F}^{}_{\!M,k}  \,\ =\,\  -F^{}_{M,k}
 \al 
 holds except on a countable set. More precisely, the second equality in~\eqref{Eq:h_M,k=-F_M,k} 
 holds at $x\in\R$ except when $k=1$ and $M(\{x\})\neq 0$, 
 and the first then also holds except perhaps when $x=0$. 
\item                    \label{part:int_g_dM_formal}
 Let $\alpha\in[0,\infty[$ and $g\in\cG_{k,\alpha}$. Then we have
 \la                      \label{Eq:int_g_dM_formal}  
  \int g\,\dd M &=& \sum_{j=0}^{k-1}\frac{g^{(j)}(0)}{j!}\mu_j(M) 
    + \int g^{(k)} h^{}_{M,k}\,\dd \leb \quad\text{ if }M\in\cM_{k+\alpha}\,. 
 \al
\item                            \label{part:F_M,k+ell}
 Let also $\ell\in\N$, and $M\in\cM_{k+\ell-1,k-1}$. Then we have
 \la                                               \label{Eq:h_M,k+ell}
  h_{M,k+\ell}(x) &=& \left( (x>0)\!\int\limits_x^\infty -\ (x<0)\!\!\int\limits_{-\infty}^x \right) 
       \frac{(y-x)^{\ell-1}}{(\ell-1)!} h_{M,k}(y)\,\dd y  \quad\text{ for } x\in\R\,.
 \al
\end{parts}
\end{Lem}

This is proved in just a few lines starting on page~\pageref{page:Proof_of_Lemma_F_M,k_and_int_g_dM}
in section~\ref{sec:zeta_distances}.

For $r\in\N$, we consider now the following subsets of $\cG_{r,0}\,$:   
\la                         
 \cF_r &\coloneqq&             \label{Eq:Def_cF_r}
     \big\{ g\in\C^\R : g^{(r-1)}\text{ absolutely continuous},\ 
            \left\|g^{(r)}\right\|_\infty \le 1 \big\}\,, \\
 \cF_r^\infty             \label{Eq:Def_cF_r^infty} 
   &\coloneqq& \big\{ g\in \cF_r : \left\|g\right\|_\infty < \infty \big\}\,, \\
 \cF_{r,r-1} &\coloneqq&                            \label{Eq:Def_cF_r,r-1}   
   \big\{ g \in \cF_r : g(0)=g'(0)=\ldots=g^{(r-1)}(0)=0\big\}\, , \\ 
 \cF_{r,r-1}^\infty &\coloneqq&                            \label{Eq:Def_cF_r,r-1^infty}
    \cF_r^\infty \cap \cF_{r,r-1} \,\ =\,\ 
    \big\{ g \in \cF^\infty_r : g(0)=g'(0)=\ldots=g^{(r-1)}(0)=0\big\}\, .
\al
To include for later convenience also the case of $r=0$, we further put
\la                         \label{Eq:Def_cf_0,0} \label{Eq:Def_cF_0,0}
  \cF^\infty_{0,-1} &\coloneqq&    \cF^\infty_{0}  \,\ \coloneqq\,\   \cF_0 \,\ \coloneqq\,\
  \big\{ g\in\C^\R : g\text{ Borel and } 
    \sup_{x\in\R}|g(x)| \le 1 \big\}\,,
\al
with in~\eqref{Eq:Def_cF_0,0} 
a true rather than a merely $\leb$-essential supremum bound required,
and also $\cM_{0,-1}\coloneqq\cM$.

The following definition achieves the 
desire to replace $\varkappa_r$ from~\eqref{Eq:Def_varkappa} in case of $r\in\N$
by a weaker eqnorm, at least on some large subspace of $\cM$,
as shown by Lemma~\ref{Lem:Comparison_of_various_eqnorms} below, 
which is based on \eqref{Eq:int_g_dM_formal}.

\begin{Def}[$\zeta$ eqnorms, two variants]
Let $r\in\N_0$. For $M\in\cM$, we put 
\la                                   \label{Eq:Def_zeta_r}
 \zeta_r(M) &\coloneqq& \sup_{g\in\cF_r^\infty}\left| \int g\,\dd M\right| \,, \qquad
 \underline{\zeta}_r(M) \,\ \coloneqq\,\    \label{Eq:Def_zeta_underbar_r}
     \sup_{g\in\cF_{r,r-1}^\infty}\left| \int g\,\dd M\right| \,.
\al
\end{Def}

\begin{Lem}[Representations of $\zeta$ and $\underline{\zeta}$,
comparison with other eqnorms on $\cM$]    \label{Lem:Comparison_of_various_eqnorms}
                                           \label{page:Lemma_Comparison_of_various_eqnorms}
In parts 
{\rm \ref{part:zeta_underbar_le_zeta},\ref{part:zeta_underbar_on_cM_r},\ref{part:zeta_on_cM_r_r-1}} 
below, let $r\in\N$.

\begin{parts}
\item                                          \label{part:zeta_underbar_le_zeta} 
 On $\cM$, $\zeta_r$ is an enorm, $\underline{\zeta}_r$ is an eqnorm, and we have
 \la 
  && \underline{\zeta}_r \ \le \ \zeta_r \quad\text{everywhere,} 
                               \label{Eq:zeta_vs_zeta_underbar_on_cM_and_oncM_r,r-1}
    \qquad   \underline{\zeta}_r \ = \ \zeta_r \,\ <\,\ \infty \quad\text{on }\cM_{r,r-1}\,, \\
  &&  \underline{\zeta}_r \ <\ \zeta_r    \label{Eq:zeta_vs_zeta_underbar_on_cM_and_on_cM_r_setminus_cM_r,r-1}                   
       \ =\  \infty \quad\text{on } \cM_r\!\setminus\! \cM_{r,r-1} \,. 
 \al
In particular, on the vector space $\cM_{r,r-1}$, $\zeta_r$ and $\underline{\zeta}_r$ 
are identical norms.
\item Let $M\in\cM_r$. Then                   \label{part:zeta_underbar_on_cM_r} 
 \la
  && \zeta_r(M) \,\ = \,\ \sup_{g\in\cF_r}\left| \int\! g\,\dd M\right| \,,
                                                                     \label{Eq:zeta_on_cM_r_g_unbounded} \\ 
  &&\underline{\zeta}_r(M)   \label{Eq:zeta_underbar_g_unbounded}  \label{Eq:zeta_underbar_on_cM_r}   
    \,\ =\,\ \sup_{g\in\cF_{r,r-1}}\left| \int\! g\,\dd M\right| \label{Eq:Four_zeta_underbar_expressions}
    \,\ =\,\ \int \big| h_{M,r}\big| \,\dd\leb   \,\ =\,\ \int g\,\dd M \\  
  && \nonumber \phantom{ \underline{\zeta}_r(M) \,\ = \,\ }
      \text{ with }g\in \cF_{r,r-1}\text{ defined by }
       \ g^{(r)} = \sgn\circ\,h_{M,r}\ \text{ $\leb$-a.e.}, \\        
  &&                                     \label{Eq:zeta_underbar_vs_varkappa_etc} 
   \tfrac{1}{r!}\max\big\{\,|\mu_r(M)|\,,\,\big|\int\!|x|^r\dd M(x)\big| \,\big\}
   \,\ \le\,\ \underline{\zeta}_r(M) \,\ \le \,\ \tfrac{1}{r!}\varkappa_r(M)
     \,\ \le \,\ \tfrac{1}{r!}\nu_r(M) \,.
 \al
\item Let $M\in\cM_{r,r-1}$. Then               \label{part:zeta_on_cM_r_r-1}  
 \la                                  \label{Eq:Four_zeta_expressions}
  && \zeta_r(M)                                        
    \,\ =\,\ \sup_{g\in\cF_r}\left| \int\! g\,\dd M\right| 
    \,\ =\,\ \int \big| F_{M,r}\big| \,\dd\leb   \,\ =\,\ \int g\,\dd M \\
  && \nonumber \phantom{ \zeta_r(M) \,\ = \,\ }
      \text{ if }g\in \cF_{r}\text{ satisfies }
       \ g^{(r)} = -\sgn \circ\, F_{M,r}\ \text{ $\leb$-a.e.}, \\
  &&                                            \label{Eq:zeta_vs_varkappa_etc} 
   \tfrac{1}{r!}\max\big\{\,|\mu_r(M)|\,,\,\big|\int\!|x|^r\dd M(x)\big| \,\big\}
   \,\ \le\,\ \zeta_r(M) \,\ \le \,\ \tfrac{1}{r!}\varkappa_r(M)
     \,\ \le \,\ \tfrac{1}{r!}\nu_r(M) \,.
 \al    
\item For $r=1$ we have equality in the central inequalities  \label{part:zeta_1=_varkappa_1} 
 in {\rm(\ref{Eq:zeta_underbar_vs_varkappa_etc},\ref{Eq:zeta_vs_varkappa_etc})\,:}
 \la                              \label{Eq:zeta_1=kappa_1}
   && \underline{\zeta}_1 \ =\  \varkappa_1  \quad\text{on }\,\cM_1\,,    \label{Eq:zeta_1=kappa_1_new}
    \qquad \zeta_1 \ =\ \varkappa_1  \quad\text{on }\,\cM_{1,0}\,.
 \al  
\item For the case of $r=0$ excluded 
 in {\rm \ref{part:zeta_underbar_le_zeta},\ref{part:zeta_underbar_on_cM_r},\ref{part:zeta_on_cM_r_r-1}},
 we have
 \la
   && \left\|\,\cdot\,\right\|_{\mathrm{K}} \ \le \ \nu_0     \label{Eq:Kolmogorov_vs_nu_0_zeta_0}
       \ =\ \zeta_0  \ =\ \underline{\zeta}_0 \quad\text{on } \,\cM,
   \qquad\left\|\,\cdot\,\right\|_{\mathrm{K}} 
    \ \le \ \tfrac{1}{2}\nu_0 \quad\text{on } \,\cM_{0,0}\,.
 \al  
\end{parts}
\end{Lem}

This is also proved in section~\ref{sec:zeta_distances}, starting there on 
page~\pageref{page:proof_of_Lem_Comparison_of_various_eqnorms}. 
Here the claims~\eqref{Eq:zeta_on_cM_r_g_unbounded} and
the first identity in~\eqref{Eq:zeta_underbar_g_unbounded},
about removing boundedness assumptions in~\eqref{Eq:Def_zeta_r},
slightly generalise a ``well-known'' result actually proved in
\citet[p.~498, Theorem~1.7(d)]{MattnerShevtsova},
and are perhaps not completely trivial.

The third expression $\int \big| h_{M,r}\big| \,\dd\leb$
in~\eqref{Eq:Four_zeta_underbar_expressions} was,
for $M=P-Q$ with $P,Q\in\Prob(\R)$, but without assuming $M\in\cM_r$,
proposed by \citet[(4) with $r\in\N$, hence $\delta=0$, $\frac{1}{r!}\tau_r$]{Christoph1979}.
Following \citet[\S1.5, in particular p.~386]{Zolotarev1976},
it is however customary and usually convenient to use~\eqref{Eq:Def_zeta_r}
to define an enorm $\zeta_r$ on all of $\cM$.
According to \citet[p.~351]{Senatov1998}, $\zeta_r$ was thus
first introduced by Zolotarev in a seminar at the Steklov Institute of Mathematics  
in November 1975.
In \citet[p.~113]{BogachevDoledenokShaposhnikov2017} and in a few earlier references cited there,
$\zeta_r$ is defined to be $\underline{\zeta}_r$ from~\eqref{Eq:Def_zeta_underbar_r}.
We decided to distinguish here the two variants notationally.

For simplicity we have here not defined $\zeta_r$ or $\underline{\zeta}_r$ 
also for $r\in\mathopen]0,\infty\mathclose[\setminus\N$, 
or even more generally as exemplified by \citet[p.~515, the definition (3)]{Tyurin2012}.

Comparing $\zeta_r$ to $\underline{\zeta}_r\,$, we note that 
the important Lemma~\ref{Lem:Regular_eqnorms_on_cM} applies to $\zeta_r$
but not to $\underline{\zeta}_r$. For example, the so-called 
regularity~\eqref{Eq:regularity_on_cM} holds for $\|\cdot\|\coloneqq \zeta_r$ on $\cM$,
whereas the analogue~\eqref{Eq:zeta_underbar_regularity}
for~$\underline{\zeta}_r$ is more complicated.
However, $\underline{\zeta}_r$ yields finite values
on $\cM_r$ and not merely on $\cM_{r,r-1}$, as stated
in~(\ref{Eq:zeta_vs_zeta_underbar_on_cM_and_oncM_r,r-1},\ref{Eq:zeta_vs_zeta_underbar_on_cM_and_on_cM_r_setminus_cM_r,r-1}),
and this is useful even if one is just interested in $\zeta_r$ on $\cM_{r,r-1}$.
For example, the asymptotic relation~\eqref{Eq:zeta_k_P-P_standardised_rounded_asymp}
is for $k\in\{2,3\}$ by~\eqref{Eq:zeta_vs_zeta_underbar_on_cM_and_oncM_r,r-1} a result about $\zeta_k$, 
but in its proof occurs $\underline{\zeta}_k(M)$ 
for a certain $M\coloneqq P_\mathrm{rd}-P\in\cM_k$ not necessarily belonging to $\cM_{k,k-1}$. 

In this paper, we use $\zeta_r$ on $\cM_{r,r-1}$ for $r\in\{1,3\}$
in our main result Theorem~\ref{Thm:BE_K_Z} and in its proof,
$\zeta_0$ in the proof of Zolotarev's Theorem~\ref{Fact:Zolotarev1986},
$\zeta_4$ in Example~\ref{Example:Zolotarev_1973},
and $\underline{\zeta}_r$ in effect only for $r\in\{2,3\}$ 
in Lemma~\ref{Lem:Roundings_and_histograms} 
to prepare for Example~\ref{Example:Discretised_normal_laws}.

That for $2\le r\in \N$ and on $\cM_{r,r-1}$ the norm $\zeta_r$ is 
strictly weaker than $\varkappa_r$,
and not merely weaker as stated in~\eqref{Eq:zeta_vs_varkappa_etc}, 
is in case of $r=3$ shown by the symmetric binomial central limit theorem 
convergence rates
\la                                \label{Eq:symmetric_Bernoulli_CLT_zeta_3_varkappa_3} 
 \zeta_3\big(\widetilde{\mathrm{B}_{\frac{1}{2}}^{\ast n}}-\mathrm{N}\big)\asymp\frac{1}{n}
 &\text{ but }&
 \varkappa_3\big(\widetilde{\mathrm{B}_{\frac{1}{2}}^{\ast n}}-\mathrm{N} \big)\asymp\frac{1}{\sqrt{n}}
 \quad\text{ for }n\in\N\,,
\al
of which the first holds for example by \citet[p.~500, Theorem~1.10]{MattnerShevtsova},
and the second by Remark~\ref{Rem:CLT_rate_varkappa_r} below.

\subsection{Known solutions to Zolotarev's problem with $\zeta$ distances}
                                    \label{subsec:Known_solutions_with_zeta}

Coming now back to providing solutions to Problem~\ref{Zolotarev-Problem}, 
we observe that it follows from \citet{Zolotarev1973} that \eqref{Eq:Problem_with_d_n} holds with each of
\la                                                                  
 d_n(\widetilde{P},\mathrm{N}) &\coloneqq&            \label{Eq:Zolotaref1973-zeta_1_3-distance}
  \left( (\zeta_1\!\vee\!\zeta_3)^{\frac{n}{n+1}}\!\vee\!\zeta_3\right)(\widetilde{P}-\mathrm{N}) \,, \\
 d_n(\widetilde{P},\mathrm{N}) &\coloneqq&             \label{Eq:Zolotaref1973-zeta_3-distance}
  \left( \zeta_3^{\frac{n}{3n+1}}\!\vee\!\zeta_3\right)(\widetilde{P}-\mathrm{N}) \,.
\al
More precisely, on the one hand we are not aware of any explicit statement of \eqref{Eq:Problem_with_d_n} 
with either~\eqref{Eq:Zolotaref1973-zeta_1_3-distance} or~\eqref{Eq:Zolotaref1973-zeta_3-distance}
in the previous literature up to now, but on the other hand,
given the defining representation of $\zeta_1$ and $\zeta_3$ from \eqref{Eq:Def_zeta_r},
the present claims obviously follow 
using \citet[p.~540, the proof of Lemma~2]{Zolotarev1973}.
Results similar to \eqref{Eq:Problem_with_d_n} with~\eqref{Eq:Zolotaref1973-zeta_3-distance}
are a weaker one of~\citet[Theorem 1 for $\alpha=2$ and $r=3$, 
with $\rho_{[r]}=0$ requiring $\mu_3(\widetilde{P})=0$ apparently accidentally]{Christoph1979},
and an incomparable one of \citet[p.~64, the inequality involving $\tau_3$]{Paditz1988},
with both authors using the third expression $ \int \big| F_{M,r}\big| \,\dd\leb$
in~\eqref{Eq:Four_zeta_expressions} for $\zeta_3$,
without mentioning \citeposs{Zolotarev1976} definition~\eqref{Eq:Def_zeta_r}.
See also \citet[p.~65, Theorem 3.11]{ChristophWolf1992}
for further related references.

We should also mention here the inequality
\la                              \label{Eq:zeta_1_vs_zeta_3}
 \zeta_1(\widetilde{P}-\mathrm{N}) &\le& 3\!\cdot\!2^{\frac{1}{3}}\big(\zeta_3(\widetilde{P}-\mathrm{N}) \big)_{}^\frac{1}{3}
   \quad\text{ for }P\in \cP_2\,,
\al
a special case of \citet[p.~29, \foreignlanguage{russian}{Teorema}~3 with $n=r=1$, $s=2$; 
in the English version p.~2227, Theorem 3]{Zolotarev1979}, 
which allows to upper bound R.H.S.\eqref{Eq:Zolotaref1973-zeta_1_3-distance}
i.c.f.~by $\left( \zeta_3^{\frac{n}{3n+3}}\!\vee\!\zeta_3\right)\,(\widetilde{P}-\mathrm{N})$, 
which is however a bit worse than R.H.S.\eqref{Eq:Zolotaref1973-zeta_3-distance}.

While (\ref{Eq:Zolotaref1973-zeta_1_3-distance},\ref{Eq:Zolotaref1973-zeta_3-distance})
improve on (\ref{Eq:Ulyanov1976-distance},\ref{Eq:Ulyanov1978-distance})
by weakening $\varkappa_3$ to $\zeta_3$, their exponents $\frac{n}{n+1}$ and $\frac{n}{3n+1}$
are worse than $1-2^{-n}$ and $\frac12(1-2^{-n})$.
Succeeding in replacing $\frac{n}{n+1}$ in~\eqref{Eq:Zolotaref1973-zeta_1_3-distance} by 
the ideal exponent~$1$, 
but at the cost of introducing the Kolmogorov norm in addition, 
\citet[Theorem 1, for dimension $k=1$ and $g(u)=u$]{Senatov1980} proved that~\eqref{Eq:Problem_with_norm}
holds with
\la                                    \label{Eq:Senatov_1980_with_K_on_RHS}
   \|\cdot\| &\coloneqq& \zeta_1\!\vee\!\zeta_3\vee \left\|\,\cdot\,\right\|_{\mathrm{K}}\, 
\al
and he improved this in \citet[p.~161, Theorem~4.3.1]{Senatov1998} to: For every $\gamma \in\mathopen]0,\infty\mathclose[$
there is 
a $c_\gamma\in \mathopen]0,\infty\mathclose[$ with
\la                                             \label{Eq:Senatov_1998_with_K_on_RHS}
  \left\|\widetilde{P^{\ast n}} -\mathrm{N}\right\|_{\mathrm{K}}   
     &\le& c_\gamma \left( \frac{\zeta_1\!\vee\!\zeta_3}{\sqrt{n}} + \frac{\left\|\,\cdot\,\right\|_{\mathrm{K}}}{n^\gamma} \right)
            \!(\widetilde{P}-\mathrm{N})    \quad\text{ for $P\in\cP_3$ and $n\in\N$}\,. 
\al
As the term with $\|\cdot\|^{}_{\mathrm{K}}$ can not be omitted in~\eqref{Eq:Senatov_1998_with_K_on_RHS}
in case of $n=1$, by~\eqref{Eq:zeta_vs_varkappa_etc} and the optimality of the exponent
$1\wedge\frac{1}{1+1}$ in~\eqref{Eq:Shiganov_combined}, 
or directly by  (\ref{Eq:||_||_K_Zolotarev-example},\ref{Eq:Zolotarev_example_r-norm_errors_n=1})
in Example~\ref{Example:Zolotarev_1973},
\citet[p.~174]{Senatov1998} asked whether it nevertheless can be so for $n\ge n_0$ 
with some $n_0\ge2$.

To sum up: Of the solutions to Zolotarev's 
Problem~\ref{Zolotarev-Problem} reviewed above, and assuming here $n\ge4$ for simplicity,
the six i.c.f.~jointly best ones are~\eqref{Eq:Problem_with_d_n} 
with any of~(\ref{Eq:Shiganov1987-distance_3},\ref{Eq:Ulyanov1976-distance},%
\ref{Eq:Zolotaref1973-zeta_1_3-distance},\ref{Eq:Zolotaref1973-zeta_3-distance}),
\eqref{Eq:Senatov_1998_with_K_on_RHS},
and, if we do not insist on $d_n$ in~\eqref{Eq:Problem_with_d_n} 
being decreasing in~$n$, also~\eqref{Eq:simple_Kolomogorov-telescoping},
namely with $d_n(\widetilde{P},\mathrm{N})\coloneqq 
n^{\frac{3}{2}}\|\widetilde{P}-\mathrm {N}\|^{}_\mathrm{K}$\,. 
More precisely, each of the four solutions~\eqref{Eq:Problem_with_d_n}  
with any of~(\ref{Eq:Trivial_distance},\ref{Eq:Paulauskas1969-distance},%
\ref{Eq:Ulyanov1978-distance}), and~\eqref{Eq:Senatov_1980_with_K_on_RHS},
is i.c.f.~strictly worse than one of the six indicated solutions, 
and it seems to us - admittedly without having checked it in detail - 
that none of the latter be worse than any of the remaining five.
\label{page:It_seems_to_us}

Theorem~\ref{Thm:BE_K_Z} below  answers affirmatively Senatov's question 
mentioned a few lines  above, and reduces the list of jointly best solutions 
to Problem~\ref{Zolotarev-Problem},
among the ones considered here and in case of $n\ge 4$, 
to the following three:
\eqref{Eq:Berry-Esseen_for_Z-close_to_normal},
\eqref{Eq:Problem_with_d_n} with~\eqref{Eq:Shiganov1987-distance_3},
and~\eqref{Eq:simple_Kolomogorov-telescoping}. 

\subsection{An improved solution, Theorem~\ref{Thm:Berry-Esseen_for_Z-close_to_normal}, 
to Problem \ref{Zolotarev-Problem}}

We recall the meaning of ``i.c.f.'' from~\eqref{Eq:i.c.f.}, and 
some notation introduced in subsection~\ref{subsec:Some_notation}.
So $\|\cdot\|^{}_{\mathrm{K}}$ is the Kolmogorov norm from~\eqref{Eq:Def_Kolmogorov-norm},
$\widetilde{P}$ denotes the standardisation of a law $P$,
below assumed to be non-Dirac and with a finite third moment by~\eqref{Eq:Def_cP_r},
$\ast$ indicates convolution,
$\mathrm{N}$~is the standard normal law, 
and $\zeta_r$ is defined by~(\ref{Eq:Def_cF_r},\ref{Eq:Def_cF_r^infty},\ref{Eq:Def_zeta_r}),
with alternative representations provided 
by~\eqref{Eq:Four_zeta_expressions} and (\ref{Eq:zeta_1=kappa_1},\ref{Eq:varkappa_1=L^1-enorm}),
and $\big(\zeta_1\!\vee\!\zeta_3\big)(M) \coloneqq \zeta_1(M)\!\vee\!\zeta_3(M)
\coloneqq \max\{ \zeta_1(M) ,\zeta_3(M)\}$.

\begin{Thm}[Berry-Esseen for summands Zolotarev-close to normal]
                               \label{Thm:Berry-Esseen_for_Z-close_to_normal} \label{Thm:BE_K_Z} 
                                                    \label{page:Thm_BE_K_Z}
There exists a constant $c\in\mathopen]0,\infty\mathclose[$ satisfying 
 \la                                                \label{Eq:Berry-Esseen_for_Z-close_to_normal}
  \left\|\widetilde{P^{\ast n}} -\mathrm{N}\right\|_{\mathrm{K}}
    &\le& \frac{c}{\sqrt{n}}\, \big(\zeta_1\!\vee\!\zeta_3\big)(\widetilde{P}-\mathrm{N})
     \quad\text{ for $P\in\cP_3$ and $n\ge 2$\,.}
 \al
One may take here $c=7.2$\;.
\end{Thm}

I.c.f.~this strictly improves the classical Berry-Esseen Theorem~\eqref{Eq:Berry-Esseen_inequality},
say in view of (\ref{Eq:zeta_1_distance_to_normal_bounded_by_nu_3},%
\ref{Eq:zeta_3_distance_to_normal_bounded_by_nu_3})
and Example~\ref{Example:Discretised_normal_laws}, as explained in more detail below.

Theorem~\ref{Thm:BE_K_Z} is proved in section~\ref{sec:3}, using the main technical result
of this paper, Theorem~\ref{Thm:F_star_G_vs_H_star_H}, in combination with 
Zolotarev's Theorem~\ref{Thm:Zolotarev's_zeta_1-B-E-Thm}
and, to obtain the stated value of the constant~$c$,
the Berry-Esseen Theorem~\eqref{Eq:Berry-Esseen_inequality} with Shevtsova's constant
from~\eqref{Eq:Def_cSH}. That proof actually yields the following result, which is
a bit more complicated, i.c.f.~equivalent, but numerically sharper.

\begin{Rem}                 \label{Rem:BE_K_Z_sharper}
Let $P\in\cP_3$ and $n\ge 2$. Then we have
\la                                \label{Eq:Berry-Esseen_for_Z-close_to_normal_sharper}
  \left\|\widetilde{P^{\ast n}} -\mathrm{N}\right\|_{\mathrm{K}}
    &\le& \frac{1}{\sqrt{n}}
      \min \left\{ c_1\frac{\zeta_1+\alpha\zeta_3}{(1-\lambda\zeta_3)^{}_+}(\widetilde{P}-\mathrm{N})
       \,,\, \cSh\nu_3(\widetilde{P}) \right\}
\al
with the constants
$c_1=1.1708\ldots$,
$\alpha= 0.9678\ldots$,
$\lambda=3.9447\ldots$,
$\cSh=0.469$
defined in {\rm(\ref{Eq:c_1_proof_main_theorem},\ref{Eq:Def_alpha_beta_gamma},%
\ref{Eq:xi_near_zero},\ref{Eq:Def_cSH})},
and with the convention $\frac{x}{0}\coloneqq\infty$ if $x>0$,
and $\frac{0}{0}\coloneqq0$.
We have $\text{\rm R.H.S.\eqref{Eq:Berry-Esseen_for_Z-close_to_normal_sharper}}
\le \text{\rm R.H.S.\eqref{Eq:Berry-Esseen_for_Z-close_to_normal}}$ with
$c=7.16553\ldots$ from~\eqref{Eq:c_approx_7.2_in_main_thm},
with strict inequality unless
$\zeta_1(\widetilde{P}-\mathrm{N}) =\zeta_3(\widetilde{P}-\mathrm{N})
 = \frac{1}{6}(\nu_3(\widetilde{P})-\nu_3(\mathrm{N}))= \zeta^\ast = 0.171989\ldots$
 from~\eqref{Eq:zeta^ast}, or $\widetilde{P}=\mathrm{N}$.

Let $\omega \coloneqq \frac{c_1}{\cSh} = 2.49647\ldots$,
hence $\frac{1}{\omega}= 0.400565\ldots$.
Then $\text{\rm R.H.S.\eqref{Eq:Berry-Esseen_for_Z-close_to_normal_sharper}}$
is attained at the first minimand
iff  $\zeta_3(\widetilde{P}-\mathrm{N})\le
\frac{\nu_3(\widetilde{P})\,-\,\omega\zeta_1(\widetilde{P}-\mathrm{N})}
{\lambda\nu_3(\widetilde{P})\,+\,\omega \alpha}$.
\end{Rem}

We finish this already long first section of the present paper by addressing 
the sharpness of R.H.S.\eqref{Eq:Berry-Esseen_for_Z-close_to_normal}, 
the computability of $\zeta_3(\widetilde{P}-\mathrm{N})$, 
the question of lower bounds for L.H.S.\eqref{Eq:Berry-Esseen_for_Z-close_to_normal}, 
and the possibility of improving~\eqref{Eq:Berry-Esseen_for_Z-close_to_normal} by
generalisation or by increasing L.H.S.\eqref{Eq:Berry-Esseen_for_Z-close_to_normal}.
Let us start with the Examples~\ref{Example:Discretised_normal_laws}
and~\ref{Examples:zeta_3=mu_3}, which show in particular that either 
of the two terms $\zeta_1(\widetilde{P}-\mathrm{N})$
and $\zeta_3(\widetilde{P}-\mathrm{N})$
in~\eqref{Eq:Berry-Esseen_for_Z-close_to_normal}
may dominate the other in interesting cases with both of them small.

Example~\ref{Example:Discretised_normal_laws} also shows that~\eqref{Eq:Berry-Esseen_for_Z-close_to_normal},
even if considered only for $\big(\zeta_1\!\vee\!\zeta_3\big)(\widetilde{P}-\mathrm{N})$
arbitrarily small and $n$ arbitrarily large, necessitates $c\ge \frac{2}{\sqrt{2\pi}} = 0.797884\ldots\;$,
since~\eqref{Eq:Discretised_normal_zeta_asymptotics}
and~\eqref{Eq:zeta_1_vee_zeta_3_over_Esseen1956_discretised_normal} yield
\la                                           \label{Eq:Asymp_quality_BE_K_Z_discretised_normal}
 \lim_{\epsilon\downarrow0}\, \sup\left\{
      \frac{ \text{\rm R.H.S.\eqref{Eq:Esseen1956_asymptotics}} }{ \big(\zeta_1\!\vee\!\zeta_3\big)(\widetilde{P}-\mathrm{N})}
    \,:\, P\in\cP_3\!\setminus\!\{\mathrm{N}\} \,,\,
     \big(\zeta_1\!\vee\!\zeta_3\big)(\widetilde{P}-\mathrm{N}) <\epsilon
    \right\}  &\ge & \frac{2}{\sqrt{2\pi}} \,.
\al
However, \eqref{Eq:Berry-Esseen_for_Z-close_to_normal} just for $n=2$ and still $\big(\zeta_1\!\vee\!\zeta_3\big)(\widetilde{P}-\mathrm{N})$
arbitrarily small even necessitates
$c\ge \text{R.H.S.}\eqref{Eq:n=2_sharpness_BE_K_Z_Zolotarev_example}= 1.1020\ldots$
from Example~\ref{Example:Zolotarev_1973},
and this is the best lower bounding for~$c$ in~\eqref{Eq:Berry-Esseen_for_Z-close_to_normal}
presently known to us.

\begin{Example}[Discretised normal laws]  \label{Example:Discretised_normal_laws}
For $\mu\in\R$,  $\sigma,\eta\in\mathopen]0,\infty\mathclose[\,$, and
$\alpha\in\mathopen]0,1\mathclose[\,$, let 
\la                          \label{Eq:_Def_discretised_normal}
  P &\coloneqq\,\ 
  P_{\mu,\sigma,\eta,\alpha}
    \,\ \coloneqq \,\ \sum_{j\in\Z} \mathrm{N}_{\mu,\sigma^2}
       \big(\,\mathopen](\alpha+j-\tfrac{1}{2})\eta\,,(\alpha+j+\tfrac{1}{2})\eta\mathclose]\big) 
              \delta_{(\alpha+j)\eta}\,.
\al
For $\eta\rightarrow0$ with $\mu,\sigma,\alpha$ fixed we then have,
recalling the lattice span notation~\eqref{Eq:Def_lattice_span_h(P)},
\la
 && \zeta_1(\widetilde{P}-\mathrm{N}) \,\ \sim\,\ \frac{\eta}{4\sigma}\,,\qquad
    \zeta_3(\widetilde{P}-\mathrm{N}) \,\ \llcurly \,\ \eta\,,   \label{Eq:Discretised_normal_zeta_asymptotics} \\
 && h(\widetilde{P}) \,\ \sim\,\ \frac{\eta}{\sigma}\,, \qquad   \nonumber
    \mu_3(\widetilde{P}) \,\ \llcurly \,\ \eta 
\al
and hence
\la                                            \label{Eq:zeta_1_vee_zeta_3_over_Esseen1956_discretised_normal}
 \frac{ \text{\rm R.H.S.\eqref{Eq:Esseen1956_asymptotics}} }
    {\big(\zeta_1\!\vee\!\zeta_3\big)(\widetilde{P}-\mathrm{N})}
  &\sim & \frac{ \frac{h(\widetilde{P})}{2\sqrt{2\pi}}}{\zeta_1(\widetilde{P}-\mathrm{N})}
  \,\ \rightarrow\,\ \frac{2}{\sqrt{2\pi}}  \,,
\al
by using {\rm(\ref{Eq:zeta_1_P-P_standardised_rounded_asymp},\ref{Eq:zeta_k_P-P_standardised_rounded_asymp})}
from Lemma~\ref{Lem:Lattice_and_histogram_approximations},
with the present $(\mathrm{N}_{\mu,\sigma^2} , P_{\mu,\sigma,\eta,\alpha})$ 
in the role of $(P,P_\mathrm{rd})$ there, and also
$|\mu_3(\widetilde{P})| = |\mu_3(\widetilde{P}-\mathrm{N})| \le 6\,\zeta_3(\widetilde{P}-\mathrm{N})$
due to~\eqref{Eq:zeta_vs_varkappa_etc}.

In particular for $\mu=\alpha=0$, symmetry of $P$ yields $\mu_3(\widetilde{P})=0$
and hence $\text{\rm R.H.S.\eqref{Eq:Esseen1956_asymptotics}}= \frac{h(\widetilde{P})}{2\sqrt{2\pi}}\,$,
and taking then also $\sigma=1$ yields numerically, for example,

\vspace{\baselineskip}
\centerline{
$   \renewcommand{\arraystretch}{2}
 \begin{array}{|l|l|l|l|l|l|l|} \hline
     & \zeta_1(\widetilde{P}-\mathrm{N}) & \zeta_3(\widetilde{P}-\mathrm{N})
          & \text{\rm R.H.S.\eqref{Eq:Esseen1956_asymptotics}}
              & c^{}_\mathrm{E}\nu_3(\widetilde{P})
                 & 7.2\,\big(\zeta_1\!\vee\!\zeta_3\big)(\widetilde{P}-\mathrm{N})
                 & \sqrt{n}\,\mathrm{R.H.S.\eqref{Eq:Berry-Esseen_for_Z-close_to_normal_sharper}} \\ \hline
  \eta =1             & 0.2417\ldots   & 0.0051\ldots & 0.1916\ldots& 0.6562\ldots  & 1.740\ldots & 0.294\ldots\\ \hline
  \eta =\frac{1}{10}  & 0.0249\ldots  & < 10^{-5} & 0.01993\ldots& 0.6538\ldots  & 0.179\ldots  & 0.0292\ldots\\ \hline
  \eta =\frac{1}{100} & 0.00249\ldots & < 10^{-5}& 0.001994\ldots & 0.6538\ldots  & 0.0179\ldots& 0.00292\ldots\\ \hline
  \end{array}
$}

\vspace{\baselineskip}\noindent
so that here, for $n\ge2$,
already for $\eta =\frac{1}{10}$
the error bound~\eqref{Eq:Berry-Esseen_for_Z-close_to_normal}
with $c=7.2$ is better than the
Berry-Esseen theorem~\eqref{Eq:Berry-Esseen_inequality} with the
hypothetical constant $c^{}_\mathrm{E}$ from~\eqref{Eq:Def_c_E},
and~\eqref{Eq:Berry-Esseen_for_Z-close_to_normal_sharper}
is better even for $\eta=1$.
\end{Example}

To actually compute $\big(\zeta_1\!\vee\!\zeta_3\big)(\widetilde{P}-\mathrm{N})$ 
for a given $P\in\cP_3$ by straightforward integrations, 
we have to compute $\mu(P),\sigma(P)$ and then, if $F_P$ is at hand,
need one integration for 
$\zeta_1(\widetilde{P}-\mathrm{N})=\int | F_{\widetilde{P}-\mathrm{N}}|\,\dd\leb$ 
using~(\ref{Eq:zeta_1=kappa_1},\ref{Eq:varkappa_1=L^1-enorm}), and two further ones 
for $\zeta_3(\widetilde{P}-\mathrm{N})=\int | F_{\widetilde{P}-\mathrm{N}\,,\,3}|\,\dd\leb$
using~(\ref{Eq:Def_F_M,k},\ref{Eq:Four_zeta_expressions}).
The latter two integrations simplify if the following known Lemma~\ref{Lem:zeta_3_computation}
is applicable, as in Examples~\ref{Examples:zeta_3=mu_3} below.
We let here $S^{-}(h)$ denote the number of sign changes of a function $h:\R\rightarrow\R$,
as defined more precisely in~\eqref{Eq:Def:S^-(f)}, 
immediately after a definition of initial positivity or negativity.

\begin{Lem}[Sufficient conditions for $\zeta_3(\widetilde{P}-\mathrm{N})
  = \frac{1}{6}|\mu_3(\widetilde{P})|\,$]                   \label{Lem:zeta_3_computation}
Let $P\in\cP_3$ with distribution function $F$.  \label{page:Lemmma_zeta_3_computation}

\begin{parts}
\item Let $S^-(\widetilde{F}-\Phi)\le2$.     \label{part:Phi-F_tilde_2_sign_changes} 
 Then $\zeta_3(\widetilde{P}-\mathrm{N}) = \frac{1}{6}|\mu_3(\widetilde{P})|$,
  $S^-(\widetilde{F}-\Phi)=2$ unless $\widetilde{P}=\mathrm{N}$, and
  \[
   \mu_3(\widetilde{P})\begin{Bmatrix}\ge \\ \le \end{Bmatrix} 0
   &\iff& \widetilde{F}-\Phi\text{ initially }\begin{Bmatrix}\text{negative}\\ \text{positive} \end{Bmatrix}
   \,\ \iff \,\ \widetilde{F}-\Phi\text{ finally }\begin{Bmatrix}\text{negative}\\ \text{positive} \end{Bmatrix}.
  \]
\item Let $\widetilde{P}=\widetilde{f}\leb$   \label{part:phi-f_tilde_3_sign_changes}  
 for some $\leb$-density $\widetilde{f}$ 
 with $S^-(\widetilde{f}-\phi)\le 3$. Then the assumption of part~{\rm\ref{part:Phi-F_tilde_2_sign_changes}} 
 is fulfilled, $S^-(\widetilde{f}-\phi)=3$ unless $\widetilde{P}=\mathrm{N}$, and
 \[ \quad
   \mu_3(\widetilde{P}) \begin{Bmatrix}\ge \\ \le \end{Bmatrix} 0
   \,\iff \, \widetilde{f}-\phi\text{ 
                                      initially }
     \begin{Bmatrix}\text{negative}\\ \text{positive} \end{Bmatrix}  
   \, \iff \, \widetilde{f}-\phi\text{ 
                                       finally }\begin{Bmatrix}\text{positive}\\ \text{negative} \end{Bmatrix}.
  \]
\end{parts}
\end{Lem}
\begin{proof}
 Theorem~\ref{Thm:Cut_criteria} with $M\coloneqq\pm(\widetilde{P}-\mathrm{N})$, $r\coloneqq 3$,
 the implications $(B_0)\Rightarrow(B_1)\Rightarrow \zeta_3(M)=\frac{1}{3!}\mu_3(M)$,
 here $\mu_3(M) = \pm\mu_3(\widetilde{P})$, 
  and~\eqref{Eq:S^-_F_k_vs_F_k-1} for $k\in\{0,1\}$.
\end{proof}

\begin{Examples}                 \label{Examples:zeta_3=mu_3} 
In each of the following three parts we have:
The assumption of Lemma~{\rm\ref{Lem:zeta_3_computation}\ref{part:Phi-F_tilde_2_sign_changes}}
is fulfilled, with $\widetilde{P}\neq \mathrm{N}$ and $\widetilde{F}-\Phi$ initially negative.
With the exception of the present part~{\rm\ref{part:Left-winsorised_normal_laws}},
even the assumption of Lemma~{\rm\ref{Lem:zeta_3_computation}\ref{part:phi-f_tilde_3_sign_changes}}
is fulfilled, with $\widetilde{f}-\phi$ initially negative. Hence 
\la                                                \label{Eq:zeta_3=mu_3_in_Examples}
 \zeta_3(\widetilde{P}-\mathrm{N}) &=& \tfrac{1}{6}\mu_3(\widetilde{P}) \,\ >\,\ 0\,.
\al
Further, here 
$\text{\rm R.H.S.\eqref{Eq:Esseen1956_asymptotics}}=\frac{1}{6\sqrt{2\pi}}\mu_3(\widetilde{P})$
due to $h(\widetilde{P})=0$, and, 
at least under the parameter restrictions as indicated in each part below, we have 
\la                                                 \label{Eq:zeta_1_le_zeta_3_in_Examples_zeta_3=mu_3}
 \zeta_1(\widetilde{P}-\mathrm{N}) &\le& \zeta_3(\widetilde{P}-\mathrm{N})
\al 
and hence then
\la                   \label{Eq:zeta_1_vee_zeta_3_over_Esseen1956_truncated_normal_etc}
  \frac{ \text{\rm R.H.S.\eqref{Eq:Esseen1956_asymptotics}} }
       { \big(\zeta_1\!\vee\!\zeta_3\big)(\widetilde{P}-\mathrm{N}) }  
  &=& \frac{1}{\sqrt{2\pi}}   \,.     
 \al 
\begin{parts}
\item {\rm Left-truncated normal laws.}   \label{part:Left-truncated_normal_laws} 
 For $t\in\R$ let $P\coloneqq P_t\coloneqq \mathrm{N}\bed{\cdot}{]-t,\infty[}$, that is, 
 $P=f\leb$ with $f\coloneqq f_t\coloneqq \1_{]-t,\infty[}\frac{\phi}{\mathrm{N}(]-t,\infty[)}$,
 and $\widetilde{P} = \widetilde{f}\leb$ with
 $\widetilde{f}(x)\coloneqq\widetilde{f_t}(x)\coloneqq \sigma(P)f_t\big(\sigma(P)x+\mu(P)\big)$ for $x\in\R$.
 
 Here we have, with asymptotics referring to $t\rightarrow\infty$,
 \la                                       \label{Eq:zeta_3_1_truncated_normal_standardised}
   && \zeta_3(\widetilde{P}-\mathrm{N})       
    \, \sim\, \tfrac{1}{6} t^2\phi(t)\,,  \qquad 
  \zeta_1(\widetilde{P}-\mathrm{N})        
    \, \sim\, \tfrac{1}{\sqrt{2\pi}}t\phi(t)
    \, \llcurly \, \zeta_3(\widetilde{P}-\mathrm{N})\,,
 \al
 and hence~\eqref{Eq:zeta_1_vee_zeta_3_over_Esseen1956_truncated_normal_etc}  holds for $t$ sufficiently large. 
 The second asymptotic equality in~\eqref{Eq:zeta_3_1_truncated_normal_standardised} 
 follows~from      
 \la            \label{Eq:zeta_1_three_asymptotics_truncated_normal}
  && \zeta_1(P-\mathrm{N}) \, \sim\,  \phi(t)\,,  \qquad
     \zeta_1(\bdot{P}-P) \, \sim\, \phi(t)\,,  \qquad
     \zeta_1(\widetilde{P} - \bdot{P}) \, \sim\, \tfrac{1}{\sqrt{2\pi}}t\phi(t)\,,
 \al
 and this also shows that here we have, perhaps surprisingly,
 $\zeta_1(P-\mathrm{N}) \llcurly \zeta_1(\widetilde{P}-\mathrm{N})$.  
 
\item {\rm Left-winsorised normal laws.}   \label{part:Left-winsorised_normal_laws}  
 For $t\in\R$, let 
 $ P\coloneqq P_t \coloneqq \Phi(-t)\delta_{-t}
  +\mathrm{N}(\cdot\cap\mathopen]-t,\infty\mathclose[)$.
 Here we have, with asymptotics referring to $t\rightarrow\infty$, 
 \la            \label{Eq:zeta_3_1_truncated_normal_winsorised}
  && \zeta_3(\widetilde{P}-\mathrm{N}) \,\sim\, \tfrac{1}{2}\phi(t)\,, \qquad
     \zeta_1(\widetilde{P}-\mathrm{N}) \,\sim\, \tfrac{2}{\sqrt{2\pi}} \frac{\phi(t)}{t}
        \, \llcurly\, \zeta_3(\widetilde{P}-\mathrm{N})\,,
 \al
 and hence~\eqref{Eq:zeta_1_vee_zeta_3_over_Esseen1956_truncated_normal_etc}
 holds for $t$ sufficiently large.
 
\item {\rm Gamma laws and some of their power transforms.}  \label{part:Gamma_laws_and_power_transforms}
 For $\alpha,\lambda\in\mathopen]0,\infty\mathclose[$ and $\beta\in\R\setminus\{0\}$, 
 let $P\coloneqq\Gamma_{\alpha,\lambda,\beta}$ be the law on $\R$ with the 
 $\leb$-density given by  
 \[
  f(x) &\coloneqq&f_{\Gamma_{\alpha,\lambda,\beta}}(x) \,\ \coloneqq\,\ \frac{\lambda^\alpha|\beta|}{\Gamma(\alpha)} 
    x^{\alpha\beta-1}\exp(-\lambda x^\beta)\cdot(x>0)\quad \text{ for } x\in\R\,,
 \]
 so that in case of $\beta=1$ we have a usual gamma law 
 $\Gamma_{\alpha,\lambda}\coloneqq \Gamma_{\alpha,\lambda,1}\,$, 
 and in general a power transformed gamma law
 $\Gamma_{\alpha,\lambda,\beta} = (x\mapsto x^\frac{1}{\beta})\im \Gamma_{\alpha,\lambda}\,$, 
 and let here the parameter pair $(\alpha,\beta)$ be restricted by 
 \la               \label{Eq:restriction_beta_power_transformed_gamma}
   \beta &\in& \mathopen]-\infty, - \tfrac{3}{\alpha}  \mathclose[  \ \cup\ \mathopen]0,2\mathclose] \,.
 \al
 
 Here, for $(\alpha,\beta)$ unrestricted,
 the condition $\nu_r(P)<\infty$ is for $r\in\mathopen]0,\infty\mathclose[$
 equivalent to  $\beta >0$ or $\beta < -\frac{r}{\alpha}$, 
 which may be rewritten as $\alpha + \frac{r}{\beta}>0$,
 and under this condition 
 \la                                    \label{Eq:nu_r_power_transformed_gamma}
  && \nu_r(P) \,=\, \lambda^{-\frac{r}{\beta}}G(\tfrac{r}{\beta},\alpha)
   \quad \text{ with } \quad 
    G(a,x) \, \coloneqq\,  \frac{\Gamma(x+a)}{\Gamma(x)}  \,,
 \al
 and, in the presence of the assumption $\nu_3(P)<\infty$,
 the further 
 condition $\beta\le2$ in~\eqref{Eq:restriction_beta_power_transformed_gamma}
 is equivalent to the assumption $S^-(\widetilde{f}-\phi)\le 3$ 
 in Lemma~{\rm\ref{Lem:zeta_3_computation}\ref{part:phi-f_tilde_3_sign_changes}}.
 
 We have here 
\la
  \zeta_3(\widetilde{P}-\mathrm{N})     \label{Eq:mu_3_Gamma_power_transformed_standardised}
   &=& \frac{1}{6}\,\frac{G(\tfrac{3}{\beta},\alpha)-3\,G(\tfrac{2}{\beta},\alpha)G(\tfrac{1}{\beta},\alpha) 
              +2\,G^{3}_{}(\tfrac{1}{\beta},\alpha)}
      {\left(G(\tfrac{2}{\beta},\alpha) -G^2_{}(\tfrac{1}{\beta},\alpha)\right)^\frac{3}{2}} \\
   & =& \sgn(\beta)\big(\tfrac{1}{\beta}-1+\tfrac{3}{8}\beta-\tfrac{1}{24}\beta^2\big)
             \frac{1}{\sqrt{\alpha}} \,+\, O\big(\alpha^{-\frac{3}{2}} \big)\,, \nonumber
    \\
   \zeta_3(\widetilde{P}-\mathrm{N})                  \label{Eq:mu_3_Gamma_standardised}
     &=& \zeta_3(\widetilde{\Gamma_{\alpha,\lambda}}-\mathrm{N}) 
     \,\ =\,\  \frac{1}{3\sqrt{\alpha}}   \quad \text{ if }\beta=1 \,,
\al
with the $O(\ldots)$-claim in~\eqref{Eq:mu_3_Gamma_power_transformed_standardised}
valid at least if $\beta$ is fixed and $\alpha\ge 1\!\vee\!(-\frac{3}{\beta}+1)$,
and we have~\eqref{Eq:zeta_1_le_zeta_3_in_Examples_zeta_3=mu_3} at least in the gamma case of $\beta=1$ 
 with $\alpha$ sufficiently large, 
 since 
 \la                                   \label{Eq:zeta_1/zeta_3->...<1_for_gamma_laws}
  \lim_{\alpha\rightarrow\infty} \frac{\zeta_1}{\zeta_3}( \widetilde{\Gamma_{\alpha,\lambda,1}}  -\mathrm{N}) 
  &=& \frac{4}{\sqrt{2\pi\mathrm{e}}} \,\ =\,\ 0.967882\ldots \,\ <\,\ 1\,.
 \al 
\end{parts}
\end{Examples}

This is proved in section~\ref{sec:special_laws}, 
starting on page~\pageref{page:Proof_of_Examples:zeta_3=mu_3}.

Laws of sums of i.i.d.~truncated normal random variables occur naturally in certain 
statistical problems, see for example \citet[Chapters 2 and 3]{Cohen1991}
and also \citet[``gestutzte Normalverteilungen'']{Rasch1995},
and have been studied at least since 
\citet{Francis1946} in the one-sided case as in 
Example~\ref{Examples:zeta_3=mu_3}\ref{part:Left-truncated_normal_laws}, 
and \citet{BirnbaumAndrews1949} in the symmetric two-sided case. 

The identity in~\eqref{Eq:zeta_3=mu_3_in_Examples} was previously obtained
for Erlang laws, that is, with $\alpha\in\N$ and $\beta=1$ 
in Example~\ref{Examples:zeta_3=mu_3}\ref{part:Gamma_laws_and_power_transforms},
and conjectured also for Weibull 
laws ($\alpha=1$ in Example~\ref{Examples:zeta_3=mu_3}\ref{part:Gamma_laws_and_power_transforms}),
by \citet[p.~1264, (55), p.~1255, line 4]{Boutsikas2011}.
General power transformed gamma laws are useful for unifying certain computations, 
and are hence introduced under various names for example in 
\citet[pp.~348--353, ``generalized gamma'' or ``gamma-Weibull'', ``extended'' if $\beta<0$]{MarshallOlkin2007}
and, for $\beta>0$ only, in
\citet[pp.~299--301, ``one-sided hyper-exponential'']{Hoffmann-JorgensenI}
and \citet[pp.~655--666, $\gamma_{\alpha,\nu,\lambda}$]{StorchWiebeIII}.

For certain other laws,  
$\zeta_3(\widetilde{P}-\mathrm{N})$ is a simple function 
of the third absolute moment of $\widetilde{P}\,$:

\begin{Example}[\citeposs{Subbotin1923} generalisations of normal, bilaterally exponential, and uniform laws]  
                                                                  \label{Example:Subbotin_laws}
Let $\beta\in\mathopen]0,\infty\mathclose]$,      \label{page:Example_Subbotin_laws}
\la                                              \label{Eq:Def_Subbotin_f_beta} 
  f_{\beta}(x) &\coloneqq& \left\{\begin{array}{ll} 
    \frac{\beta}{2\,\Gamma(\frac{1}{\beta})}\exp(-|x|^\beta) & \text{ if }\beta\in\mathopen]0,\infty\mathclose[, \\
      \frac12 \1_{[-1,1]}(x) & \text{ if }\beta=\infty 
   \end{array}\right\}  \quad\text{ for }x\in\R,
\al
$\alpha\in\mathopen]0,\infty\mathclose[\,$,
$f^{}_{\beta,\alpha}(x)\coloneqq\frac{1}{\alpha}f^{}_\beta(\frac{x}{\alpha})$ for $x\in\R$,
and $P_{\beta,\alpha}\coloneqq f_{\beta,\alpha}\leb$. Then 
\la                                             \label{Eq:zeta_3_Subbotin_-_N} 
 \zeta_3(\widetilde{P_{\beta,\alpha}} -\mathrm{N}) 
   &=& \frac16\left| \nu_3(\widetilde{P_{\beta,\alpha}})  - \nu_3(\mathrm{N}) \right|
   \,\ =\,\ \frac{\sgn(2-\beta)}{6} 
     \left( \frac{\Gamma(\tfrac{4}{\beta})\Gamma(\tfrac{1}{\beta})^\frac{1}{2}}{\Gamma(\tfrac{3}{\beta})^\frac{3}{2}} 
            - \frac{4}{\sqrt{2\pi}} \right)   \\  \nonumber 
  & & \!\!\!\!\!\!\!\!
   \left\{\begin{array}{ll} 
             \uparrow\,\ \infty & \text{ for } 2\ge \beta\downarrow 0\,, \\
             =\,\ \frac{1}{6}\big(\frac{3}{\sqrt{2}}-\frac{4}{\sqrt{2\pi}}\big) \,\ =\,\ 0.0875918\ldots  & \text{ if }\beta=1\,,\\
             =\,\ 0                                                                        & \text{ if }\beta=2\,,\\  
             \uparrow\,\ \frac{1}{6}\big(\frac{4}{\sqrt{2\pi}} - \frac{3\sqrt{3}}{4}\big) \,\ =\,\ 0.0494551\ldots  
              & \text{ for }
                2\le \beta \uparrow \infty \,, 
            \end{array} \right.          
\al
with the $\Gamma$ quotient in the third expression being decreasing 
in $\beta\in \mathopen]0,\infty\mathclose[\,$, and defined at $\beta=\infty$ to be its limit.
\end{Example} 
This is proved in section~\ref{sec:special_laws}, 
starting on page~\pageref{page:Proof_of_Example_Subbotin}.

In Example~\ref{Example:Subbotin_laws}, however,
$\text{R.H.S.}\eqref{Eq:Esseen1956_asymptotics}=0$ due to symmetry and 
absolute continuity, 
and the convergence rate in \eqref{Eq:Lindeberg-Levy-Thm}
is then, using also finiteness of~$\mu_4(P)$, in fact $\frac{1}{n}$ 
by \citet[p.~173, Theorem~5.21 with $k=4$]{Petrov1995}
or \citet[p.~54, Corollary~3.1]{Yaroslavtseva2008b}.

Example~\ref{Example:Discretised_normal_laws} together with either of 
Example~\ref{Examples:zeta_3=mu_3}\ref{part:Left-truncated_normal_laws} 
or~\ref{part:Left-winsorised_normal_laws}
shows that, simultaneously,
the two distances
\mbox{$\zeta_1(\widetilde{P}-\mathrm{N})$} $=$ $\varkappa_1(\widetilde{P}-\mathrm{N})$ 
and $\zeta_3(\widetilde{P}-\mathrm{N})$ occurring in~\eqref{Eq:Berry-Esseen_for_Z-close_to_normal}
may be 
both arbitrarily small,
with either one being arbitrarily large compared to the other,
and such that~\eqref{Eq:Berry-Esseen_for_Z-close_to_normal}
is of asymptotically correct order as made more precise
by~\eqref{Eq:zeta_1_vee_zeta_3_over_Esseen1956_discretised_normal} 
or~\eqref{Eq:zeta_1_vee_zeta_3_over_Esseen1956_truncated_normal_etc}.
In particular, none of the two distances may simply be omitted 
in~\eqref{Eq:Berry-Esseen_for_Z-close_to_normal}.
Since, however, $\zeta_3(\widetilde{P}-\mathrm{N})$ may not always be easy to 
compute, or to bound accurately from above, one may consider the following
simple consequence of Theorem~\ref{Thm:Berry-Esseen_for_Z-close_to_normal}.

\begin{Cor}[Berry-Esseen for summands $\varkappa$-close to normal]
             \label{Cor:Berry-Esseen_for_varkappa-close_to_normal} \label{Cor:BE_K_varkappa} 
There exists a constant $c\in\mathopen]0,\infty\mathclose[$ satisfying 
 \la                   \label{Eq:Berry-Esseen_for_varkappa-close_to_normal}
    \left\|\widetilde{P^{\ast n}} -\mathrm{N}\right\|_{\mathrm{K}}   
     &\le& \frac{c}{\sqrt{n}}\, \big(\varkappa_1\!\vee\!\varkappa_3\big)(\widetilde{P}-\mathrm{N})
   \quad\text{ for $P\in\cP_3$ and $n\ge 2$.}
 \al
One may take here $c=7.2$, or more precisely take $(7.2\,\varkappa_1)\!\vee\!(1.2\,\varkappa_3)$
in place of $c\,\varkappa_1\!\vee\!\varkappa_3$\,.
\end{Cor}
\begin{proof}
Inequality~\eqref{Eq:Berry-Esseen_for_Z-close_to_normal}
combined with (\ref{Eq:zeta_1=kappa_1},\ref{Eq:zeta_vs_varkappa_etc}). 
\end{proof}

Now the example~\eqref{Eq:symmetric_Bernoulli_CLT_zeta_3_varkappa_3} 
shows that the distance $\varkappa_3(\widetilde{P}-\mathrm{N})$
occurring in~\eqref{Eq:Berry-Esseen_for_varkappa-close_to_normal}
may be arbitrarily large compared to the distance $\zeta_3(\widetilde{P}-\mathrm{N})$
occurring in~\eqref{Eq:Berry-Esseen_for_Z-close_to_normal}, 
but this does not yet rule out the possibility 
of~\eqref{Eq:Berry-Esseen_for_varkappa-close_to_normal}
being i.c.f.~equivalent to~\eqref{Eq:Berry-Esseen_for_Z-close_to_normal}. 
That this is in fact not so is shown by the following example, 
which is admittedly a bit artificial compared to Examples~\ref{Example:Discretised_normal_laws},
\ref{Examples:zeta_3=mu_3}\ref{part:Left-truncated_normal_laws},
and~\ref{Examples:zeta_3=mu_3}\ref{part:Left-winsorised_normal_laws}.

\begin{Example}[Tail-discretised normal laws] \label{Example:Tail-discretised_normal_laws}
For $t,\eta\in\mathopen]0,\infty\mathclose[$\,, 
let $I\coloneqq I_t\coloneqq \mathopen]-t,t\mathclose[$ and 
\[
 P &\coloneqq& P_{t,\eta} \,\ \coloneqq \,\ \mathrm{N}( \cdot \cap I)
    + \sum_{j\in\Z} \mathrm{N}
       \big(\,\mathopen](j-\tfrac{1}{2})\eta\,,(j+\tfrac{1}{2})\eta\mathclose] \cap I^\mathrm{c} \big) 
              \delta_{j\eta}\,.
\]
Then $\lim_{\eta\rightarrow0} \big(\zeta_1\!\vee\!\zeta_3\big)(\widetilde{P}-\mathrm{N}) =0$ 
for each $t$, 
and 
\[
 \lim_{t\rightarrow\infty}\limsup_{\eta\rightarrow0}
    \frac {\zeta_1\!\vee\!\zeta_3}{\varkappa_1\!\vee\!\varkappa_3 }(\widetilde{P}-\mathrm{N}) &=& 0\,.
\]
\end{Example}
This is proved in section~\ref{sec:Roundings_and_histograms}, 
starting on page~\pageref{page:Proof_of_Example_Tail-discretised_normal_laws}.

In inequality~\eqref{Eq:Berry-Esseen_for_Z-close_to_normal}
we cannot just omit the assumption ``$n\ge 2$'', 
by~\eqref{Eq:zeta_vs_varkappa_etc} and the optimality of the exponent 
$1\wedge\frac{1}{1+1}$ in~\eqref{Eq:Shiganov_combined},
or  more 
directly by (\ref{Eq:||_||_K_Zolotarev-example},\ref{Eq:Zolotarev_example_r-norm_errors_n=1})
in \citeposs{Zolotarev1972,Zolotarev1973} Example~\ref{Example:Zolotarev_1973}.
But using the simple inequality~\eqref{Eq:Kolmogorov_vs_varkappa_1_near_N} from
Lemma~\ref{Lem:Kolmogorov_bounded_by_varkappa_and_Lipschitz}, we can alternatively state~\eqref{Eq:Berry-Esseen_for_Z-close_to_normal}
in the form 
\la           \label{Eq:Berry-Esseen_for_Z-close_to_normal_also_for_n=1}
    \left\|\widetilde{P^{\ast n}} -\mathrm{N}\right\|_{\mathrm{K}}   
     &\le& \frac{c}{\sqrt{n}}\,\big(\zeta_1^{1\wedge\frac{n}{2}}
     \!\vee\!\zeta^{}_3\big) (\widetilde{P}-\mathrm{N}) 
   \quad\text{ for $P\in\cP_3$ and $n\in \N$}\,,   
\al
still with $c=7.2$\;.

Inequality~\eqref{Eq:Berry-Esseen_for_Z-close_to_normal_also_for_n=1}
i.c.f.~improves the classical Berry-Esseen theorem~\eqref{Eq:Berry-Esseen_inequality},
since it improves earlier improvements of~\eqref{Eq:Berry-Esseen_inequality}
as discussed in the next paragraph. Independently of this argument, 
and now also considering constant factors, let us first note the inequalities
\la
 \zeta_1(\widetilde{P}-\mathrm{N})    \label{Eq:zeta_1_distance_to_normal_bounded_by_nu_3}
   &\le&(1\!+\!\tfrac{2}{\sqrt{2\pi}})\wedge\nu_3(\widetilde{P}) \quad\text{ for }P\in\cP_2\,,\\
 \zeta_3(\widetilde{P}-\mathrm{N})   \label{Eq:zeta_3_distance_to_normal_bounded_by_nu_3}
  &\le&\tfrac{1}{6}\nu_3(\widetilde{P})  \quad\text{ for }P\in\cP_2\,. 
\al
Here the bound with just the first minimand in~\eqref{Eq:zeta_1_distance_to_normal_bounded_by_nu_3}
is due to~(\ref{Eq:zeta_1=kappa_1},\ref{Eq:zeta_1_le}), and the remaining two bounds 
are trivial in case of $P\in\cP_2\setminus\cP_3$, 
and are else the
Goldstein-Tyurin theorem~\eqref{Eq:Goldstein-Tyurin} with $n=1$
in case of~\eqref{Eq:zeta_1_distance_to_normal_bounded_by_nu_3},
and \citet[Theorem 4 with $n=1$]{Tyurin2010}
in case of~\eqref{Eq:zeta_3_distance_to_normal_bounded_by_nu_3}.
We obtain
\la                       \label{Eq:zeta_1_vee_zeta_3_distance_to_normal_bounded_by_nu_3}
 \big(\zeta_1^{1\wedge\frac{n}{2}}
     \!\vee\!\zeta^{}_3\big) (\widetilde{P}-\mathrm{N}) 
 &\le& \nu_3(\widetilde{P})\quad \text{ for }P\in\cP_3
\al
by (\ref{Eq:zeta_1_distance_to_normal_bounded_by_nu_3},%
\ref{Eq:zeta_3_distance_to_normal_bounded_by_nu_3})
in case of $n\ge 2$ or $\zeta_1(\widetilde{P}-\mathrm{N})\ge 1$,
and in the remaining case of $n=1>\zeta_1(\widetilde{P}-\mathrm{N})$
by $\text{L.H.S.}\eqref{Eq:zeta_1_vee_zeta_3_distance_to_normal_bounded_by_nu_3}
\le 1 \!\vee\! \zeta_3(\widetilde{P}-\mathrm{N}) \le \nu_3(\widetilde{P})$
due to~(\ref{Eq:nu_3(widetilde(P))_ge_1},\ref{Eq:zeta_3_distance_to_normal_bounded_by_nu_3}).
Hence~\eqref{Eq:Berry-Esseen_for_Z-close_to_normal_also_for_n=1}
with $c=7.2$ is always better than~\eqref{Eq:Berry-Esseen_inequality} with $c=7.2$\,.

As claimed in the first paragraph of subsection~\ref{subsec:Aim},
inequality~\eqref{Eq:Berry-Esseen_for_Z-close_to_normal_also_for_n=1}
i.c.f.~improves, usually strictly, inequality~\eqref{Eq:Problem_with_d_n}  for
each choice of (\ref{Eq:Shiganov1987-distance_1},%
\ref{Eq:Ulyanov1976-distance},%
\ref{Eq:Zolotaref1973-zeta_1_3-distance},\ref{Eq:Zolotaref1973-zeta_3-distance}),
and also inequality~\eqref{Eq:Senatov_1998_with_K_on_RHS},
except that we have to assume $n\ge2$ 
in case of \eqref{Eq:Zolotaref1973-zeta_3-distance} 
or~\eqref{Eq:Senatov_1998_with_K_on_RHS}. 
More precisely, 
by~(\ref{Eq:zeta_vs_varkappa_etc},\ref{Eq:zeta_1=kappa_1},%
\ref{Eq:zeta_1_le},\ref{Eq:zeta_1_vs_zeta_3})
or trivially, inequality~\eqref{Eq:Berry-Esseen_for_Z-close_to_normal_also_for_n=1}
i.c.f.~improves, perhaps nonstrictly, 
each of the other inequalities under the stated restriction on $n$, 
with in case of~\eqref{Eq:Problem_with_d_n} with~\eqref{Eq:Zolotaref1973-zeta_3-distance}
an intermediate improvement being the following:

\begin{Cor}                  \label{Cor:Berry-Esseen_for_just_zeta_3-close_to_normal}
There exists a constant $c\in\mathopen]0,\infty\mathclose[$ satisfying 
\la               \label{Eq:Berry-Esseen_for_just_zeta_3-close_to_normal}
  \left\|\widetilde{P^{\ast n}} -\mathrm{N}\right\|_{\mathrm{K}}   
    &\le& \frac{c}{\sqrt{n}}\,\big(\zeta_3^{\frac{1}{3}} \!\vee\!\zeta_3\big)(\widetilde{P}-\mathrm{N}) 
      \quad\text{ for $P\in\cP_3$ and $n\ge 2$.}
\al
One may take here $c=28$. Here, for each $n\ge2$,
the exponent $\frac{1}{3}$ is not increasable beyond $\frac{1}{2}$.
\end{Cor}
\begin{proof}
Inequality~\eqref{Eq:Berry-Esseen_for_Z-close_to_normal} with $c=7.2$ combined
with~\eqref{Eq:zeta_1_vs_zeta_3} gives~\eqref{Eq:Berry-Esseen_for_just_zeta_3-close_to_normal}
with $c=27.2142\ldots$\,.

The final claim follows from considering $P\coloneqq \mathrm{B}_{\frac{1}{2}}^{\ast k}$ 
with $k\in\N$ arbitrarily large, since there are constants $c_1,c_2<\infty$ 
with then 
$\text{L.H.S.\eqref{Eq:Berry-Esseen_for_just_zeta_3-close_to_normal}} \ge \frac{c_1}{\sqrt{kn}}$
by~\eqref{Eq:Esseen1956_asymptotics} with the present $(\mathrm{B}_\frac{1}{2},kn)$ in the role
of $(P,n)$ there, and $\zeta_3(\widetilde{P}-\mathrm{N})\le \frac{c_2}{k}$ 
by~\eqref{Eq:symmetric_Bernoulli_CLT_zeta_3_varkappa_3}.
\end{proof}

Corollary~\ref{Cor:Berry-Esseen_for_just_zeta_3-close_to_normal}
also i.c.f.~improves the result of
\citet[p.~64, the inequality involving~$\tau_3$]{Paditz1988}
already mentioned above in connection
with~\eqref{Eq:Problem_with_d_n} with~\eqref{Eq:Zolotaref1973-zeta_3-distance}.

Further, to justify the ``usually strictly'' above, we note:
Inequality~\eqref{Eq:Berry-Esseen_for_Z-close_to_normal_also_for_n=1}
is for each $n\in\N$ i.c.f.~strictly better
than~\eqref{Eq:Problem_with_d_n} with~\eqref{Eq:Shiganov1987-distance_1}
for discrete examples with
$\big(\zeta_1^{}
     \!\vee\!\zeta^{}_3\big) (\widetilde{P}-\mathrm{N})$
arbitrarily small, as in Example~\ref{Example:Discretised_normal_laws}.
And assuming now  $n\ge 2$,
Example~\ref{Example:Tail-discretised_normal_laws} shows that
\eqref{Eq:Berry-Esseen_for_Z-close_to_normal_also_for_n=1}
is i.c.f.~strictly better
than~\eqref{Eq:Berry-Esseen_for_varkappa-close_to_normal},
and hence i.c.f.~strictly better than~\eqref{Eq:Problem_with_d_n}
with~\eqref{Eq:Ulyanov1976-distance},
while Example~\ref{Example:Zolotarev_1973}
with (\ref{Eq:||_||_K_Zolotarev-example},\ref{Eq:Zolotarev_example_r-norm_errors_n=1})
shows that~\eqref{Eq:Berry-Esseen_for_Z-close_to_normal_also_for_n=1}
is also i.c.f.~strictly better
than~\eqref{Eq:Problem_with_d_n} with~\eqref{Eq:Zolotaref1973-zeta_1_3-distance},
than~\eqref{Eq:Senatov_1998_with_K_on_RHS},
and than~\eqref{Eq:Berry-Esseen_for_just_zeta_3-close_to_normal},
which is in turn i.c.f.~strictly better
than~\eqref{Eq:Problem_with_d_n} with~\eqref{Eq:Zolotaref1973-zeta_3-distance}.

Restricting attention to $n\ge4$ for simplicity, we see that of the previously
best solutions to Problem~\ref{Zolotarev-Problem}, as summarised at the end of
subsection~\ref{subsec:Known_solutions_with_zeta}, only~\eqref{Eq:Problem_with_d_n}
with~\eqref{Eq:Shiganov1987-distance_3}, that is,
\la              \label{Eq:Berry-Esseen_for_nu_3-close_to normal}
 \left\|\widetilde{P^{\ast n}} -\mathrm{N}\right\|_{\mathrm{K}}
  &\le& \frac{35}{\sqrt{n}}\,\big(\nu_3^{1\wedge\frac{n}{4}}
     \!\vee\!\nu^{}_3\big) (\widetilde{P}-\mathrm{N})
   \quad\text{ for $P\in\cP_3$ and $n\in \N$}\,,
\al
and~\eqref{Eq:simple_Kolomogorov-telescoping} have not been shown to
be i.c.f.~strictly worse than~\eqref{Eq:Berry-Esseen_for_Z-close_to_normal_also_for_n=1},
and the three bounds now under consideration are in fact mutually incomparable.
For incomparability
of~\eqref{Eq:Berry-Esseen_for_Z-close_to_normal_also_for_n=1}
and~\eqref{Eq:Berry-Esseen_for_nu_3-close_to normal},
we note that in Example~\ref{Example:Zolotarev_1973}
and in case of $n\ge3$,
\eqref{Eq:Berry-Esseen_for_Z-close_to_normal_also_for_n=1} is i.c.f.~strictly worse than~\eqref{Eq:Berry-Esseen_for_nu_3-close_to normal},
by~\eqref{Eq:Zolotarev_example_r-norm_errors_n=1} with $r\in\{1,3\}$,
while in discrete cases like Example~\ref{Example:Discretised_normal_laws}
and for every $n\in\N$,
\eqref{Eq:Berry-Esseen_for_Z-close_to_normal_also_for_n=1} is
i.c.f.~strictly better than~\eqref{Eq:Berry-Esseen_for_nu_3-close_to normal}.

So Theorem~\ref{Thm:BE_K_Z}
can be improved by replacing there $\zeta_1\!\vee\!\zeta_3$ 
by the on $\widetilde{\cP_3}-\mathrm{N}$ also i.c.f.~strictly smaller
functional $(\zeta_1\!\vee\!\zeta_3)\wedge  
(\nu_3^{1\wedge\frac{n}{4}} \!\vee\!\nu^{}_3)$, and so in case of $n\ge 4$ by
$(\zeta_1\!\vee\!\zeta_3)\wedge   \nu^{}_3
= (\zeta_1\wedge\nu_3)\vee \zeta_3$, 
using~\eqref{Eq:zeta_vs_varkappa_etc} in the last step.
But, apart from considering~\eqref{Eq:simple_Kolomogorov-telescoping},  
are there any further and perhaps nicer improvements?     
This we presently do not know, but we can rule out the following idea: 
 
A natural try for improving Theorem~\ref{Thm:BE_K_Z} is to consider 
replacing $\zeta_1$ in~\eqref{Eq:Berry-Esseen_for_Z-close_to_normal} 
by the so-called dual bounded Lipschitz norm $\beta$ defined,
recalling~\eqref{Eq:Def_Lip-norm}, through
\la
 \left\| g \right\|_{\mathrm{BL}} &\coloneqq& \left\|g\right\|_{\mathrm{L}} + \left\|g\right\|_{\infty} 
            \quad\text{ for }g\in \C^\R, \qquad  
  \cG  \,\ \coloneqq \,\   
    \left\{ g\in \C^\R :  \left\| g \right\|_{\mathrm{BL}} \le 1\right\}, \nonumber \\  
 \beta(M) &\coloneqq& 
      \sup_{g\in\cG} \left|\int g\,\dd M\right| \quad\text{ for }M\in\cM\,. \label{Eq:Def_beta}
\al
Here $\beta$ is indeed a norm on $\cM$, 
was introduced by \citet[pp.~277--278]{FortetMourier1953}, and popularised
by R.M.~Dudley in particular, as in \citet[Chapter 11]{Dudley2003} and in the references given there.
Recalling~(\ref{Eq:Def_cF_r^infty},\ref{Eq:Def_zeta_r}), 
we observe that $\cG\subseteq \cF^\infty_1$, and hence $\beta\le \zeta_1$ on $\cM$.
In fact $\beta$ is i.c.f.~strictly smaller than $\zeta_1$ 
even on $\widetilde{\cP_3}-\mathrm{N}$, 
by the following example proved in section~\ref{sec:special_laws}, 
starting there on page~\pageref{page:Proof_of_Example_beta_llcurly_zeta_1}.
We recall that $\phi=\Phi'$ denotes the standard normal density.

\begin{Example}                          \label{Example:beta_llcurly_zeta_1}
For $t\in\mathopen]0,\infty\mathclose[\,$, there are unique
$p=p_t\in\mathopen]0,1\mathclose[$ and $s=s_t\in\mathopen]0,\infty\mathclose[$
with 
\la                                \label{Eq:Def_example_beta_llcurly_zeta_1}
 P &\coloneqq& P_t \,\ \coloneqq \,\ \big(\phi-\phi(t)\big)\1_{]-t,t[}\leb
      + \tfrac{p}{2}(\delta_{-s}+\delta_s) 
   \,\ \in\,\ \widetilde{P_3}\,,
\al
and then, with asymptotics referring to $t\rightarrow\infty$,  
\[
 \beta(P-\mathrm{N}) &\le& \nu_0(P-\mathrm{N}) \,\ \sim\,\ 2t\,\phi(t),\\
 \zeta_1 (P-\mathrm{N}) &\ge& \nu_1(P)-\nu_1(\mathrm{N}) 
   \,\ \sim\,\ ( \tfrac{2}{\sqrt{3}} - 1)t^2\phi(t), 
\]
and hence $\beta(P-\mathrm{N}) \llcurly \zeta_1 (P-\mathrm{N})$. 
\end{Example}

On the other hand we nevertheless have 
\la                      \label{Eq:varkappa_vee_zeta_bounded_by_beta_vee_zeta_in_Intro}
  \beta\!\vee\!\zeta_3 &\le& \zeta_1\!\vee\!\zeta_3 
    \,\ \le\,\ (2+3^\frac{1}{3})\, \beta\!\vee\!\zeta_3  \quad \text{ on }\cM,
\al
with the right hand inequality being~\eqref{Eq:varkappa_vee_zeta_bounded_by_beta_vee_zeta}
from Lemma~\ref{Lem:zeta_1_bounded_by_beta_and_zeta_3}.
Hence we obtain the following  i.c.f.~equivalent 
version of Theorem~\ref{Thm:BE_K_Z}:

\begin{Cor}                                    \label{Cor:B-E-Z_with_beta}
Theorem~\ref{Thm:BE_K_Z} remains true if $\zeta_1$ is decreased to $\beta$,
and $7.2$ increased to $25$.
\end{Cor}
\begin{proof} Inequality~\eqref{Eq:Berry-Esseen_for_Z-close_to_normal} combined 
with~\eqref{Eq:varkappa_vee_zeta_bounded_by_beta_vee_zeta_in_Intro}, and
$7.2\,(2+3^\frac{1}{3})= 24.784\ldots  \le  25$.
\end{proof}

Apparently not much is known about
lower bounds for $\text{L.H.S.\eqref{Eq:Berry-Esseen_for_Z-close_to_normal}}$.
The following nontrivial but presumably improvable result merely addresses the case of $n=2$.
It is in its more interesting first part a reformulation and specialisation,
and in its second part an improvement by elimination of logarithmic factors,
of \citet[Theorems 1.2 and 1.3]{BCG2012}.

\begin{Thm}[mainly Bobkov, Chistyakov, and G\"otze \citeyear{BCG2012}]  \label{Thm:BoChiGoe}
There exist constants $c,C\in\mathopen]0,\infty\mathclose[$ such that the following holds:
For a function $h:[0,1]\toitself$ to satisfy
\la                                                \label{Eq:BCG_lower_bound}
  h\!\left(\, \left\|\widetilde{P}-\mathrm{N} \right\|_{\mathrm{K}} \,\right)
   &\le& 
    \left\|\widetilde{P^{\ast 2}} -\mathrm{N}\right\|_{\mathrm{K}}   
   \quad\text{ for $P\in\cP_2$}\,,
\al
it is sufficient that we have 
\la                                    \label{Eq:BCG_our_h_=_inverse of theirs}
  h(t) &=& c\,\frac{t^{\frac{5}{2}}}{\,1\!\vee\!\log(\frac{1}{t})\,} \quad \text{ for }t\in[0,1],
\al
and necessary, even if the Kolmogorov norm $\left\|\cdot\right\|_{\mathrm{K}}$
on the right in \eqref{Eq:BCG_lower_bound} is replaced by the 
total variation norm $\nu_0$, that we have
\la
   h(t) &\le& C\,t^2 \quad \text{ for }t\in[0,1].
\al
\end{Thm}

This is proved in section~\ref{sec:Proof_of_Theorem_about_lower_bounds}. 
As a side remark, which one might take into account when trying to improve 
Theorem~\ref{Thm:BoChiGoe}, let us mention the following sharpening 
of the implication ``$\Rightarrow$'' in~\eqref{Eq:Elementary_Cramer_Levy}: 
If $\left\|\widetilde{P}-\mathrm{N}\right\|_{\mathrm{K}}>0$, 
then not only $\text{R.H.S.}\eqref{Eq:BCG_lower_bound}>0$, 
but even $\sup\limits_{x\in\mathopen]-\infty,x_0\mathclose]} \left|
 F_{\widetilde{P^{\ast n}}}(x)-\Phi(x)  \right| >0$
for every $x_0\in\R$ and every $n\ge2$.
This is a result of~\citet[\foreignlanguage{russian}{Teorema}~1]{Titov1981},
also presented by \citet[p.~85, Corollary~4.7.6]{RossbergJesiakSiegel1985}.
 
Looking at Theorem~\ref{Thm:BE_K_Z}, one should of course ask for extensions to questions
as those indicated in the final paragraph of subsection~\ref{subsec:Zolotarev-Problem} above.
Let us pose here just one specific such question, and a further one as
Question~\ref{Question:BE_K_Z_for_simple_random_sampling} in the next section below.

\begin{Question}
Do we have
\la                                       \label{Eq:BE_Nag_Z?}
  \qquad\sup_{x\in\R}\big(1+|x|^3\big)\left| F_{\widetilde{P^{\ast n}}}(x)-\Phi(x) \right| &\le& 
   \frac{c}{\sqrt{n}}\big(\zeta_1\!\vee\!\zeta_3\big)(\widetilde{P}-\mathrm{N})  
   \quad\text{ for $P\in\cP_3$ and $n\ge n_0$}
\al
for any constants $c\in\mathopen]0,\infty\mathclose[$ and $n_0\in\N$?
\end{Question}

A positive answer would, for $n\ge n_0$, improve 
i.c.f.~\citeposs[p.~570, Theorem 3.1 with dimension $k=1$]{Sazonov1972}
and \citeposs[Theorem~2, $g(x)=|x|$, Remark B]{Ulyanov1976}
incomparable improvements of 
\citeposs[the second theorem]{Nagaev1965}
improvement $\text{L.H.S.\eqref{Eq:BE_Nag_Z?}}\le\text{R.H.S.\eqref{Eq:Berry-Esseen_inequality}}$
of the Berry-Esseen theorem~\eqref{Eq:Berry-Esseen_inequality}.
By the following example, 
$n_0=2$ will not do in \eqref{Eq:BE_Nag_Z?},
in contrast to~\eqref{Eq:Berry-Esseen_for_Z-close_to_normal}.

\begin{Example}[Left-winsorised normal laws]\label{Example:Left-winsorised_normal_re_Nagaev}
For $P=P_t$ as in 
Example~{\rm\ref{Examples:zeta_3=mu_3}\ref{part:Left-winsorised_normal_laws}}, 
we have
\[
 \lim_{t\rightarrow\infty} 
   \frac{\sup\limits_{x\in\R}\big(1+|x|^3\big)\left| F_{\widetilde{P^{\ast 2}}}(x)-\Phi(x)\right|}
        {\big(\zeta_1\!\vee\!\zeta_3\big)(\widetilde{P}-\mathrm{N})}
 &= & \infty\,.       
\]
\end{Example}

This is proved in section~\ref{sec:special_laws},
starting on page~\pageref{page:Proof_of_Example_Left-winsorised_normal_re_Nagaev}. 
 
\section{Theorem~\ref{Thm:BE_K_Z} applied to sums of simple random samples from a finite population}
This section is not logically necessary for understanding the rest of the present paper, 
and may hence be skipped. Its purpose is to illustrate by 
Corollary~\ref{Cor:Normal_approximation_for_sampling_from_finite_populations_Z-close_to_normal}
the importance of having error bounds 
of the form (\ref{Eq:Problem_with_d_n},\ref{Eq:Problem_with_norm})
with metrics or norms strictly weaker on discrete laws than~$\nu_3$. 
We use below the customary letter $N$ for a population size, which should not lead to any confusion
with the upright letter $\mathrm{N}$ denoting the standard normal law.
\begin{Cor}[A normal approximation error bound for sums of samples from a finite
population Zolotarev-close to normal]
             \label{Cor:Normal_approximation_for_sampling_from_finite_populations_Z-close_to_normal}
Let $M$ be a set of cardinality $N\coloneqq\#M\in \N$, $x\in\R^M$ a ``population'' 
with the ``value range'' 
$\cX\coloneqq\{x_i:i\in M\}$ and the ``diversity'' $d\coloneqq\#\cX\ge 2$,
and $P\coloneqq \frac{1}{N}\sum_{i\in M}\delta_{x_i} \in\Prob(\R)$.
Let further $n\in\{1,\ldots,N\}$ and let the random variable~$S$, on some probability space $(\Omega,\cA,\P)$,
be a simple random sample of size $n$ from $M$, that is, with $\cS\coloneqq\{s\subseteq M:\#s=n\}$
we require $S:\Omega\rightarrow\cS$ to be uniformly distributed, namely $\P(S=s)=\binom{N}{n}^{-1}$
for $s\in\cS$.
Then the real-valued random variable 
\[
   Z &\coloneqq& \frac{\sum
   _{i\in S}x_i-n\mu(P)}{\sqrt{n}\sigma(P)}
\]
satisfies
\la    \label{Eq:Normal_approximation_for_sampling_from_finite_populations_Z-close_to_normal}
 \big\| \P(Z\in\cdot\,) - \mathrm{N}\big\|^{}_{\mathrm{K}} 
 &\le& \frac{7.2}{\sqrt{n}}
  \big(\zeta_1\!\vee\!\zeta_3\big)(\widetilde{P}-\mathrm{N}) 
    \ +\ \left(\tfrac{n-1}{2}\wedge d \right)\tfrac{n}{N} \quad\text{ if 
    $n\ge2$.}
\al
\end{Cor}          
\begin{proof} 
Let, on a possibly different probability space again denoted by $(\Omega,\cA,\P)$,
$T=\lpp T_1,\ldots,T_n\rpp$ and $U=\lpp U_1,\ldots,U_n\rpp$ be random variables with 
$T$ uniformly distributed on $M^n_{\neq}\coloneqq \{t\in M^n : t_i\neq t_j\text{ for }i\neq j\}$,
and $U$ uniformly distributed on $M^n$ (``successive random samples of size $n$ from $M$,
without, respectively with, replacement''). 
Here we use the double parentheses notation $\lpp \ldots \rpp$
for tuple-valued functions, in order to avoid abusing the notation $(\ldots)$ for tuples of functions.
With ``$\sim$'' here to be read as ``is distributed as'', 
and with an abuse of notation analogous to the one just avoided,
we then have $S \sim \{T_1,\ldots,T_n\}$,
and hence 
\[
    Z &\sim&   \frac{\sum_{j=1}^nx^{}_{T_j} -\mu(P)}{\sqrt{n}\sigma(P)}\,.
\]
With 
\[
   W &\coloneqq& \frac{\sum_{j=1}^nx^{}_{U_j} -\mu(P)}{\sqrt{n}\sigma(P)}
\]
we have $\P(W\in\cdot\,)=\widetilde{P^{\ast n}}$ and get
\[
 \big\| \P(Z\in\cdot\,) - \mathrm{N}\big\|^{}_{\mathrm{K}} 
  &\le& \big\| \P(W\in\cdot\,) - \mathrm{N}\big\|^{}_{\mathrm{K}}
     +  \big\| \P(Z\in\cdot\,) - \P(W\in\cdot\,) \big\|^{}_{\mathrm{K}}
\]
with on the right the first summand $\le$ $\text{R.H.S.\eqref{Eq:Berry-Esseen_for_Z-close_to_normal}}$ 
if $n\ge2$.
Writing now $D(X,Y) \coloneqq \sup_{B\in\cB}|\P(X\in B)-\P(Y\in B)|$ 
for the supremum distance of the laws of any $(\cX,\cB)$-valued random variables $X,Y$
(which, logically unnecessary to state here, but perhaps helpful to avoid the 
usual confusion,
is  $\frac12\nu_0\big( \P(X\in\cdot\,)-\P(Y\in\cdot\,)\big)$ according to~\eqref{Eq:Def_nu_r} 
with $r=0$ if $\cX=\R$),
we have
\[
 \big\| \P(Z\in\cdot\,) - \P(W\in\cdot\,) \big\|^{}_{\mathrm{K}}
  &\le& D(T,U) \,\ \le\,\ \tfrac{(n-1)n}{2N}
\]
with the last bound noted by \citet{Freedman1977} and by \citet[p.~84]{Stam1978},
and with
\[
     A \,\ \coloneqq\,\ \left\lpp \sum_{j=1}^n \1_{\{x^{}_{T_j}=\xi\}} : \xi \in \cX  \right\rpp,  
  && B \,\ \coloneqq\,\ \left\lpp \sum_{j=1}^n \1_{\{x^{}_{U_j}=\xi\}} : \xi \in \cX  \right\rpp
\]
we also have 
\[
 \left\| \P(Z\in\cdot\,) - \P(W\in\cdot\,) \right\|_{\mathrm{K}}
  &\le& D(A,B) \,\ \le \,\ \tfrac{dn}{N}
\]  
by \citet[p.~746, Theorem (4)]{DiaconisFreedman1980}. Hence the claim.
\end{proof}             
             
Let us compare~\eqref{Eq:Normal_approximation_for_sampling_from_finite_populations_Z-close_to_normal}
with a classical and recently improved Berry-Esseen type theorem
for standardised sums of samples from a finite population: In the situation of 
Corollary~\ref{Cor:Normal_approximation_for_sampling_from_finite_populations_Z-close_to_normal},
we have the well-known variance formula
\[
 \sigma^2\left(\sum_{i\in S}x_i\right) &=& n\tfrac{N-n}{N-1}\sigma^2(P)\,, 
\]
so that, recalling the notation $\widetilde{X}$ for the standardisation of a nondegenerate finite variance 
real-valued random variable $X$, and assuming from now on $n\le N-1$, 
\la                                                         \label{Eq:Def_tilde(Z)}
 \widetilde{Z} &=& \sqrt{\tfrac{N-1}{N-n}}Z\,,
\al
and the theorem in question yields
\la                      \label{Eq:Hoeglund1976}
  \left\| \P(\widetilde{Z}\in\cdot\,) - \mathrm{N}\right\|_{\mathrm{K}} 
 &\le& c\frac{\nu_3(\widetilde{P})}{\sqrt{n\frac{N-n}{N-1}}}
\al
with $c=82.4$. A result equivalent to~\eqref{Eq:Hoeglund1976} with the universal constant 
$c<\infty$ unspecified was apparently first obtained by~\citet{Hoeglund1976}, namely 
\la                         \label{Eq:Hoeglund1976_really}
 \left\| \P\big(\sqrt{\tfrac{N-1}{N}}\widetilde{Z}\in\cdot\,\big) - \mathrm{N}\right\|_{\mathrm{K}} 
    &\le& c\frac{\nu_3(\widetilde{P})}{\sqrt{n\frac{N-n}{N}}}\,,
\al
with the equivalence becoming clear by observing 
$\text{L.H.S.\eqref{Eq:Hoeglund1976_really}}
 =\| \P(\widetilde{Z}\in\cdot\,) - \mathrm{N}_{\sqrt{\frac{N}{N-1}}}\|^{}_{\mathrm{K}}$
by 
scale invariance of $\|\cdot\|^{}_{\mathrm{K}}$ 
and therefore, 
with $c_0\coloneqq \frac{1}{\sqrt{2\pi\mathrm{e}}}$ and using Lemma~\ref{Lem:Kolmogorov_distance_of_centred_normals}
in the second step,
$|\text{L.H.S.\eqref{Eq:Hoeglund1976_really}} -\text{L.H.S.\eqref{Eq:Hoeglund1976}}|$
 $\le$ $\| \mathrm{N} - \mathrm{N}_{\sqrt{\frac{N}{N-1}}}\|_{\mathrm{K}}
 \,\le\, c_0|\sqrt{\frac{N}{N-1}}-1|
 \,\le\, \frac{c^{}_0}{2(N-1)} \le \text{R.H.S.\eqref{Eq:Hoeglund1976}
          with $c=\frac{c^{}_0}{2}$}$, 
and also  $\text{R.H.S.\eqref{Eq:Hoeglund1976}}
 \,\le\, \text{R.H.S.\eqref{Eq:Hoeglund1976_really}}
 \,\le\, \sqrt{2}\,\text{R.H.S.\eqref{Eq:Hoeglund1976}}$.
Inequality~\eqref{Eq:Hoeglund1976} with $c=451$
was obtained by \citet[p.~337, Corollary~1.4 in the special case of $\sigma_i=0$]{ChenFang2015}, 
follows from \citet[Corollary 1.2 and the line before it]{Thanh2013} with $c=90$
as a consequence of~\eqref{Eq:Hoeglund1976} improved by the additional factor 
$\frac{N-1}{N}\big((\frac{N-n}{N})^2+(\frac{n}{N})^2\big) \in [\frac{N-1}{2N},1[$ on the right, 
and follows with $c=82.4$ as claimed above from the bound in
\citet[Corollary 1.1]{Roos2020}, which is in fact strictly better than~\eqref{Eq:Hoeglund1976}
i.c.f.~thanks to an additional minimum operation.
Better admissible constants for~\eqref{Eq:Hoeglund1976} appear to be known for special cases only,
namely $c=1.1166$  
for the hypergeometric case of $d=2$, 
and $ c= \frac{1}{\sqrt{2\pi}}$ for the symmetric hypergeometric subcase of 
$\widetilde{P}=\frac12(\delta_{-1}+\delta_{1})$, by 
results reviewed or proved by \citet[first paragraph on p.~733, 
$h(1)$ on p.~729]{MattnerSchulz2018}.  

A lower bound for the unknown optimal constant in~\eqref{Eq:Hoeglund1976} is $c^{}_{\mathrm{E}}$
from~\eqref{Eq:Def_c_E}, since H\"oglund's theorem with any constant $c$ in~\eqref{Eq:Hoeglund1976}
yields as a limiting case the classical Berry-Esseen theorem~\eqref{Eq:Berry-Esseen_inequality} with 
the same $c$, 
see 
for example~\citet[pp.~728--729, in particular Lemma 1.1]{MattnerSchulz2018}.
Now, in view of the scale invariance of $\left\|\cdot\right\|_{\mathrm{K}}$ again, 
either of~\eqref{Eq:Normal_approximation_for_sampling_from_finite_populations_Z-close_to_normal}
and~\eqref{Eq:Hoeglund1976} provides an error bound for some approximation of, say, 
$\P(\widetilde{Z}\in\cdot)$, namely for the approximation $\mathrm{N}$ in 
case of~\eqref{Eq:Hoeglund1976}, and for $\mathrm{N}_{\sqrt{(N-1)/(N-n)}}$ in case 
of~\eqref{Eq:Normal_approximation_for_sampling_from_finite_populations_Z-close_to_normal},
and the error bound 
in~\eqref{Eq:Normal_approximation_for_sampling_from_finite_populations_Z-close_to_normal} 
can easily be smaller than the one in~\eqref{Eq:Hoeglund1976}, not only with the presently 
best admissible value $82.4$ for $c$, but even if we most optimistically 
assume~\eqref{Eq:Hoeglund1976} to be true with $c=c^{}_{\mathrm{E}}$.

A comparison of the classical Berry-Esseen theorem~\eqref{Eq:Berry-Esseen_inequality}, 
its improvement Theorem~\ref{Thm:BE_K_Z},  
and H\"oglund's generalisation~\eqref{Eq:Hoeglund1976} of~\eqref{Eq:Berry-Esseen_inequality},
suggests to us:

\begin{Question}                        \label{Question:BE_K_Z_for_simple_random_sampling}
In the situation of 
Corollary~\ref{Cor:Normal_approximation_for_sampling_from_finite_populations_Z-close_to_normal}
and with $\widetilde{Z}$ as in~\eqref{Eq:Def_tilde(Z)}, do we have
\[ 
 \text{L.H.S.\eqref{Eq:Hoeglund1976}} &\le&  
  c \frac{ \big(\zeta_1\!\vee\!\zeta_3\big)(\widetilde{P}-\mathrm{N})}
     { \sqrt{n\frac{N-n}{N-1}}}\quad\text{ if }n\ge 2
\] 
for some universal constant $c$? 
\end{Question}

Coming now finally to the main point of this section within the present paper, 
we observe that neither 
Corollary~\ref{Cor:Normal_approximation_for_sampling_from_finite_populations_Z-close_to_normal}
nor a positive answer to Question~\ref{Question:BE_K_Z_for_simple_random_sampling}
would yield any improvement i.c.f.~over H\"oglund's~\eqref{Eq:Hoeglund1976}
if $\zeta_1\!\vee\!\zeta_3$ were replaced by any bound involving~$\nu_3$, 
as~\eqref{Eq:Problem_with_d_n} with~\eqref{Eq:Shiganov1987-distance_3},
since here each $P$ is discrete and hence satisfies 
$\nu_3(\widetilde{P}-\mathrm{N})=\nu_3(\widetilde{P})+\nu_3(\mathrm{N})\ge \nu_3(\widetilde{P})$.

\section{Zolotarev's $\zeta_1\!\vee\!\zeta_3$ Theorem~\ref{Thm:Zolotarev's_zeta_1-B-E-Thm}, 
the convolution inequality Theorem~\ref{Thm:F_star_G_vs_H_star_H}, 
and a proof of Theorem~\ref{Thm:BE_K_Z}}                            \label{sec:3}
In this section we state 
Theorems~\ref{Thm:Zolotarev's_zeta_1-B-E-Thm} 
and~\ref{Thm:F_star_G_vs_H_star_H}, 
postpone their proofs to sections~\ref{sec:zeta_distances} 
and~\ref{sec:proof_of_convolution_inequality},
but already apply them here to prove Theorem~\ref{Thm:BE_K_Z}.

In Theorem~\ref{Thm:Zolotarev's_zeta_1-B-E-Thm} below, the triple use of the symbol $\zeta$, 
namely to denote with index~$1$ or~$3$ a Zolotarev norm on $\cM_{3,2}$, 
with no index and no argument a variable, and with no index and the argument~$\frac{3}{2}$ 
a value of the Riemann zeta function, should not cause any confusion. 

\begin{Thm}[essentially Zolotarev \citeyear{Zolotarev1986}, \citeyear{Zolotarev1997}]
                         \label{Fact:Zolotarev1986}\label{Thm:Zolotarev's_zeta_1-B-E-Thm}
There exists a constant $c\in\mathopen]0,\infty\mathclose[$ satisfying 
 \la                                                            \label{Eq:Zolotarev's_zeta_1-B-E-Thm}
    \zeta_1(\widetilde{P^{\ast n}} -\mathrm{N})   
     &\le& \frac{c}{\sqrt{n}}\big(\zeta_1\!\vee\!\zeta_3\big)(\widetilde{P}-\mathrm{N}) 
   \quad\text{ for $P\in\cP_3$ and $n\ge 1$.}
 \al
One may take here $c=14$. More precisely, we have
\la                                                              \label{Eq:Zolotarev's_zeta_1-B-E-Thm_sharper}
  \zeta_1(\widetilde{P^{\ast n}} -\mathrm{N}) &\le& \frac{1}{\sqrt{n}}
   \xi\left( \zeta_1(\widetilde{P}-\mathrm{N}) , \zeta_3(\widetilde{P}-\mathrm{N}) \right) 
   \quad\text{ for $P\in\cP_3$ and $n\ge 1$},
\al
where the function $\xi:[0,\infty[^2\rightarrow[0,\infty[$ is defined through
\la
  &&
 \alpha\coloneqq \frac{4\mathrm{e}^{-1/2}}{\sqrt{2\pi}}= 0.9678\ldots,\quad
 \beta\coloneqq \frac{4}{\sqrt{2\pi}} = 1.5957\ldots,  \quad 
   \gamma\coloneqq  \frac{2+8\mathrm{e}^{-3/2}}{\sqrt{2\pi}}=  1.5100\ldots\,,  \label{Eq:Def_alpha_beta_gamma} \\
 && g(\eta) \,\ \coloneqq\,\ \sum_{j=1}^\infty \frac{1}{(j+\eta^2)^{3/2}}
  \,\ < \,\ \frac{2}{\eta}
 \quad \text{ for } \eta\in[0,\infty[\,,            \label{Eq:Def_g(eta)} \\
 && \xi(\varkappa,\zeta)\,\ \coloneqq\,\                               \label{Eq:Def_xi(kappa,zeta)}
   \inf\left\{\frac{\varkappa+\alpha\zeta+\beta\eta}{1-\gamma g(\eta)\zeta} : 
     \eta\in\mathopen[0,\infty\mathclose[,\ \gamma g(\eta)\zeta <1  \right\}         
     \quad\text{ for }(\varkappa,\zeta)\in[0,\infty[^2,
\al
which easily yields \eqref{Eq:Zolotarev's_zeta_1-B-E-Thm} with $c=23.21\ldots$, 
but we also have 
\la                                                                            \label{Eq:xi_near_zero}
 \xi(\varkappa,\zeta)&\le& \frac{\varkappa +\alpha\zeta }{1- \lambda\zeta} \quad\text{ for } 
    \zeta < \frac1\lambda = 0.2535\ldots
    \text{ with }\lambda \coloneqq \gamma\zeta(\tfrac{3}{2})=3.9447\ldots,
\al 
and, from just~\eqref{Eq:Zolotarev's_zeta_1-B-E-Thm_sharper} and~\eqref{Eq:xi_near_zero} combined
with the \citet{Goldstein2010}-\citet{Tyurin2010} theorem~\eqref{Eq:Goldstein-Tyurin} below, 
the validity of \eqref{Eq:Zolotarev's_zeta_1-B-E-Thm} with $c=13.3803\ldots$\;.
\end{Thm}

This is proved in section~\ref{sec:Proof_Zolotarev-Thm}.

Let us note that for $g(\eta)$ in~\eqref{Eq:Def_g(eta)} 
we have  $g(\eta) = \zeta(\frac{3}{2},\eta^2) - \eta^{-3} $
with the Hurwitz zeta function $\zeta(\cdot,\cdot)$ from \citet[pp.~607--610, section 25.11]{Olver2010NIST},
so that various representations of the latter might be used
- apparently uninterestingly for our present purposes - 
to refine the inequality in~\eqref{Eq:Def_g(eta)} and hence improve the constant  $c=23.21\ldots$
a bit.

In 
\citet[pp.~365--368, in particular~(6.5.41)]{Zolotarev1997}
inequality~\eqref{Eq:Zolotarev's_zeta_1-B-E-Thm}
is stated with the constant $c=8.35$, but the proof presented there yields, after 
correcting the trite error there of having $\frac{1}{2\eta}$ rather than $\frac{2}{\eta}$ in~\eqref{Eq:Def_g(eta)}, 
only a somewhat larger value of $c$.
The present essentially self-contained version of Zolotarev's proof,  given below on 
pages~\pageref{Start_of_proof_Zolotarev's_zeta_1-B-E-Thm}--\pageref{End_of_step_2_in_proof_Zolotarev's_zeta_1-B-E-Thm}
in the steps~1 and~2, improves a bit on this latter constant by using a better and actually simpler choice of a 
parameter ($m$ in Zolotarev's notation), but is then followed in step~4  
by a use of the Goldstein-Tyurin theorem to arrive at the constant $13.3803\ldots$\,. 
While this still seems to be rather large,  we observe that the factors $1$ and $\alpha$ 
of $\varkappa$ and $\zeta$ in the numerator in~\eqref{Eq:xi_near_zero}, and 
for $\zeta\rightarrow0$ only this is asymptotically relevant, 
are quite small; in particular the factor $1$ is optimal, as can be seen by taking $n=1$ 
in~\eqref{Eq:Zolotarev's_zeta_1-B-E-Thm_sharper} and
any examples where $\zeta_3(\widetilde{P}-\mathrm{N})$ is small compared to $\zeta_1(\widetilde{P}-\mathrm{N})$, 
as in Example~\ref{Example:Discretised_normal_laws}.

One may easily ``improve'' the error bounding~\eqref{Eq:Zolotarev's_zeta_1-B-E-Thm}
by combining it with
the simpler fact~\eqref{Eq:Simple_zeta_3-CLT_error_bound} below, yielding
\la                                        \label{Eq:Zolotarev's_zeta_1-B-E-Thm_combined_with_zeta_3}
  \big(\zeta_1\!\vee\!\zeta_3\big)(\widetilde{P^{\ast n}} -\mathrm{N})   
     &\le& \frac{c}{\sqrt{n}}\big(\zeta_1\!\vee\!\zeta_3\big)(\widetilde{P}-\mathrm{N}) 
   \quad\text{ for $P\in\cP_3$ and $n\ge 1$},
\al
with the same norm $\zeta_1\!\vee\!\zeta_3$ occurring on both sides, justifying to some extent the description
of Theorem~\ref{Thm:Zolotarev's_zeta_1-B-E-Thm} in the title of the present section.


We now state the main technical result of the present paper, using here
standard analytical notation as in and directly before~\eqref{Eq:Def_Lip-norm},
the perhaps not so standard definition of $\star$ in~\eqref{Eq:Def_conv_distr_fcts},
and the ad hoc notation
\la
 \ltrivert f \rtrivert &\coloneqq&
  \left\{
   \begin{array}{ll}
    \tfrac{1}{2}\left( \left\| f\right\|_1+ \big|\int \!f\,\dd\leb\big|\,\right) &\text{ if } \left\| f\right\|_1<\infty,\\
    \infty  &\text{ if } \left\| f\right\|_1=\infty,
   \end{array}
  \right.
\al
so that we have $\frac{1}{2} \left\| f\right\|_1 \le \ltrivert f \rtrivert \le  \left\| f\right\|_1$.
\begin{Thm}                                                                \label{Thm:F_star_G_vs_H_star_H}
                                       \label{page:Thm_convolution_inequality}
Let $F_1,F_2,H_1,H_2$ be probability distribution functions on $\R$.
Then we have
\la                                 \label{Eq:The_convolution_inequality_new}
  \quad \left\| F_1\star F_2 -H_1\star H_2 \right\|_\infty &\le &
   \left( \sqrt{ \left\| H_2 \right\|_{\mathrm{L}}    \ltrivert F_1-H_1\rtrivert} +
      \sqrt{ \left\| H_1\right\|_{\mathrm{L}} \ltrivert F_2-H_2\rtrivert}\right)^2 .
\al
\end{Thm}
Obviously, inequality~\eqref{Eq:The_convolution_inequality_new} always implies
the weaker and slightly simpler inequality
\la           \label{Eq:F_star_G_vs_H_star_H} \label{Eq:The_convolution_inequality_old}
  \quad \left\| F_1\star F_2 -H_1\star H_2 \right\|_\infty &\le & 
   \left( \sqrt{ \left\| H_2 \right\|_{\mathrm{L}}    \left\|F_1-H_1\right\|_1} +
      \sqrt{ \left\| H_1\right\|_{\mathrm{L}} \left\|F_2-H_2\right\|_1}\right)^2 ,
\al
but reduces to~\eqref{Eq:The_convolution_inequality_old} improved by a factor $\frac{1}{2}$ on the right
in case of $\int (F_i-H_i)\,\dd\leb =0$ for $i\in\{1,2\}$.
This latter case is used for obtaining Corollary~\ref{Cor:Main_for_n=2} below.

Inequality~\eqref{Eq:The_convolution_inequality_new} is true even if, say,
$\left\| H_1\right\|_{\mathrm{L}}=\infty$,
for then either $ \ltrivert F_2-H_2\rtrivert >0$ and then~\eqref{Eq:F_star_G_vs_H_star_H}
is trivial due to $\text{R.H.S.}=\infty$,
or $F_2=H_2$ and then~\eqref{Eq:F_star_G_vs_H_star_H} reduces by the convention $\infty\cdot 0=0$
to $\left\| (F_1-H_1)\star H_2\right\|_\infty
 \le \left\|H_2\right\|_\mathrm{L} \ltrivert F_1-H_1\rtrivert$\,,
which is contained in, and is by the homogeneity of enorms actually equivalent to,
\eqref{Eq:||cdot||_K-smoothing_sharper}
in the simple Lemma~\ref{Lem:||cdot||_K-smoothing} below.
For the much more interesting case of $\left\| H_1\right\|_{\mathrm{L}}$
and $\left\| H_2\right\|_{\mathrm{L}}$ being both finite,
our proof of~\eqref{Eq:F_star_G_vs_H_star_H}
is a bit complicated and takes essentially all of
section~\ref{sec:proof_of_convolution_inequality}.

We now switch back to notation as introduced or explained 
in (\ref{Eq:Def_N_sigma_N_1},\ref{Eq:Def_F_M_and_overline{F}_M},\ref{Eq:Def_Kolmogorov-norm}),
in the paragraph around~\eqref{Eq:Def_cP_r}, and 
in~(\ref{Eq:Def_zeta_r},\ref{Eq:zeta_1=kappa_1},\ref{Eq:Def_varkappa},\ref{Eq:varkappa_1=L^1-enorm}).

\begin{Cor} Let $P,Q\in\cP_2$ with standard deviations $\sigma,\tau$. Then we have     \label{Cor:Main_for_n=2} 
\la
  \left\|\widetilde{P\ast Q} -\mathrm{N}\right\|_{\mathrm{K}}  
   &\le& \frac{1}{2\sqrt{2\pi}}\left(\sqrt{\tfrac{\sigma}{\tau}\zeta_1(\widetilde{P}-\mathrm{N})}+
    \sqrt{\tfrac{\tau}{\sigma}\zeta_1(\widetilde{Q}-\mathrm{N})}  \right)^2, \label{Eq:PQ-N}   \\
  \left\|\widetilde{P^{\ast2}}-\mathrm{N}\right\|_{\mathrm{K}}  
             &\le& \frac{2}{\sqrt{2\pi}}\zeta_1(\widetilde{P}-\mathrm{N}) .  \label{Eq:PP-N}
\al 
\end{Cor}
\begin{proof}[Proof of Corollary~\ref{Cor:Main_for_n=2} assuming Theorem~\ref{Thm:F_star_G_vs_H_star_H}]
Assuming w.l.o.g.~$\mu(P)=0=\mu(Q)$
and $\sigma^2+\tau^2=1$, we get
\[
 \text{L.H.S.\eqref{Eq:PQ-N}} &=& \left\| P\ast Q-\mathrm{N}_\sigma\ast\mathrm{N}_\tau\right\|_{\mathrm{K}} \\
   &\le& \left(\sqrt{\tfrac{1}{\tau\sqrt{2\pi}}\zeta_1(P-\mathrm{N}_\sigma)} 
    + \sqrt{\tfrac{1}{\sigma\sqrt{2\pi}}\zeta_1(Q-\mathrm{N}_\tau)} \right)^2 \,\ =\,\ \text{R.H.S.\eqref{Eq:PQ-N}}
\]
by applying in the second step Theorem~\ref{Thm:F_star_G_vs_H_star_H} to 
$F_1\coloneqq F_P$, $F_2 \coloneqq F_Q$, $H_1\coloneqq \Phi(\frac{\cdot}{\sigma})$, 
$H_2\coloneqq \Phi(\frac{\cdot}{\tau})$,
and in the third~\eqref{Eq:homogeneity_on_cM_conclusion}
for $\|\cdot\|=\zeta_r$ with $r=1$.
Specialising~\eqref{Eq:PQ-N} to $Q=P$ yields~\eqref{Eq:PP-N}.
\end{proof}

For inequality~\eqref{Eq:F_star_G_vs_H_star_H} to be nontrivial, 
we must for each $i$ have $\left\|F_i-H_i\right\|_1<\infty$, 
but the laws $P_i,R_i$ corresponding to $F_i,H_i$ need not have finite first moments.
Hence analogues of~(\ref{Eq:PQ-N},\ref{Eq:PP-N}) for general stable  
laws in place of $\mathrm{N}$ follow similarly from Theorem~\ref{Thm:F_star_G_vs_H_star_H}.

While in~\eqref{Eq:The_convolution_inequality_old} we have asymptotic equality
for appropriate $F_1=F_2$ close to but different from $H_1\coloneqq H_2\coloneqq \Phi$,
by Example~\ref{Example:Convolution_inequality_sharp}\ref{part:Example_Convolution_inequality_sharp_a},
and analogously in~\eqref{Eq:The_convolution_inequality_new} in a case
with $\ltrivert F_i-H_i\rtrivert=\frac{1}{2}\left\| F_i-H_i\right\|_1$
by Example~\ref{Example:Convolution_inequality_sharp}\ref{part:Example_Convolution_inequality_sharp_b},
it seems likely that the constant factors on the right in~\eqref{Eq:PQ-N}
and~\eqref{Eq:PP-N} might be improved
by exploiting equal variances in addition to equal means.
But certainly the constant $2$ in the numerator in~\eqref{Eq:PP-N} can
not be replaced by any number strictly smaller
than $\frac{15+6\sqrt{3}}{13}=1.9532\ldots$, 
by~\eqref{Eq:n=2_sharpness_BE_K_Z_Zolotarev_example} 
in Example~\ref{Example:Zolotarev_1973}.

Coming now to the proof of Theorem~\ref{Thm:Berry-Esseen_for_Z-close_to_normal},
let us first present its simple idea by restricting attention to $n=2k$ even
and by not caring about constant factors. We then have
\[
 \left\| \widetilde{P^{\ast n}} -\mathrm{N} \right\|_{\mathrm{K}}
 &\le& \frac{2}{\sqrt{2\pi}}\zeta_1\left(\widetilde{P^{\ast k}}-\mathrm{N}\right)
 \,\ \le\,\ \frac{2c_0}{\sqrt{2\pi} \sqrt{k}} \big(\zeta_1\!\vee\!\zeta_3\big)(\widetilde{P}-\mathrm{N})
 \,\ =\,\ \frac{c}{\sqrt{n}} \big(\zeta_1\!\vee\!\zeta_3\big)(\widetilde{P}-\mathrm{N})
\]
by applying in the first step~\eqref{Eq:PP-N} with $P^{\ast k}$ in place of  $P$,
in the second~\eqref{Eq:Zolotarev's_zeta_1-B-E-Thm} 
with an admissible value
$c_0$ in place of $c$,
and by taking in the third step
$c\coloneqq \sqrt{2} \frac{2c_0}{\sqrt{2\pi}}$,
so $c = 15.098\ldots$ in case of $c_0=13.3803\ldots$ as admissible by Theorem~\ref{Thm:Zolotarev's_zeta_1-B-E-Thm}.
The actual proof of Theorem~\ref{Thm:Berry-Esseen_for_Z-close_to_normal}
is a bit longer due to allowing also $n\ge 3$ odd and aiming at a better value of $c$.

\begin{proof}[Proof of Theorem~\ref{Thm:Berry-Esseen_for_Z-close_to_normal}
and Remark~\ref{Rem:BE_K_Z_sharper}
assuming Theorem~\ref{Thm:Zolotarev's_zeta_1-B-E-Thm} and Corollary~\ref{Cor:Main_for_n=2}] 
Let $P\in\cP_3$ and $n\in\N$ with $n\ge 2$.

\smallskip
1. Let  $\Xi : \cP_3\rightarrow[0,\infty[$ be a functional such that we have 
\la           \label{Eq:Abstract_Zolotarev_zeta_1-theorem}
  \zeta_1(\widetilde{P^{\ast k}} -\mathrm{N})
     &\le& \frac{ \Xi(P) }{\sqrt{k}}    \quad\text{ for $P\in\cP_3$ and $k\in \N$}.
\al
Since $n=2k +\varrho$ with $k\in\N$ and $\varrho\in\{0,1\}$,
we get,
using below the convolution inequality~\eqref{Eq:PQ-N} from Corollary~\ref{Cor:Main_for_n=2}
with $P^{\ast k}$ and $P^{\ast(k+\varrho)}$ instead of $P$ and $Q$ in the second step,
and assumption~\eqref{Eq:Abstract_Zolotarev_zeta_1-theorem}
once for $k$ and once for $k+\varrho$
in the third,
\la
 \left\| \widetilde{P^{\ast n}} -\mathrm{N} \right\|_{\mathrm{K}}  \label{Eq:Ineq_chain_BE_for_Z-close_to_normal}
  &=&  \left\| \widetilde{ P^{\ast k}\ast P^{\ast(k+\varrho)}    } -\mathrm{N} \right\|_{\mathrm{K}}       \\   
  &\le& \frac{1}{2\sqrt{2\pi}}\left(\sqrt{ \sqrt{\tfrac{k}{k+\varrho}}\,\zeta_1(\widetilde{P^{\ast k}}-\mathrm{N})}+
   \sqrt{\sqrt{\tfrac{k+\varrho}{k}}\,\zeta_1(\widetilde{P^{\ast(k+\varrho)}}-\mathrm{N})} \right)^2\nonumber  \\
  &\le& \frac{\Xi(P)}{2\sqrt{2\pi}} \left( (k+\varrho)^{-1/4} + k^{-1/4} \right)^2
  \,\ =\,\    \frac{ h(k,\varrho) }{2 \sqrt{2\pi n} }  \Xi(P)   \nonumber
\al
where, for $x\in\mathopen]0,\infty\mathclose[$ and $\varrho\in[0,\infty[$,
\[ 
 h(x,\varrho) &\coloneqq& \sqrt{2x+\varrho}\left( (x+\varrho)^{-1/4} + x^{-1/4} \right)^2
  \,\ =\,\ \left( \left(2-\tfrac{\varrho}{x+\varrho}\right)^{1/4}
   +  \left(2+\tfrac{\varrho}{x}\right)^{1/4}\right)^2
\] 
satisfies
\[
 4\frac{\dd}{\dd x}\sqrt{h(x,\varrho)}
 &=&  \left(2-\tfrac{\varrho}{x+\varrho}\right)^{-3/4}\tfrac{\varrho}{(x+\varrho)^2}
     -  \left(2+\tfrac{\varrho}{x}\right)^{-3/4}  \tfrac{\varrho}{x^2}  \\
 &=& (2x+\varrho)^{-3/4}\left(\tfrac{\varrho}{(x+\varrho)^{5/4}} 
     - \tfrac{\varrho}{x^{5/4}} \right) \,\ \le \,\ 0     
\]
and hence 
\la   \quad                                                                \label{Eq:h(k,rho)_bound}
  \sup_{k\in\N,\varrho\in\{0,1\}} h(k,\varrho) &=& \max_{\varrho\in\{0,1\}} h(1,\varrho)
    \,\ = \,\ h(1,1) \,\ =\,\  \sqrt{3}\left(2^{-1/4}+1\right)^2 \,\ =\,\  5.86974\ldots,
\al
in the second step by comparison with $h(1,0)=4\sqrt{2}=5.656854\ldots$,
so that the inequality chain~\eqref{Eq:Ineq_chain_BE_for_Z-close_to_normal} yields
\la                                                       \label{Eq:c_1_proof_main_theorem}
 \left\| \widetilde{P^{\ast n}} -\mathrm{N} \right\|_{\mathrm{K}} 
  &\le& c^{}_1\frac{\Xi(P)}{\sqrt{n}} \quad 
   \text{ with } \quad c^{}_1 \,\ \coloneqq\,\ \frac{ \text{R.H.S.\eqref{Eq:h(k,rho)_bound}}}{2\sqrt{2\pi}}
   \,\ =\,\ 1.1708\ldots\,. 
\al

2. Theorem~\ref{Thm:Zolotarev's_zeta_1-B-E-Thm} yields~\eqref{Eq:Abstract_Zolotarev_zeta_1-theorem}
with $\Xi(P)\coloneqq\xi\left( \zeta_1(\widetilde{P}-\mathrm{N}) , \zeta_3(\widetilde{P}-\mathrm{N}) \right)$
with $\xi$ satisfying~\eqref{Eq:xi_near_zero}, and hence~\eqref{Eq:c_1_proof_main_theorem} yields
\la              \label{Eq:Berry-Esseen_for_Z-close_to_normal_sharper_when_close}
 \left\| \widetilde{P^{\ast n}} -\mathrm{N} \right\|_{\mathrm{K}}
  &\le& \frac{c^{}_1}{\sqrt{n}} \frac{\zeta_1+\alpha\zeta_3}{(1-\lambda\zeta_3)_+}(\widetilde{P}-\mathrm{N})\,.
\al
Combining~\eqref{Eq:Berry-Esseen_for_Z-close_to_normal_sharper_when_close}
with the Berry-Esseen inequality~\eqref{Eq:Berry-Esseen_inequality}
with $c=\cSh$ from \eqref{Eq:Def_cSH}
yields~\eqref{Eq:Berry-Esseen_for_Z-close_to_normal_sharper}.

By~\eqref{Eq:zeta_vs_varkappa_etc} we have
\[          \label{page:step_2_of_proof_of_thm_2.1}
 \nu_3(\widetilde{P})
   & \le &  6\zeta_3(\widetilde{P}-\mathrm{N})+\nu_3(\mathrm{N})
   \,\ =\,\  6\zeta_3(\widetilde{P}-\mathrm{N})+\beta
\]
with $\beta \coloneqq \frac{4}{\sqrt{2\pi}}= 1.5957\ldots$ from~\eqref{Eq:nu_1_and_nu_3_of_N}.
Hence
$\text{\rm R.H.S.\eqref{Eq:Berry-Esseen_for_Z-close_to_normal_sharper}}
\le \text{\rm R.H.S.\eqref{Eq:Berry-Esseen_for_Z-close_to_normal}}$ with
\[
 c &\coloneqq& \sup_{\varkappa,\zeta>0}\frac{1}{\varkappa\!\vee\!\zeta}\left(
  \Big(c^{}_1\frac{\varkappa+\alpha\zeta}{1-\lambda\zeta}\Big)_{\!+} \wedge    \big(\cSh \left(6\zeta+\beta\right)\big)  \right)
  \,\ =\,\  \sup_{\zeta>0} \left(  \Big(c^{}_1\frac{1+\alpha}{1-\lambda\zeta}\Big)_{\!+} \wedge
      \Big(\cSh \, \big(6+\frac{\beta}{\zeta}\big)\Big) \right).
\]
The last supremum above is uniquely attained at the
positive solution $\zeta^\ast$,
which by monotonicity obviously exists uniquely
and is automatically $<$~$\frac1\lambda$,
of the quadratic equation
\[
  c^{}_1\frac{1+\alpha}{1-\lambda\zeta} &=& \cSh \, \big(6+\frac{\beta}{\zeta}\big)
\]
for $\zeta$. 
With  $\omega = \frac{c_1}{\cSh}$ we get
\la \quad                                               \label{Eq:zeta^ast}
 \zeta^\ast &=&  -\frac{\omega\left(1+\alpha\right)+\beta\lambda-6}{12\lambda}
 +\sqrt{ \left(\frac{\omega\left(1+\alpha\right)+\beta\lambda-6}{12\lambda} \right)^2 +\frac{\beta}{6\lambda} }
\,\ = \,\ 0.171989\ldots\,, \\  
  c &=& \cSh \,\big(6+\frac{\beta}{\zeta^\ast}\big) \label{Eq:c_approx_7.2_in_main_thm}
    \,\ = \,\  7.16553\ldots\,. 
\al
This proves Remark~\ref{Rem:BE_K_Z_sharper} up to ``from~\eqref{Eq:c_approx_7.2_in_main_thm}'', and hence
in particular Theorem~\ref{Thm:BE_K_Z}.

The remaining claims of Remark~\ref{Rem:BE_K_Z_sharper} can be checked easily.
\end{proof}

One may of course improve upon the above value of $c=7.16553\ldots$ a tiny bit by restricting first attention to $n$ even,
replacing in this case $c_1$ from~\eqref{Eq:c_1_proof_main_theorem} by 
 $c_{1,\text{even}}\coloneqq h(1,0)/(2\sqrt{2\pi})=1.12837\ldots$ and correspondingly  getting
$\zeta^\ast_{\text{even}}= 0.174306\ldots <\frac{1}{\lambda}$, hence
$c_{\text{even}} =  \cSh \,\big(6+\frac{\beta}{\zeta^\ast_{\text{even}}}\big)=7.10768\ldots$,
then using for small odd $n=2k+1$ with $1\le k\le k_0$ just, for example,
\[
 \left\| \widetilde{P^{\ast (2k+1)}} -\mathrm{N} \right\|_{\mathrm{K}}
  &\le&  \text{the third term in inequality chain~\eqref{Eq:Ineq_chain_BE_for_Z-close_to_normal} 
             with $\varrho=1$} \\
  &=& \frac{1}{2\sqrt{2\pi}}
   \left(
     \sqrt{\tfrac{1}{\sqrt{k+1}}\,\zeta_1(\widetilde{P}_{}^{\ast k}-\mathrm{N}^{\ast k}_{})}
    +\sqrt{\tfrac{1}{\sqrt{k}}\,  \zeta_1(\widetilde{P}_{}^{\ast(k+1)}-\mathrm{N}^{\ast(k+1)}_{})} 
   \,\right)^2\nonumber  \\
  &\le& \frac{1}{2\sqrt{2\pi}}
   \left(\sqrt{\tfrac{k}{\sqrt{k+1}}} + \sqrt{\tfrac{k+1}{\sqrt{k}}}\right)^2\zeta_1(\widetilde{P}-\mathrm{N})  \\
  &\le& \frac{1}{2\sqrt{2\pi}}
   \left(\sqrt{\tfrac{k_0}{\sqrt{k_0+1}}} + \sqrt{\tfrac{k_0+1}{\sqrt{k_0}}}\right)^2\zeta_1(\widetilde{P}-\mathrm{N}) 
   \vee \zeta_3(\widetilde{P}-\mathrm{N})
\]
by~\eqref{Eq:semiadditivity} below with $\|\cdot\|=\zeta_1$ in the third step,
and for the remaining odd $n$ the modification of~\eqref{Eq:h(k,rho)_bound} obtained by adding the condition $k>k_0$ 
in the supremum. But trying to optimise such an approach does not appear to be worthwhile with the present still
rather high value of~$c_{\text{even}}$.

\section{Proof of the convolution inequality Theorem~\ref{Thm:F_star_G_vs_H_star_H}}
                                                     \label{sec:proof_of_convolution_inequality}

In the proof of the basic Lemma~\ref{Lem:cK_and_its_extreme_points}, we roughly speaking 
use in~\eqref{Eq:Use_of_lambda_1} and~\eqref{Eq:Use_of_lambda_2}
the ``mean $\mu(Q-P)$'' for certain $P,Q\in\Prob(\R)$, although 
we could perhaps have,
say, $\int x_+\,\dd(Q-P)_+(x)= \int x_+\,\dd(Q-P)_-(x)=\infty$, and then 
$\mu(Q-P)= \int x\,\dd(Q-P)(x)$ were actually undefined. 
Hence we define here an appropriate extension $\lambda$ of $\mu$, namely  
the special case of $\lambda_r$ with $r=1$ in the following Lemma~\ref{Lem:Generalised_signed_moments},
which is given here in a generality reusable in section~\ref{sec:varkappa_distances}.
We recall the notation $\cM,F_M,\nu_r,\cM_r,h_M,\varkappa_r$ from 
(\ref{Eq:Def_cM},\ref{Eq:Def_F_M_and_overline{F}_M},%
\ref{Eq:Def_nu_r},\ref{Eq:Def_cM_r},\ref{Eq:Def_h_M},\ref{Eq:Def_varkappa});
in particular $\cM$ denotes the vector space of all bounded signed measures on the 
Borel-$\sigma$-algebra on~$\R$.

\begin{Lem}[Generalised signed moments]                     \label{Lem:Generalised_signed_moments}
Let $r\in\mathopen]0,\infty\mathclose[$ and
\[
 \cM_{\varkappa_r} &\coloneqq& \left\{ M\in\cM : \varkappa_r(M)<\infty\right\}\,.
\]
Then $\cM_{\varkappa_r}$ is a vector space with $\cM_r \subsetneq \cM_{\varkappa_r}\subsetneq \cM$.
For $M\in\cM_{\varkappa_r}$ we have
\la                                \label{Eq:Def_lambda_etc}
 \lambda_r(M) &\coloneqq & \int rx^{r-1} h_M(x) \,\dd x \\ 
  &=& 
   \left\{\begin{array}{ll}   \displaystyle
    - \int  r|x|^{r-1} F_M(x)\,\dd x &  \quad\text{ if } 
                                                      M(\R)=0, \\
    \displaystyle \int \sgn(x)|x|^r\dd M(x) & \quad \text{ if  } 
   \nu_r(M)<\infty
   \end{array}\right\} ,  \nonumber \\
  |\lambda_r(M)| &\le& \varkappa_r(M) \,\ \le\,\ \nu_r(M)\,,\\
  \lambda_r(M)&=&                   \label{Eq:lambda_r=mu_r_if_r_odd}
         \mu_r(M) \quad\text{ if } r\in\N\text{ is odd and } \nu_r(M) < \infty\,,
\al
and $\lambda_r$ thus defined is a linear functional on $\cM_{\varkappa_r}$. 

With $\lambda \coloneqq \lambda_1$ we have for $P,Q\in\Prob(\R)$ in particular 
\la                        \label{Eq:lambda_Q-P_via_F-G}
 \lambda(Q-P) &=& \int (F_P-F_Q)\,\dd\leb \quad
  \text{ if }F_P-F_Q \in \mathrm{L}^1(\R), 
\al
and $\lambda(Q-P)=\mu(Q-P)$ if $P,Q\in\Prob(\R)$ with $\int |x|\,\dd|Q-P|(x)<\infty$. 
\end{Lem}
\begin{proof}[Proof of Lemma~\ref{Lem:Generalised_signed_moments}
  and of claims in~\eqref{Eq:Def_varkappa} and~\eqref{Eq:|mu|_le_varkappa_le_nu}]
Integrating 
\la                 \label{Eq:|y|^r=int...}
 |y|^r &=& \int_0^yr|x|^{r-1}\sgn(x)\,\dd x
   \,\ =\,\ \int_\R r|x|^{r-1} \big( (0<x<y) + (y<x<0) \big) \,\dd x
\al
w.r.t.~$|M|$ and applying Fubini, justified by positivity, yields the first equality in~\eqref{Eq:Def_varkappa},
for arbitrary $M\in\cM$.
Since the map $\cM\ni M\mapsto h_M$ is linear with $|h_M|\le |h_{|M|}|$, we indeed
have the inequality in~\eqref{Eq:Def_varkappa}, and $\cM_{\varkappa_r}$ is a vector subspace of $\cM$
with $\cM_r\subseteq\cM_{\varkappa_r}$.
We have, for example, $\sum_{j\in\N}j^{-r-1}\delta_j\in\cM\setminus\cM_{\varkappa_r}$
and $\sum_{j\in\N}j^{-r-1}(\delta_{2j}-\delta_{2j-1})\in\cM_{\varkappa_r}\!\!\setminus\cM_r$,
and hence the stated strict inclusions hold.

Obviously, $\lambda_r$ is well-defined and linear, 
the first alternative representation in~\eqref{Eq:Def_lambda_etc} holds, 
the second follows from integrating the analogue of~\eqref{Eq:|y|^r=int...} for $\sgn(y)|y|^r$,
and the remaining claims follow.
\end{proof}

In this section, we also abbreviate $\varkappa \coloneqq \varkappa_1$, 
so that we have in particular
\la                     \label{Eq:varkappa_via_F-G}
 \varkappa(P-Q) &=& \left\| F-G \right\|_1 \,\ =\,\ \int_\R |F-G|\,\dd\leb 
 \,\ \in\,\ [0,\infty]
\al
for $P,Q\in\Prob(\R)$ with distribution functions $F,G$, 
by  \eqref{Eq:varkappa_1=L^1-enorm} applied to $M\coloneqq P-Q$.

Theorem~\ref{Thm:F_star_G_vs_H_star_H} is proved below using the Bauer maximum principle, 
combined with a simple stochastic ordering argument. To this end, 
let, in Lemmas~\ref{Lem:kappa-balls_tight} and~\ref{Lem:cK_and_its_extreme_points},
and in the proof of Theorem~\ref{Thm:F_star_G_vs_H_star_H}, 
the space $\cM$
be equipped with the (probabilist's) 
weak topology of convergence of integrals of bounded continuous functions, 
so that $\cM$ becomes a Hausdorff locally convex vector space.
Let further~$\stle$ denote the usual stochastic order 
on $\Prob(\R)$, so, for $P,Q \in\Prob(\R)$ with distribution functions $F,G$, 
\la                                 \label{Eq:usual_stochastic_order}
   P \stle Q &:\iff& F\ge G\,. 
\al
With this notation, (\ref{Eq:lambda_Q-P_via_F-G},\ref{Eq:varkappa_via_F-G}) then yield
\la                              \label{Eq:kappa(P-Q)_if_stle}
 \quad && \varkappa(P-Q) \ =\ \int_\R (F-G)\,\dd \leb  
   \ =\ \lambda(Q-P) \ \in\ [0,\infty[
   \quad\text{ if $P\stle Q$ and $F-G\in\mathrm{L}^1(\R)$}  \,.
\al

We recall that $q$ is a $u$-quantile of the probability distribution function $H$ if
$H(q-)\le u\le H(q)$ in case of $u\in\mathopen]0,1\mathclose[$,  $q=\inf\{y\in\R:H(y)>0\}$ if $u=0$,
and $q=\sup\{y\in\R: H(y)<1\}$ if $u=1$.

\begin{Lem}                                   \label{Lem:kappa-balls_tight}
Let $H$ be the distribution function of a law $R\in\Prob(\R)$.

\begin{parts}
 \item Let $F$ be a further probability distribution function on $\R$. Then we have
 \la                                                                                   \label{Eq:(F(x)-u)(q-x)}
  (F(x)-u)(q-x)&\le&\left\|F-H\right\|_1  \quad\text{ for $x\in\R$, $u\in[0,1]$, $q$ any $u$-quantile of $H$}.
 \al
 \item Let $\epsilon\in[0,\infty[$. Then the (possibly degenerate) Kantorovich ball \label{part:kappa-balls_tight}
 $B\coloneqq\{P\in\Prob(\R): \varkappa(P,R) \le \epsilon\}$ is weakly compact in $\cM$.
\end{parts}
\end{Lem}
\begin{proof} (a) If $x\le q$, then (even if $q=\infty$, with the usual conventions of measure 
theory) we have 
\[
 \text{L.H.S.\eqref{Eq:(F(x)-u)(q-x)}} &=&\int\limits_{[x,q[}(F(x)-u)\,\dd y
   \,\ \le \,\ \int\limits_{[x,q[}|F(y)-H(y)|\,\dd y \,\ \le \,\  \text{R.H.S.\eqref{Eq:(F(x)-u)(q-x)}}
\]
by using in the second step $F(x)\le F(y)$ for $x\le y$, and $H(y)\le u$ for $y<q$.
If $x\ge q$, then analogously 
$ \text{L.H.S.\eqref{Eq:(F(x)-u)(q-x)}}  = \int_{]q,x]}(u-F(x))\,\dd y
\le  \int_{]q,x]}|H(y)-F(y)|\,\dd y \le \text{R.H.S.\eqref{Eq:(F(x)-u)(q-x)}}$
by $F(y)\le F(x)$ for $y\le x$, and $H(y)\ge u$ for $y\ge q$.

\smallskip(b) Given $\delta>0$, let $u_1\coloneqq\frac{\delta}{4}$, $u_2\coloneqq1-\frac{\delta}{4}$,
and let $q_i$ be a $u_i$-quantile of $H$ and $x_i$ be chosen such that $x_1\le x_2$,
$x_1<q_1$, $q_2< x_2$, and $\frac{\epsilon}{|q_i-x_i|}<\frac{\delta}{4}$ for $i\in\{1,2\}$.
For $P\in B$ with distribution function $F$ we then get
\[
 P(\R\!\setminus\mathopen]x_1,x_2]) &=& F(x_1)+1-F(x_2)
  \,\ \le \,\ u_1 + \frac{\epsilon}{q_1-x_1} +1-u_2 +\frac{\epsilon}{x_2-q_2} \,\ <\,\ \delta
\]
by using~\eqref{Eq:(F(x)-u)(q-x)} twice in the second step. Hence $B$ is uniformly tight in $\Prob(\R)$, 
in the usual sense that \citet[p.~293, Theorem~9.3.3]{Dudley2003} applies.
If $P_\sdot$ is any sequence in $B$ converging weakly to  $P\in \Prob(\R)$, 
with corresponding distribution functions $F_n$ and $F$, 
then Fatou's Lemma yields 
$\int |F(x)-H(x)|\,\dd x \le \liminf\limits_{n\rightarrow\infty}\int|F_n(x)-H(x)|\,\dd x \le \epsilon$,
and hence $P\in B$. Hence, with respect to weak convergence, $B$ is compact in $\Prob(\R)$ and, equivalently, in $\cM$. 
\end{proof}

For the remainder of this section we recall that an element $x$ of a subset $\cK$ of
a vector space over $\R$ is called an \textit{extreme} point of $\cK$ if
$x=\frac{1}{2}(y+z)$ with $y,z\in\cK$ implies $x=y=z$.

\begin{Lem}                                             \label{Lem:cK_and_its_extreme_points} 
Let $R\in\Prob(\R)$ with a continuous distribution function,
and let $\epsilon\in\mathopen[0,\infty\mathclose[$. Then 
\la                                                                \label{Eq:Def_cK_epsilon}
 \cK \,\coloneqq \,\ \cK^{}_{R,\epsilon} &\coloneqq& \{ P\in\Prob(\R) :
   P \stle R, \,  \varkappa(P-R)\le\epsilon  \}                  \\
   &=& \{ P\in\Prob(\R) :    P \stle R, \,  \lambda(R-P)\le \epsilon\}  \nonumber
\al
is a convex and weakly compact subset of $\cM$.
A law $P\in\Prob(\R)$ 
is an extreme point of $\cK$ 
iff there exists a countable (possibly finite, possibly even empty) pairwise
disjoint family $(\,\mathopen[a_i,b_i\mathclose[ : i\in I)$ of nonempty half-open intervals 
such that with 
\la                                                     \label{Eq:Def_p_i}
  p_i &\coloneqq& R(\,\mathopen[a_i,b_i\mathclose[\,) \quad\text{ for }i\in I
\al  
we have 
\la                                                     \label{Eq:Decomp_P}
  P &=& R + \sum_{i\in I} S_i
\al
with the family $(S_i:i\in I)$ in $\cM$ satisfying
\la
 && S_i \,\ =\,\ p_i\delta_{a_i} -R(\cdot\cap\mathopen[a_i,b_i\mathclose[)
   \quad\text{ for every $i\in I$ with at most one exception},   \label{Eq:Def_S_i_normal_case}  \\
 && S_i \,\ =\,\  s\delta_{a_i} +(p_i-s)\delta_\xi -R(\cdot\cap\mathopen[a_i,b_i\mathclose[)
  \quad\text{ for some $\xi\in\mathopen]a_i,b_i\mathclose[$, 
  $s\in\mathopen]R(\mathopen[a_i,\xi\mathclose[), p_i\mathclose[$}  \label{Eq:Def_S_i_exceptional_case} \\  \nonumber
 &&    \phantom{S_i \,\ =\,\   q\delta_{a_i} +(p_i-q)\delta_\xi -R(\cdot\cap\mathopen[a_i,b_i\mathclose[)\quad}
  \text{ if $i\in I$ is  exceptional  in \eqref{Eq:Def_S_i_normal_case}},    \\
 && \varkappa(P-R) \,\ \le \,\ \epsilon,  \label{Eq:kappa(P-R)_le_epsilon} \\
 && \varkappa(P-R) \,\ = \,\ \epsilon \quad\text{ if an exception actually 
       occurs in \eqref{Eq:Def_S_i_normal_case}}.       \label{Eq:kappa(P-R)=epsilon}
\al
\end{Lem}

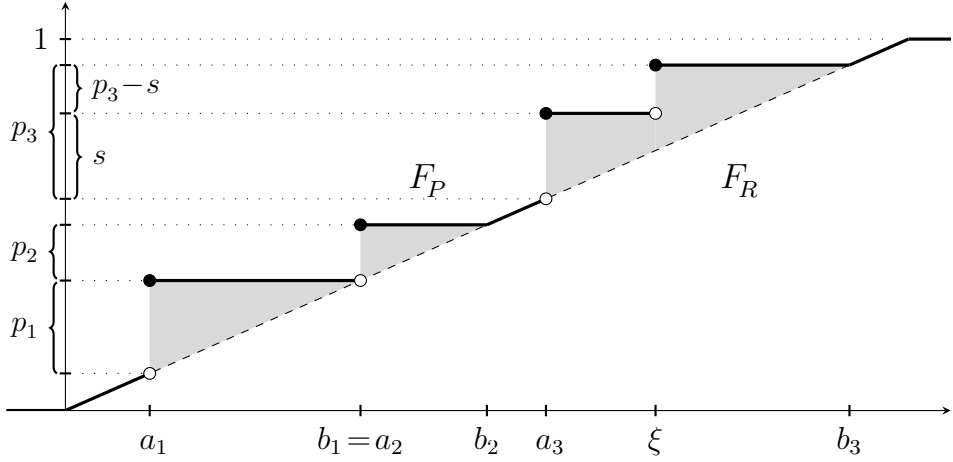
\medskip
\begin{center}

\begin{tikzpicture}

\begin{axis}[xmin = -0.7, xmax = 10.5, ymin = 0, ymax = 1.1, 
xtick = {1,3.5,5,5.7,7,9.3}, 
xticklabels = {\color{white} b \color{black} \; $a_1$ \color{white} b \color{black}, $b_1 \! = \! a_2$, $b_2$, \color{white} b \color{black} \; $a_3$  \color{white} b \color{black},
$\xi$, $b_3$}, 
ytick = {0.1, 0.35, 0.5, 0.57, 0.8, 0.93, 1}, 
yticklabels = { , , , , , , 1}, 
axis lines=middle, 
x tick style={black,thick},
y tick style={black,thick},
height = 7cm, width = 14cm] 

\addplot[domain=0:10, color = black, dashed, name path = A]{0.1*x}; 
\addplot[domain=0:1, color = black, very thick]{0.1*x};
\addplot[color=black, very thick] coordinates {
		(-0.7,0)
		(0,0) 
	};

\addplot[color=black,name path = B, very thick] coordinates {
		(1,0.35)
		(3.5,0.35) 
	};

\addplot[color=black,name path = C, very thick] coordinates {
		(3.5,0.5) 
		(5,0.5) 
	};
	
\addplot[domain=5:5.7, color = black, very thick]{0.1*x};	

\addplot[color=black,name path = D, very thick] coordinates {
		(5.7,0.8) 
		(7,0.8) 
	};
\addplot[color=black,name path = E, very thick] coordinates {
		(7,0.93) 
		(9.3,0.93) 
	};
\addplot[domain=9.3:10, color = black, very thick]{0.1*x};
	
\addplot[black!15] fill between[of=A and B, soft clip={domain=1:3.5}];
\addplot[black!15] fill between[of=A and C, soft clip={domain=3.5:5}];
\addplot[black!15] fill between[of=A and D, soft clip={domain=5.7:7}];
\addplot[black!15] fill between[of=A and E, soft clip={domain=7:9.3}];
%

\draw[-, very thick] (10,1) -- (11,1);

\node at (axis cs:4.3,0.62)[color = black, scale = 1.1]{$F_{\!P}^{}$};
\node at (axis cs:8,0.62)[color = black, scale = 1.1]{$F_{\!R}^{}$};

\draw [decorate, decoration={brace}, color = black, thick](-0.1,0.1) -- (-0.1,0.345) node [midway,left, xshift = -0.1cm, scale = 0.95] {$p^{}_1$};
\draw [decorate, decoration={brace}, color = black, thick](-0.1,0.355) -- (-0.1,0.5) node [midway,left, xshift = -0.1cm, scale = 0.95] {$p^{}_2$};
\draw [decorate, decoration={brace}, color = black, thick](-0.1,0.57) -- (-0.1,0.93) node [midway,left, xshift = -0.1cm, scale = 0.95] {$p^{}_3$};
\draw [decorate, decoration={brace, mirror}, color = black, thick](0.1,0.57) -- (0.1,0.795) node [midway,right, xshift = 0.1cm, scale = 0.95] {$s$};
\draw [decorate, decoration={brace, mirror}, color = black, thick](0.1,0.805) -- (0.1,0.93) node [midway,right, xshift = 0.1cm, scale = 0.95] {$p^{}_3\!-\!s$};


\addplot[color=black, loosely dotted, color = black!100] coordinates  {
		(0.2,0.1)
		(1,0.1) 
	};
\addplot[color=black, loosely dotted, color = black!100] coordinates {
		(0.2,0.35)
		(1,0.35) 
	};	
\addplot[color=black, loosely dotted, color = black!100] coordinates {
		(0.2,0.5)
		(3.5,0.5) 
	};
\addplot[color=black, loosely dotted, color = black!100] coordinates {
		(0.2,0.57)
		(5.7,0.57) 
	};
\addplot[color=black, loosely dotted, color = black!100] coordinates {
		(0.2,0.8)
		(5.7,0.8) 
	};
\addplot[color=black, loosely dotted, color = black!100] coordinates {
		(0.2,0.93)
		(7,0.93) 
	};
	\addplot[color=black, loosely dotted, color = black!100] coordinates {
		(0.2,1)
		(9.8,1) 
	};

\node at (axis cs:1,0.1)[draw = black, fill = white, circle, scale = 0.4](above){};
\node at (axis cs:1,0.35)[draw = black, fill = black, circle, scale = 0.4](above){};
\node at (axis cs:3.5,0.35)[draw = black, fill = white, circle, scale = 0.4](above){};
\node at (axis cs:3.5,0.5)[draw = black, fill = black, circle, scale = 0.4](above){};
\node at (axis cs:5.7,0.57)[draw = black, fill = white, circle, scale = 0.4](above){};
\node at (axis cs:5.7,0.8)[draw = black, fill = black, circle, scale = 0.4](above){};
\node at (axis cs:7,0.8)[draw = black, fill = white, circle, scale = 0.4](above){};
\node at (axis cs:7,0.93)[draw = black, fill = black, circle, scale = 0.4](above){};

\end{axis}
\end{tikzpicture}
\end{center}
\centerline{An extremal $P$, with $I=\{1,2,3\}$, $i\!=\!3$ exceptional, grey area $=$ $\epsilon$,
 $F^{}_{\!R}$ dashed, $F^{}_{\!P}$ solid.}
\medskip

\begin{proof} Let $H\coloneqq F_R$.  
We will use the continuity assumption on $H$ in this proof only 
for the more important ``only if'' part, 
when deriving~\eqref{Delta(x-)_neq_0} below.
Let $\cE$ denote the set of extreme points of $\cK$, 
and let in this proof $\cP$ denote the set of all $P\in\Prob(\R)$ as described after 
the ``iff'' in the claim (with no relation to the notation~\eqref{Eq:Def_cP_r}).

\smallskip 
1. The alternative representation of $\cK$ follows 
from~(\ref{Eq:varkappa_via_F-G},\ref{Eq:kappa(P-Q)_if_stle}).
Since $\{P\in\Prob(\R): P\stle R\}$ is a convex and weakly closed subset of $\cM$,
and since $B$ from Lemma~\ref{Lem:kappa-balls_tight}\ref{part:kappa-balls_tight} is convex and weakly compact,
the stated convexity and compactness of $\cK$ follows.

\smallskip
2. $\cP \subseteq \cE$: Let $P\in\cP$, with (\ref{Eq:Def_p_i}--\ref{Eq:kappa(P-R)=epsilon}), and let $F\coloneqq F_P$.
Then (\ref{Eq:Def_p_i},\ref{Eq:Def_S_i_normal_case},\ref{Eq:Def_S_i_exceptional_case}) yield 
\la                                                                          \label{Eq:F_S_i_ge_or_=0}
  \quad&& F_{S_i} \ge 0 \text{ everywhere}, \quad F_{S_i}(x-)=S_i(\mathopen]-\infty,x\mathclose[)  = 0  
   \text{ for }x\in \R\!\setminus\!\mathopen]a_i,b_i\mathclose[\,,
   \quad\text{ for every $i\in I$}.  
\al
With~\eqref{Eq:Decomp_P} we hence get $F=F_R+\sum_{i\in I} F_{S_i} \ge H$ and thus $P\stle R$, 
and~\eqref{Eq:kappa(P-R)_le_epsilon} then yields 
$P\in\cK$.
To continue, let us write
\[
  I_0 &\coloneqq& \{i \in I : i \text{ nonexceptional in \eqref{Eq:Def_S_i_normal_case}}\}, \\
  U   &\coloneqq& \bigcup_{i\in I}\, \mathopen]a_i,b_i\mathclose[\,, \, \qquad A \,\ \coloneqq \R\setminus U,  \qquad
       U_0 \,\ \coloneqq\,\  \bigcup_{i\in I_0}]a_i,b_i[\,, \qquad A_0 \,\ \coloneqq\,\ \R\setminus U_0\,.        
\]
By \eqref{Eq:Decomp_P} and \eqref{Eq:F_S_i_ge_or_=0} we then have 
\la                                         \label{Eq:F=H_on_...}
 F(x-) &=& H(x-)+\sum_{i\in I} F_{S_i}(x-)   \,\ = \,\ H(x-)   \quad \text{ for }x\in A\,.
\al

Let now $P_0,P_1\in\cK$ with $P=\frac12(P_0+P_1)$, and with corresponding distribution functions $F_0,F_1$.
Then $F,F_0,F_1\ge H$ everywhere, and $F=\frac12(F_0+F_1)$, and with~\eqref{Eq:F=H_on_...} we obtain
\la                                                     \label{Eq:f_0=F_F_0=F_on_A-}
 F_0(x-)=F_1(x-)=F(x-) \quad\text{ for }x\in A\,.
\al

If $i\in I$, then we have $a_i,b_i \in A$ by the disjointness assumption,
and for ${t}\in\{0,1\}$ then
\la                                                    \label{Eq:P_j(a_i,b_i[)}
       P_{t}([a_i,b_i[)=F_{t}(b_i-)-F_{t}(a_i-)=F(b_i-)-F(a_i-)=p_i\,.
\al
If in addition $i$ is nonexceptional in~\eqref{Eq:Def_S_i_normal_case}, then 
$P = p_i\delta_{a_i}$ on $[a_i, b_i[$, hence $P=\frac12(P_0+P_1)$ yields
$P_{t}=c_{t}\delta_{a_i}$ on $[a_i, b_i[$ for ${t}\in\{0,1\}$
with 
some $c_{t}\in[0,1]$, but then $c_0=c_1=p_i$ by~\eqref{Eq:P_j(a_i,b_i[)}, 
and thus $P_0=P_1=P$ on $[a_i,b_i[$ in the present case.
Thus if $x\in U_0$, so $x\in \mathopen]a_i,b_i\mathclose[$ for some $i\in I_0$,
then for ${t}\in\{0,1\}$
\[
 F_{t}(x-) &=& 
   P_{t}([a_i,x[) +F_{t}(a_i-)  
  \,\ = \,\ P([a_i,x[) + F(a_i-) \,\ =\,\ F(x-)\,,
\]
using in the central step also~\eqref{Eq:f_0=F_F_0=F_on_A-} with there $x=a_i$.
So up to now we have obtained
\la                                                                 \label{Eq:OK_on_A_cup_U_0}
 F_0(x-) &=& F_1(x-)\,\  =\,\ F(x-) \quad\text{ for }x\in A\cup U_0\,.
\al
Hence if $I = I_0$, then we get  $P_0=P_1=P$ as desired.

If finally $i \in I\setminus I_0$ and ${t}\in\{0,1\}$, then~\eqref{Eq:OK_on_A_cup_U_0} yields  
$P_{t} = P$ on $\R\setminus[a_i,b_i[\,$, and $\frac{1}{2}(P_0+P_1)=P$ now yields
$P_{t}=c_{t}\delta_{a_i} + d_{t}\delta_\xi$ 
on $[a_i, b_i[$ with 
$\frac{c_0+c_1}{2}=s$, $\frac{d_0+d_1}{2}= p_i-s$, and, by~\eqref{Eq:P_j(a_i,b_i[)}, $c_{t}+d_{t}=p_i\,$, 
and hence we get, recalling Lemma~\ref{Lem:Generalised_signed_moments} and 
using in the third step 
in particular~\eqref{Eq:kappa(P-R)=epsilon},
\la                                   \label{Eq:Use_of_lambda_1}
 \epsilon & \ge &\lambda(R-P_{t}) \,\ =\,\ \lambda(R-P) + \mu(P-P_{t})
  \,\ =\,\ \epsilon + (s-c_{t})a_i+ (p_i-s-d_{t})\xi \\ \nonumber
  &=&  \epsilon + (c_{t}-s)(\xi-a_i) \,.
\al
Hence $\frac{c_0+c_1}{2}=s$ and $\xi>a_i$ yield 
$c_0=c_1=s$, 
and hence again $P_0=P_1=P$.

\smallskip 3. $\cE\subseteq\cP$: 
Let $P\in\cE$, with distribution function $F$. Then $\Delta \coloneqq F-H\ge 0$. 
Using the rightcontinuity of $\Delta$, the leftcontinuity of 
$x\mapsto\Delta(x-)$, and $\Delta(x)=\Delta((x+)-)$
yield the openness, and hence the representation, of 
\[
 U &\coloneqq& \{ x\in \R : \Delta(x-)>0\text{ and }\Delta(x)>0\} 
  \,\ =\,\  \bigcup_{i\in I}\, \mathopen]a_i,b_i\mathclose[
\]
with some countable (possibly finite or even empty) pairwise disjoint family $(\,\mathopen]a_i,b_i\mathclose[:i\in I)$ 
of open 
intervals 
with $a_i<b_i$. We have 
\la                                \label{Eq:a_i_and_b_i_notin_U_i}
   a_i,b_i &\notin& U \quad\text{ for }i\in I,
\al
and in particular the corresponding family $(\mathopen[a_i,b_i\mathclose[ : i\in I)$
of  half-open intervals is also pairwise disjoint
(though we might have $b_i=a_j$ for some $i,j\in I$).
We further have the implication 
\la                              \label{Delta(x)_neq_0}
  \Delta(x) \neq 0 \,&\Rightarrow&\, x \in  \bigcup_{i\in I}\, \mathopen[a_i,b_i\mathclose[ \eqqcolon V
\al
just by $\Delta\ge 0$, the definition of $U$, and rightcontinuity of $\Delta$,  and even
\la                                \label{Delta(x-)_neq_0}
  \Delta(x-) \neq 0 \,&\Rightarrow&\, x \in U
\al
since leftcontinuity yields the conclusion with 
$\bigcup_{i\in I}\, \mathopen]a_i,b_i\mathclose]$ in place of $U$, 
and for each $i\in I$ 
we have $0\le \Delta(b_i-)= F(b_i-)-H(b_i)\le \Delta(b_i)$
by the continuity of $H$, and hence $\Delta(b_i-)=0$ due to~\eqref{Eq:a_i_and_b_i_notin_U_i}. 

For $i\in I$ we define $p_i$ by~\eqref{Eq:Def_p_i}, and
\[
   S_i &\coloneqq& \big(P-R\big)(\cdot\cap\mathopen[a_i,b_i\mathclose[\,)\,,
\]
and then obtain, using~(\ref{Eq:a_i_and_b_i_notin_U_i},\ref{Delta(x-)_neq_0})
in the second step below,
\la
 S_i(\mathopen[a_i,b_i\mathclose[\,)              \label{Eq:S_i([a_i,b_i[)=0}
   &=& \Delta(b_i-)-\Delta(a_i-)\,\ =\,\ 0\,, \\
 P(\mathopen[a_i,b_i\mathclose[\,)                \label{Eq:P([a_i,b_i[)=p_i}
  &=& R(\mathopen[a_i,b_i\mathclose[\,) 
   + S_i(\mathopen[a_i,b_i\mathclose[\,) \,\ =\,\ p_i\,.
\al
For $x\in\R$ we then get
\[
 F^{}_{\sum\limits_{i\in I}S_i}(x) &=& \sum\limits_{i\in I}F^{}_{S_i}(x)
 \,\ =\,\ \sum\limits_{i\in I}
 \left\{\begin{array}{ll} 0  &\text{ if } x\notin\mathopen[a_i,b_i\mathclose[ \\
     \Delta(x)-\Delta(a_i-)  &\text{ if } x\in\mathopen[a_i,b_i\mathclose[  \end{array} \right\} \\
 &=&  \left\{\begin{array}{ll} 0  &\text{ if } x\notin V \\
     \Delta(x)  &\text{ if } x\in V  \end{array} \right\} 
  \,\ =\,\ F(x)-H(x)
\]
using~\eqref{Eq:S_i([a_i,b_i[)=0} for $x\ge b_i$ in the second step,
$\Delta(a_i-)=0$  by~(\ref{Eq:a_i_and_b_i_notin_U_i},\ref{Delta(x-)_neq_0})
and the disjointness of the $[a_i,b_i[$ in the third,
and~\eqref{Delta(x)_neq_0} in the final fourth step.
Hence~\eqref{Eq:Decomp_P} holds.

Let $i\in I$ be fixed in this and in the next two  paragraphs. We have 
\la                                                    \label{Eq:P([a_i,xi[)>0}
 P([a_i,\xi[) &>& 0 \quad\text{ for every }\xi\ >\ a_i\,,
\al
for else we would have $P([a_i,\xi[)=0$ for some $\xi \in \mathopen]a_i,b_i\mathclose[$, 
and then $F(a_i)=F(a_i-) = F(\xi-)>H(\xi-)\ge H(a_i) \ge H(a_i-)$, in contradiction 
to $a_i \notin U$. 

If there exist $\xi,\eta$ with $a_i<\xi<\eta<b_i$ and 
$P([a_i,\xi[\,) P([\xi,\eta[\,)P([\eta,b_i[\,) >0 $,
then we can find $\widetilde{a_i}\in[a_i,\xi[$
and $\widetilde{b_i}\in\mathopen[\eta,b_i\mathclose[$ with 
\[
 \alpha &\coloneqq& P(\,[\widetilde{a_i},\xi[\,)\,\ >\,\ 0\,,\quad
    \beta \,\ \coloneqq\,\ P(\,[\xi,\eta[\,)\,\ > \,\ 0\,, \quad
    \gamma \,\ \coloneqq \,\ P(\,[\eta,\widetilde{b_i}]\,) \,\ > \,\ 0\,, \nonumber \\
 \varrho &\coloneqq& \inf_{x\in[\widetilde{a_i},\widetilde{b_i}]}
   \big( F(x)- H(x)\big) \,\ > \,\ 0\,,
\]
taking $\widetilde{a_i}\coloneqq a_i$ in case of $P(\,]a_i,\xi[\,)=0$, for then
$F(a_i)=F(\xi-) > H(\xi-) \ge H(a_i)$, and else $\widetilde{a_i}> a_i$ small 
enough to have $P(\,[\widetilde{a_i},\xi[\,)>0$.
Then the three conditional probability measures
\[
 A &\coloneqq& P\bed{\cdot}{[\widetilde{a_i},\xi[\,}\,,  \qquad
 B \,\ \coloneqq\,\ P\bed{\cdot}{[\xi, \eta[\,}\,,  \qquad
 C \,\ \coloneqq\,\ P\bed{\cdot}{ [\eta , \widetilde{b_i}]\,}
\]
have means 
\[
  a &\coloneqq& \mu(A) \,\ <\,\ b\,\ \coloneqq\,\ \mu(B) \,\ <\,\ c \,\ \coloneqq\,\ \mu(C)\,,
\]
and for $t,u\in\R$ with $|t|+|u|\le \alpha\wedge\beta\wedge\gamma\wedge\varrho$, we have
\[
 P_{t,u} &\coloneqq& P+tA+uB-(t+u)C \,\ \in\,\ \Prob(\R)\,,\\
 P_{t,u} &=& P \,\ \text{ iff }\,\ (t,u)=(0,0)\,,\\
 F_{t,u}&\coloneqq& F_{P_{t,u}} \,\ \ge H\,,
\]
so in particular $P_{t,u}\stle R$, 
and, recalling Lemma~\ref{Lem:Generalised_signed_moments},
\la                                            \label{Eq:Use_of_lambda_2}    
 \lambda(R-P_{t,u}) &=& 
   \lambda(R-P)  - \mu(P_{t,u}-P) 
  \,\ =\,\  \lambda(R-P) -ta-ub+(t+u)c \\ \nonumber
  &=& \lambda(R-P)  \quad\text{ if } u=-\frac{c-a}{c-b}t\,,
\al
so that $P_{t,u}\,,P_{-t,-u}  \in\cK\setminus\{P\}$ for some $(t,u)\neq(0,0)$,
and then $P = \frac{1}{2}(P_{t,u}+P_{-t,-u})$,
which is incompatible with the assumption $P\in\cE$.

The above contradiction shows that $P$ is on $[a_i,b_i[$ a measure supported in at most two 
points, and in view of~\eqref{Eq:P([a_i,xi[)>0} we then must have $P(\{a_i\})>0$,
and hence $P(\{\xi\})>0$ for at most one $\xi\in\mathopen]a_i,b_i\mathclose[\,$,
and by~\eqref{Eq:P([a_i,b_i[)=p_i} then
either $P=p_i\delta_{a_i}$ on $[a_i,b_i[\,$,
and then the equality in~\eqref{Eq:Def_S_i_normal_case}  holds for the present $i$,
or 
\la \quad                                            \label{Eq:P_exceptional_i}
  P&=& s_i\delta_{a_i} +(p_i-s_i)\delta_{\xi_i} \text{ on  $[a_i,b_i[\,$, \quad
   with some $\xi_i\in\mathopen]a_i,b_i\mathclose[$ and 
   $s_i\in\mathopen]R(\mathopen[a_i,\xi_i\mathclose[),p_i\mathclose[$}\,,
\al
where the lower bound on $s_i$ results from
\[
 H(\xi_i-) &<& F(\xi_i-) \,\ =\,\ F(a_i-)+ P([a_i,\xi_i[) \,\ =\,\ H(a_i-) + s_i\,,
\]
and then we have the equality in~\eqref{Eq:Def_S_i_exceptional_case}, 
with $\xi=\xi_i$ and $s=s_i$.

If we now had two different indices $j,k\in I$ with $P$ 
satisfying~\eqref{Eq:P_exceptional_i} for $i\in\{j,k\}$, 
then for $t\in\R\setminus\{0\}$ with $|t|$ sufficiently small and for 
$u\coloneqq \frac{a_j-\xi_j}{a_k-\xi_k}t$, the law $P_{t,u}$ defined by 
\[
 P_{t,u} &\coloneqq & \begin{Bmatrix} P  \\
     (s_j+t)\delta_{a_j} +(p_j-s_j-t)\delta_{\xi_j} \\ 
     (s_k-u)\delta_{a_k} +(p_k-s_k+u)\delta_{\xi_k}
     \end{Bmatrix} \text{ on }
     \begin{Bmatrix}
     \R\!\setminus\!\big(\mathopen[a_j,b_j\mathclose[ \cup \mathopen[a_k,b_k\mathclose[ \big) \\
    [a_j,b_j[ \\   [a_k,b_k[
 \end{Bmatrix}
\]
would satisfy $P_{t,u}\stle R$ and $\mu(P_{t,u}- P) 
= t\left(a_j-\xi_j\right)-u\left( a_k - \xi_k\right) =0$, and hence $P_{t,u}\in\cK$, 
but then $P=P_{0,0}=\frac12(P_{t,u} +P_{-t,-u})$ would contradict the assumption $P\in\cE$.
Hence~\eqref{Eq:Def_S_i_normal_case} holds.

If we finally had an index $i\in I$ satisfying the condition 
in~\eqref{Eq:Def_S_i_exceptional_case}, 
that is, \eqref{Eq:P_exceptional_i} with $\xi_i=\xi$ and $ s_i=s$, 
but $\varkappa(P-R)<\epsilon$, then we would have 
\[
 P_t &\coloneqq& \begin{Bmatrix} P  \\
     (s +t)\delta_{a_i} +(p_i-t)\delta_{\xi}
     \end{Bmatrix} \text{ on }
     \begin{Bmatrix}
     \R\setminus\mathopen[a_i,b_i\mathclose[  \\
    [a_i,b_i[ 
 \end{Bmatrix}
 \,\ \in \,\ \cK\setminus\{P\}
\]
for $t \neq 0$ with
$ s-R([a_i,\xi[) \le  t \le (p_i-s) \wedge\frac{\epsilon - \varkappa(P-R)}{\xi-a_i}$,
and then $P=\frac{1}{2}(P_t+P_{-t})$ would contradict $P\in\cE$.
Hence~\eqref{Eq:Def_S_i_exceptional_case} holds.
\end{proof}

\begin{Lem}                                                  \label{Lem:S_i_ast_T_j_bounds}
Let $a,b,c,d\in\R$ with $a\le b$ and $c\le d$, and let $U,V$ be bounded positive Borel measures on $\R$, 
 respectively supported in $[a,b], [c,d]$ with total masses $p,q$, that is, 
 \[
   p \,\ \coloneqq\,\  U([a,b]) \,\ =\,\ U(\R),&& q \,\ \coloneqq \,\ V([c,d])\,\ =\,\ V(\R)\,.
 \]
 Then for $z\in\R$ we have
 \la                                                                 \label{Eq:F_{p_delta_a-U_etc}}
   F^{}_{(p\delta_a-U)\ast(q\delta_c-V)}(z) &\begin{Bmatrix} = 0 \\ \le pq \\ \le 0 \end{Bmatrix} &
    \text{ if } \quad \begin{Bmatrix} z< a+c \,\text{ or }\, z \ge b+d \\
    a+c \le z <  (a+d)\wedge(b+c) \\ 
     (a+d)\wedge(b+c)\le z < b+d \end{Bmatrix}.   
 \al
\end{Lem}
\begin{proof} 
We may normalize to $p=q=1$, for we can write 
$\text{L.H.S.\eqref{Eq:F_{p_delta_a-U_etc}}}= pq F^{}_{(\delta_a-\frac{1}{p}U)\ast(\delta_c-\frac{1}{q}V)}(z)$
if $pq>0$, and in case of $pq=0$ the claim is trivial.
The signed measure $(\delta_a-U)\ast(\delta_c-V)$ then 
has total mass zero and its support contained in $[a+c,b+d]$\,;
hence the first of the three claims in~\eqref{Eq:F_{p_delta_a-U_etc}} is obvious. 
For arbitrary $z\in\R$, we have 
\[
 \text{L.H.S.\eqref{Eq:F_{p_delta_a-U_etc}}} 
  &=& F_{\delta_{a+c}}(z) -  F_{\delta_a\ast V}(z)  +F_{U\ast V}(z) - F_{U\ast \delta_c}(z) \\
  &\le& F^{}_{\delta_{a+c}-\delta_a\ast V}(z)
  \,\ \le \,\ F^{}_{\delta_{a+c}-\delta_{a+d}}(z)
   \,\ \begin{Bmatrix} \le 1 =pq \text{ always}  \\   = 0 \text{ if }z\ge a+d\  \end{Bmatrix},
\]
by $\delta_c\stle V$ in the second step, and $V \stle \delta_d$ in the third.
The above, together with the analogous result for
$(a,b,U)$ and $(c,d,V)$ interchanged, yields~\eqref{Eq:F_{p_delta_a-U_etc}}.
\end{proof}

We recall the notation~(\ref{Eq:Def_Kolmogorov-norm},\ref{Eq:Def_Lip-norm},%
\ref{Eq:Def_M_Lipschitz},\ref{Eq:varkappa_via_F-G}).

\begin{Lem}        \label{Lem:(FG-HH)(z)} \label{Lem:CS_for_extreme_points}
Let $P_1,P_2,R_1,R_2\in\Prob(\R)$ with $P_i\stle R_i$ for $i\in\{1,2\}$.
Then
\la     \label{Eq:Kolmogorov_distance_of_product_of_differences_le}
 \left\| (P_1-R_1)\ast(P_2-R_2)\right\|_\mathrm{K}
  &\le& 
  2\sqrt{\left\|R_1\right\|_\mathrm{L}
    \left\|R_2\right\|_\mathrm{L}\varkappa(P_1-R_1)\varkappa(P_2-R_2)}\,.
\al
\end{Lem}
\begin{proof} Let $H_i$ denote the distribution function of $R_i$ for $i\in\{1,2\}$,
so $\left\|R_i\right\|_\mathrm{L}=\left\|H_i\right\|_\mathrm{L}$.

\smallskip 1. We may assume
the $\left\|R_i\right\|_\mathrm{L}$ and the $\varkappa(P_i-R_i)$
to be finite, since by $\left\|R_i\right\|_\mathrm{L}>0$
we else either have
$\text{R.H.S.\eqref{Eq:Kolmogorov_distance_of_product_of_differences_le}}=\infty$,
or~$\varkappa(P_i-R_i)=0$ for some $i$ and then also
$\text{L.H.S.\eqref{Eq:Kolmogorov_distance_of_product_of_differences_le}}=0$.

\smallskip 2.
For $z\in\R$ we have
 \[
  -F^{}_{(P_1-R_1)\ast(P_2-R_2)}(z)
    &=& \int F^{}_{P_1-R_1}(z-y)\,\dd\big(R_2-P_2)(y) \\
   &\le& \int F^{}_{P_1-R_1}(z-y)\,\dd R_2(y)
   \,\ \le \,\ \varkappa(P_1-R_1)\left\|R_2\right\|_\mathrm{L}
 \]
 by using $F^{}_{P_1-R_1}\ge 0$ in the second step,
 and, say, Lemma~\ref{Lem:||cdot||_K-smoothing} in the third.
 Hence
 \la   \label{Eq:Lower_bound_complementing_Lemma_4.5}
  \sup_{z\in\R}\left( -F^{}_{(P_1-R_1)\ast(P_2-R_2)}(z)\right)
   &\le& \big( \varkappa(P_1-R_1)\left\|R_2\right\|_\mathrm{L}\big)
          \wedge \big( \varkappa(P_2-R_2)\left\|R_1\right\|_\mathrm{L} \big) \\
  &\le& \tfrac{1}{2}\,
   \text{R.H.S.\eqref{Eq:Kolmogorov_distance_of_product_of_differences_le}} \,.\nonumber
\al

3. It remains to prove
\la                                          \label{Eq:CS_for_extreme_points_new}
 \sup_{z\in\R} F^{}_{(P_1-R_1)\ast(P_2-R_2)}(z)
   &\le& \text{R.H.S.\eqref{Eq:Kolmogorov_distance_of_product_of_differences_le}}\,,
\al
which in view of~\eqref{Eq:Lower_bound_complementing_Lemma_4.5}
then yields~\eqref{Eq:Kolmogorov_distance_of_product_of_differences_le}.
For $i\in\{1,2\}$ let $\epsilon_i\coloneqq\varkappa(P_i-R_i)$, so that with
the notation~\eqref{Eq:Def_cK_epsilon} we have $P_i\in\cK_{R_i,\epsilon_i}$.
Hence it suffices to prove that
\la                  \label{Eq:Def_Psi_new}
   \Psi(P_1',P_2',z) &\coloneqq& F^{}_{(P_1'-R_1)\ast(P_2'-R_2)} (z)
     \,\ \le \,\  2\sqrt{ \left\|R_1\right\|_\mathrm{L} \left\|R_2\right\|_\mathrm{L}\epsilon_1\epsilon_2}
     \quad \text{ for }z\in\R
\al
holds for $(P_1',P_2')\in\cK_{R_1,\epsilon_1}\times\cK_{R_2,\epsilon_2}$,
and the condition~\eqref{Eq:Def_Psi_new} remains equivalent if
there $z$ is replaced by $z-$, indicating left hand limits.

Let now $z\in\R$ be fixed. Then 
$\Psi(P_1',P_2',z-)$ is
separately in each of its two variables $P_1',P_2'\in\Prob(\R)$
a function affine-linear, and hence convex,  and
(even jointly, but this is not needed here)
weakly upper semi-continuous.
The latter claim becomes clear by writing
\[
 \Psi(P_1',P_2',z-) &=& F^{}_{P_1'\ast P_2' - P_1'\ast R_2  - R_1\ast P_2' + R_1\ast R_2}   (z-)  \\
  &=&  F^{}_{P_1'\ast P_2'}(z-) - \int H_2(z-x)\,\dd P_1'(x)- \int H_1(z-x)\,\dd P_2'(x)
      + F^{}_{R_1\ast R_2} (z-)
\]
and noting first that $(P_1',P_2')\mapsto F^{}_{P_1'\ast P_2'}(z-)$ is weakly upper semicontinuous
by continuity of convolution in $\Prob(\R)$ and by the
portmanteau theorem applied to the open set $]-\infty,z[\,$,
see for example~\citet[pp.~45--48, Theorem~3.1 and Corollary~3.4]{BCR1984},
and second that the remaining summands are weakly continuous by the
continuity and boundedness of the $H_i$.

Hence, using also the convexity and compactness of the $\cK_i\coloneqq \cK_{R_i,\epsilon_i}$
established in Lemma~\ref{Lem:cK_and_its_extreme_points},
two applications of the \citet[p.~392, Korollar]{Bauer1958} maximum principle,
presented also by \citet[p.~102, Theorem 25.9]{Choquet1969} and by \citet[p.~298]{AliprantisBorder2006},
yield
\[
 \sup_{P_1'\in\cK_1,P_2'\in\cK_2} \Psi(P_1',P_2',z-)
 &=& \sup_{P_1'\in\cK_1,P_2'\in\cE_2} \Psi(P_1',P_2',z-)
 \,\ =\,\  \sup_{P_1'\in\cE_1,P_2'\in\cE_2} \Psi(P_1',P_2',z-)\,.
\]
Therefore it suffices to prove~\eqref{Eq:Def_Psi_new}
in case of each $P_i'$ being
an extreme point of $\cK_{R_i,\epsilon_i}$,
and for such $P_i'$ we have $\varkappa(P_i'-R_i)\le \epsilon_i$.

Using the isotonicity of R.H.S.\eqref{Eq:Def_Psi_new} in the $\epsilon_i$,
we conclude that for proving~\eqref{Eq:CS_for_extreme_points_new} for the given $P_1,P_2$,
we may from now on assume that,
for {\em some} finite $\epsilon_i\ge\varkappa(P_i-R_i)$ not necessarily with equality,
$P_1$ is an extreme point of $\cK_{R_1,\epsilon_1}$
and $P_2$ is an extreme point of $\cK_{R_2,\epsilon_2}$.

\smallskip 4.
Under the assumptions as just stated at the end of the previous step,
let us now change notation from $P_1,P_2$ to $P,Q$, but inconsequentially keep
$R_1,R_2,H_1,H_2$,
in order to reuse in this proof the notation of Lemma~\ref{Lem:cK_and_its_extreme_points}
and to avoid double indices.

Let $S_i$ etc.~be as in Lemma~\ref{Lem:cK_and_its_extreme_points}
applied to $R_1,\epsilon_1, P$. Analogously,  with Lemma~\ref{Lem:cK_and_its_extreme_points} applied to $R_2,\epsilon_2,Q$,
we have $Q-R_2= \sum_{j\in J}T_j$ with $q_j = R_2([c_j,d_j])$ for $j\in J$, 
$T_j= q_j\delta_{c_j} - R_2(\cdot\cap[c_j,d_j])$ for every $j\in J$ with at most one exception,
$T_j=t\delta_{c_j}+(q_j-t)\delta_{\eta}-R_2(\cdot\cap [c_j,d_j])$ for some $\eta\in\mathopen]c_j,d_j\mathclose[$ and
$t \in\mathopen]R(\mathopen[c_j,\eta\mathclose[),q_j\mathclose[$ if $j$ actually is exceptional.

Let also $z\in\R$ be fixed in the rest of this proof.

\smallskip
5. The unexceptional case: Let us assume here that neither an exceptional $i$ in~\eqref{Eq:Def_S_i_normal_case}
nor an analogous exceptional $j$ occurs. We then observe that the set of  pairs
\[
   A &\coloneqq& A_z \, \coloneqq \,\ \{(i,j)\in I\times J : a_i+c_j\le z <
    (a_i+d_j)\wedge(b_i+c_j)\},  
\]
is (the graph of) an injective function $j(\cdot):I_0\rightarrow J$ for some $I_0\subseteq I$,
for if $(i,j_1),(i,j_2)\in A$, then we have
$z-a_i\in\mathopen[c_{j_1},d_{j_1}\mathclose[ \cap \mathopen[c_{j_2},d_{j_2}\mathclose[$
and hence $j_1=j_2$ by the pairwise disjointness of $(\mathopen[c_j,d_j\mathclose[ : j\in J)$,
and if $(i_1,j),(i_2,j)\in A$, then similarly 
$z-c_j\in\mathopen[a_{i_1},b_{i_1}\mathclose[\cap\mathopen[a_{i_2},b_{i_2}\mathclose[$
and hence $i_1=i_2$. Thus we get                            
\la                                                            \label{Eq:CS_in_the_proof}                                                                                      
 F^{}_{(P-R_1)\ast(Q-R_2)}(z)  &=& \sum_{(i,j)\in I\times J}  F^{}_{S_i\ast T_j}(z) 
  \,\ \le \,\  \sum_{(i,j)\in A} p_iq_j \,\ =\,\ \sum_{i\in I_0} p_iq_{j(i)} \\ 
  &\le& \sqrt{\sum_{i\in I_0} p_i^2 \sum_{i\in I_0} q_{j(i)}^2}   \,\ \le \,\
  \sqrt{ \sum_{i\in I} p_i^2  \sum_{j\in J} q_j^2 }               \nonumber \\
  &\le& 2\,\sqrt{\left\|H_1\right\|_\mathrm{L}
   \left\|H_2\right\|_\mathrm{L}
   \varkappa(P-R_1)\varkappa(Q-R_2)}\nonumber
\al
by applying in the second step above Lemma~\ref{Lem:S_i_ast_T_j_bounds}
to $U\coloneqq R_1(\cdot\cap[a_i,b_i])$ 
and $V\coloneqq R_2(\cdot\cap[c_j,d_j])$ for each pair $(i,j)$,
and by using in the last step $\varkappa(P-R_1) =\sum_{i\in I}  \int_{a_i}^{b_i} \left(H_1(b_i)-H_1(x)\right)\dd x$ and
\[
 \int_{a_i}^{b_i} \big(H_1(b_i)-H_1(x)\big)\,\dd x 
  &\ge& \int_{a_i}^{a_i
  + \frac{p_i}{\,\left\|H_1\right\|_\mathrm{L}}}
  \big(p_i-\left\|H_1\right\|_\mathrm{L}(x-a_i)\big)\,\dd x
 \,\ =\,\  \frac{p_i^2}{ 2\left\| H_1\right\|_\mathrm{L}} \,,
\]
and the analogous inequalities for the $q_j^2$. This
proves~\eqref{Eq:CS_for_extreme_points_new} in the present unexceptional case.

\smallskip 6. Reduction of the general case to the unexceptional one:
Let now $P$ and $Q$ be arbitrary as specified in step 4,
but without loss of generality we assume $I\neq\emptyset \neq J$. Let  $i\in I$ be fixed, 
exceptional if possible, and arbitrary else; in the latter case we also choose 
an arbitrary $\xi \in\mathopen]a_i,b_i\mathclose[$. 
In any case we then put 
\[
 P_\sigma &\coloneqq&  P-S_i + \sigma\delta_{a_i}+(p_i-\sigma)\delta_\xi
  - R_1(\cdot\cap[a_i,b_i[\,) \quad\text{ for }\sigma\in [R_1([a_i,\xi[\,),p_i].
\]
Let analogously $j\in J$ be fixed, exceptional if possible, and $\eta\in\mathopen]c_j,b_j\mathclose[$
chosen if necessary, and 
\[
 Q_\tau &\coloneqq& Q-T_j + \tau\delta_{c_j} +(q_j-\tau)\delta_{\eta} -  R_2(\cdot\cap[c_j,d_j[\,) \quad\text{ for }
  \tau \in[R_2([c_j,\eta[\,),q_j].
\]
For each $\sigma$ then $P_\sigma\stle R_1$,
say by Lemma~\ref{Lem:cK_and_its_extreme_points} applied to $R_1$
and $\epsilon'\coloneqq \varkappa(P_\sigma-R_1)$,
and hence by~\eqref{Eq:kappa(P-Q)_if_stle} the map
$\sigma\mapsto \varkappa(P_\sigma-R_1)=\int F_{P_\sigma-R_1}\dd\leb$ is
affine-linear. With the analogous property for $Q$ we conclude that
the function $\Omega$, defined by
\[
 \Omega(\sigma,\tau) &\coloneqq& F_{(P_\sigma-R_1)\ast(Q_\tau-R_2)}(z)
 - 2\,\sqrt{\left\|H_1\right\|_\mathrm{L}\left\|H_2\right\|_\mathrm{L}\varkappa(P_\sigma-R_1)\varkappa(Q_\tau-R_2)}
\]
for $(\sigma,\tau)$ belonging to the square $[R_1([a_i,\xi[),p_i]\times [R_2([c_j,\eta[),q_j]$,
is separately convex in each of its two variables, and hence assumes its maximal value at 
some of the four corners. But if $(\sigma,\tau)$ is one of the four corners, then $P_\sigma$ and $Q_\tau$
are nonexceptional extreme points (with possible different $\epsilon_1, \epsilon_2$, 
and with the index set $I$ enlarged by one element in case of $\sigma= R_1([a_i,\xi[\,)$,
analogously for $J$), and hence we then get $\Omega(\sigma,\tau)\le 0$ by step~5 of this proof.
Thus we have $\Omega\le0$ everywhere. Since $P=P_\sigma$ and $Q=Q_\tau$ for some $\sigma,\tau$, we are done.
\end{proof}

\begin{proof}[Proof of Theorem~\ref{Thm:F_star_G_vs_H_star_H}]
Let $P_1,P_2,R_1,R_2\in\Prob(\R)$ be the laws corresponding to $F_1,F_2,H_1,H_2$.
We recall~(\ref{Eq:Def_Kolmogorov-norm},\ref{Eq:Def_M_Lipschitz},\ref{Eq:varkappa_via_F-G},%
\ref{Eq:lambda_Q-P_via_F-G}).

\smallskip
1. With the notation  $\check{P}$ for the reflection of a law $P$,
introduced on page~\pageref{page:Def_reflection_of_M},
we have
\[ 
 \text{L.H.S.\eqref{Eq:The_convolution_inequality_new}}
   &=& \left\|P_1\ast P_2-R_1\ast R_2\right\|_\mathrm{K} \\
   &=& \sup_{z\in\R}\max\{ F_{P_1\ast P_2-R_1\ast R_2}(z), - F_{P_1\ast P_2-R_1\ast R_2}(z-)\}  \\
   &=& \sup_{z\in\R}\max\{ F_{P_1\ast P_2- R_1\ast R_2}(z),
       F_{\check{P_1}\ast\check{P_2}- \check{R}_1 \ast \check{R}_2}(z)\}\,,
\]
and, for $i\in\{1,2\}$, 
$\varkappa(\check{P_i}-\check{R_i})=\varkappa(P_i-R_i)$,
$|\lambda(\check{P_i}-\check{R_i})|=|\lambda(P_i-R_i)|$
in case of $\varkappa(P_i-R_i)<\infty$,
and $\left\|\check{R_i}\right\|_\mathrm{L} = \left\|R_i\right\|_\mathrm{L}$.
Hence, passing to the reflections if necessary,
it is enough to prove
\[                                              
 \sup_{z\in\R} F^{}_{P_1\ast P_2 -R_1\ast R_2}(z)   &\le &
 \text{R.H.S.\eqref{Eq:The_convolution_inequality_new}}\,.
\]

\smallskip 2. For $i\in\{1,2\}$ let $Q_i\coloneqq P_i\wedge_\mathrm{st} R_i$,
the infimum with respect to the stochastic order~$\stle$ recalled in~\eqref{Eq:usual_stochastic_order},
having as distribution function the pointwise supremum $G_i\coloneqq F_i\!\vee H_i$.
Then $ F_{Q_1\ast Q_2}(z)\ge F_{P_1\ast P_2}(z)$
for $z\in \R$ by the isotonicity w.r.t.~$\stle$ of convolution,
and, for $i\in\{1,2\}$,
\[\quad
 \left\| G_i-H_i\right\|_1 &=& \int (F_i-H_i)_+\,\dd\leb
  \,\ =\,\ \int \!\tfrac{1}{2}\big( \left|F_i-H_i\right| + F_i-H_i\big)\,\dd\leb
  \,\ \le\,\ \ltrivert F_i-H_i \rtrivert .
\]
Hence it is enough to prove
\la                                                     \label{Eq:P_ast_Q_vs_N^2_pointwise}
 \sup_{z\in\R} F^{}_{P_1\ast P_2 -R_1\ast R_2}(z)   &\le &
 \text{R.H.S.\eqref{Eq:The_convolution_inequality_old}}
\al
with the $(Q_i,G_i)$ in place of the $(P_i,F_i)$ or, equivalently,
as it stands with the $(P_i,F_i)$ but assuming $P_i\stle R_i$.

\smallskip 3. Using below the ring identity
\la                                       \label{Eq:ring_identity}
  P_1\ast P_2  - R_1\ast R_2 &=& (P_1-R_1)\ast(P_2-R_2) +(P_1-R_1)\ast R_2  + R_1\ast(P_2-R_2)
\al
in the first step,
and in the second for the first summand the assumption $P_i\stle R_i$
and Lemma~\ref{Lem:CS_for_extreme_points}
and for the remaining two summands two sumands just
Lemma~\ref{Lem:||cdot||_K-smoothing},
%
yields                                   \label{page:use_of_Lem:||cdot||_K-smoothing}
\[
 \text{L.H.S.\eqref{Eq:P_ast_Q_vs_N^2_pointwise}}
   &\le& \left\| (P_1-R_1)\ast(P_2-R_2)  \right\|_\mathrm{K}     
    + \| (P_1-R_1)\ast R_2 \|^{}_{\mathrm{K}} + \|(P_2-R_2)\ast R_1 \|^{}_{\mathrm{K}} \\
   &\le&  \text{R.H.S.\eqref{Eq:Kolmogorov_distance_of_product_of_differences_le}}
     + \varkappa(P_1-R_1)\left\|R_2\right\|_\mathrm{L}
     + \varkappa(P_2-R_2)\left\|R_1\right\|_\mathrm{L} \\
   &=& \text{R.H.S.\eqref{Eq:P_ast_Q_vs_N^2_pointwise}}\,.
\]

\vspace{-1.5\baselineskip}
\end{proof}

The following presumably known side remark suggests to us that a slight complication like 
using upper semi-continuity of the separately affine-linear function
$(P_1',P_2')\mapsto\Psi(P_1',P_2',z-)$  in step~3 of the
proof of Lemma~\ref{Lem:CS_for_extreme_points}
might be unavoidable.

\begin{Rem}                                            \label{Rem:phi_discontinuous_on_K}
A discontinuous linear functional, vanishing at each extreme point of a compact
and convex subset $K$ of a topological vector space $X$, need not be bounded on $K$,
even if $X$ is a Hilbert space.
\end{Rem}
\begin{proof}
Let $X\coloneqq \ell^2$, the usual Hilbert space of all real quadratically summable sequences,
and let $K\coloneqq\{x\in X: |x_n|\le \frac{1}{n} \text{ for }n\in\N\}$, the Hilbert cube.
Then $K$ is compact and convex, and its set of extreme points is
$E\coloneqq\{x\in X: |x_n|=\frac1n\text{ for }n\in\N\}$.
If $x\in\mathrm{span\,}   E$, that is, $x=\sum_{j=1}^k\alpha_je^j$ for some $k\in\N$,
$\alpha_j\in\R$, and  $e^j\in E$,
then $nx_n\in \{\sum_{j=1}^k \alpha_j\epsilon_j : \epsilon\in\{-1,1\}^k\}$ for each $n\in\N$, 
and so the set $\{nx_n:n\in\N\}$ is finite. Hence $b^k \coloneqq (n^{-k-1})_{n\in\N} \in K\setminus\mathrm{span\,}E$ for $k\in\N$.
Choosing $E_0\subseteq E$ maximal linearly independent, and extending the linearly independent set 
$E_0\cup\{b^k:k\in\N\}$ to an algebraic basis of $X$ by some $B_0\subseteq X$, we  may define a linear functional $\phi$
on $X$ by requiring $\phi(b)=0$ for $b\in E_0\cup B_0$ and $\phi(b^k)=k$ for $k\in\N$, 
and get $\phi=0$ on $E$ but $\sup_{x\in K} \phi(x)=\infty$.  
\end{proof}

\section{Auxiliary results for $\zeta$ and related distances} \label{sec:zeta_distances}
In this section and in the next one, we often write convolution of laws or more 
general bounded signed measures simply as juxtaposition, as in $PQ \coloneqq P\ast Q$, 
and similarly for convolution powers, $P^n\coloneqq P^{\ast n}$.
We need some 
well-known auxiliary facts about Kolmogorov and $\zeta$ distances,
and we might as well state the first few, namely variations of the so-called 
regularity~\eqref{Eq:regularity_on_cM} or its special case~\eqref{Eq:regularity_on_Prob},
and of the homogeneity~
\eqref{Eq:homogeneity_on_cM_conclusion}, 
in a more natural generality.
Below, a set $\cF$ of functions defined on $\R$ is {\em translation invariant} if
$f\in\cF$ and $a\in\R$ imply $f(\,\cdot+a)\in\cF$, 
and {\em reflection invariant} if $f\in\cF$ implies $f(-\,\cdot\,)\in\cF$.
We put 
\[
 \cL^\infty &\coloneqq&  \big\{ g \in \C^\R : g\text{ Borel and }\sup_{x\in\R}|g(x)|<\infty\big\}\,.
\]

\begin{Lem}                                                            \label{Lem:Regular_eqnorms_on_cM}
Let $\cF\subseteq\cL^\infty$ be a translation invariant subset, and let 
\[
 \|M\| &\coloneqq& \|M\|^{}_{\cF} \,\ \coloneqq\,\ \sup\{ |Mf| : f\in\cF\} \quad\text{ for }M\in\cM\,.
\]
Then $\|\cdot\|$ is an eqnorm on $\cM$, for $M,M_1,M_2,M_3\in\cM$ we have 
\la
 \|\delta_a M\| &=& \|M\| \quad\text{ for }a\in\R\,, \label{Eq:translation_invariance_on_cM}  \\
 \|M_1M_2\| &\le& \|M_1\| \nu_0(M_2)\,, \label{Eq:regularity_on_cM} \\
 \|M_1M_2\| &\le& \|M_1M_3\| + \min\big\{\,\|M_1\|\nu_0(M_2-M_3)\,,\,\nu_0(M_1)\|M_2-M_3\|\,\big\}\,,
                                                         \label{Eq:changing_the_regularizer}   
\al
and for $n\in\N_0$ and $P,Q,R,P_1,\ldots,P_n,Q_1,\ldots,Q_n\in\Prob(\R)$ we have
\la
 \|PR-QR\| &\le& \|P-Q\|\,,             \label{Eq:regularity_on_Prob}\\
 \|P-Q\| &\le& \|PR-QR\| + 2\|R-\delta_0\|\,, \label{Eq:smoothing_on_Prob}     \\
 \left\| \bigconv_{j=1}^nP_j -\bigconv_{j=1}^nQ_j \right\|   \label{Eq:semiadditivity}
    &\le& \sum_{j=1}^n\left\| P_j-Q_j \right\|. 
\al

Further, $\ltrivert \cdot \rtrivert \coloneqq \nu_0 \vee \|\cdot\| = \|\cdot\|^{}_{\cF_0\cup\cF} $
with $\cF_0$ from~\eqref{Eq:Def_cF_0,0} is an enorm on $\cM$, and is
submultiplicative in the sense of 
\la                                        \label{Eq:|||_submultiplicative}
  \ltrivert M_1M_2\rtrivert &\le& \ltrivert M_1\rtrivert\ltrivert M_2\rtrivert
  \quad\text{ for }M_1,M_2\in\cM.
\al

If $\cF$ is reflection invariant, then so is $\|\cdot\|$, that is, then 
$\|\check{M}\|=\|M\|$ for $M\in\cM$. 

If $r\in\R$ is such that the implication 
\la                                       \label{Eq:homogeneity_on_cM_assumption}
 f\in\cF,\, \lambda\in\mathopen]0,\infty\mathclose[ 
  &\Rightarrow&  \lambda^{-r}f(\lambda\,\cdot\,) \in\cF
\al
holds, then we have, for $M\in\cM$,
\la                                    \label{Eq:homogeneity_on_cM_conclusion} 
 \left\|M\!\left(\tfrac{\cdot}{\lambda}\right) \right\| 
   &=& \big\| (x\mapsto \lambda x)\im M  \big\|
   \,\ =\,\ \lambda^r \big\| M\big\|
  \quad\text{ for }\lambda\in\mathopen]0,\infty\mathclose[\,.
\al
\end{Lem}
\begin{proof} The eqnorm claim is obvious.  
For  $f\in\cF$ we have
\[
 \left|\int f\,\dd M_1M_2\right| &=& \left| \int\int f(x+y)\,\dd M_1(x)\dd M_2(y)\right|
   \,\ \le \,\ \int \left| \int f(x+y)\,\dd M_1(x)\right| \dd\left| M_2\right|(y)   \\  
 &\le& \int \left\| M_1 \right\|_\cF \dd\left| M_2\right|(y)  \,\ =\,\ \text{R.H.S.\eqref{Eq:regularity_on_cM}}\,,
\]
and this proves~\eqref{Eq:regularity_on_cM}.
The latter applied once to $(M_1,M_2)\coloneqq (M,\delta_a)$
and once to  $(M_1, M_2) \coloneqq (\delta_a M, \delta_{-a})$ 
yields~\eqref{Eq:translation_invariance_on_cM}.
Writing $M_1M_2=M_1M_3+M_1(M_2-M_3)$ and applying first subadditivity of $\|\cdot\|$, 
and then~\eqref{Eq:regularity_on_cM} in two ways, yields~\eqref{Eq:changing_the_regularizer}.  
\eqref{Eq:regularity_on_Prob} is just~\eqref{Eq:regularity_on_cM} with $M_1\coloneqq P-Q$ and $M_2\coloneqq R$.
\eqref{Eq:smoothing_on_Prob} is~\eqref{Eq:changing_the_regularizer}
with $M_1\coloneqq P-Q$, $M_2\coloneqq\delta_0$, $M_3\coloneqq R$, $\nu_0(P-Q)\le 2$, taking the second minimand.
\eqref{Eq:semiadditivity} follows from
\la                                                                      \label{Eq:telescope_P_j_Q_j}
  \bigconv_{j=1}^nP_j -\bigconv_{j=1}^nQ_j 
    &=& \sum_{k=1}^n \left(\bigconv_{j=1}^{k-1} P_j\right)
    \left(P_k-Q_k\right)
    \bigconv_{j=k+1}^nQ_j
\al
by applying subadditivity of $\|\cdot\|$ and then~\eqref{Eq:regularity_on_Prob} with $P-Q=P_k-Q_k$.

\eqref{Eq:regularity_on_cM} applied to $\cF_0$ instead of $\cF$ 
yields the well-known total variation norm inequality 
\la                                       \label{Eq:nu_0_algebra_norm}
 \nu_0(M_1M_2) &\le& \nu_0(M_1)\nu_0(M_2)\,,
\al
and this combined with~\eqref{Eq:regularity_on_cM} as it stands yields~\eqref{Eq:|||_submultiplicative}.

The remaining claims, about reflection invariance and scaling behaviour, 
are also easy to check.
\end{proof}

Lemma~\ref{Lem:Regular_eqnorms_on_cM} may of course be adapted to more general measurable monoids
in place of $(\R,+)$. As it stands it applies in particular to $\nu_0$,
as already noted in the above proof,
and to the Kolmogorov norm $\|\cdot\|^{}_{\mathrm{K}}$
as defined by~\eqref{Eq:Def_Kolmogorov-norm}.
In these two cases, \eqref{Eq:homogeneity_on_cM_conclusion} applies with $r=0$,
and we get the reflection and scale invariances 
\la                                    \label{Eq:nu_0_and_Kolmogorov_scale_invariant}
  \nu_0\big(M(\tfrac{\cdot}{\lambda})\big) &=& \nu_0(M), 
  \quad \left\| M(\tfrac{\cdot}{\lambda})\right\|^{}_{\mathrm{K}} \,\ =\,\ \left\| M\right\|_{\mathrm{K}} 
  \quad\text{ for } M\in\cM\text{ and }\lambda\in\R\!\setminus\!\{0\} \,.
\al
Lemma~\ref{Lem:Regular_eqnorms_on_cM} further applies
to each of the enorms  $\zeta_r$ with $r\in\N_0$ defined by~\eqref{Eq:Def_zeta_r}, 
and with the exception of~\eqref{Eq:homogeneity_on_cM_conclusion} also  
to the dual bounded Lipschitz norm~$\beta$ from~\eqref{Eq:Def_beta}. 
Special cases of inequality~\eqref{Eq:changing_the_regularizer}
are given by \citet[p.~365, (6.5.43) and (6.5.44) combined, p.~366, (6.5.46) and (6.5.47) 
combined]{Zolotarev1997} and also, on $\R^k$, by \citet[p.~85, (2.8.1) and (2.8.2) combined, pp.~122-123, 
(2.10.25) and next display combined]{Senatov1998}. The remainder of Lemma~\ref{Lem:Regular_eqnorms_on_cM} is 
even better known. 

Most of Lemma~\ref{Lem:Regular_eqnorms_on_cM} does not apply to 
the eqnorms $\nu_r$ from~\eqref{Eq:Def_nu_r} with $r>0$,
although we have $\nu_r(M) = \sup\{|Mf| : f\in\cF\}$ 
with $\cF \coloneqq \{f\in\cL^\infty :  |f(x)|\le |x|^r\text{ for }x\in \R\}$,
since for example~\eqref{Eq:regularity_on_cM}  
with $\|\cdot\| \coloneqq \nu_r$
would yield the absurdity 
$\nu_r(M)=\nu_r(\delta_0M) \le \nu_r(\delta_0)\nu_0(M)=0$ for every $M\in\cM$. 
This illustrates the importance of the translation invariance of $\cF$ in Lemma~\ref{Lem:Regular_eqnorms_on_cM}, 
violated by the present $\cF$. 
However, we obviously do have~\eqref{Eq:homogeneity_on_cM_conclusion} and reflection invariance 
for $\|\cdot\| \coloneqq \nu_r$, that is,
\la                                      \label{Eq:nu_r_scaling_behaviour}
 \nu_r\big(M\!\left(\tfrac{\cdot}{\lambda}\right)\big) &=& |\lambda|^r \nu_r( M)
  \quad\text{ for }r\in[0,\infty[\,,\,\lambda\in\R\!\setminus\!\{0\}\,,\,
  M\in\cM \,,
\al
and we have analogous identities for $\varkappa_r$ in~\eqref{Eq:varkappa_scaling},
for $\underline{\zeta}_r$ since~\eqref{Eq:homogeneity_on_cM_assumption}
is fulfilled for $\cF\coloneqq\cF^\infty_{r,r-1}$ from~\eqref{Eq:Def_cF_r,r-1^infty}. 
And, as an analogue of~\eqref{Eq:nu_0_algebra_norm} in the style 
of~\eqref{Eq:|||_submultiplicative}, 
used in Example~\ref{Example:Zolotarev_1973},
we have
\la\qquad                                           \label{Eq:nu_0_with_nu_r_submultiplicative}
 \big(\nu_0\!\vee\! \nu_r\big)(M_1M_2) 
  &\le& 2^{r\vee1} \big(\nu_0\!\vee\! \nu_r\big)(M_1)\,\big(\nu_0\!\vee\! \nu_r\big)(M_2)
  \quad\text{ for } r\in[0,\infty[\,,\, M_1,M_2\in\cM\,,
\al
since we have $|x+y|^r\le (|x|+|y|)^r \le 2^{(r-1)\vee0}(|x|^r+|y|^r)$ for $x,y\in\R$,
and hence indeed also 
$ \nu_r(M_1M_2) = \iint|x+y|^r\,\dd|M_1|(x)\dd|M_2|(y)
 \le 2^{(r-1)\vee0}(\nu_r(M_1)\nu_0(M_2) +\nu_0(M_1)\nu_r(M_2))
\le \text{R.H.S.\eqref{Eq:nu_0_with_nu_r_submultiplicative}}$. 

The scaling behaviour \eqref{Eq:homogeneity_on_cM_assumption} $\Rightarrow$  \eqref{Eq:homogeneity_on_cM_conclusion},
also called {\em homogeneity}, 
somewhat in conflict with the absolute homogeneity of just any 
eqnorm, yields, 
using also the translation invariance~\eqref{Eq:translation_invariance_on_cM}, in particular 
\la                                                  \label{Eq:homogeneity_of_zeta}
  \zeta_r(P-Q) &=& \lambda^r\zeta_r(\widetilde{P}-\widetilde{Q})
  \quad \text{ for }P,Q\in\cP_2 \text{ with }\sigma(P)=\sigma(Q)=\lambda,\ r\in\N_0\,.
\al
This is used, for example,
in the proofs of Theorem~\ref{Thm:Zolotarev's_zeta_1-B-E-Thm} and 
of the following simple and well-known result.

\begin{Cor}
\la                          \label{Eq:Simple_zeta_3-CLT_error_bound}
  \zeta_3\!\left(\widetilde{P^{\ast n}}-\mathrm{N}\right) &\le& 
   \frac{\zeta_3\!\left(\widetilde{P}-\mathrm{N}\right)}{\sqrt{n}}
   \quad\text{ for }P\in\cP_3\text{ and }n\in\N\,.
\al
\end{Cor}
\begin{proof} We have 
$\text{L.H.S.\eqref{Eq:Simple_zeta_3-CLT_error_bound}} 
 =  \zeta_3\!\left(\widetilde{\widetilde{P}_{}^{\ast n}}-\widetilde{\mathrm{N}^{\ast n}}\right)
 = n^{-\frac{3}{2}}  \zeta_3\!\left(\widetilde{P}^{\ast n}_{}-\mathrm{N}_{}^{\ast n}\right) 
 \le \text{R.H.S.\eqref{Eq:Simple_zeta_3-CLT_error_bound}}$,
using~\eqref{Eq:homogeneity_of_zeta} with $r=3$ and $\lambda=\sqrt{n}$ in the second step,
and~\eqref{Eq:semiadditivity} in the third.
\end{proof}

For $\|\cdot\|=\|\cdot\|^{}_{\mathrm{K}}$ and $M_2$ sufficiently regular, 
the following simple alternative to~\eqref{Eq:regularity_on_cM}
might be preferable, and is used on page~\pageref{page:use_of_Lem:||cdot||_K-smoothing}
in the proof of Theorem~\ref{Thm:F_star_G_vs_H_star_H}.
We recall (\ref{Eq:Def_Kolmogorov-norm},\ref{Eq:Def_M_Lipschitz},\ref{Eq:Def_cM_r,k}),
in particular $ \cM_{0,0}=\{M\in\cM : M(\R)=0\}$,
\eqref{Eq:varkappa_1=L^1-enorm},
for~\eqref{Eq:||cdot||_K-smoothing_sharper} also~\eqref{Eq:Def_lambda_etc},
and~\eqref{Eq:zeta_1=kappa_1} for the possibility
of replacing $\varkappa_1(M_1)$ by $\zeta_1(M_1)$
 in case of $M_1\in\cM_{1,0}$.
\begin{Lem}                                        \label{Lem:||cdot||_K-smoothing}
Let $M_1\in\cM_{0,0}$ and $M_2\in\cM$. Then
$\left\|M_1M_2\right\|_{\mathrm{K}}\le\varkappa_1(M_1)\left\|M_2\right\|_{\mathrm{L}}$ always,
and
\la                        \label{Eq:||cdot||_K-smoothing_sharper}
  \left\|M_1M_2\right\|_{\mathrm{K}} &\le&
    \tfrac{1}{2}\big(\varkappa_1(M_1)+\lambda_1(M_1)\big)\left\|M_2\right\|_{\mathrm{L}}
  \quad \text{ if $M_2\ge 0$ and $\varkappa_1(M_1)<\infty$}\,.
\al
\end{Lem}
\begin{proof}
 We may assume  $\left\|M_2\right\|_{\mathrm{L}}<\infty$,
 and hence $M_2=f\leb$ for some
 measurable function~$f$ with $\left\|f\right\|_\infty = \left\|M_2\right\|_{\mathrm{L}}\,$
 by~\eqref{Eq:Lip=sup_|f'|}.
 Then  $F_{M_1M_2}(x)=\int F_{M_1}(y)f(x-y)\dd y$ for $x\in\R$ yields
 $\left\|M_1M_2\right\|_{\mathrm{K}}\le \int \left|  F_{M_1}(y)  \right|\left\|f\right\|_\infty \dd y
  =  \varkappa_1(M_1)\left\|M_2\right\|_{\mathrm{L}}$
 by using $M_1M_2\in\cM_{0,0}$ and~(\ref{Eq:Def_Kolmogorov-norm},\ref{Eq:varkappa_1=L^1-enorm}).
 For \eqref{Eq:||cdot||_K-smoothing_sharper} we write
 $F_{M_1M_2}(x)=\int F_{M_1}(y)\big(f(x-y)-\frac{1}{2}\left\|f\right\|_\infty\big) \dd y
 +  \frac{1}{2}\left\|f\right\|_\infty \int F_{M_1}(y)\,\dd y$
 and use positivity of $f$.
\end{proof}

The following perhaps not completely trivial norm comparison lemma 
is used in the proof of Corollary~\ref{Cor:B-E-Z_with_beta}.
We recall~(\ref{Eq:Def_zeta_r},\ref{Eq:Lip=sup_|f'|},\ref{Eq:Def_beta}).
%

\begin{Lem}                                    \label{Lem:zeta_1_bounded_by_beta_and_zeta_3}
On $\cM$ we have 
\la                              
 \zeta_1 &\le&                              \label{Eq:varkappa_linearly_bounded_by_beta_and_zeta}
       2\beta+3_{\phantom{3}}^\frac{1}{3}\beta_{\phantom{3}}^\frac{2}{3}\zeta_3^\frac{1}{3}\,, \\
 \zeta_1\!\vee\!\zeta_3              \label{Eq:varkappa_vee_zeta_bounded_by_beta_vee_zeta}
    &\le & (2+3^\frac{1}{3})\, \beta\!\vee\! \zeta_3 \,.         
\al 

The pair $(\frac{2}{3},\frac{1}{3})$ of exponents in~\eqref{Eq:varkappa_linearly_bounded_by_beta_and_zeta}
is in the following sense i.c.f.~optimal even on~$\widetilde{\cP_3}-\mathrm{N}${\,\rm:}
There are no constants $c<\infty$ and $\alpha \in[0,\frac{1}{3}[$ with 
$\zeta_1\le c\,(\beta\vee (\beta^{1-\alpha}\zeta_3^\alpha))$ on $\widetilde{\cP_3}-\mathrm{N}$.
\end{Lem}
\begin{proof}
\eqref{Eq:varkappa_vee_zeta_bounded_by_beta_vee_zeta}
follows trivially from~\eqref{Eq:varkappa_linearly_bounded_by_beta_and_zeta}
by using $\beta\le \beta\!\vee\! \zeta_3$ and $\zeta_3 \le  \beta\!\vee\! \zeta_3$.

To prove~\eqref{Eq:varkappa_linearly_bounded_by_beta_and_zeta}, let $t\in\mathopen]0,\infty\mathclose[$
and $\psi(x)\coloneqq\frac{1}{t}{\big(1-\frac{|x|}{t}\big)}_+$ for~$x\in\R$.
For $f\in\cF_1^\infty$ as defined in~\eqref{Eq:Def_cF_r^infty}, we then put $f_1\coloneqq f- f\ast \psi$ 
and $f_2 \coloneqq f\ast \psi$, get 
\[
 \left\|f_1\right\|_\infty 
  &=& \sup_{x\in\R} \left| \int \big(f(x)-f(x-y)\big)\psi(y) \dd y  \right| 
  \,\ \le \,\  \int  |y|\psi(y) \dd y 
  \,\ = \,\ \tfrac{t}{3}\,, \\
 \left\| f_1 \right\|_{\mathrm{L}} 
  &\le & \left\| f \right\|_{\mathrm{L}} + \left\| f\ast\psi \right\|_{\mathrm{L}}
   \,\ \le \,\ 2\,,
\]
and with $\psi'' =\frac{1}{t^2}(\delta_{-t}-2\delta_0 +\delta_t) $ 
in the sense of distributions and then $f_2''=f\ast\psi''$, 
say by \citet[17.5.12.2, 17.5.7.1, 17.11.11.1, 17.11.1.1]{Dieudonne3}, also  
\[
 \left\| f_2''\right\|_{\mathrm{L}} &= &  \left\| f\ast\psi''\right\|_{\mathrm{L}}
  \,\ \le \,\  \left\| f \right\|_{\mathrm{L}} \nu_0(\psi'') 
  \,\ \le \,\ \tfrac{4}{t^2}
\] 
(with, we recall, $\nu_0$ denoting the usual total variation norm of a signed measure), 
and for $M\in\cM$  therefore  
\[
 \left| \int f \dd M\right| &\le&  \left|\int f_1 \dd M\right| + \left|\int f_2 \dd M\right|
  \,\ \le\,\ (\tfrac{t}{3}+2)\beta(M) + \tfrac{4}{t^2}\zeta_3(M)\,. 
\]
Minimising the right hand side above,  unless it is zero or infinite anyway, 
at $t=(24\frac{\zeta_3}{\beta}(M))^\frac{1}{3}$ 
yields~\eqref{Eq:varkappa_linearly_bounded_by_beta_and_zeta}.

The final claim is indeed an optimality claim, since for fixed 
values $\beta,\zeta_3\in\mathopen]0,\infty\mathclose[$
and then with $f(\alpha)\coloneqq \beta^{1-\alpha}\zeta_3^\alpha$ for $\alpha\in\R$,
we have $\beta\vee (\beta^{1-\alpha}\zeta_3^\alpha)= f(0)\vee f(\alpha)$
decreasing in $\alpha\in\mathopen]-\infty,0\mathclose]$
and increasing in $\alpha\in[0,\infty[$, by convexity of $f$.
If now $\zeta_1\le c\,(\beta\vee (\beta^{1-\alpha}\zeta_3^\alpha))$ on $\widetilde{\cP_3}-\mathrm{N}$,
with some $c<\infty$ and, say, $\alpha\in[0,1]$, 
then recalling the asymptotic bounds for $P=P_t$ from Example~\ref{Example:beta_llcurly_zeta_1},
and also $\zeta_3(P-\mathrm{N})\asymp t^4\phi(t)$ for $t\rightarrow\infty$
from~\eqref{Eq:zeta_3_for_Example_beta_llcurly_zeta_1}, yields 
\[
 t^2\phi(t) \,\ \preccurlyeq\,\ \zeta_1(P-\mathrm{N}) 
  & \preccurlyeq & \big(\beta\vee (\beta^{1-\alpha}\zeta_3^\alpha)\big)(P-\mathrm{N}) \\
  & \preccurlyeq & \big(t\phi(t)\big) \vee \Big( \big(t\phi(t)\big)^{1-\alpha}\big(t^4\phi(t)\big)^\alpha \Big)
  \,\ \sim\,\ t^{1+3\alpha}\phi(t)
\]
and hence $\alpha\ge \frac{1}{3}$.
\end{proof}

We recall from~\eqref{Eq:Def_N_sigma_N_1}
that $\mathrm{N}_\sigma$ denotes the centred normal law on $\R$ with standard deviation 
$\sigma\in\mathopen[0,\infty\mathclose[\,$. 
A specialisation of the so-called smoothing inequality~\eqref{Eq:smoothing_on_Prob}
yields:

\begin{Lem}                                           \label{Lem:zeta_1-smoothing}
We have 
$\zeta_1(P-Q) \le \zeta_1(P\mathrm{N}_\epsilon-Q\mathrm{N}_\epsilon) + \frac{4}{\sqrt{2\pi}}\epsilon$
for $P,Q\in\Prob(\R)$ and $\epsilon\ge0$.
\end{Lem}
\begin{proof}
\eqref{Eq:smoothing_on_Prob}   with $\|\cdot\| \coloneqq\zeta_1$, 
$R\coloneqq\mathrm{N}_\epsilon$, 
$\zeta_1(\mathrm{N}_\epsilon-\delta_0)=\epsilon\,\nu_1(\mathrm{N})=\frac{2\epsilon}{\sqrt{2\pi}}$
by (\ref{Eq:zeta_1_under_x_maspto_ax},\ref{Eq:nu_3_von_N}).
\end{proof}

\begin{Lem}                           \label{Lem:zeta_s_vs_zeta_(s+k)}
Let $M\in\cM$, 
$s,k\in\N_0$, 
and $\sigma\in\mathopen]0,\infty\mathclose[$.
Then we have
\la                                                 \label{Eq:zeta_s_smoothed_vs_zeta_{s+k}}
   \zeta_s(M\mathrm{N}_\sigma) &\le& 
   \left\|\phi^{(k)}\right\|_1\frac{\zeta_{s+k}(M)}{\sigma^k}
\al
where 
\[
 && \left\|\phi^{(0)}\right\|_1 \,=\,1,\qquad 
  \left\|\phi^{(1)}\right\|_1 \,=\,\frac{2}{\sqrt{2\pi}}= 0.797884\ldots,\qquad
 \left\|\phi^{(2)}\right\|_1 \,=\, \frac{4\mathrm{e}^{-1/2}}{\sqrt{2\pi}}
   \,=\,   0.967882\ldots, \\
 && \left\|\phi^{(3)}\right\|_1 \,=\, \frac{2+8\mathrm{e}^{-3/2}}{\sqrt{2\pi}} \,=\, 1.510013\ldots,
    \qquad  \\ 
 &&   \left\|\phi^{(4)}\right\|_1 
     \,=\, 4\frac{ \sqrt{18-6\sqrt{6}} \mathrm{e}^{ -\frac{3-\sqrt{6}}{2}   }
              +\sqrt{18+\sqrt{6}}\mathrm{e}^{ -\frac{3+\sqrt{6}}{2}  }  }
             {\sqrt{2\pi}}
     \,=\, 2.800600
       \ldots\,.  
\]
\end{Lem}
\begin{proof} 
The stated values of the $\left\|\phi^{(k)}\right\|_1$ are well-known and easily checked. 
So only~\eqref{Eq:zeta_s_smoothed_vs_zeta_{s+k}} remains to be considered:

The case of $k=0$ is contained in~\eqref{Eq:regularity_on_cM} 
of Lemma~\ref{Lem:Regular_eqnorms_on_cM}, 
and may hence be excluded here. In case of $s>0$, then, 
inequality~\eqref{Eq:zeta_s_smoothed_vs_zeta_{s+k}} is proved, 
assuming but not using $M=P-Q$ with $P,Q\in\Prob(\R)$, and otherwise 
more generally, 
in \citet[p.~47, Theorem 1.4.5]{Zolotarev1997} with $\mathrm{N}_\sigma$ replaced by any law with a 
$k$ times differentiable density and with $s\in\mathopen]0,\infty\mathclose[$ 
not necessarily an integer, 
and in \citet[p.~108, Lemma~2.10.1]{Senatov1998} with $s,k\in\mathopen]0,\infty[$ not necessarily integers 
and with 
a multivariate generalisation. Essentially the latter proof is, in the univariate case,  given in a bit more detail 
in~\citet[pp.~513--515, Lemma 4.1]{MattnerShevtsova}. Of these references, each  gives the definition of 
$\zeta_s$ for $s\in\mathopen]0,\infty\mathclose[$, but unfortunately none treats the case $s=0$.

For the case of $s=0$, and with a $k$ times differentiable probability density $f$ in place of~$\phi$,
\citet[pp.~1257-1258, Lemma~17]{Boutsikas2011} gives a sketch of a proof and provides related references.
Here we wrote ``sketch'' since there the necessary integrability properties of $f'$ 
are not addressed, and no reference to a fact like \citet[p.~149, Theorem~7.21]{Rudin1987} occurs.
So let us give a short alternative proof for the special normal case considered here: 

To prove \eqref{Eq:zeta_s_smoothed_vs_zeta_{s+k}}, for arbitrary $s,k\in\N_0$, we may assume $\sigma=1$,
as 
for arbitrary $\sigma\in\mathopen]0,\infty\mathclose[$ then
$\text{L.H.S.\eqref{Eq:zeta_s_smoothed_vs_zeta_{s+k}}}
 =  \zeta_s\big( \big(M(\sigma\cdot)\mathrm{N}\big) (\frac{\cdot}{\sigma})\big) 
 =\sigma^s \zeta_s\big(M(\sigma\cdot)\mathrm{N}\big)
 \le \sigma^s  \left\|\phi^{(k)}\right\|_1 \zeta_{s+k}\big(M(\sigma\cdot)\big)
 = \text{R.H.S.\eqref{Eq:zeta_s_smoothed_vs_zeta_{s+k}}}
$.

Let now $s=0$, $k\in\N$,  and $\sigma=1$. 
Given any function $f\in\cF_0$ from~\eqref{Eq:Def_cF_0,0} 
and writing $g(x)\coloneqq\int f(x-y)\phi(y)\,\dd y$ and
$h(x)\coloneqq g(x)/\left\|\phi^{(k)}\right\|_1$ for $x\in\R$,
it is sufficient to prove that $h\in\cF_k^\infty=\cF^\infty_{s+k}$, for then we would get
\[
 \left|(M\mathrm{N}) f\right| &=& \left| Mg\right|
 \,\ =\,\ \left\|\phi^{(k)}\right\|_1 \left|M h\right| \,\ \le\,\ \text{R.H.S.\eqref{Eq:zeta_s_smoothed_vs_zeta_{s+k}}}
\]
as desired. So let $f,g,h$ be as above. Then $h$ is 
bounded. We have 
$g(x)=\int f(y)\phi(x-y)\,\dd y$ and hence $g^{(k)}(x)=\int f(y)\phi^{(k)}(x-y)\,\dd y$ for $x\in\R$,
say by the well-known differentiability of Laplace transforms under the integral as in \citet[Example]{Mattner2001},
and we hence get $\left\|g^{(k-1)}\right\|_{\mathrm{L}}=\left\|g^{(k)}\right\|_\infty \le \left\|\phi^{(k)}\right\|_1  
= \sigma^{-k}\left\|\phi^{(k)}\right\|_1$, and thus $h\in\cF_k^\infty$. 
\end{proof}

In the following Lemma~\ref{Lem:nu_r_smoothed_vs_zeta_k}, 
the presumably rather imperfect inequality~\eqref{Eq:nu_r_smoothed_vs_zeta_k}
supplements the case of $s=0$ in~\eqref{Eq:zeta_s_smoothed_vs_zeta_{s+k}},
and is used in Example~\ref{Example:Zolotarev_1973}.
\begin{Lem}                   \label{Lem:nu_r_smoothed_vs_zeta_k}
Let $r\in[0,\infty[\,$. Then we have
\la                      \label{Eq:nu_r_vs_nu_0}
 \nu_r(M) &\le& t^r\nu_0(M)\quad \quad\text{ for }t\in[0,\infty[ \text{ and }
  M\in\cM\text{ with }M(\cdot\!\setminus\![-t,t])=0\,.
\al
If further $k\in\N_0$, then there is a constant 
$c=c^{}_{r,k}\in\mathopen]0,\infty\mathclose[$ with
\la\qquad                              \label{Eq:nu_r_smoothed_vs_zeta_k}
 \nu_r(M\mathrm{N}_\sigma) 
    &\le& c\left(\sigma\!\vee\! t\right)^r\frac{\zeta_k(M)}{\sigma^k}
    \quad\text{ for }\sigma,t\in[0,\infty[
    \text{ and }M\in\cM\text{ with }M(\cdot\!\setminus\![-t,t])=0\,.
\al
\end{Lem}
\begin{proof} \eqref{Eq:nu_r_vs_nu_0} is obvious.
In case of $\sigma=0$ and $\text{R.H.S.}\eqref{Eq:nu_r_smoothed_vs_zeta_k}<\infty$,
and even with the standard convention $\frac{0}{0}\coloneqq0$,
we have $k=0$ and then~\eqref{Eq:nu_r_smoothed_vs_zeta_k} with $c=1$ 
by $\zeta_0=\nu_0$ and~\eqref{Eq:nu_r_vs_nu_0},
or $t=0<r$ and then $\big(M\mathrm{N}_\sigma\big)(\cdot\!\setminus\!\{0\})=0$
and hence $\text{L.H.S.}\eqref{Eq:nu_r_smoothed_vs_zeta_k}=0$,
or $\zeta_k(M)=0$  and then $M=0$ and hence again
$\text{L.H.S.}\eqref{Eq:nu_r_smoothed_vs_zeta_k}=0$.

Hence we may assume $\sigma>0$, but then w.l.o.g.~$\sigma=1$, 
since~\eqref{Eq:nu_r_smoothed_vs_zeta_k} in the special case of $\sigma=1$ yields the general 
case through
$\text{L.H.S.}\eqref{Eq:nu_r_smoothed_vs_zeta_k}
 = \nu_r\big( \big( M(\sigma\cdot)\mathrm{N}\big)(\frac{\cdot}{\sigma})\big) 
 = \sigma^r  \nu_r\big( M(\sigma\cdot)\mathrm{N}\big)
 \le \sigma^r c\left(1\!\vee\!\frac{t}{\sigma}\right)^r\zeta_k(M(\sigma\cdot)) 
 = \text{R.H.S.}\eqref{Eq:nu_r_smoothed_vs_zeta_k} $.

Let $f\in\cL^\infty$ with $|f(x)|\le |x|^r$ for $x\in\R$.
With $g(x)\coloneqq \int f(x-y)\phi(y)\,\dd y = \int f(y)\phi(x-y)\,\dd y$ for $x\in\R$,
we then have $g\in\cL^\infty$ and 
\la                            \label{Eq:g_derivative_bounds_Lem_5.7}
 |g^{(j)}(x)| &=&  \left|\int f(x-y)\phi^{(j)}(y)\,\dd y\right|
   \,\ \le \int|x-y|^r\left|\phi^{(j)}(y)\right|\,\dd y
   \,\ \le \,\ c^{}_0\left(1\!\vee\!|x|\right)^r
\al
for $j\in\{0,\ldots,k\}$ and $x\in\R$, 
where $c^{}_0 \in \mathopen]0,\infty\mathclose[$ depends only on $r$ and $k$.
Let now also $t\in[0,\infty[$. 
We let $h:\R\rightarrow\C$ be the 
$k$ times continuously differentiable function which extrapolates
$g|_{[-t,t]}$\,, vanishes on $\mathopen]-\infty,-t-1\mathclose]\cup[t+1,\infty[$\,,
is on $]t,t+1[$ the Hermite interpolation polynomial
for the two interpolation points $t$ and $t+1$ and with 
there the derivatives of orders $0$ to $k$ as already determined, 
and is analogously defined on $\mathopen]-t-1,-t[$\,.
Then, using~\eqref{Eq:g_derivative_bounds_Lem_5.7} 
and \citet[pp.~502--503, Lemma~2.1(d)]{MattnerShevtsova},
we get $\left\|h^{(k)}\right\|_\infty\le c \left(1\vee t\right)^r$
for some 
$c\in \mathopen]0,\infty\mathclose[$ depending only on $r$ and $k$,
and then 
\[
 \left|(M\mathrm{N})f \right| &=& \left|M g \right| \,\ =\,\ \left|M h \right|
  \,\ \le \,\ c \left(1\vee t\right)^r \zeta_k(M)\,.
\]
This proves~\eqref{Eq:nu_r_smoothed_vs_zeta_k} in case of $\sigma=1$. 
\end{proof}

\begin{Lem}[$\underline{\zeta}$ distances of distorted images] \label{Lem:zeta_distorted_images}
Let $M\in\cM$ and $r\in\N$. If $S,T:\R\rightarrow \R$ are measurable functions, then
\la                                                \label{Eq:zeta_r_under_distortions} 
 \underline{\zeta}_r(T\im M - S\im M) 
  &\le& \tfrac{1}{(r-1)!}
   \int \big( |S|\!\vee\! |T| \big)^{r-1}\, \big|T-S\big|\,\dd|M|\,. 
\al
Let further $a,b,c,d\in\R$. Then 
\la                                 
 \underline{\zeta}_r\big((x\mapsto bx)\im M  \label{Eq:zeta_r_under_x_mapsto_ax} 
       -  (x\mapsto ax)\im     M  \big)
   &\le&   |b-a|\tfrac{(\,|a|\vee|b|\,)^{r-1}}{(r-1)!}\nu_r(M)\,, 
\al
\begin{multline}                         
 \quad \underline{\zeta}_r\big((x\mapsto cx+d)\im M     \label{Eq:zeta_r_under_x_mapsto_ax+b}
                         - (x\mapsto ax+b)\im M\big)  \\
  \,\ \le\,\ (\,|c-a|\!+\!|d-b|\,)\tfrac{(\,|a|\vee|c| +|b|\vee|d|\,)^{r-1}}{(r-1)!}\big(\nu_0\!\vee\!\nu_r\big)(M)
                    \,.  \quad 
\end{multline}

If $r=1$, then $\underline{\zeta}_1=\zeta_1$ 
in~{\rm(\ref{Eq:zeta_r_under_distortions},\ref{Eq:zeta_r_under_x_mapsto_ax},%
\ref{Eq:zeta_r_under_x_mapsto_ax+b})}, 
and we further have 
\la
 \zeta_1\big((x\mapsto bx)\im M -(x\mapsto ax)\im M\big) \label{Eq:zeta_1_under_x_maspto_ax} 
   &=& |b-a|\nu_1(M) \quad\text{ if }M\ge 0 \text{ and }ab\ge0\,, \\   
 \zeta_1(\delta_b\ast M -\delta_a\ast M)         \label{Eq:zeta_1_under_x_maspto_a+x} 
     &=& |b-a| \nu_0(M) \quad\text{ if }M\ge 0\,, \\
 \zeta_1(\delta_b-\delta_a) &=& |b-a|\,. \label{Eq:zeta_s_delta_b-delta_a} 
\al
\end{Lem}
\begin{proof} If $g\in \cF_{r,r-1}^\infty$, then, using 
$|g'(\xi)|= |g'(\xi)-\sum_{j=0}^{r-2}g^{(1+j)}(0)\frac{\xi^j}{j!}| 
\le \frac{|\xi|^{r-1}}{(r-1)!}$ $\leb$-a.e.~in the third step, 
\[
 \left| \int g\,\dd (T\im M - S\im M)  \right|
   &=& \left| \int \big(g\circ T -g\circ S \big)  \,\dd M  \right| \\
   &\le& \tfrac{1}{(r-1)!}\int \esssup_{\xi\in[S(x),T(x)]\cup[T(x),S(x)]}
     |g'(\xi)|\,\big|T(x)-S(x)\big| \,\dd|M|(x) \\
   &\le& \text{R.H.S.\eqref{Eq:zeta_r_under_distortions}}\,;
\]
hence~\eqref{Eq:zeta_r_under_distortions} holds.
If in particular $S(x)=ax+b$ and $T(x)=cx+d$ for $x\in\R$, 
then 
\la  \label{Eq:zeta_r_under_x_mapsto_ax+b_proof}   
 \text{R.H.S.\eqref{Eq:zeta_r_under_distortions}}
  &=& \tfrac{1}{(r-1)!}\int \left(|ax+b|^{r-1}\vee|cx+d|^{r-1}\right) \big|(c-a)x+d-b\big|\,\dd|M|(x), 
\al
which in case of $b=d=0$ 
equals~R.H.S.\eqref{Eq:zeta_r_under_x_mapsto_ax} 
with $c$ in place of $b$,
and is in any case at most 
\[
 \tfrac{1}{(r-1)!}\int \sum_{j=0}^{r-1}\textstyle{\binom{r-1}{j}} 
    (\,|a|\!\vee\!|c|\,)^{j}(\,|b|\!\vee\!|d|\,)^{r-1-j}|x|^j \big(\, |c-a|\, |x|+|d-b|\,\big)\,\dd|M|(x)
  &\le&  \text{R.H.S.\eqref{Eq:zeta_r_under_x_mapsto_ax+b}}
\]
by using in the final step
$\nu_j\!\vee\!\nu_{j+1} \le \nu_0\!\vee\!\nu_r$
from~\eqref{Eq:Lyapunov_three_moments_inequality}.

We have $\big(S\im M\big)(\R)=\big(T\im M\big)(\R)$, hence here $\underline{\zeta}_1=\zeta_1$.

Specialising~\eqref{Eq:zeta_r_under_x_mapsto_ax} to $r=1$ 
yields ``$\le$''  in~\eqref{Eq:zeta_1_under_x_maspto_ax}.
Assuming now $M\ge 0$ and w.l.o.g. $0\le a\le b$, a consideration of
$g_n\coloneqq |\cdot|\wedge n \in\cF_1^\infty$ for $n\in\N$ yields
\[
 \text{L.H.S.\eqref{Eq:zeta_1_under_x_maspto_ax}}
  &\ge & \lim_{n\rightarrow\infty} \int \big( |bx|\!\wedge\!n - |ax|\!\wedge\!n\big)\,\dd M(x)
  \,\ = \,\ \text{R.H.S.\eqref{Eq:zeta_1_under_x_maspto_ax}}
\]
by monotone convergence.

We similarly have ``$\le$'' in~\eqref{Eq:zeta_1_under_x_maspto_a+x} by 
specialising~\eqref{Eq:zeta_r_under_x_mapsto_ax+b_proof} to $r=1$ and there $a=c=1$.
Assuming now $M\ge 0$ and w.l.o.g. $a\le b$, 
a consideration of $g_n(x) \coloneqq (-n)\vee x\wedge n$ yields 
\[
 \text{L.H.S.\eqref{Eq:zeta_1_under_x_maspto_a+x}}
  &\ge & \lim_{n\rightarrow\infty} \int \big( g_n(x+b)-g_n(x+a) \big)\,\dd M(x) 
  \,\ = \,\ \text{L.H.S.\eqref{Eq:zeta_1_under_x_maspto_a+x}}
\]
by dominated convergence. 

Finally, \eqref{Eq:zeta_s_delta_b-delta_a} 
is~\eqref{Eq:zeta_1_under_x_maspto_a+x} in the special case of $M=\delta_0$. 
\end{proof}

\begin{Lem}[$\underline{\zeta}$ norms and convolutions]   \label{Lem:zeta_convolutions}
Let $M_1,M_2\in \cM$ and $r\in\N$. Then we have
\la                       \label{Eq:zeta_convolutions} \label{Eq:zeta_underbar_regularity}
 \underline{\zeta}_r(M_1M_2) &\le& \underline{\zeta}_r(M_1)\nu_0(M_2)
   + \sum_{j=0}^{r-1}\tfrac{1}{j!}\,|\mu_j(M_1)|\,\underline{\zeta}_{r-j}(M_2)
   \quad\text{ if } \nu_{r}(M_1) <\infty\,.
\al
\end{Lem}
\begin{proof} Let $g\in\cF^\infty_{r,r-1}$. 
We then have $\|g^{(j)}\|^{}_\infty<\infty$ for $j\in\{0,\ldots,r\}$, 
by, for example, \citet[p.~9, Theorem~1.2]{KwongZettl1992}, 
and hence with $T_y(x)\coloneqq \sum_{j=0}^{r-1}g^{(j)}(y)\frac{x^j}{j!}$,
and using just the assumption $\nu_{r-1}(M_1)<\infty$, no integrability problems
arise in verifying the first two steps below in
\[
 \left|\int g \,\dd (M_1M_2)\right| &=& \left|\iint T_y(x)\,\dd M_1(x)\dd M_2(y)
              + \iint \big(g(x\!+\!y)-T_y(x)\big)\,\dd M_1(x)\dd M_2(y) \right|  \\
  &\le& \sum_{j=0}^{r-1}\tfrac{1}{j!}\,|\mu_j(M_1)|\,\big|\int g^{(j)}\,\dd M_2\big|
    +  \int \big| \int\! \big(g(x\!+\!y)-T_y(x)\big)\dd M_1(x)    \big|\,\dd|M_2|(y) \\
  &\le& \text{R.H.S.\eqref{Eq:zeta_convolutions}}\,.
\]
In the final step we use $g^{(j)}\in\cF^\infty_{r-j,r-j-1}$
to bound the sum $\sum_{j=0}^{r-1}$\,,
and $ g(\cdot+y)-T_y \in \cF_{r,r-1}$ and the full 
assumption $\nu_r(M)<\infty$ in order to apply~\eqref{Eq:zeta_underbar_on_cM_r}
from Lemma~\ref{Lem:Comparison_of_various_eqnorms}
to bound the rest.
\end{proof}

We next provide proofs of Lemmas~\ref{Lem:F_M,k_and_int_g_dM} and 
\ref{Lem:Comparison_of_various_eqnorms} from the introduction.

\begin{proof}[Proof of Lemma~\ref{Lem:F_M,k_and_int_g_dM} from page \pageref{page:Lemma_F_M,k_and_int_g_dM}]
                                               \label{page:Proof_of_Lemma_F_M,k_and_int_g_dM}
The parts \ref{part:F_M,k_with_k=1} and \ref{part:F_M_overline_complementary} are obvious.

\smallskip\ref{part:int_g_dM_formal} 
 Identity~\eqref{Eq:int_g_dM_formal} follows from 
 integrating the Taylor formula 
\la                      \label{Eq:Taylor_formula}
 g(y)- \sum_{j=0}^{k-1}\frac{g^{(j)}(0)}{j!}y^j
  &=& \int_0^y g^{(k)}(x)\frac{(y-x)^{k-1}}{(k-1)!}\dd x \quad\text{ for }y\in\R
\al
with respect to $M$ and using Fubini, 
which is justified, with both integrals in~\eqref{Eq:int_g_dM_formal} finite,
since 
$\int \int_0^{|y|}|g^{(k)}(x)\frac{(y-x)^{k-1}}{(k-1)!}|\,\dd x \,\dd|M|(y)
 \le  \left\|g^{(k)}/(1\!+\!|\cdot|^\alpha)\right\|_\infty
 \int\int_0^{|y|}(1\!+\!|x|^\alpha)\frac{(|y|-x)^{k-1}}{(k-1)!}\,\dd x \,\dd|M|(y)
 <\infty$, as the inner integral is $\frac{|y|^k}{k!}+\frac{\mathrm{B}(\alpha+1,k)}{(k-1)!}|y|^{k+\alpha}$.  

\smallskip\ref{part:F_M,k+ell} The case of $x=0$ is trivial.
 If $x\neq 0$, then we apply \eqref{Eq:int_g_dM_formal} with $\alpha\coloneqq \ell-1$
 and with $g(y)\coloneqq \frac{(y-x)^{k+\ell-1}}{(k+\ell-1)!}(y>x)$ if $x>0$,  
 and with $g(y)\coloneqq \frac{(y-x)^{k+\ell-1}}{(k+\ell-1)!}(y<x)$ if $x<0$,
 and get~\eqref{Eq:h_M,k+ell}  by observing that in either case $g(0)=\ldots=g^{(k-1)}(0)=0$.
\end{proof}

\begin{proof}[Proof of Lemma~\ref{Lem:Comparison_of_various_eqnorms}
 from page \pageref{page:Lemma_Comparison_of_various_eqnorms}]  
                                         \label{page:proof_of_Lem_Comparison_of_various_eqnorms}
Let $r\in\N$ in steps 1--7 below.

\smallskip1. For $g\in\cF_r$, let us put 
$g_0(x)\coloneqq g(x)-\sum_{j=0}^{r-1} \frac{g^{(j)}(0)}{j!}x^j$ for $x\in\R$,
so that $g=g_0$ iff $g\in\cF_{r,r-1}$\,, and in any case $|g_0|\le \frac{1}{r!}|\cdot|^r$
by~\eqref{Eq:Taylor_formula}.

\smallskip2. 
 If $g\in\cF_r$, then there exists a sequence $(g_n)$ in $\cF^\infty_r$ with $g_n\rightarrow g$
 pointwise and, for some constants $a,b\in[0,\infty[\, $, $|g_n|\le a+b|\cdot|^r$ for each $n$,
 by \citet[p.~504, Lemma 2.2]{MattnerShevtsova}. 
 If even $g\in\cF_{r,r-1}$, then we may also take $g_n\in\cF_{r,r-1}$,
 since in the proof of the lemma just cited, where the present $g,g_n$ 
 are called $f,f_n$, any condition $f^{(k)}(0)=0$ with $k\in\N_0$ 
 obviously implies $f_n^{(k)}(0)=0$ for each~$n$.  
 
 Let now $M\in\cM_r$. If $g\in\cF_r$, then 
 with $\cF_r\ni g_n\rightarrow g$ as above, dominated convergence yields
 $|\int g\,\dd M| = \lim_{n\rightarrow\infty} |\int g_n\,\dd M| \le \zeta_r(M)$;
 hence we get $\zeta_r(M)\le \sup_{g\in\cF_r}\left| \int\! g\,\dd M\right|\le \zeta_r(M)$,
 that is, \eqref{Eq:zeta_on_cM_r_g_unbounded} holds.
 If $g\in\cF_{r,r-1}$, then, by step~1, $|g|\le \frac{1}{r!}|\cdot|^r$
 and hence $|\int g\,\dd M| \le \frac{1}{r!}\nu_r(M)$, and with 
 $\cF_{r,r-1}\ni g_n\rightarrow g$ as above we now get 
 $|\int g\,\dd M| = \lim_{n\rightarrow\infty} |\int g_n\,\dd M| \le \underline{\zeta}_r(M)$;
 hence we get the first identity in~\eqref{Eq:zeta_underbar_on_cM_r}, 
 and finiteness of $\underline{\zeta}_r\le\frac{1}{r!}\nu_r$ on $\cM_r$.
 
\smallskip 3. Let $M\in\cM$. Then trivially $\underline{\zeta}_r(M)\le\zeta_r(M)$.
 If $M\in\cM_{r,r-1}$, then we also get 
 \[
  \zeta_r(M) &\le& \sup_{g\in\cF_r} \left|\int g\,\dd M \right|
    \,\ =\,\ \sup_{g\in\cF_{r,r-1}}\left|\int g\,\dd M \right|
    \,\ =\,\ \underline{\zeta}_r(M)
 \]
 trivially in the first step (or less trivially actually with equality 
 by~\eqref{Eq:zeta_on_cM_r_g_unbounded} proved in step~2),
 using $\int g\,\dd M = \int g_0\,\dd M$ in the second,
 and by step~2 in the last. 
 This proves~\eqref{Eq:zeta_vs_zeta_underbar_on_cM_and_oncM_r,r-1}. 
 
\smallskip 4. Let $M\in\cM_r\!\setminus\!\cM_{r,r-1}$. Then $\mu_j(M)\neq 0$
 for some $j\in\{0,\ldots,r-1\}$, and with $g_t(x)\coloneqq t x^j$ for $t,x\in\R$ 
 we 
 get $\zeta_r(M) \ge \sup_{t\in\R}|\int g_t \,\dd M|=\infty$,
 using~\eqref{Eq:zeta_on_cM_r_g_unbounded} and $g_t\in\cF_r\,$.
 Hence~\eqref{Eq:zeta_vs_zeta_underbar_on_cM_and_on_cM_r_setminus_cM_r,r-1} holds. 
 
\smallskip 5. Obviously $\underline{\zeta}_r$ and $\zeta_r$ are eqnorms on $\cM$.
 If $M\in\cM$ with $\zeta_r(M)=0$, then $M=0$ for example by 
 the uniqueness theorem for Fourier transforms, considering the functions
 $g_t\coloneqq (x\mapsto t^{-r}\mathrm{e}^{\mathrm{i}tx}) \in \cF^\infty_r$ for $t\in\R\!\setminus\!\{0\}$;
 hence $\zeta_r$ is an enorm. This completes the proof of part~\ref{part:zeta_underbar_le_zeta}.

\smallskip 6. The second and the third identity in~\eqref{Eq:zeta_underbar_on_cM_r}
 follow from~\eqref{Eq:int_g_dM_formal} with $(k,\alpha)\coloneqq(r,0)$.

 The first inequality in~\eqref{Eq:zeta_underbar_vs_varkappa_etc} 
 follows from considering in the second term in~\eqref{Eq:zeta_underbar_on_cM_r} 
 the functions $\frac{1}{r!}(\cdot)^r\, ,\, \frac{1}{r!}|\cdot|^r \in\cF_{r,r-1}$. 
 
 The second inequality in~\eqref{Eq:zeta_underbar_vs_varkappa_etc} 
 is obtained via $\underline{\zeta}_r(M)= \int|h_{M,r}|\,\dd\leb$ from~\eqref{Eq:zeta_underbar_on_cM_r}:
 In case of $r=1$, the last integral is just the one defining $\varkappa_1(M)$ in~\eqref{Eq:Def_varkappa}, 
 and we hence obtain even equality. In case of $r\ge2$, we obtain
 \[
  \underline{\zeta}_r(M) &=& \int \big| h_{M,1+(r-1)}\big|\,\dd\leb \\
   &\le& \int  
       \left( (x>0)\!\int\limits_x^\infty +\ (x<0)\!\!\int\limits_{-\infty}^x \right) 
       \frac{|y-x|^{r-2}}{(r-2)!} \big|h_{M,1}(y)\big|\,\dd y 
       \,\dd x \\
   &=&  \iint \big( (0<x<y)+(y<x<0)\big) \frac{|y-x|^{r-2}}{(r-2)!}\,\dd x\  \big|h_{M}(y)\big|\,\dd y    
   \,\ = \,\ \varkappa_r(M)     
 \]
 by using in the second step~\eqref{Eq:h_M,k+ell} with $(k,\ell)\coloneqq(1,r\!-\!1)$.
 
 The final inequality in~\eqref{Eq:zeta_underbar_vs_varkappa_etc} is known from~\eqref{Eq:Def_varkappa}.
 
 This proves part~\ref{part:zeta_underbar_on_cM_r} and, using~\eqref{Eq:h_M,k=-F_M,k}, 
 also part~\ref{part:zeta_1=_varkappa_1}.
 
\smallskip 7. Part \ref{part:zeta_on_cM_r_r-1} follows from~\ref{part:zeta_underbar_on_cM_r}, 
 using~\eqref{Eq:zeta_vs_zeta_underbar_on_cM_and_oncM_r,r-1},
 and~\eqref{Eq:h_M,k=-F_M,k} with $k\coloneqq r$.

\smallskip 8. The inequalities in~\eqref{Eq:Kolmogorov_vs_nu_0_zeta_0} 
 are rather obvious and well-known, in case of the last one due to 
 $\| M \|_\mathrm{K} = \sup_{x\in\R} |\int (\1_{]-\infty,x]}-\frac{1}{2})\,\dd M|$ 
 for $M\in\cM_{0,0}$\,.
\end{proof}       
Parts of the above proof could have been replaced, less naturally, 
by references to \citet[p.~498, Theorem 1.7, 
with $P\coloneqq \frac{1}{|M|(\R)}M_+$ and $Q\coloneqq\frac{1}{|M|(\R)}M_-$]{MattnerShevtsova}.

The following Theorem~\ref{Thm:Cut_criteria} is essentially a reformulation of known results
collected or refined in \citet[Theorem 4.2, Lemma~2.8]{MattnerShevtsova}, 
and some earlier relevant references are given below after the proof.
Here the formulation is in terms of signed measures, rather than in pairs $(P,Q)$ 
corresponding to the case of $M=Q-P$, and thus seems more natural.
Also statements involving $S^-(F_0)$ are directly included 
in~(\ref{Eq:S^-_F_k_vs_F_k-1},\ref{Eq:S^-_F_k_ge_r-k}) and in the conditions $(B_k),(C_k)$.
If Theorem~\ref{Thm:Cut_criteria} is specialised to $M\coloneqq\widetilde{P}-\mathrm{N}$
and $r\coloneqq 3$, 
then its parts~\ref{part:S^-(F_k)} and~\ref{part:zeta_r=mu_r/r!} yield in particular 
Lemma~\ref{Lem:zeta_3_computation}, while part~\ref{part:zeta_r=nu_r/r!} 
is used in Examples~\ref{Example:Subbotin_laws} and \ref{Example:Zolotarev_1973}.

In the next two paragraphs, we define ``initially positive'' and the notation $S^-(f)$,
both needed for the applications of Theorem~\ref{Thm:Cut_criteria} in the present paper.
The then following three paragraphs up to the definition of ``$M \ge^{}_{r-\mathrm{cx}} 0$'' 
may be skipped here. Let $D\subseteq \R$ and let $f:D\rightarrow\R$ be a function.

$f$ is called {\em initially positive}
if either $f=0$ on $D$ or there exists an $x_0\in D$ with $f(x_0)>0$ 
and $f\ge 0$ on $D\cap \mathopen]-\infty, x_0\mathclose[$. Initial negativity 
and final positivity or negativity of $f$ are defined analogously.

The so-called {\em number of sign changes} of $f$ is defined to be
\la                                                                  \label{Eq:Def:S^-(f)}
 S^-(f) &\coloneqq& \sup\big\{ n\in\N_0 : \exists\, x \in D^{n+1}\text{ with }x_i<x_{i+1} \text{ and } \\
        &         & \phantom{\sup\big\{ n\in\N_0 :}\ f(x_i)f(x_{i+1})<0 \text{ for }i\in\{1,\ldots,n\}\big\}, \nonumber
\al
and this is either $\infty$ or is more accurately called the 
{\em maximal number of inequivalent sign change points} of $f$. 
For example, with the definition of the next paragraph, 
the function $f\coloneqq \1_{[1,\infty[} - \1_{]-\infty,-1]}$ on $\R$ has $S^-(f)=1$, 
but each $z\in[-1,1]$ is a sign change point of $f$; hence the qualifier ``inequivalent'' 
in the preceding sentence.

If now for simplicity $D$ is assumed to be a nonempty interval, then we have $S^-(f)=n\in\N_0$
iff there exists a decomposition $D=\bigcup_{j=0}^n I_j$ into nonempty
(put possibly one-point, as in the example $f\coloneqq -1+2\1_{\{0\}}$ on $\R$)
intervals $I_j$ with, for $j\in\{0,\ldots,n\}$, $f(x)f(y)\ge0$ for $x,y\in I_j$,
but in case of $j\ge 1$ also $\sup I_{j-1}\eqqcolon z_j = \inf I_j$
and $f(x)f(y)<0$ for some $x\in I_{j-1}$ and some $y\in I_j$.
In this case, such a $(z_1,\ldots,z_n)$ is called a {\em sign change tuple} of $f$,
and any of its entries a {\em sign change point}.
Two sign change points of $f$ are called {\em inequivalent} if they occur 
in a same sign change tuple.

For $r\in\N_0$,  a function $g:\R\rightarrow\R$ is called {\em $r$-convex}
if $g\ge 0$ in case of $r=0$, $g$ is increasing in case of $r=1$, 
and $g$ is $r-2$ times differentiable with $g^{(r-2)}$ convex in case if $r\ge 2$.
Standard examples are the polynomials of degree at most $r-1$
and the functions given by $g(x)=x^r$ and $g(x)=|x|^r$.
We refer to \citet{PinkusWulbert2005} and also 
\citet[p.~505]{MattnerShevtsova} for a more detailed introduction and some appropriate references. 

For $r\in\N$ and 
$M\in\cM_{r-1}$, we define $M \ge^{}_{r-\mathrm{cx}} 0$
to mean $\int g\,\dd M\ge 0$ for every $r$-convex function $g$ with $\int|g|\,\dd|M|(x)<\infty$.
In the case of $r\ge 1$ and  $M=Q-P$ with $P,Q\in \Prob_r(\R)$,
this condition is easily checked to be equivalent to the so-called $r$-convex ordering 
$P \le^{}_{r-\mathrm{cx}} Q$ considered by \citet{DenuitLefevreShaked1998}, 
as defined for example by \citet[p.~515]{MattnerShevtsova}.
Considering the polynomials of degree at most $r-1$, one observes that $M \ge^{}_{r-\mathrm{cx}} 0$
implies $\mu_k(M)=0$ for $k\in\{0,\ldots,r-1\}$, that is, $M\in\cM_{r-1,r-1}$.

We recall the notation~\eqref{Eq:Def_F_M,k}, in particular
$F_{M,1}(x)=F_M(x)=M(\,]-\infty,x])$ for $M\in\cM$ and $x\in\R$.
We also recall that, by the Radon-Nikod\'ym theorem, the assumption 
$M=f\mu$ below is always fulfilled with, for example, $\mu\coloneqq|M|$.

\begin{Thm}[Cut criteria for computing $\zeta$ norms]   \label{Thm:Cut_criteria}
Let $r\in\N$ and $M\in\cM_{r-1,r-1}\,$, and let $F_k\coloneqq F_{M,k}$ be defined 
by~\eqref{Eq:Def_F_M,k} for $k\in\{1,\ldots,r\}$. 
Let further $M=f\mu$ for some positive measure~$\mu$ and a $\mu$-integrable
$\R$-valued function~$f$,
and let us write here $F_0\coloneqq -f$. 

\begin{parts}
\item We have    \label{part:S^-(F_k)}
 \la
  S^-(F_k) &\le& S^{-}(F_{k-1}) - 1 \quad\text{ for }\,k\in\{1,\ldots, r\},\label{Eq:S^-_F_k_vs_F_k-1} \\
  S^-(F_k) &\ge& r\!-\!k \ \text{ or }\ M = 0 \quad\text{ for }\,k\in\{0,\ldots, r\}. \label{Eq:S^-_F_k_ge_r-k}
 \al
\item  For $k\in\{0,\ldots, r\}$ let $(B_k)$ be the condition defined by  \label{part:zeta_r=mu_r/r!}
 \[
  (B_k) &:\iff& S^-(F_k) \,\le\, r\!-\!k \ \text{ and }\ (-1)^kF_k\text{ is initially positive}.
 \]
 Then we have the implications 
 \la
  (B_0) &\Rightarrow& (B_1) \,\ \Rightarrow\,\  \ldots \,\ \Rightarrow\,\ (B_r)
   \,\ \iff \,\ (-1)^rF_r \ge 0 \,\ \iff\,\ M  \ge^{}_{r-\mathrm{cx}} 0\,.  
 \al
 If even $M\in\cM_{r,r-1}\,$, then we further have
 \la
  M \ge^{}_{r-\mathrm{cx}} 0 &\iff& \zeta_r(M) =\tfrac{1}{r!}\mu_r(M)\,.
 \al
\item For $k\in\{0,\ldots, r\}$ let $(C_k)$ be the  \label{part:zeta_r=translated_nu_r/r!}
  condition defined by
 \[
 (C_k) &:\iff& S^-(F_k) \,=\, r\!-\!k\!+\!1 \ \text{ and }\ (-1)^kF_k\text{ is initially positive}.
 \]
 Then
 \la
  (C_0) &\Rightarrow& (C_1) \,\ \Rightarrow\,\  \ldots \,\ \Rightarrow\,\ (C_r) .  
 \al
 If even $M\in\cM_{r,r-1}\,$, then we further have 
 \la
  (C_r) &\iff& \zeta_r(M) =\tfrac{1}{r!}\int|x-x_0|^r\,\dd M(x)
     \quad \text{for some sign change point $x_0$ of }F_r\,,
 \al
 and this remains true with ``some'' replaced by ``some and every''.
 Further, if $(C_k)$ holds for some $k\in\{0,\ldots,r-1\}$,
 then each sign change point of $F_r$ belongs to the interior 
 of the convex hull of the set of the entries of each sign change tuple of $F_k$.
\item  Suppose that even $M \in\cM_{r,r-1}$,  $M$ is symmetric, \label{part:zeta_r=nu_r/r!}
 and $(C_r)$ from part~\ref{part:zeta_r=translated_nu_r/r!} holds. 
 Then $r$ is odd, and $\zeta_r(M)= -\frac{1}{r!}\int|x|^r\dd M(x)$.
\end{parts}
\end{Thm}
\begin{proof} Follows from \citet[Lemma~2.8 and Theorem 4.2(b,c,d)]{MattnerShevtsova}
and Remark~\ref{Rem:Correction_of_Mattner_Shevtsova_4.2(d)} below,
with some obvious modifications.
For example, if $M\neq0$, then the cited theorem may be applied to 
$P\coloneqq \frac{M_-}{|M|(\R)}$ and $Q\coloneqq \frac{M_+}{|M|(\R)}$. 
\end{proof}

\begin{Rem}                  \label{Rem:Correction_of_Mattner_Shevtsova_4.2(d)}
In \citet[pp.~517--518, Theorem 4.2(d)]{MattnerShevtsova}, the assumption
``$\overline{H}_s$ lastly positive'' is missing in the statement
but used in the proof.
\end{Rem}

Theorem~\ref{Thm:Cut_criteria} is an instance of
refined \citet{KarlinNovikoff1963} type cut criteria 
as presented by \citet{DenuitLefevreShaked1998}, \citet{BoutsikasVaggelatou2002},
and \citet{MattnerShevtsova}. 
We have to note here that in these papers partial priority should have been acknowledged 
to \citet[in particular section~2]{vonMises1937}.

\section{A proof of Zolotarev's $\zeta_1\!\vee\!\zeta_3$
    Theorem~\ref{Thm:Zolotarev's_zeta_1-B-E-Thm}}   \label{sec:Proof_Zolotarev-Thm}
We proceed to proving Theorem~\ref{Thm:Zolotarev's_zeta_1-B-E-Thm},
following \citet[pp.~365--368]{Zolotarev1997} in
using just the simple properties
of $\zeta_0,\zeta_1,\zeta_3$ 
from 
Lemmas~\ref{Lem:Regular_eqnorms_on_cM},\,\ref{Lem:zeta_1-smoothing},\,\ref{Lem:zeta_s_vs_zeta_(s+k)}
in a not very complicated inductive argument.
Merely for obtaining the value $c=13.3803\ldots$ defined in~\eqref{Eq:Def_c=13.3803...} below, 
we also use the following nontrivial result:

\begin{Thm}[Goldstein, Tyurin, 2010]                            \label{Thm:Goldstein-Tyurin}
We have 
\la            \label{Eq:Goldstein-Tyurin}
  \zeta_1\big(\widetilde{P^{\ast n}}-\mathrm{N}\big)  &\le& \frac{1}{\sqrt{n}}\nu_3(\widetilde{P})
   \quad\text{ for }P\in\cP_3 \text{ and }n\in\N\,,
\al
with the constant~$1$ on the right hand side not reducible beyond
\la                        \label{Eq:Goldstein-Tyurin_0.53}
  \zeta_1\big(\widetilde{\mathrm{B}_{\frac12}}-\mathrm{N}\big) 
    &=& 4\Phi(1)+4\phi(1) -2\phi(0) - 3 
    \,\ =\,\   0.535377\ldots\,.
\al
\end{Thm}
\begin{proof} Inequality~\eqref{Eq:Goldstein-Tyurin} is a special case, for identical
convolution factors, of apparently independently obtained theorems
of \citet[Theorem 1.1]{Goldstein2010}
and~\citet[\foreignlanguage{russian}{Teorema} 4]{Tyurin2010}.
The former paper also contains the remark involving~\eqref{Eq:Goldstein-Tyurin_0.53}.
\end{proof}

In the above ``proof'' we have cited the first peer-reviewed publications,
of their respective authors, containing complete proofs of the result in question,
thus justifying in some sense the 2010 in our caption of Theorem~\ref{Thm:Goldstein-Tyurin}.
For prepublications and submission dates one may consult \citet[p.~1688]{Goldstein2010}
and \citet[p.~1]{Tyurin2009arxiv}. The latter paper actually contains improvements
compared to \citet{Tyurin2010}, but apparently not so with respect to 
Theorem~\ref{Thm:Goldstein-Tyurin}.
For completeness let us mention that 
\citet[p.~1674, line~4, the claim ``$c_\infty=1/2$'']{Goldstein2010}
appears there without any justification of an apparent interchange of a limit with a supremum. 

\begin{proof}[Proof of Theorem~\ref{Thm:Zolotarev's_zeta_1-B-E-Thm}] 1.     \label{Start_of_proof_Zolotarev's_zeta_1-B-E-Thm} 
Let $P\in\cP_3$ with $\widetilde{P}\neq\mathrm{N}$,
and let $\xi_0\coloneqq\xi(\zeta_1(\widetilde{P}-\mathrm{N}),
\zeta_3(\widetilde{P}-\mathrm{N}))$. Let $n\in\N$ be such that we have
\la                                                    \label{Eq:inductive_hypothesis}
  \sqrt{k}\zeta_1(\widetilde{P^k}-\mathrm{N}) &\le& \xi_0 \quad\text{ for }k\in\{1,\ldots,n-2\}
\al
(yes, in this inductive proof, the validity of the inequality in~\eqref{Eq:inductive_hypothesis}
for $k=n-1$ is not used for establishing it for $k=n$). 
We are going to prove that we then also have 
\la                                                    \label{Eq:inductive_conclusion}
  \sqrt{n}\zeta_1(\widetilde{P^n}-\mathrm{N}) &\le& \xi_0\,.
\al
To this end, we may assume w.l.o.g.~$\mu(P)=0$ and $\sigma(P)=\frac{1}{\sqrt{n}}$. 
We put $Q\coloneqq\mathrm{N}_{\frac{1}{\sqrt{n}}}$. 
If further $\epsilon\in[0,\infty[$,  then we get, using about $\xi_0$ 
initially only that it is some number satisfying~\eqref{Eq:inductive_hypothesis},
\[
 \zeta_1(\widetilde{P^n}-\mathrm{N}) &=& \zeta_1(P^n-Q^n) \,\ \le \,\  
   \zeta_1(P^n\mathrm{N}_\epsilon-Q^n\mathrm{N}_\epsilon)  + \beta\epsilon  \\
 &\le&  \beta\epsilon + \zeta_1(P^n\mathrm{N}_\epsilon-P^{n-1}Q\mathrm{N}_\epsilon) 
    + \sum_{j=1}^{n-1} \zeta_1(P^{n-j}Q^j\mathrm{N}_\epsilon-P^{n-j-1}Q^{j+1}\mathrm{N}_\epsilon) \\
 &\le& \beta\epsilon +\zeta_1(P-Q)  \\
  && +\ \sum_{j=1}^{n-1}\left( 
     \zeta_1(PQ^{n-1}\mathrm{N}_\epsilon-Q^n\mathrm{N}_\epsilon) 
      +\zeta_1(P^{n-j-1}\!-Q^{n-j-1})\zeta_0(PQ^j\mathrm{N}_\epsilon-Q^{j+1}\mathrm{N}_\epsilon)   \right) \\
 &\le& \beta\epsilon + \zeta_1(P-Q) + (n-1)\alpha\frac{\zeta_3(P-Q)}{\frac{n-1}{n}+\epsilon^2} 
  +\sum_{j=1}^{n-2}\frac{\xi_0}{\sqrt{n}} \, \gamma\,\frac{\zeta_3(P-Q)}{(\frac{j}{n}+\epsilon^2)^{3/2}} \\
 &\le& \beta\epsilon +\frac{\zeta_1(\widetilde{P}-\mathrm{N})}{\sqrt{n}}
      + \alpha\frac{\zeta_3(\widetilde{P}-\mathrm{N})}{\sqrt{n}}
      + \gamma\frac{\xi_0\cdot\zeta_3(\widetilde{P}-\mathrm{N})}{\sqrt{n}}
      \sum_{j=1}^{n-2}\frac{1}{(j+n\epsilon^2)^{3/2}} 
\]
by using in the second step Lemma~\ref{Lem:zeta_1-smoothing}, in the third just the triangle inequality for $\zeta_1$
applied to~\eqref{Eq:telescope_P_j_Q_j} times $\mathrm{N}_\epsilon$ with $P_j\coloneqq P$ and $Q_j\coloneqq Q$,
in the fourth from Lemma~\ref{Lem:Regular_eqnorms_on_cM} the regularity~\eqref{Eq:regularity_on_Prob} applied 
to $R\coloneqq P^{n-1}\mathrm{N}_\epsilon$ and~\eqref{Eq:changing_the_regularizer}  
applied to $M_1\coloneqq  PQ^j\mathrm{N}_\epsilon- Q^{j+1}\mathrm{N}_\epsilon$, 
$M_2\coloneqq P^{n-j-1}$, $M_3\coloneqq Q^{n-j-1}$
and taking the second minimand,
in the fifth Lemma~\ref{Lem:zeta_s_vs_zeta_(s+k)} with $(s,k)=(1,2)$ and with $(s,k)=(0,3)$,
and also the homogeneity~\eqref{Eq:homogeneity_of_zeta} of $\zeta_1$ and the inductive 
hypothesis~\eqref{Eq:inductive_hypothesis} in order to get
\[
  \zeta_1(P^{n-j-1}-Q^{n-j-1}) &=& \sqrt{\frac{n-j-1}{n}} \zeta_1(\widetilde{P^{n-j-1}}- \mathrm{N})
  \,\ \le \,\ \frac{\xi_0}{\sqrt{n}}\quad\text{ for }j\in\{1,\ldots,n-2\},
\]
and in the final sixth step the homogeneity~\eqref{Eq:homogeneity_of_zeta}
of $\zeta_1$ and of $\zeta_3$.

Hence, if $\eta\in[0,\infty[$, we get by applying the above to $\epsilon\coloneqq\frac{\eta}{\sqrt{n}}$,
and by recalling the definition of $g$ from~\eqref{Eq:Def_g(eta)}, 
\[
 \sqrt{n} \zeta_1(\widetilde{P^n}-\mathrm{N}) 
  &\le& \zeta_1(\widetilde{P}-\mathrm{N}) + \alpha\zeta_3(\widetilde{P}-\mathrm{N}) + \beta\eta 
  +\gamma g(\eta)  \zeta_3(\widetilde{P}-\mathrm{N})\xi_0 \\
  & \eqqcolon& A(\eta)+B(\eta)\xi_0\,,
\]
and hence, with $\Eta\coloneqq\{\eta\in[0,\infty[\,:\,B(\eta)<1\}$ and now using the definition
of $\xi_0$ through the function $\xi$ from~\eqref{Eq:Def_xi(kappa,zeta)}, we get 
$\xi_0 = \inf_{\eta'\in\Eta} \frac{A(\eta') }{1-B(\eta')}$ and hence
\[
 \sqrt{n} \zeta_1(\widetilde{P^n}-\mathrm{N})  
  &\le& \inf_{\eta,\eta'\in\Eta} \left(A(\eta)+B(\eta)\frac{A(\eta')}{1-B(\eta')}\right) \\
  &\le& \inf_{\eta\in\Eta}  \left(A(\eta)+B(\eta)\frac{A(\eta)}{1-B(\eta)}\right)         \,\ =\,\ \xi_0\,,
\]
that is, \eqref{Eq:inductive_conclusion}.
This proves \eqref{Eq:Zolotarev's_zeta_1-B-E-Thm_sharper}.

\smallskip 2. 
The inequality in~\eqref{Eq:Def_g(eta)} of course follows from $(j+\eta^2)^{-3/2}<\int_{j-1}^j(x+\eta^2)^{-3/2}\,\dd x$.
Using it in the first step below, and $\eta\coloneqq 4\gamma\zeta$ in the second, yields
\[                                           
  \xi(\varkappa,\zeta) &\le& \inf\left\{\frac{\varkappa+\alpha\zeta+\beta\eta}{1-\frac{2\gamma\zeta}{\eta}} : \eta\in[0,\infty[,\
   \frac{2\gamma\zeta}{\eta} <1 \right\}  \\
   &\le& 2\varkappa + 2(\alpha+4\beta\gamma)\zeta 
    \,\ =\,\  2\varkappa + 21.212827\ldots \zeta   
\]
for $(\varkappa,\zeta)\in[0,\infty[^2$, and hence~\eqref{Eq:Zolotarev's_zeta_1-B-E-Thm} with 
$c=2+2(\alpha+4\beta\gamma) = 23.212827\ldots$\;.     \label{End_of_step_2_in_proof_Zolotarev's_zeta_1-B-E-Thm}

\smallskip 3. Considering $\eta\coloneqq0$ in~\eqref{Eq:Def_xi(kappa,zeta)} and using   
$g(0)=\zeta(\frac{3}{2})=2.612375\ldots$ yields~\eqref{Eq:xi_near_zero}.

\smallskip 4.
Using below in the first step the Goldstein-Tyurin inequality~\eqref{Eq:Goldstein-Tyurin},
in the second the first inequality in~\eqref{Eq:zeta_vs_varkappa_etc}
for $M\coloneqq \widetilde{P}-\mathrm{N}$,
and in the third~\eqref{Eq:nu_1_and_nu_3_of_N},
we get 
\la &&
 \zeta_1\big(\widetilde{P^{\ast n}},\mathrm{N}\big) \,\ \le\,\ \frac{\nu_3(\widetilde{P})}{\sqrt{n}}
  \,\ \le \,\   \frac{6\zeta_3(\widetilde{P}-\mathrm{N})+ \nu_3(\mathrm{N})}{\sqrt{n}}
   \,\ = \,\   \frac{6\zeta_3(\widetilde{P}-\mathrm{N})+ \beta}{\sqrt{n}}  \quad\text{ for }P\in\cP_3,    
\al
and hence  \eqref{Eq:Zolotarev's_zeta_1-B-E-Thm} holds with 
\la                                                \label{Eq:Def_c=13.3803...}
   c &\coloneqq& \sup_{\varkappa,\zeta>0} \frac{1}{\varkappa\!\vee\!\zeta}
   \left( \frac{\varkappa +\alpha\zeta}{\left(1-\lambda\zeta\right)_+} \wedge \Big(6\zeta + \beta \Big)   \right)
   \,\ =\,\ \sup_{\zeta>0} \frac{1+\alpha}{\left(1-\lambda\zeta\right)_+} \wedge \left(6+\frac{\beta}{\zeta}\right).
\al
The last supremum above is uniquely attained at the positive solution $\zeta^\ast$,
automatically $<$~$\frac{1}{\lambda}$,
of the quadratic equation
\[
          \frac{1+\alpha}{1-\lambda\zeta} &=& 6+\frac{\beta}{\zeta}
\]
for $\zeta$.
We hence get $c$ as claimed through
\[
  \zeta^\ast &=& -\frac{\lambda\beta+\alpha-5}{12\lambda}
  + \sqrt{  \left(\frac{\lambda\beta+\alpha-5}{12\lambda}\right)^2 + \frac{\beta}{6\lambda}}
  \,\ = \,\ 0.216219\ldots  \,, \\
 c &=& 6+\frac{\beta}{\zeta^\ast} \,\ =\,\ 13.3803\ldots\,.
\]                                  \label{End_of_proof_Zolotarev's_zeta_1-B-E-Thm} 
\vspace{-2\baselineskip}

\end{proof}

In steps 3 and 4 above, step 2 was not used, and step 1 only in the slightly simpler special case of $\epsilon=0$, 
but the general case of step 1 makes the proof of Theorem~\ref{Thm:Zolotarev's_zeta_1-B-E-Thm} just up 
to~\eqref{Eq:Zolotarev's_zeta_1-B-E-Thm} self-contained. A plot of 
$\eta\mapsto(\varkappa + \alpha\zeta +\beta\eta)/(1-\gamma g(\eta)\zeta)$ with $\varkappa \coloneqq \zeta \coloneqq\zeta^\ast$
suggest that no improvement upon $c=13.3803\ldots$ seems possible using just the present ideas. 
In particular, it does not seem to help to modify the definition of $c$ in~\eqref{Eq:Def_c=13.3803...} taking into account 
that $\cP_2\ni P \mapsto \zeta_1(\widetilde{P}-\mathrm{N})$ is actually bounded, 
since the obvious bound 
$\zeta_1(\widetilde{P}-\mathrm{N})  \le
  1+ \frac{2}{\sqrt{2\pi}} = 1.79788\ldots
  $
from~\eqref{Eq:zeta_1_le}
is irrelevant due to $\zeta^\ast$ being much smaller, 
and since here $1+\frac{2}{\sqrt{2\pi}}$ can surely not be improved 
beyond  $\zeta_1(\widetilde{\mathrm{B}_p}-\mathrm{N})$ for any $p$,
and the simplest choice of $p=\frac{1}{2}$
yields by~\eqref{Eq:Goldstein-Tyurin_0.53} the value $0.535377\ldots > \zeta^\ast$.

\section{Auxiliary results for $\varkappa$ distances}                          \label{sec:varkappa_distances}
In this section, which is admittedly of only marginal importance 
in the present paper, we first provide the simple
Lemma~\ref{Lem:varkappa_scale_power_transformations}, which is 
used in Example~\ref{Example:Tail-discretised_normal_laws}.
We recall our notation~\eqref{Eq:Def_image_measure} for image measures.

\begin{Lem}[$\varkappa_r$ and scale or power transformations] 
                                                    \label{Lem:varkappa_scale_power_transformations}
Let $M\in\cM$ and $r\in\mathopen]0,\infty\mathclose[\,$,
and let us write $T_s(x)\coloneqq \sgn(x)|x|^s $ for $s\in\mathopen]0,\infty\mathclose[$ and $x\in\R$.

\begin{parts}
\item{\rm Scalings.}            \label{part:varkappa_scalings}
Let $a,b\in\R$. Then
\la
  \varkappa_r\big( (x\mapsto ax)\im M\big) &\!=\!& |a|^r\varkappa_r(M)\,, \label{Eq:varkappa_scaling} \\
  \varkappa_r\big((x\mapsto bx)\im M -  (x\mapsto ax)\im M \big) \label{Eq:varkappa_different_scalings}
    &\!\le\! & \big|T_r(b)-T_r(a)\big|\,\nu_r(M) \,,      \\
 \qquad\ \varkappa_r\big((x\mapsto bx)\im M -  (x\mapsto ax)\im M \big) \label{Eq:varkappa_different_scalings_=}
    & \!=\! &  \big| |b|^r-|a|^r\big|\,\nu_r(M) \quad\text{if }ab\ge0\text{ and }0\le M\in\cM_r \,. 
\al
\item{\rm Power transformations.}                \label{part:varkappa_power_transformed}
Let also $s\in\mathopen]0,\infty\mathclose[\,$. Then
\la
  \varkappa_r(T_s\im M) &=& \varkappa_{rs}(M) \,, \label{Eq:varkappa_r_power_s}  \\ 
  \varkappa_r(M) &=& \varkappa_1(T_r\im M)  \,.  \label{Eq:varkappa_1_power_r}
\al
\end{parts}
\end{Lem}
\begin{proof} 
\ref{part:varkappa_power_transformed} Recalling the definitions (\ref{Eq:Def_h_M},\ref{Eq:Def_varkappa}),
we observe that $h_{T_s\im M} = h_M\!\circ\!T_s^{-1}$ and 
hence, by the obvious change of variables in the second step below,
\[
 \text{L.H.S.\eqref{Eq:varkappa_r_power_s}}
  &=& \int r|x|^{r-1}\big|h_M(T^{-1}_s(x))\big|\,\dd x
   \,\ =\,\ \int rs|x|^{rs-1}\big|h_M(x)\big|\,\dd x
   \,\ =\,\ \text{R.H.S.\eqref{Eq:varkappa_r_power_s}}\,.
\]
Identity~\eqref{Eq:varkappa_1_power_r} follows by letting  $(1,r)$ 
play the role of $(r,s)$ in~\eqref{Eq:varkappa_r_power_s}. 

\smallskip\ref{part:varkappa_scalings}
Let us write here $S_c(x)\coloneqq cx$ for $c,x\in\R$.

If $a=0$, then $S_a\im M = M(\R)\delta_0$, hence $h_{S_a\im M} =0$,
and hence $\text{L.H.S.\eqref{Eq:varkappa_scaling}}=0=\text{R.H.S.\eqref{Eq:varkappa_scaling}}$.
If $a\neq0$, then $h_{S_a\im M}=\sgn(a) h_M\circ S_a^{-1}$, and 
hence~\eqref{Eq:varkappa_scaling} follows by a scale change of variables.

If $\nu_r(M)=\infty$, then $\text{R.H.S.\eqref{Eq:varkappa_different_scalings}}$
is finite only if $a=b$, in which case $\text{R.H.S.\eqref{Eq:varkappa_different_scalings}}=0$.
Hence we may assume $M\in\cM_r$ in proving~\eqref{Eq:varkappa_different_scalings},
and with $T_r\circ S_c = S_{T_r(c)}\circ T_r$ for $c\in\R$ we then get
\[
 \text{L.H.S.\eqref{Eq:varkappa_different_scalings}}
  &=& \underline{\zeta}_1\big( T_r\im ( S_b\im M- S_a\im M)\big)   
  \,\ =\,\ \underline{\zeta}_1\big( S_{T_r(b)}\im(T_r\im M) - S_{T_r(a)}\im(T_r\im M)\big) \\
  &\le & \big| T_r(b)- T_r(a)\big|\,\nu_1(T_r\im M) 
  \,\ = \,\ \text{R.H.S.\eqref{Eq:varkappa_different_scalings}}
\]
by (\ref{Eq:varkappa_1_power_r},\ref{Eq:zeta_1=kappa_1}) in the first step, 
the linearity and associativity properties of forming image measures in the second,
and~\eqref{Eq:zeta_r_under_x_mapsto_ax} in the third.
In case of $ab\ge 0$ and $M\ge 0$, we have equality everywhere in the above, 
by~(\ref{Eq:zeta_1=kappa_1},\ref{Eq:zeta_1_under_x_maspto_ax}),
and hence~\eqref{Eq:varkappa_different_scalings_=} holds. 
\end{proof}

Of the above, at least~\eqref{Eq:varkappa_1_power_r} is not only simple but also
well known, namely stated by \citet[p.~70, $\kappa_s(X,Y)$]{Zolotarev1997}.

Next our goal is to derive Lemmas~\ref{Lem:Kolmogorov_bounded_by_varkappa_and_Lipschitz}--\ref{Lem:varkappa_1_vs_varkappa_r_and_Lip}. 
Lemma~\ref{Lem:Kolmogorov_bounded_by_varkappa_and_Lipschitz} is used for
deducing~\eqref{Eq:Berry-Esseen_for_Z-close_to_normal_also_for_n=1} 
from~\eqref{Eq:Berry-Esseen_for_Z-close_to_normal}, 
and, 
with Lemma~\ref{Lem:varkappa_r_vs_varkappa_s_and_Kolmogorov}, 
for proving Lemma~\ref{Lem:varkappa_1_vs_varkappa_r_and_Lip}.
The latter is used for showing that~\eqref{Eq:Problem_with_d_n} with~\eqref{Eq:Ulyanov1978-distance}
is i.c.f.~worse than~\eqref{Eq:Problem_with_d_n} with~\eqref{Eq:Ulyanov1976-distance}.
Lemma~\ref{Lem:varkappa_r_vs_varkappa_s_and_Kolmogorov} is finally used
in the proof of Remark~\ref{Rem:CLT_rate_varkappa_r}.

While we only need in the present paper inequalities for $\varkappa_r(M)$ with $M=P-Q$ where $P,Q\in\Prob(\R)$,
and in fact $Q=\mathrm{N}$, it appears natural to consider more generally $Q\in\cM$ with $Q(M)=1$, 
called {\em signed laws} in the captions of Lemmas~\ref{Lem:Kolmogorov_bounded_by_varkappa_and_Lipschitz}
and~\ref{Lem:varkappa_1_vs_varkappa_r_and_Lip},
since this might become useful in connection with error bounds for Edgeworth approximations as provided,
although with the strong norm distances $\nu_r(\widetilde{P}-\mathrm{N})$ but not yet 
instead with $\varkappa_r(\widetilde{P}-\mathrm{N})$
or better quantities, by \citet{Yaroslavtseva2008a,Yaroslavtseva2008b}.

For Lemmas~\ref{Lem:Kolmogorov_bounded_by_varkappa_and_Lipschitz} and~\ref{Lem:varkappa_1_vs_varkappa_r_and_Lip}
we recall the definition~\eqref{Eq:Def_M_Lipschitz} of the Lipschitz constant
$\|M\|^{}_{\mathrm{L}}$ of an $M\in\cM$.
We further recall the generalised signed moments $\lambda_r$ from
Lemma~\ref{Lem:Generalised_signed_moments}, and also
write $\lambda_r(F_M)\coloneqq \lambda_r(M)$
for $M\in\cM$ with $\varkappa_r(M)<\infty$.

\begin{Lem}[Kolmogorov bounded by $\varkappa$ distances of laws 
             to signed Lipschitz laws]
                   \label{Lem:Kolmogorov_bounded_by_varkappa_and_Lipschitz}
Let $r\in\mathopen]0,\infty\mathclose[$ and 
\la                                                 \label{Eq:c_r_for_Kolmogorov_bounded_by_varkappa_and_Lipschitz}
  && c_r\ \coloneqq\ \left(r\min_{a\in\R}\int_0^1|x-a|^{r-1}2\,(1-x)\,\dd x\right)^{-\frac{1}{r+1}} 
   \ =\  \left(r\min_{a\in[0,1]} 2\left(a^r+\tfrac{(1-a)^{r+1}-a^{r+1}}{r+1}\right)\right)^{-\frac{1}{r+1}}  ,
\al
so that $2^\frac{r-1}{r+1}\vee\left(\frac{r+1}{2r}\right)^\frac{1}{r+1}\le c_r <\infty$, 
$c_1=1$, $c_2 =(\frac32+\frac{3}{4}\sqrt{2})^\frac13=1.36809\ldots$, 
$c_3=6^\frac{1}{4}=1.56508\ldots\;$. 
Then for 
$Q\in \cM$ Lipschitz with 
$Q(\R)=1$, we have
\la                                              \label{Eq:Kolmogorov_bounded_by_varkappa_and_Lipschitz}                                             
 \qquad \left\|P-Q\right\|_{\mathrm{K}} 
  &\le&  c_r \|Q\|_{\mathrm{L}}^\frac{r}{r+1}\left(\varkappa_r(P-Q)
  + \big|\lambda_r(P-Q)\big| \,\right)^\frac{1}{r+1} \\
  & \le & 2^\frac{1}{r+1}c_r \|Q\|_{\mathrm{L}}^\frac{r}{r+1} \varkappa_r(P-Q)^\frac{1}{r+1}  
   \quad\text{ for }P \in\Prob(\R) \text{ with } \varkappa_r(P-Q)<\infty. \nonumber
\al
In particular, 
\la                              \label{Eq:Kolmogorov_vs_varkappa_1_near_N}
  \left\|P-\mathrm{N}\right\|_{\mathrm{K}}
  &\le&       (2\pi)^{-1/4}\sqrt{\varkappa_1(P-\mathrm{N})} 
   \quad\text{for }P\in\Prob_1(\R) \text{ with }\mu(P)=0. 
\al
\end{Lem}
\begin{proof}
 Let $Q$ and $P$ be as stated, and $L\coloneqq\|Q\|_{\mathrm{L}}$.
 We may assume $P\neq Q$ and choose $x_0\in\R$ and distribution functions $F,G$ 
 of $P,Q$ or of their reflections in such a way  that  
 $\varrho\coloneqq\left\|P-Q\right\|_{\mathrm{K}} = F(x_0)-G(x_0)$.
 With  $a\coloneqq-\frac{Lx_0}{\varrho}$ we then get
 \[
  \varkappa_r(P-Q) &=& \int r|x|^{r-1}\big|F(x)-G(x)\big|\,\dd x           \\  
   &=&  2r\int |x|^{r-1}\big(F(x)-G(x)\big)_+\dd x   + \lambda_r(F-G)  \\
   &\ge& 2r \int_{x_0}^{x_0+\frac{\varrho}{L}}|x|^{r-1}\big(\varrho-L\,(x-x_0)\big)\dd x - |\lambda_r(F-G)|  \\
   &=& \frac{\varrho^{r+1}}{L^r} r \int_0^1|x-a|^{r-1} 2(1-x)\dd x - \left|  \lambda_r(P-Q) \right|
 \]
 by using~\eqref{Eq:varkappa_if_M(IR)=0} in the first step,  $|y|=2y_+-y$ in the second, 
 isotonicity of $F$ and the Lipschitz property of $G$ in the third,
 and the change of variables $x\mapsto \frac{\varrho}{L}x+x_0$ 
 and the reflection invariance of $|\lambda_r|$ in the fourth.
 This yields the first inequality in~\eqref{Eq:Kolmogorov_bounded_by_varkappa_and_Lipschitz},
 and the second follows from $|\lambda_r|\le\varkappa_r$.
 
 The alternative representation of $c_r$ in~\eqref{Eq:c_r_for_Kolmogorov_bounded_by_varkappa_and_Lipschitz}
 results from computing the integral in the definition, say starting with an integration by parts.
 The stated lower bound for $c_r$ follows from considering $a=\frac12$ and $a=1$ in, say, 
 the alternative representation.

 For $r\in\{1,2,3\}$, the minimum in 
 the definition of $c_r$ is, 
 respectively, the integral $1$, the mean absolute deviation $\frac{2-\sqrt{2}}{3}$ from the median $1-\frac{1}{\sqrt{2}}$, 
 and the variance $\frac1{18}$ of the probability density $[0,1]\ni x\mapsto 2\,(1-x)$, which yields the stated 
 values for $c_1,c_2,c_3$. 
 
 For $Q=\mathrm{N}$, we have $L=\frac{1}{\sqrt{2\pi}}$, and 
 so~\eqref{Eq:Kolmogorov_bounded_by_varkappa_and_Lipschitz} with $r=1$
 yields~\eqref{Eq:Kolmogorov_vs_varkappa_1_near_N}.
\end{proof}

In Lemma~\ref{Lem:Kolmogorov_bounded_by_varkappa_and_Lipschitz} with $r=1$,
the first bound in~\eqref{Eq:Kolmogorov_bounded_by_varkappa_and_Lipschitz}
improves \citet[p.~528, proof of Theorem~C]{Erickson1974}
and \citet[p.~2811, Corollary, (1)]{Shiganov1987},
and even just the second one improves 
\citet[p.~353, Proposition~2(iii) with $s=1$]{BoutsikasVaggelatou2002}
attributed there to \citet{RachevRueschendorf1991} or \citet{Rachev1991}, 
%
 and \citet[p.~47, Theorem~3.3]{ChenGoldsteinShao2011}. 

\begin{Lem}[$\varkappa_r$ versus $\left\|\cdot\right\|_\mathrm{K}$ and $\varkappa_s$ on $\cM$]
                                        \label{Lem:varkappa_r_vs_varkappa_s_and_Kolmogorov}
Let $0<r<s<\infty$. Then
\la                                               \label{Eq:varkappa_r_vs_varkappa_s_and_Kolmogorov}
 \varkappa_r(M) &\le& 2^{1-\frac{r}{s}}
   \left\|M\right\|_{\mathrm{K}}^{1-\frac{r}{s}} \varkappa_s(M)^{\frac{r}{s}}
   \quad\text{ for }M\in\cM,
\al
with finite equality iff for some $c,t\in[0,\infty[$ we have 
$|h_M|= c\,\1_{[-t,t]\setminus\{0\}}$ with~\eqref{Eq:Def_h_M},
in case of $M(\R)=0$ equivalently  $|F_M|=c\,\1_{[-t,t[}$. 
\end{Lem}
\begin{proof}
Let $M\in\cM$, without loss of generality neither 
a multiple of $\delta_0$ nor  
$\varkappa_s(M)=\infty$. Then, for every $t\in\mathopen]0,\infty\mathclose[$, we get,
using in the third step the positivity of the integrand and 
$|h_M(x)| \le \left\|M\right\|_{\mathrm{K}}$ 
by~(\ref{Eq:Def_Kolmogorov-norm},\ref{Eq:Def_h_M}),
\[                                                          
  \varkappa_r(M) &\le& \int_{|x|\le t} r|x|^{r-1} \left|h_M(x)\right|\dd x
         +\int_{|x|>t} \frac{s|x|^{s-1}}{\frac{s}{r}t^{s-r}} \left|h_M(x)\right|\dd x \\
   &=&  \int_{|x|\le t} \left(r|x|^{r-1} -\tfrac{r}{s}t^{r-s}s|x|^{s-1}\right) |h_M(x)|\,\dd x 
     + \tfrac{r}{s}t^{r-s} \varkappa_s(M) \\
   &\le& 
         \left\|M\right\|_{\mathrm{K}} 2 \tfrac{s-r}{s}t^r + 
          \tfrac{r}{s}t^{r-s} \varkappa_s(M) 
\] 
with equality throughout iff, using the one-sided continuity properties of $h_M$, 
we have $|h_M|=c\,\1_{[-t,t]\setminus\{0\}}$ for some $c\in\mathopen]0,\infty\mathclose[$. 
Minimising the bound at 
$t= \left(\frac{\varkappa_s(M)}{2\left\|M\right\|_{\mathrm{K}}}\right)^{\frac{1}{s}}$
yields the claim.
\end{proof}

Lemma~\ref{Lem:varkappa_r_vs_varkappa_s_and_Kolmogorov} improves 
\citet[Theorem 1]{Zolotarev1979}, which is misstated in 
\citet[p.~74, Remark 1.5.4]{Zolotarev1997} where $(\rho\kappa_r,\rho\kappa_s)$
should be $(\frac{\kappa_r}{\rho},\frac{\kappa_s}{\rho})$
and \citet{Zolotarev1978} should be \citet{Zolotarev1979}, 
and also a result of \citet{MitalauskasStatulevicius1976}
as presented in \citet[p.~30, Lemma 2.10]{ChristophWolf1992}.

We get equality in~\eqref{Eq:varkappa_r_vs_varkappa_s_and_Kolmogorov} for $M=P-Q$
with $P,Q\in\Prob(\R)$ and arbitrarily given \linebreak[4]   
$\left( \left\|M\right\|_{\mathrm{K}}, \varkappa_s(M)\right)=(\varrho,\lambda)
 \in\mathopen]0,1\mathclose]\times\mathopen]0,\infty\mathclose[$ 
by taking $P\coloneqq \varrho\delta_{-t}+(1-\varrho)R$ and  $Q\coloneqq \varrho\delta_{t}+(1-\varrho)R$
with $t\coloneqq\left(\frac{\lambda}{2\varrho}\right)^{\frac{1}{s}}$ and $R\in\Prob(\R)$ arbitrary,
and analogously so under the additional condition $\mu(P)=\mu(Q)=0$ for $(\varrho,\lambda)
 \in\mathopen]0,\frac12\mathclose]\times\mathopen]0,\infty\mathclose[$
by taking $P\coloneqq\varrho\left(\delta_{-t}+\delta_t\right)+(1-2\varrho)R$ and 
$Q\coloneqq 2 \varrho\delta_0 +(1-2\varrho)R$ with $t$ and $R$ as above  and also $\mu(R)=0$,
and even with $P$ and $Q$ standardised for $0<\varrho<\frac14$ 
and $0< \lambda \le (2\varrho)^{1-\frac{s}{2}}$
by taking $P\coloneqq\varrho\delta_{-t}+2\varrho\delta_{0}+\varrho\delta_t+(1-4\varrho)R$
and $Q\coloneqq2\varrho\left( \delta_{-\frac{t}{\sqrt{2}}} + \delta_{\frac{t}{\sqrt{2}}}\right)+(1-4\varrho)R$
with $t$ and $R$ as above and also $\mu(R)=0$ and 
$\sigma^2(R)=  (1-2\varrho t^2)/( 1-4\varrho)$.
In the proof of Lemma~\ref{Lem:varkappa_1_vs_varkappa_r_and_Lip} below, however, 
inequality~\eqref{Eq:varkappa_r_vs_varkappa_s_and_Kolmogorov} is applied to $M=P-Q$ with 
$Q$ Lipschitz, excluding finite nonzero equality in~\eqref{Eq:varkappa_r_vs_varkappa_s_and_Kolmogorov},
which explains why we thus only get there inequalities~\eqref{Eq:varkappa_1_vs_varkappa_r_and_Lip}
and~\eqref{Eq:varkappa_1_vs_varkappa_r_and_Lip_centred}  of merely unimprovable order
in $\varkappa_3(M)$, for bounded $\varkappa_3(M)$, but with presumably improvable constants.

\begin{Lem}[$\varkappa_1$ versus $\varkappa_r$ distances of laws to signed Lipschitz laws]
                                                            \label{Lem:varkappa_1_vs_varkappa_r_and_Lip}
Let $r\in\mathopen[1,\infty\mathclose[$ and let 
$Q\in \cM$ be Lipschitz with $Q(\R)=1$.
For $P\in\Prob(\R)$ then
\la                                                        \label{Eq:varkappa_1_vs_varkappa_r_and_Lip}
 \varkappa_1(P-Q)&\le& \left(\min\left\{2^\frac{r+2}{r}c^{}_r,8\right\}\|Q\|_{\mathrm{L}}  \right)^\frac{r-1}{r+1}
   \varkappa_r(P-Q)^{\frac{2}{r+1}}  , \\ 
 \varkappa_1(P-Q)&\le& \left(4\|Q\|_{\mathrm{L}}\right)^\frac{r-1}{r+1}  \label{Eq:varkappa_1_vs_varkappa_r_and_Lip_centred}
  \varkappa_r(P-Q)^{\frac{2}{r+1}} \quad
      \text{ if }\,\mu(P-Q)=0                                                               
\al
with $c_r$ from~\eqref{Eq:c_r_for_Kolmogorov_bounded_by_varkappa_and_Lipschitz}.
The upper bound in~\eqref{Eq:varkappa_1_vs_varkappa_r_and_Lip_centred} is strictly smaller than the one
in~\eqref{Eq:varkappa_1_vs_varkappa_r_and_Lip} unless $P=Q$ or $r=1$. 
The minimum in~\eqref{Eq:varkappa_1_vs_varkappa_r_and_Lip} is equal to its first term in case 
of $r\in\{1,2,3\}$, namely
\[
  2^3c_1 = 8,\qquad 2^2c_2 =5.47239\ldots,\qquad 2^\frac{5}{3}c_3 = 4.96883\ldots\,.
\]

Even for $Q=\mathrm{N}$ and $P$ standardised and arbitrarily $\varkappa_r$-close to $\mathrm{N}$, 
the upper bound in~\eqref{Eq:varkappa_1_vs_varkappa_r_and_Lip_centred} can not be improved by any constant factor 
strictly less than
\la                                                              \label{Eq:varkappa_1_vs_varkappa_r_sharp}
 \lim_{\epsilon\downarrow0} \frac{\varkappa_1(P_\epsilon-\mathrm{N})}
 { \left(\frac{4}{\sqrt{2\pi}}\right)^\frac{r-1}{r+1} \varkappa_r\!\left(P_\epsilon-\mathrm{N}\right)^\frac{2}{r+1}} 
  & =& \frac{5-2\sqrt{3}}
            { 12 \Big( \frac{1}{2} \big( \frac{1}{r+1} 
                   + (\frac{2r}{r+1} -\sqrt{3}\,)3^{-\frac{r+1}{2}} \big) \Big)^\frac{2}{r+1} }
\al
with $P_\epsilon$ as in Zolotarev's Example~\ref{Example:Zolotarev_1973}.
For $\,r=3$, 
$\text{\rm R.H.S.\eqref{Eq:varkappa_1_vs_varkappa_r_sharp}} 
 = \sqrt{\frac{105}{118} - \frac{76}{177}\sqrt{3}}
 = 0.382263\ldots\,$. 
\end{Lem}
\begin{proof}
 Let $P\in\Prob(\R)$,  $M\coloneqq P-Q$, and  $L\coloneqq \|Q\|_{\mathrm{L}}$.
 We will apply Lemma~\ref{Lem:varkappa_r_vs_varkappa_s_and_Kolmogorov}
 with the pair
 $(r,s)$ there being the present $(1,r)$.
 We may assume $\varkappa_r(M)<\infty$, and 
 by~\eqref{Eq:varkappa_r_vs_varkappa_s_and_Kolmogorov} then 
 have $\varkappa_1(M)<\infty$.
 
 Replacing $\left\|M\right\|_{\mathrm{K}}$ in~\eqref{Eq:varkappa_r_vs_varkappa_s_and_Kolmogorov}
 by its final upper bound from~\eqref{Eq:Kolmogorov_bounded_by_varkappa_and_Lipschitz},
 taking there $r$ once equal to the present $r$, and once equal to $1$, yields
 \[
  \varkappa_1(M) &\le& 2^\frac{r-1}{r}\varkappa_r(M)_{}^\frac{1}{r}
   \left(\min\left\{2^\frac{1}{r+1}c_rL^\frac{r}{r+1}\varkappa_r(M)^\frac{1}{r+1},  
   \sqrt{2L\varkappa_1(M)}\right\}\right)^\frac{r-1}{r}
 \]
 and solving this inequality for $\varkappa_1(M)$ yields~\eqref{Eq:varkappa_1_vs_varkappa_r_and_Lip}.
 
 If $\mu(M)=0$, we get as above, taking now the first upper bound 
 from~\eqref{Eq:Kolmogorov_bounded_by_varkappa_and_Lipschitz} just with there $r=1$, 
 and hence $\lambda_r(M)=0$, to get 
 \[
  \varkappa_1(M) &\le& 2^\frac{r-1}{r}\varkappa_r(M)_{}^\frac{1}{r}\sqrt{L\varkappa_1(M)}
 \] 
 and hence~\eqref{Eq:varkappa_1_vs_varkappa_r_and_Lip_centred}.
 
 \vspace{-.4\baselineskip}
 The claim comparing~\eqref{Eq:varkappa_1_vs_varkappa_r_and_Lip_centred} 
 with~\eqref{Eq:varkappa_1_vs_varkappa_r_and_Lip}     
 follows from 
 $2^\frac{r+2}{r}c_r\ge2^\frac{r+2}{r}2^\frac{r-1}{r+1}=4^{\frac{r^2+r +1}{r(r+1)}}>4$.

 Claim~\eqref{Eq:varkappa_1_vs_varkappa_r_sharp} follows from the asymptotics
 (\ref{Eq:zeta_1_Zolotarev-example},\ref{Eq:varkappa_r_Zolotarev-example})
 for $\varkappa_r(P_\epsilon-\mathrm{N})$ given in Example~\ref{Example:Zolotarev_1973}.
\end{proof}

The following surely imperfect remark is used in justifying parts 
of~(\ref{Eq:Binomial_CLT_varkappa_r_but_not_nu_r},\ref{Eq:symmetric_Bernoulli_CLT_zeta_3_varkappa_3}).
We use here the standard notation $\preccurlyeq$ explained in
section~\ref{sec:Asymptotic_comparison_terminology},
and the lattice span notation~\eqref{Eq:Def_lattice_span_h(P)}.

\begin{Rem}[CLT convergence rate with respect to $\varkappa_r$]       \label{Rem:CLT_rate_varkappa_r}
Let $P\in\cP_3$ and $r\in[1,\infty[$ be fixed.
                          
\begin{parts}
\item If $P\in\cP_{r+\epsilon}$ for some $\epsilon>0$, then   \label{part:Osipov_applied_to_varkappa_r}
 $\varkappa_r(\widetilde{P^{\ast n}}-\mathrm{N}) \preccurlyeq \frac{1}{\sqrt{n}}\,$.
\item If $h(P)>0$ or $\mu_3(\bdot{P})\neq 0$, then  \label{part:Esseen_1958_applied_to_varkappa_r}
 $\varkappa_r(\widetilde{P^{\ast n}}-\mathrm{N}) \succcurlyeq \frac{1}{\sqrt{n}}\,$.
\end{parts}        
\end{Rem}
\begin{proof}
\ref{part:Osipov_applied_to_varkappa_r}
 Osipov's theorem as in \citet[p.~167, Theorem~5.15]{Petrov1995} yields 
 $|F_{\widetilde{P^{\ast n}}}(x)-\Phi(x)| 
  \le \frac{ C_{r.P}}{\sqrt{n} (1+|x|)^{r+\epsilon}} $
 with some constant $C_{r,P}<\infty$, and hence a simple integration yields the claim.
 
\smallskip\ref{part:Esseen_1958_applied_to_varkappa_r} For $r=1$, 
\citet[pp.~21-22, Theorem 4.2]{Esseen1958} yields more precisely the existence, with an explicit
formula, of
$ c_P \coloneqq \lim_{n\rightarrow\infty} \sqrt{n}\varkappa_1(\widetilde{P^{\ast n}}-\mathrm{N})>0$, 
in analogy to~\eqref{Eq:Esseen1956_asymptotics}.
For $r>1$ we use Lemma~\ref{Lem:varkappa_r_vs_varkappa_s_and_Kolmogorov},
with the present $(1,r)$ in the role of $(r,s)$ there,
to get $\varkappa_r(M)\ge 2^{1-r}\varkappa_1^r(M)\|M\|_{\mathrm{K}}^{1-r}$ for $M\in\cM$, 
and hence 
$ \varliminf_{n\rightarrow\infty} \sqrt{n} \varkappa_r(\widetilde{P^{\ast n}}-\mathrm{N})
 \ge 2^{1-r} c_P^r\, (\text{R.H.S.\eqref{Eq:Esseen1956_asymptotics}})^{1-r}>0$. 
\end{proof}

It seems likely to us that in
Remark~\ref{Rem:CLT_rate_varkappa_r}\ref{part:Osipov_applied_to_varkappa_r}
the assumption on $P$ can be weakened to $P\in\cP_r$ even if $r\ge3$, and that in
any case $\lim_{n\rightarrow\infty} \sqrt{n} \varkappa_r(\widetilde{P^{\ast n}}-\mathrm{N})$
exists, with a more or less explicit formula analogous to Esseen's special case of $r=1$.

\section{Proof of Theorem~\ref{Thm:BoChiGoe} about lower bounds}
                                                 \label{sec:Proof_of_Theorem_about_lower_bounds} 
\begin{proof}[Proof of Theorem \ref{Thm:BoChiGoe}]
1.~Sufficiency:
\citet[Theorem 1.2]{BCG2012} is actually equivalent to the existence of a constant
$c\in\mathopen]0,\infty\mathclose[$ such that \eqref{Eq:BCG_our_h_=_inverse of theirs} implies,
more generally than~\eqref{Eq:BCG_lower_bound}, that 
\la                                                \label{Eq:BCG_lower_bound_nonstandardised}
  h\!\left( \left\|P-\mathrm{N}\right\|_{\mathrm{K}} \right)
   &\le& 
    \left\|P^{\ast 2} -\mathrm{N}^{\ast2}\right\|_{\mathrm{K}}   
   \quad\text{ for $P\in\Prob(\R)$}
\al
holds. To spell out a proof for the direction of this equivalence actually needed here,
let $h$ be defined by~\eqref{Eq:BCG_our_h_=_inverse of theirs} with $c$ indicated below,
$P\in\Prob(\R)$, $\epsilon\coloneqq\left\|P^{\ast 2} -\mathrm{N}^{\ast2}\right\|_{\mathrm{K}}$,
and $t\coloneqq\left\|P-\mathrm{N}\right\|_{\mathrm{K}}$. 
Then $h(t)\le\epsilon$ holds trivially if $t=0$, if $\epsilon=0$ 
by~\eqref{Eq:Elementary_Cramer_Levy} for $n=2$, 
and if $\epsilon\ge\frac{1}{\mathrm{e}}$ if we choose $c\le \frac{1}{\mathrm{e}}$.  
So let us now assume $\epsilon,t>0$ and  $\epsilon<\frac{1}{\mathrm{e}}$. 
Then the cited theorem states, for some absolute constant here denoted
by $B$, that we have $t\le B\left(\epsilon\log(\frac{1}{\epsilon})\right)^\frac{2}{5}$.
 
If real numbers $x,y$ satisfy 
\la               \label{Eq:x/log(x)_le_y} 
   x\ >\ 1\,,&& \frac{x}{\log x} \ \le\ y\,,
\al
then $y\ge e$, $x\le y\log x$, $\log x \le\log(y\log x)= \log y  +  \log\log x$, and hence 
\la                \label{Eq:x_le_c_y_log_y}
  y\ \ge\ \mathrm{e}, &&x  \le\  \left(\log y +\log\log x\right) \,\ \le\,\ \frac{\mathrm{e}}{\,\mathrm{e}\!-\!1\,}y\log y
\al
by using in the last step $\log\log x\le 0$ in case of $x\le\mathrm{e}$, and else 
\[
 \frac{\,\log\log x\,}{\log y} &\le & \frac{\,\log\log x\,}{\log\left(\frac{x}{\log x}\right)}
 \,\ = \,\ \frac{z}{1-z}  \,\ \le \,\ \frac{1}{\mathrm{e}-1}
\]
with $z \coloneqq  \frac{\log\log x}{\log x}$, being maximal for $\log x=\mathrm{e}$. 

Applying in case of $\epsilon\in\mathopen]0,\frac{1}{\mathrm{e}}[$ and $t>0$ the 
implication $\eqref{Eq:x/log(x)_le_y} \Rightarrow\eqref{Eq:x_le_c_y_log_y}$  to $x\coloneqq\frac{1}{\epsilon}$  
and $y\coloneqq (\frac{B}{t})^{\frac{5}{2}}$  yields
$\frac{1}{\epsilon}\le\frac{\mathrm{e}}{\mathrm{e}-1}B^\frac{5}{2}t^{-\frac{5}{2}}\frac{5}{2}
 (\log(B)+\log\frac{1}{t}) 
 \le\frac{\mathrm{e}}{\mathrm{e}-1}B^\frac{5}{2}t^{-\frac{5}{2}}\frac{5}{2}
  (\log(B)+1)(1\vee \log\frac{1}{t})$, 
hence again $h(t)\le\epsilon$ if $c$ is small enough.  
 
\smallskip 2.~Necessity:
Let $h:[0,1]\toitself$ be such that~\eqref{Eq:BCG_lower_bound} holds. 
For $\epsilon>0$, let $P_\epsilon$ be the (standardised) law as in \citeposs{Zolotarev1973} 
Example~\ref{Example:Zolotarev_1973}, and $t_\epsilon\coloneqq  \left\|P_\epsilon-\mathrm{N}\right\|_{\mathrm{K}}$.
Then, for $\epsilon\downarrow0$ and 
using~\eqref{Eq:||_||_K_Zolotarev-example},
we get $h(t_\epsilon)\le \left\|\widetilde{P_\epsilon^{\ast2}}-\mathrm{N}\right\|_{\mathrm{K}}
\sim \frac{\epsilon^2}{2\pi} \sim \frac{1}{2\pi}(\sqrt{3}\sqrt{2\pi}t_\epsilon)^2=3t_\epsilon^2$.
Hence, using the continuity of $\epsilon\mapsto t_\epsilon$, 
we get, say, $h(t)\le4t^2$ for $t\le t_0$ with some $t_0\in\mathopen]0,1]$,
hence $h(t)\le (4\vee t_0^{-2})t^2$ for every $t\in[0,1]$.

If $\left\|\cdot\right\|_{\mathrm{K}}$
on the right in \eqref{Eq:BCG_lower_bound} is replaced by $\nu_0$, then with~(\ref{Eq:n=2_nu_0_error_Zolotarev_example},\ref{Eq:n=2_Kolmogorov_error_Zolotarev_example}) we get
$h(t_\epsilon)\le \nu_0( \widetilde{P_\epsilon^{\ast2}}-\mathrm{N})
\sim \frac{16-2\sqrt{3}}{3\pi}(\sqrt{3}\sqrt{2\pi}t_\epsilon)^2
=(32-4\sqrt{3})t_\epsilon^2$, and we finish as before.
\end{proof}

Instead of the paragraph above establishing $\eqref{Eq:x/log(x)_le_y} \Rightarrow\eqref{Eq:x_le_c_y_log_y}$, 
we could alternatively have used \citet[p.~88, exemple (8.5.1)]{DieudonneCI}.

\section{Asymptotic comparison terminology and notation} \label{sec:Asymptotic_comparison_terminology}
The purpose of this section is to recall briefly some standard terminology and notation 
for ``local'' or ``asymptotic'' comparisons  of 
functions as presented in \citet[Chapter V, \S1, sections 1 and 2]{Bourbaki2004},
to define our use of phrases like 
``inequality~\eqref{Eq:Berry-Esseen_for_Z-close_to_normal_also_for_n=1}
is i.c.f.~strictly better than the Berry-Esseen inequality~\eqref{Eq:Berry-Esseen_inequality}'',
and to provide some simple facts used in Example~\ref{Example:Zolotarev_1973}.

Let $\mathfrak{F}$ be a filter base. Then for functions $f,g$ defined
{\em along\,} $\mathfrak{F}$, that is, defined on some $F\in\mathfrak{F}$, and with values in a normed vector space $(V,\|\cdot\|)$,
one writes $f \preccurlyeq g$ $:\iff$ there exist an $F\in\mathfrak{F}$ and a $c\in[0,\infty[$
with $\|f\|\le c\|g\|$ on $F$,
$f \asymp g$ $:\iff$ $f\preccurlyeq g$ and $ g\preccurlyeq f$,
$f \llcurly g$ $:\iff$ for every $\epsilon\in\mathopen]0,\infty\mathclose[$ there is an 
$F\in \mathfrak{F}$ with $\|f\|\le \epsilon\|g\|$ on $F$, 
and $f \sim g  :\iff f-g\llcurly g $. 
Analogous definitions of $\preccurlyeq, \asymp, \llcurly$ for $[0,\infty]$-valued functions.
We read $\asymp$ as ``is of the same order as'',
and~$\sim$ as ``is asymptotically equal to'', along $\mathfrak{F}$,
and we recall that not only $\asymp$ but also $\sim$ is an equivalence relation.

Without explicit reference to a filter base, writing ``$f\preccurlyeq g$ on $X$'', for 
functions $f,g$ defined on the nonempty set $X$, 
means: $f\preccurlyeq g$ along the filter base $\{ X\}$.

For $[0,\infty]$-valued functions $f,g_1,g_2$ defined along $\mathfrak{F}$, 
an inequality $f\le g_1$
is called {\em i.c.f.}~better (or sharper, or stronger) 
than $f\le g_2$ if $g_1\preccurlyeq g_2$.
If also 
$g_2\not\preccurlyeq g_1$,
then $f\le g_1$ is i.c.f.~{\em strictly} better, 
and else the two inequalities are i.c.f.~{\em equivalent}.
Examples: Inequality~\eqref{Eq:Berry-Esseen_for_Z-close_to_normal_also_for_n=1}
is i.c.f.~strictly better than~\eqref{Eq:Berry-Esseen_inequality}, 
the filter base 
being $\{ \cP_3\times\N\}$, 
by~\eqref{Eq:zeta_1_vee_zeta_3_distance_to_normal_bounded_by_nu_3}
and by, say,   Example~\ref{Example:Discretised_normal_laws} proving strictness.
Inequality~\eqref{Eq:Berry-Esseen_for_varkappa-close_to_normal} 
is i.c.f.~better than~\eqref{Eq:Berry-Esseen_for_Z-close_to_normal},
referring 
to the filter base 
$\mathfrak{F}_1\coloneqq \{\cP_3\times \N_{\ge2}\}$,
and i.c.f.~strictly better even for $\widetilde{P}$ arbitrarily 
$\zeta_1\!\vee\!\zeta_3$-close to $\mathrm{N}$ by 
Example~\ref{Example:Tail-discretised_normal_laws}, 
referring to being i.c.f.~better w.r.t.~$\mathfrak{F}_1$
and i.c.f.~strictly better 
w.r.t.~$\mathfrak{F}_2\coloneqq \big\{\{P\in\cP_3:
 \big(\zeta_1\!\vee\!\zeta_3\big)(\widetilde{P} -\mathrm{N}) 
<\epsilon\}\times \N_{\ge2} : \epsilon >0\big\}$.
Theorem~\ref{Thm:BE_K_Z}
is i.c.f.~equivalent to Corollary~\ref{Cor:B-E-Z_with_beta}.

In Example~\ref{Example:Zolotarev_1973} we use the following quite trivial but useful complements to
\citet[Chapter V, \S1, section 2, in particular Propositions 8 and 6]{Bourbaki2004}.

\begin{Lem}                                  \label{Lem:sim_relations}
Let $\mathfrak{F}$ be a filter base.

\begin{parts}
\item  Let $V$ be a normed vector space, and let $f_1,f_2,g_1,g_2$ \label{part:sim_relations_added}
 be $V$-valued functions defined along~$\mathfrak{F}$. Then we have the implication
 \[            
  f_1 \sim g_1\,,\ f_2\sim g_2\,,\ \|g_1 \| +\|g_2\| \preccurlyeq \|g_1+g_2\|
  &\Rightarrow& f_1+f_2\sim g_1+g_2\,.
 \]
\item Let $V_1,V_2,V$ be normed vector spaces,       \label{part:sim_relations_multiplied}
 $V_1\times V_2 \ni (x,y)\mapsto xy \in V$ a continuous bilinear map,
 and for $i\in\{1,2\}$ let $f_i,g_i$ be $V_i$-valued functions defined along~$\mathfrak{F}$.
 Then 
 \[                                   
  f_1 \sim g_1\,,\ f_2\sim g_2\,,\  \|g_1 \|\!\cdot\!\|g_2\| \preccurlyeq \|g_1g_2\|
  &\Rightarrow& f_1f_2\sim g_1g_2\,.   
 \]
\item Let $V$ be a vector space with two norms $\|\cdot\|_1$ and $\|\cdot\|_2$\,,
                                                          \label{part:sim_relations_with_two_norms}
 and let $f,g$ be $V$-valued functions defined along $\mathfrak{F}$. Then
 \[
  f\sim g\text{ w.r.t.\,}\left\|\cdot\right\|_1\,,\ \left\|\cdot\right\|_2\preccurlyeq\left\|\cdot\right\|_1\text{ on V}\,,\
  \left\|g\right\|_1\preccurlyeq\left\|g\right\|_2
  &\Rightarrow&  f\sim g\text{ w.r.t. }\left\|\cdot\right\|_2\,.
 \]
\item Let $V$ be a normed vector space, and let $f,g$ be $V$-valued  \label{part:sim_relation_implies_same_for_norms} 
 functions defined along $\mathfrak{F}$. Then
\[
      f\sim g &\Rightarrow& \|f\| \sim\|g\|\,.
\]
\end{parts}
\end{Lem}
\begin{proof} (a) 
$\|(f_1+f_2)-(g_1+g_2)\| \le \|f_1-g_1\| + \|f_2-g_2\|
\llcurly \|g_1\| + \| g_2\| \preccurlyeq \|g_1+g_2\|$.

(b) 
$\|f_1f_2-g_1g_2\| \le  \|f_1(f_2-g_2)\| +\| (f_1-g_1)g_2\|
 \preccurlyeq \|f_1\| \!\cdot\!\|f_2-g_2\| +  \|f_1-g_1\| \!\cdot\!\|g_2\|
 \llcurly \|f_1\|\!\cdot\!\|g_2\| + \|g_1\|\!\cdot\!\|g_2\|  \asymp \|g_1\|\!\cdot\!\|g_2\| 
 \preccurlyeq \|g_1g_2\|$.

(c) $\left\|f-g\right\|_2\preccurlyeq \left\| f-g\right\|_1 \llcurly \left\|g\right\|_1\preccurlyeq\left\|g\right\|_2$\:\!. 
(d) $\big|\|f\|-\|g\|\big|  \le \|f-g\|\llcurly \|g\|$.   
\end{proof}

In Example~\ref{Example:Zolotarev_1973}, 
the above is applied to ``$\epsilon\downarrow 0$'', that is, to 
$\mathfrak{F}\coloneqq\left\{\,\mathopen]0,\epsilon_0\mathclose]
:\epsilon_0\in\mathopen]0,\infty\mathclose[\,\right\}$, 
to subspaces of the space $\cM$ of bounded signed measures on $\R$, 
with various norms,
and with the bilinear map in~\ref{Lem:sim_relations}\ref{part:sim_relations_multiplied}
being convolution. 
We also use the following, 
in particular for $\|\cdot\| = \nu_r$ on $\cM'=\cM_{r,0}$ as defined 
in~(\ref{Eq:Def_cM_r},\ref{Eq:Def_cM_r,k}), 
which may here serve as an example of Lemma~\ref{Lem:sim_relations}\ref{part:sim_relations_multiplied}.

\begin{Lem}                                    \label{Lem:sim_on_cM_scale_invariant}
Let $\cM'$ be a vector subspace of $\cM$, $r\in\R$, and $\|\cdot\|$ a norm on $\cM'$
with the scaling property~\eqref{Eq:homogeneity_on_cM_conclusion} for $M\in\cM'$.
Let $M_1,M_2,t$ be functions defined along a filter base $\mathfrak{F}$, 
with $M_1,M_2$ being $\cM'$-valued, and $t$ being $\mathopen]0,\infty\mathclose[\,$-valued.
Then, with respect to the norm $\|\cdot\|$  on $\cM'$, we have the equivalence
\la                       \label{Eq:sim_on_cM_scale_invariant}
  M_1\sim M_2 &\iff& M_1(\tfrac{\,\cdot\,}{t}) \sim M_2(\tfrac{\,\cdot\,}{t})\,.
\al
\end{Lem}
\begin{proof} 
It is enough to prove ``$\Rightarrow$'', since we then get ``$\Leftarrow$''
by considering $\frac{1}{t}$.

First proof of ``$\Rightarrow$'':
If we have L.H.S.\eqref{Eq:sim_on_cM_scale_invariant}, then 
$\|  M_1(\tfrac{\,\cdot\,}{t}) - M_2(\tfrac{\,\cdot\,}{t})\|
 = t^r \|M_1-M_2 \|  \llcurly t^r \|M_2\| = \|  M_2(\tfrac{\,\cdot\,}{t})\| $.
 
Second proof of ``$\Rightarrow$'': Apply Lemma~\ref{Lem:sim_relations}\ref{part:sim_relations_multiplied}
to $V_2=V\coloneqq\cM'$, $V_1$ the space of all bounded linear endomorphisms 
of~$V_2$, $f_2\coloneqq M_1$, $g_2\coloneqq M_2$, 
and 
$f_1(\epsilon)\coloneqq g_1(\epsilon)$ being, for $\epsilon\in F_0$ with some $F_0 \in\mathfrak{F}$, 
the map which sends any $M\in\cM'$ to $M(\frac{\,\cdot\,}{t(\epsilon)})$, so that 
$\|g_1g_2\| =\|  M_2(\tfrac{\,\cdot\,}{t})\| = t^r \|M_2\| = \|g_1\|\cdot\|g_2\|$.
\end{proof}

\section{Monotonicity of the variance under contraction, in particular
 under winsorisation} \label{sec:Variance_under_contraction} 
In the proof of Example~\ref{Examples:zeta_3=mu_3}\ref{part:Left-winsorised_normal_laws},
we use the rather obvious Corollary~\ref{Cor:Variance_under_winsorisation}  below, 
which, except for the strictness of the inequality needed by us,
is well-known as the special case of exponent~$2$ of 
\citet[Corollary~3]{ChowStudden1969} = \citet[p.~104, Corollary~2]{ChowTeicher1997}.

\begin{Lem}[Contraction decreases variance]                  \label{Lem:Variance_under_contraction}
Let $T:\R\toitself$ be a contraction, in the sense of
\la                                            \label{Eq:T_contraction}
   |T(y)-T(x)| &\le& |y-x| \quad\text{ for }(x,y)\in\R^2,
\al
and let $P\in\Prob_2(\R)$. Then $\sigma^2(T\im P)\le \sigma^2(P)$, 
with equality iff equality holds in~\eqref{Eq:T_contraction} $P^{\otimes2}$-a.e. 
\end{Lem}
\begin{proof}
$\sigma^2(T\im P)=\iint\frac{1}{2}\big(T(y)-T(x)\big)^2\,\dd P(x)\dd P(y) 
 \le \iint\frac{1}{2}(y-x)^2\,\dd P(x)\dd P(y) =\sigma^2(P)$.
\end{proof}

\begin{Cor}[Winsorisation decreases variance]                   \label{Cor:Variance_under_winsorisation}
Let $P\in\Prob_2(\R)$, $-\infty\le a\le b \le \infty$, 
and $Q \coloneqq P(\,\cdot\,\cap\, ]a,b[\,) + P(\,\mathopen]-\infty,a\mathclose]\,)\delta_a$ 
 $+$ $P(\,\mathopen[b,\infty\mathclose[\,)\delta_b$\,.
Then $\sigma^2(Q) < \sigma^2(P)$ or $Q=P$ or $\sigma^2(P)=0$.
\end{Cor}
\begin{proof}
 Lemma~\ref{Lem:Variance_under_contraction} applied to 
 $T(x)\coloneqq a\!\vee\! x\!\wedge\! b$ for $x\in\R$.
\end{proof}

\section{Roundings and histograms of laws on $\R$} \label{sec:Roundings_and_histograms} 

This is a classical if somewhat marginal topic in probability and statistics, 
going back at least to \citet{Sheppard1898}. Treatments known to us
are usually deliberately incomplete
and not always mathematically precise, with the latter exemplified by
\citet[145, p.~362]{Cramer1945} writing about the Sheppard corrections quite tautologically:
``These relations hold under the assumption that the remainder $R$ in (27.5.2) may be neglected''. 
A good entry into the relevant literature is 
\citet{SchneeweissKomlosAhmad2010}, providing 57 references,
with a mathematically precise and comparatively 
recent one among these being \citet{Janson2006}.

The present section is auxiliary to Examples~\ref{Example:Discretised_normal_laws}
and~\ref{Example:Tail-discretised_normal_laws}.
For a bounded nondegenerate interval $I\subseteq \R$, we let below
$\mathrm{U}_I\coloneqq \frac{1}{\leb(I)}\leb(\cdot\cap I)$ denote the uniform law on~$I$.

\begin{Def}          \label{Def:Rounding_and_histograms}
Let $P\in\Prob(\R)$, $\eta\in\mathopen]0,\infty\mathclose[\,$, and
$\alpha\in\mathopen]0,1\mathclose[\,$. With
\[
 I_j &\coloneqq& I_{\eta,\alpha,j} 
  \,\ \coloneqq\,\ ] (\alpha+j-\tfrac{1}{2})\eta\,, (\alpha+j+\tfrac{1}{2})\eta[\,,\\
 I_j^\natural &\coloneqq&  I_{\eta,\alpha,j}^\natural
 \,\ \coloneqq\,\ \1_{I_j} 
  + \tfrac{1}{2} \1_{\{ (\alpha+j-\tfrac{1}{2})\eta\,, (\alpha+j+\tfrac{1}{2})\eta\}}\,, \\
 p_j &\coloneqq& p_{\eta,\alpha,j} \,\ \coloneqq\,\ \int   I_j^\natural\,\dd P
\]
for $j\in \Z$, we call 
\[
 P_{\mathrm{rd}} &\coloneqq& P_{\mathrm{rd},\eta,\alpha} 
    \,\ \coloneqq\,\ \sum_{j\in\Z}p_j\delta_{(\alpha+j)\eta} 
\] 
the {\em rounding}, and 
\[
 P_{\mathrm{hist}} &\coloneqq&  P_{\mathrm{hist},\eta,\alpha} 
  \,\ \coloneqq \,\ \sum_{j\in\Z}p_j\mathrm{U}_{I_j} 
\]
the {\em histogram law}, of $P$,  with respect to the {\em rounding lattice} 
$\{(\alpha+j)\eta : j\in \Z\}$, with the {\em width}~$\eta$, the {\em shift} $\alpha\eta$,
and the {\em shift parameter} $\alpha$.
\end{Def}

In the above situation, we obviously have 
$ (P_\mathrm{hist})_\mathrm{rd}=P_\mathrm{rd} =(P_\mathrm{rd})_\mathrm{rd}$
and $ (P_\mathrm{rd})_\mathrm{hist}=P_\mathrm{hist} =(P_\mathrm{hist})_\mathrm{hist}$.

We are 
in particular interested, for $\eta$ close to zero,
in the zeta distances  
\mbox{$\zeta_r(\widetilde{(\mathrm{N}_{\mu,\sigma^2})_{\mathrm{rd},\eta,\alpha}} \!-\! \mathrm{N})$} 
for $r\in\{1,3\}$, the stan\-dard\-ised lattice span
$h( \widetilde{(\mathrm{N}_{\mu,\sigma^2})_{\mathrm{rd},\eta,\alpha}})
 =  \eta / \sigma( (\mathrm{N}_{\mu,\sigma^2})_{\mathrm{rd},\eta,\alpha}) $, 
and the standardised third moment
$\mu_3(\widetilde{(\mathrm{N}_{\mu,\sigma^2})_{\mathrm{rd},\eta,\alpha}})$, 
in order to compare R.H.S.\eqref{Eq:Esseen1956_asymptotics} with 
R.H.S.\eqref{Eq:Berry-Esseen_for_Z-close_to_normal}
in case of $P\coloneqq (\mathrm{N}_{\mu,\sigma^2})_{\mathrm{rd},\eta,\alpha}$
as in~\eqref{Eq:_Def_discretised_normal}.

\begin{Lem}[Lattice and histogram approximations in $\Prob(\R)$]  
                     \label{Lem:Lattice_and_histogram_approximations}
                     \label{Lem:Roundings_and_histograms} 
In the situation 
of Definition~\ref{Def:Rounding_and_histograms},
let $k\in\N_0$, and let any asymptotic relation $\preccurlyeq, \asymp, \llcurly,\sim$
refer to $\eta\rightarrow 0$ 
with $P,\alpha,k$ fixed.

\smallskip
\begin{parts}
\item We have, assuming $\nu_j(P)< \infty$ in every relation where $\mu_j$ occurs, 
\la \qquad\quad
 \mu_k(P_\mathrm{rd}-P_\mathrm{hist})          \label{Eq:mu_k_histogram-rounded}
  &=& \tfrac{-1}{k+1}\sum_{\ell=1}^{\lfloor\frac{k}{2}\rfloor}
   \textstyle{\binom{k+1}{2\ell+1}}\left(\tfrac{\eta}{2}\right)^{2\ell}\mu_{k-2\ell}(P_\mathrm{rd})
   \,, \\
 \mu_0(P_\mathrm{rd}-P_\mathrm{hist})    \label{Eq:mu_0_1_histogram-rounded}   
 &=& \mu_1(P_\mathrm{rd}-P_\mathrm{hist}) \,\ = \,\ 0\,,\quad
   \mu_2(P_\mathrm{rd}-P_\mathrm{hist}) \,\ =\,\ -\tfrac{\eta^2}{12}\,, \\
 \mu_3  (P_\mathrm{rd}-P_\mathrm{hist})    \label{Eq:mu_3_histogram-rounded}   
  &=& -\tfrac{\eta^2}{4}\mu_1(P_\mathrm{rd}) 
  \,\ =\,\ -\tfrac{\eta^2}{4}\mu_1(P_\mathrm{hist}) \,, \\ 
 \zeta_1 (P_\mathrm{rd}-P_\mathrm{hist}) \label{Eq:zeta_1_rounded_histogram} 
   &=& \tfrac{\eta}{4}\,,       \\
 \underline{\zeta}_k (P_\mathrm{rd} - P_\mathrm{hist} )    \label{Eq:zeta_k_rounded-histogram}      
    &\le & \tfrac{\eta^2}{8} \sum_{\ell=0}^{k-2} \tfrac{\eta^\ell}{\ell!(k-2-\ell)!}
            \nu_{k-2-\ell} (P)   \quad\text{ if } k\ge2  \,, 
\al \la                
 \underline{\zeta}_k (P_\mathrm{hist} -P)  &\le & \tfrac{\eta}{2}  \label{Eq:zeta_k_histogram_original} 
   \sum_{\ell=0}^{k-1} \tfrac{\eta^\ell}{\ell!(k\!-\!1\!-\!\ell)!}\nu_{k-1-\ell}(P_\mathrm{hist}-P)
                                                                \quad\text{ if } k\ge1\,,   \\                                                            
  |\mu_1(P_\mathrm{rd}-P)|                 \label{Eq:mu_rounded_histogram_original_zeta_1} 
    &=& |\mu_1(P_\mathrm{hist}-P)|\,\ \le\,\     \zeta_1 (P_\mathrm{hist} -P) 
    \,\ \le\,\   \tfrac{\eta}{2}\nu_0(P_\mathrm{hist}-P)  \,,   \\
 \tfrac{1}{2}|\mu_2 ( P_\mathrm{hist} -P)|   \label{Eq:zeta_2_histogram_original}    
  &\le&   \underline{\zeta}_2( P_\mathrm{hist} -P)              
  \,\ \le \,\ \tfrac{\eta}{2}\nu_1(P_\mathrm{hist}-P)+\tfrac{\eta^2}{2}\nu_0(P_\mathrm{hist} -P) \,, \\
 \left| \nu_k(P_\mathrm{hist}) - \nu_k(P) \right|  \label{Eq:nu_hist-nu_original}   
  &\preccurlyeq& \eta \quad\text{ if }\nu_k(P)<\infty    \,. 
\al
\item If $P$ is absolutely continuous with respect to $\leb$, then we have 
\la             
 \nu_k(P_\mathrm{hist}-P)      \label{Eq:nu_k_P_hist-P_a.c.}   
  &\rightarrow& 0   \quad\text{ if }\nu_k(P)<\infty\,, \\
 \underline{\zeta}_k(P_\mathrm{hist}-P)      \label{Eq:zeta_k_P_hist-P_a.c.}   
  &\llcurly & \eta   \quad\text{ if }k\ge1\text{ and }\nu_{k-1}(P)<\infty\,, \\ 
 \zeta_1 (P_\mathrm{rd} -P)  \label{Eq:zeta_1_rounded-original_asymp} 
   & \sim & \tfrac{\eta}{4} \,, \\   
 \underline{\zeta}_k(P_\mathrm{rd}-P)&\llcurly& \eta   \label{Eq:zeta_k_rounded-orignal_small}   
    \quad   \text{ if }k\ge 2 \text{ and }\nu_{k-1}(P)<\infty\,,  \\
 \mu_1(P_\mathrm{rd})\! -\!\mu_1(P)\!\!&\llcurly&\!\!\eta\,,
                            \label{Eq:zeta_1_rounded-original_centred_asymp}
   \quad \zeta_1(\bdot{P_\mathrm{rd}}- \bdot{P})     
    \,\ \sim\,\ \tfrac{\eta}{4} \quad\text{ if } \nu_1(P)<\infty  \,, \\
 \quad\mu_2(P_\mathrm{rd})\!-\!\mu_2(P)\!\!&\llcurly&\!\!\eta \,,
                     \label{Eq:zeta_1_P-P_standardised_rounded_asymp}
  \,\ \sigma(P_\mathrm{rd})\!-\!\sigma(P) \, \llcurly\, \eta\,, 
  \,\  \zeta_1(\widetilde{P_\mathrm{rd}}\!-\! \widetilde{P})
     \, \sim\, \tfrac{\eta}{4\,\sigma(P)} \quad\text{ if } \nu_2(P)<\infty \,, \\
  \underline{\zeta}_k(\widetilde{P_\mathrm{rd}}\!-\! \widetilde{P})  
                    \label{Eq:zeta_k_P-P_standardised_rounded_asymp}
    &\llcurly& \eta    
    \quad\text{ if } k\ge 2 \text{ and } \nu_k(P)<\infty \,.  
\al 
Here $\bdot{P_\mathrm{rd}}$ denotes the centring of $P_\mathrm{rd}\,$,
which may differ from 
the rounding of $\bdot{P}$. 
And $\widetilde{P_\mathrm{rd}\,}$
denotes the standardisation of $P_\mathrm{rd}$\,, which
in~{\rm(\ref{Eq:zeta_1_P-P_standardised_rounded_asymp},\ref{Eq:zeta_k_P-P_standardised_rounded_asymp})}
is indeed defined for $\eta$ sufficiently small.
\end{parts}
\end{Lem}
\begin{proof} (a)
Let $x_j\coloneqq(\alpha+j)\eta$ for $j\in\Z$. 

To prove~\eqref{Eq:mu_k_histogram-rounded},
we calculate, using merely the assumption  
$\nu_{(k-2)_+}(P)<\infty$  in the first two steps and $\nu_k(P)<\infty$ only in the last,
\[
 \text{R.H.S.\eqref{Eq:mu_k_histogram-rounded}}
  &=& \sum _{j\in\Z}p_j\left(   x_j^{\,k} - \tfrac{2}{(k+1)\eta}\sum_{\ell=0}^{\lfloor\frac{k}{2}\rfloor}
       \binom{k+1}{2\ell+1}\left(\tfrac{\eta}{2}\right)^{2\ell+1} x_j^{\,k-2\ell} \right)    \\
  &=& \sum _{j\in\Z}p_j\Big( x_j^{\,k} - \tfrac{1}{(k+1)\eta}
       \big((x_j+\tfrac{\eta}{2})^{k+1} - (x_j-\tfrac{\eta}{2})^{k+1}\big) \Big) \\
  &=& \sum_{j\in\Z}p_j\left(  x_j^{\,k} -\tfrac{1}{\eta}\int_{I_j}x^k\,\dd x  \right)    
    \,\ =\,\ \text{L.H.S.\eqref{Eq:mu_k_histogram-rounded}}.
\]
From this (\ref{Eq:mu_0_1_histogram-rounded},\ref{Eq:mu_3_histogram-rounded}) follow easily.

To prove~\eqref{Eq:zeta_1_rounded_histogram}, we calculate
\[
 \zeta_1(P_\mathrm{hist}-P_\mathrm{rd}) 
  &=& \int\big|F_{P_\mathrm{hist}-P_\mathrm{rd}}(x)\big|\,\dd x 
  \,\ =\,\ \sum_{j\in\Z} \int_{I_j} p_j\left| \tfrac{x-(\alpha+j-\frac{1}{2})\eta}{\eta}
            -(x\ge x_j) \right|\,\dd x 
  \,\ =\,\ \tfrac{\eta}{4}   \,.
\]

To prove~\eqref{Eq:zeta_k_rounded-histogram}, recalling  
the definitions (\ref{Eq:Def_zeta_underbar_r},\ref{Eq:Def_cF_r,r-1^infty}),
let $k\ge 2$ and $g\in\cF_{k,k-1}^\infty$.
Then, using $\sum_{j\in\Z}I_j^\natural = 1$
and $\int I_j^\natural  \dd P_\mathrm{hist}= p_j = \int I_j^\natural  \dd P_\mathrm{rd}$
as well as $\int x I_j^\natural(x)  \dd P_\mathrm{hist}(x)= p_jx_j = \int xI_j^\natural(x)  \dd P_\mathrm{rd}(x)$
in the first step,
and $|g''(x)|=|g''(x)-\sum_{\ell=0}^{k-3}  g^{(2+\ell)}(0)\frac{x^\ell}{\ell!}|$ 
$\le$ $\frac{|x|^{k-2}}{(k-2)!}$ $\leb$-a.e.~in the fourth, we get
\[
 \left|\int g\,\dd(P_\mathrm{hist}-P) \right|
  &\le& \sum_{j\in \Z} \left|\int \big(g(x)-g(x_j)-g'(x_j)(x-x_j)\big) I^\natural_j
   \,\dd\big(P_\mathrm{rd}-P_\mathrm{hist}\big)(x) \right| \\
  &=& \sum_{j\in \Z} \left|\int \big(g(x)-g(x_j)-g'(x_j)(x-x_j)\big) I^\natural_j
   \,\dd P_\mathrm{hist}(x) \right| \\
  &\le& \sum_{j\in\Z} \tfrac{\eta^2}{8}  \esssup_{x\in I_j}|g''(x)|\, p_j  \\
  &\le& \tfrac{\eta^2}{8}\sum_{j\in\Z}  \int \tfrac{(\,|x|+\eta)^{k-2}}{(k-2)!}I_j^\natural(x) \,\dd P(x)
  \,\ =\,\ \text{R.H.S.\eqref{Eq:zeta_k_rounded-histogram}}.
\]  

To prove~\eqref{Eq:zeta_k_histogram_original}, let again $g\in\cF_{k,k-1}^\infty$.
Then, using $\sum_{j\in\Z}I_j^\natural = 1$
and $\int I_j^\natural  \dd P_\mathrm{hist}= \int I_j^\natural  \dd P$ in the first step,
and $|g'(x)|=|g'(x)-\sum_{\ell=0}^{k-2}  g^{(1+\ell)}(0)\frac{x^\ell}{\ell!}|$ 
$\le$ $\frac{|x|^{k-1}}{(k-1)!}$ $\leb$-a.e.~in the third, we get
\[
 \left|\int g\,\dd(P_\mathrm{hist}-P) \right|
  &\le& \sum_{j\in \Z} \left|\int \big(g-g(x_j)\big) I^\natural_j\,\dd(P_\mathrm{hist}-P) \right| \\
  & \le & \sum_{j\in \Z}\tfrac{\eta}{2} \esssup_{x\in I_j}|g'(x)| \int I^\natural_j\,\dd\left|P_\mathrm{hist}-P \right|\\
  & \le & \tfrac{\eta}{2}\sum_{j\in \Z}  \int 
     \tfrac{(\,|x|+\eta\,)^{k-1}}{(k-1)!} I^\natural_j(x)\,\dd\big|P_\mathrm{hist}-P \big|(x)
  \,\ =\,\ \text{R.H.S.\eqref{Eq:zeta_k_histogram_original}}\,.   
\]  

\eqref{Eq:mu_rounded_histogram_original_zeta_1} follows 
from~(\ref{Eq:mu_0_1_histogram-rounded},\ref{Eq:zeta_vs_varkappa_etc},\ref{Eq:zeta_k_histogram_original}).

The first claim in~\eqref{Eq:zeta_2_histogram_original} is contained 
in~\eqref{Eq:zeta_underbar_vs_varkappa_etc}, 
the second in~\eqref{Eq:zeta_k_histogram_original}.

For $k=0$, \eqref{Eq:nu_hist-nu_original} is trivial due to $\text{L.H.S.}=0$.
For $k\ge 1$, we use 
\[
 \tfrac{1}{k!} \left| \big(P_\mathrm{hist}-P\big)|\cdot|^k  \right| 
  &\le& \text{L.H.S.\eqref{Eq:zeta_k_histogram_original}}
  \,\ \le \,\  \text{R.H.S.\eqref{Eq:zeta_k_histogram_original}} \\
  &\le& \tfrac{\eta}{2}\sum_{\ell=0}^{k-1} \tfrac{\eta^\ell }{\ell!(k\!-\!1\!-\!\ell)!}
  \big( \nu_{k-1-\ell}(P_\mathrm{hist}) +   \nu_{k-1-\ell}(P)  \big) \,,
\]
with finiteness of the sum above following inductively.

\smallskip(b) Let $f$ be a $\leb$-density of  $P$. 
Then $f_\eta\coloneqq \sum_{j\in\Z} p_j\1_{I_j}$ is a $\leb$-density of $P_\mathrm{hist}$, 
with $f_\eta\rightarrow f$ $\leb$-a.e.\ by the fundamental theorem of calculus. 
If $\nu_k(P)<\infty$, then  we have
$\int |\cdot|^kf_\eta\,\dd\leb \rightarrow \int |\cdot|^kf\,\dd\leb$ by~\eqref{Eq:nu_hist-nu_original},  
and hence 
$\text{L.H.S.\eqref{Eq:nu_k_P_hist-P_a.c.}}  
 =  \int\big| |\cdot|^kf_\eta - |\cdot|^kf\big|\,\dd\leb \rightarrow 0$
by ``Scheff\'e's theorem'' 
as in \citet[p.~135, Theorem~2.8.9]{Bogachev_MT_I}.

The claim~\eqref{Eq:zeta_k_P_hist-P_a.c.} follows 
from~(\ref{Eq:zeta_k_histogram_original},\ref{Eq:nu_k_P_hist-P_a.c.}).

To prove~\eqref{Eq:zeta_1_rounded-original_asymp}, 
we use~\eqref{Eq:zeta_1_rounded_histogram} in the first step 
and~\eqref{Eq:zeta_k_P_hist-P_a.c.} in the third to get 
$|\zeta_1(P_\mathrm{rd}-P)-\frac{\eta}{4}|
 = |\zeta_1(P_\mathrm{rd}-P)-\zeta_1(P_\mathrm{rd}-P_\mathrm{hist})|
 \le \zeta_1(P_\mathrm{hist}-P) \llcurly \eta$.

\eqref{Eq:zeta_k_rounded-orignal_small} follows 
from~(\ref{Eq:zeta_k_rounded-histogram},\ref{Eq:zeta_k_P_hist-P_a.c.}).
 
The first claim in~\eqref{Eq:zeta_1_rounded-original_centred_asymp} 
follows from \eqref{Eq:mu_rounded_histogram_original_zeta_1}
and either of~(\ref{Eq:nu_k_P_hist-P_a.c.},\ref{Eq:zeta_k_P_hist-P_a.c.}).
Further,
\[
 \bdot{P_\mathrm{rd}}-\bdot{P} &=& \delta_{-\mu(P)}\ast(P_\mathrm{rd}-P)  
                         + (\delta_{-\mu(P_\mathrm{rd})} - \delta_{-\mu(P)})\ast P_\mathrm{rd}
\]
with $\zeta_1( \delta_{-\mu(P)}\ast(P_\mathrm{rd}-P) ) =\zeta_1(P_\mathrm{rd}-P)\sim\frac{\eta}{4}$ 
by \eqref{Eq:translation_invariance_on_cM} and~\eqref{Eq:zeta_1_rounded-original_asymp}, 
and, using~\eqref{Eq:zeta_1_under_x_maspto_a+x} 
and the first claim in~\eqref{Eq:zeta_1_rounded-original_centred_asymp},
\[ 
 \zeta_1\big( (\delta_{-\mu(P_\eta)} - \delta_{-\mu(P)})\ast P_\mathrm{rd}\big) 
   &=&  | \mu(P_\mathrm{rd})-\mu(P)|
   \,\ \llcurly \,\ \eta \,.
\]  
Therefore, using the norm property of $\zeta_1$ on $\cM_{1,0}$\,, we get
the second claim in~\eqref{Eq:zeta_1_rounded-original_centred_asymp}.

The first claim in~\eqref{Eq:zeta_1_P-P_standardised_rounded_asymp}
follows from~(\ref{Eq:mu_0_1_histogram-rounded},\ref{Eq:zeta_2_histogram_original})
and either of~(\ref{Eq:nu_k_P_hist-P_a.c.},\ref{Eq:zeta_k_P_hist-P_a.c.}), 
and the second then follows using~\eqref{Eq:zeta_1_rounded-original_centred_asymp} 
and, for the differentiability of 
$(x,y)\mapsto \sqrt{y-x^2}$ at $(\mu_1(P),\mu_2(P))$,  also $\sigma(P)>0$. 
For $\eta$ small enough to ensure $\sigma(P_\mathrm{rd})>0$, we put 
$Q\coloneqq (x\mapsto \frac{\sigma(P)}{\sigma(P_\mathrm{rd})}x)\im \bdot{P_\mathrm{rd}}$ 
and get 
\la                                            \label{Eq:191}
 \sigma(P)\zeta_1(\widetilde{P_\mathrm{rd}}- \widetilde{P})
   &=&\zeta_1( Q- \bdot{P})
\al
by the homogeneity~\eqref{Eq:homogeneity_on_cM_conclusion} or~\eqref{Eq:homogeneity_of_zeta} 
of $\zeta_1$, and
\la                                            \label{Eq:192}
 |\zeta_1(Q-\bdot{P})- \zeta_1(\bdot{P_\mathrm{rd}}-\bdot{P})|
  &\le &  \zeta_1(Q-\bdot{P_\mathrm{rd}})  
  \,\ =\,\ \left|\tfrac{\sigma(P)}{\sigma(P_\mathrm{rd})}-1\right| \nu_1(\bdot{P_\mathrm{rd}})
  \,\ \llcurly \,\ \eta 
\al
by using in the second step~\eqref{Eq:zeta_1_under_x_maspto_ax}, 
and in the third step 
the second claim in~\eqref{Eq:zeta_1_P-P_standardised_rounded_asymp}
and the boundedness of 
$\nu_1(\bdot{P_\mathrm{rd}})\le \nu_1(\bdot{P}) + \zeta_1(\bdot{P_\mathrm{rd}}-\bdot{P})
 \preccurlyeq 1$ due to~\eqref{Eq:zeta_1_rounded-original_centred_asymp}. 
Using (\ref{Eq:191},\ref{Eq:192},\ref{Eq:zeta_1_rounded-original_centred_asymp}) we get
\[
 \sigma(P)\zeta_1(\widetilde{P_\mathrm{rd}}- \widetilde{P})
   &\sim& \zeta_1(\bdot{P_\mathrm{rd}}-\bdot{P}) \,\ \sim\,\ \tfrac{\eta}{4} 
\]  
and therefore the third claim in~\eqref{Eq:zeta_1_P-P_standardised_rounded_asymp}.

To prove finally~\eqref{Eq:zeta_k_P-P_standardised_rounded_asymp}, 
we write 
\[
 \widetilde{P_\mathrm{rd}} - \widetilde{P}
  &=& \big(x\mapsto\tfrac{x-\mu(P_\mathrm{rd})}{\sigma(P_\mathrm{rd})}\big)\im ( P_\mathrm{rd} - P )
   + \Big(\big(x\mapsto\tfrac{x-\mu(P_\mathrm{rd})}{\sigma(P_\mathrm{rd})}\big)\im P 
    - \big(x\mapsto\tfrac{x-\mu(P)}{\sigma(P)}\big)\im P \Big) \\
  &\eqqcolon& M_1+M_2  
\]
with 
\[
 \underline{\zeta}_k(M_1) &=& \underline{\zeta}_k\Big(
    \big(x\mapsto\tfrac{x}{\sigma(P_\mathrm{rd})}\big)\im \big((P_\mathrm{rd}-P)  \ast\delta_{-\mu(P_\mathrm{rd}} \big)
   \Big)  \\
 &\le& \frac{1}{\sigma(P_\mathrm{rd})^k}
   \Big( \underline{\zeta}_k(P_\mathrm{rd}-P)  
      + \sum_{j=0}^{k-1}\tfrac{(r-j)!}{j!} \big|\mu_j(P_\mathrm{rd}-P)\big|\,
       | \mu(P_\textrm{rd})|^{r-j} \Big)
 \,\ \llcurly \,\ \eta      
\]
by Lemma~\ref{Lem:zeta_convolutions} in the second step,
and by (\ref{Eq:zeta_k_rounded-orignal_small},\ref{Eq:zeta_1_rounded-original_centred_asymp}) and  
using $\frac{1}{j!} \big|\mu_j(P_\mathrm{rd}-P)\big|\le\underline{\zeta}_j(P_\mathrm{rd}-P)$ for $j\ge2$
in the third,
and with 
\[
 \underline{\zeta}_k(M_2) &\le& 
  \text{R.H.S.\eqref{Eq:zeta_r_under_x_mapsto_ax+b} with
   $a\coloneqq\tfrac{1}{\sigma(P)}\,, 
   b\coloneqq-\tfrac{\mu(P)}{\sigma(P)}\,, 
   c\coloneqq\tfrac{1}{\sigma(P_\mathrm{rd})}\,,
   d\coloneqq-\tfrac{\mu(P)_\mathrm{rd}}{\sigma(P_\mathrm{rd})}\,,  r\coloneqq k, M\coloneqq P$} \\
  &\llcurly& \eta  
\] 
by
(\ref{Eq:zeta_1_rounded-original_centred_asymp},\ref{Eq:zeta_1_P-P_standardised_rounded_asymp}).  
\end{proof}

\begin{proof}[Proof of Example~\ref{Example:Tail-discretised_normal_laws}]
                               \label{page:Proof_of_Example_Tail-discretised_normal_laws}  
1. With the conditional laws 
$\mathrm{N}\bed{\cdot}{I}$ and
$Q\coloneqq \mathrm{N}\bed{\cdot}{I^\mathrm{c}}$,
with $Q_\mathrm{rd}\coloneqq Q_{\mathrm{rd},\eta,0}$ according to 
Definition~\ref{Def:Rounding_and_histograms},
and with $\epsilon\coloneqq \mathrm{N}(I^\mathrm{c})$, 
we have $\mathrm{N}=\mathrm{N}(I)\mathrm{N}\bed{\cdot}{I} + \epsilon\, Q$
and $P=\mathrm{N}(I)\mathrm{N}\bed{\cdot}{I} + \epsilon\, Q_\mathrm{rd}$\,,
and hence 
\la                                                  \label{Eq:P-N_Example_Tail-discretised_normal_laws}
 P-\mathrm{N} &=& \epsilon\,(Q_\mathrm{rd}-Q)\,.
\al

\smallskip 2. Let in this part of the proof $t\in\mathopen]0,\infty\mathclose[$ 
and hence also $I,\epsilon, Q$
be fixed, and let any asymptotics refer to $\eta\rightarrow 0$. 
We have $\mu(P) =0$ by symmetry, 
and with $\sigma\coloneqq \sigma(P)$ hence, 
using linearity of $\mu_2$ and~\eqref{Eq:P-N_Example_Tail-discretised_normal_laws}  
in the second step, 
and~\eqref{Eq:zeta_1_P-P_standardised_rounded_asymp} in the third,
\la                                                        \label{Eq:sigma_asymptotics_Example_Tail-discretised_normal}
 \sigma^2-1 &=& \mu_2(P)-\mu_2(\mathrm{N})
  \,\ =\,\ \epsilon\,\mu_2(Q_\mathrm{rd}-Q) \,\ \llcurly \,\ \eta \,.
\al

We next get 
\la                                          \label{Eq:zeta_1_Tail-discretised_normal_-_normal}
 \zeta_1(P-\mathrm{N}) & =&  \epsilon\, \zeta_1(Q_\mathrm{rd}-Q) \,\ \sim\,\ \tfrac{\epsilon}{4}\eta  
\al
by~\eqref{Eq:P-N_Example_Tail-discretised_normal_laws} and~\eqref{Eq:zeta_1_rounded-original_asymp},
and 
\[
 \left| \zeta_1(\widetilde{P}-\mathrm{N})- \zeta_1(P-\mathrm{N} \right|
 &\le & \zeta_1(\widetilde{P}-P) \,\ =\,\ \left|\tfrac{1}{\sigma} -1\right|\nu_1(P)
 \,\ \llcurly \,\ \eta
\]
by using in the second step centredness of $P$ and~\eqref{Eq:zeta_1_under_x_maspto_ax},
and in the third \eqref{Eq:sigma_asymptotics_Example_Tail-discretised_normal}
and boundedness of $\nu_1(P) = \mathrm{N}(I)\nu_1(\mathrm{N}\bed{\cdot}{I}) +
\epsilon\,\nu_1(Q_\mathrm{rd})$ due to 
$\nu_1(Q_\mathrm{rd}) \le \nu_1(Q)+\zeta_1(Q_\mathrm{rd}-Q)\rightarrow \nu_1(Q) $.  
Combining the previous two displays yields
\la                          \label{Eq:zeta_1_standardised_Tail-discretised_normal_-_normal}
   \zeta_1(\widetilde{P}-\mathrm{N}) & = &  \varkappa_1(\widetilde{P}-\mathrm{N}) \,\ \sim\,\ \tfrac{\epsilon}{4}\eta \,.  
\al
Further $\underline{\zeta}_3(P-\mathrm{N}) =\epsilon\, \underline{\zeta}_3(Q-Q_\mathrm{rd}) \llcurly \eta$
by (\ref{Eq:P-N_Example_Tail-discretised_normal_laws},\ref{Eq:zeta_k_rounded-orignal_small}),
and $\underline{\zeta}_3(\widetilde{P}- P) \preccurlyeq \left|\frac{1}{\sigma}-1 \right| \nu_3(P) \llcurly \eta$ 
by centredness of $P$ and \eqref{Eq:zeta_r_under_x_mapsto_ax} in the first step, 
and~\eqref{Eq:sigma_asymptotics_Example_Tail-discretised_normal} 
and also boundedness of $\nu_3(P)$ 
in the second, and hence 
\la \label{Eq:zeta_3_standardised_Tail-discretised_normal_-_normal}
 \zeta_3(\widetilde{P}-\mathrm{N}) 
  &\le& \underline{\zeta}_3(P-\mathrm{N}) + \underline{\zeta}_3(\widetilde{P}- P) 
  \,\ \llcurly\,\ \eta\, .
\al 

Assuming below in the first step $\eta<2t$, recalling the definitions~(\ref{Eq:Def_h_M},\ref{Eq:Def_varkappa}), 
and  using that $h_{P-\mathrm{N}}=\epsilon\,h_{Q_\mathrm{rd}-Q}$ 
vanishes outside of $[-t+\frac{\eta}{2},t-\frac{\eta}{2}]$, we obtain
\la                    \label{Eq:varkappa_3_Tail-discretised_normal_-_normal}
 \varkappa_3(P-\mathrm{N}) &\ge & 3(t-\tfrac{\eta}{2})^2\varkappa_1(P-\mathrm{N})
   \,\ \sim \,\ 3 t^2\tfrac{\epsilon}{4}\eta\,\,
\al
using~\eqref{Eq:zeta_1_Tail-discretised_normal_-_normal}.
Further, 
\la                      \label{Eq:varkappa_3_Tail-discretised_normal_-_standardised}
 \varkappa_3(\widetilde{P}-P) &=& \left| \tfrac{1}{\sigma^3}-1 \right|\nu_3(P) \,\ \llcurly\,\ \eta
\al
by using in the first step centredness of $P$ and~\eqref{Eq:varkappa_different_scalings_=},  
and in the second~\eqref{Eq:sigma_asymptotics_Example_Tail-discretised_normal}
and boundedness of $\nu_3(P)$. 

Combining (\ref{Eq:zeta_1_standardised_Tail-discretised_normal_-_normal},%
\ref{Eq:zeta_3_standardised_Tail-discretised_normal_-_normal},%
\ref{Eq:varkappa_3_Tail-discretised_normal_-_normal},%
\ref{Eq:varkappa_3_Tail-discretised_normal_-_standardised}) yields
\la                        \label{Eq:zeta/varkappa_Tail-discretised_normal}
  \limsup_{\eta\rightarrow0}
    \frac {\zeta_1\!\vee\!\zeta_3}{\varkappa_1\!\vee\!\varkappa_3 }(\widetilde{P}-\mathrm{N}) 
    &\le& \frac{1}{3t^2} \,.
\al

\smallskip3.                                        
Convergence to zero of the right hand sides in
(\ref{Eq:zeta_1_standardised_Tail-discretised_normal_-_normal},%
\ref{Eq:zeta_3_standardised_Tail-discretised_normal_-_normal})
yields the first claim, and the second follows from 
letting $t\rightarrow\infty$ in~\eqref{Eq:zeta/varkappa_Tail-discretised_normal}.
\end{proof}

\section{Some identities, inequalities, and asymptotics for special laws}
                                        \label{sec:special_laws}

The following presumably very well-known fact is used in the discussion
of Corollary~\ref{Cor:Normal_approximation_for_sampling_from_finite_populations_Z-close_to_normal}.
\begin{Lem}[Kolmogorov distance of centred normal laws]      \label{Lem:Kolmogorov_distance_of_centred_normals}
$\|\mathrm{N}_\sigma-\mathrm{N}_\tau\|^{}_{\mathrm{K}}
 = \Phi(\omega x)-\Phi(x)
 \le \frac{1}{\sqrt{2\pi\mathrm{e}}}\frac{|\sigma-\tau|}{\sigma\wedge\tau}$ 
for $\sigma,\tau\in\mathopen]0,\infty\mathclose[$\,, $\omega\coloneqq\frac{\sigma\vee\tau}{\sigma\wedge\tau}$, 
$x\coloneqq\sqrt{\frac{2\log\omega}{\omega^{2}-1}}<1$ if $\omega>1$, $x$ $\coloneqq$ $1$ if $\omega=1$.
\end{Lem}
\begin{proof}
$\|\mathrm{N}_\sigma-\mathrm{N}_\tau\|^{}_{\mathrm{K}} 
 = \|\mathrm{N}_{\sigma\wedge\tau}-\mathrm{N}_{\sigma\vee\tau}\|^{}_{\mathrm{K}}
 = \|\mathrm{N}_\frac1{\omega}-\mathrm{N}\|^{}_{\mathrm{K}} 
 = \sup_{y > 0}\big(\Phi(\omega y)-\Phi(y)\big) = \Phi(\omega x)-\Phi(x)
 \le (\omega-1)x\phi(x) \le (\omega-1)\!\cdot\!1\!\cdot\!\phi(1)
 = \frac{1}{\sqrt{2\pi\mathrm{e}}}\frac{|\sigma-\tau|}{\sigma\wedge\tau}$,
by scale invariance of $\|\cdot\|^{}_{\mathrm{K}}$, symmetry of $\Phi$, 
and differential calculus. 
\end{proof}

The following Examples~\ref{Example:Convolution_inequality_sharp} show 
the sharpness
of Theorem~\ref{Thm:F_star_G_vs_H_star_H} in the case of $H_1=H_2=\Phi$, 
and $F_1=F_2$ close to $\Phi$. They are
simpler relatives of Zolotarev's Example~\ref{Example:Zolotarev_1973} 
below,
and the laws $P_\epsilon$ of part (a)
are 
the simplest examples of extreme points as in
Lemma~\ref{Lem:cK_and_its_extreme_points} with $R\coloneqq\mathrm{N}$.

\begin{Examples}                          \label{Example:Convolution_inequality_sharp}
We consider $\mathopen]0,\infty\mathclose[ \ni \epsilon \rightarrow 0$.

\begin{parts}
\item            \label{part:Example_Convolution_inequality_sharp_a}
 If $P \coloneqq P_\epsilon
    \coloneqq \mathrm{N}(\cdot\setminus [0,\epsilon]) + \big(\Phi(\epsilon)-\tfrac{1}{2}\big)\delta_0$,
then have
\[
  \zeta_1(P -\mathrm{N}) \,\ \sim\,\ \frac{\epsilon^2}{2\sqrt{2\pi}}\ ,
    &\quad& \left\|P_{}^{\ast 2} -  \mathrm{N}^{\ast 2} \right\|_\mathrm{K}  \,\ \sim\,\ \frac{\epsilon^2}{\pi} \ ,
\]
and with $\|\Phi\|^{}_\mathrm{L} = \frac{1}{\sqrt{2\pi}}$
hence $\text{\rm L.H.S\eqref{Eq:F_star_G_vs_H_star_H}}\sim \text{\rm R.H.S\eqref{Eq:F_star_G_vs_H_star_H}}$. 
\item           \label{part:Example_Convolution_inequality_sharp_b}
 If $P \coloneqq P_\epsilon
    \coloneqq \mathrm{N}(\cdot\setminus [-\epsilon,\epsilon]) + \big(\Phi(\epsilon)-\tfrac{1}{2}\big)
    (\delta_{-\epsilon} +\delta_\epsilon)$,
then have $\mu(P)=0$ and
\[
  \zeta_1(P -\mathrm{N}) \,\ \sim\,\ \frac{\epsilon^2}{\sqrt{2\pi}}\ ,
    &\quad& \left\|P_{}^{\ast 2} -  \mathrm{N}^{\ast 2} \right\|_\mathrm{K}  \,\ \sim\,\ \frac{\epsilon^2}{\pi} \ ,
\]
and hence $\text{\rm L.H.S\eqref{Eq:F_star_G_vs_H_star_H}}\sim \frac{1}{2}\,\text{\rm R.H.S\eqref{Eq:F_star_G_vs_H_star_H}}$.
\end{parts}
\end{Examples}
\begin{proof} \ref{part:Example_Convolution_inequality_sharp_a}
With the notation~\eqref{Eq:Def_F_M_and_overline{F}_M} we have
\[
 F_{P-\mathrm{N}}(x) &=& \big(\Phi(\epsilon)-\Phi(x)\big)\1_{[0,\epsilon]}(x) \quad \text{ for }x\in\R,
\]
and hence 
\la                                            \label{Eq:zeta_1_asymptotics_example_conv_ineq_sharp} 
 \zeta_1(P-\mathrm{N}) &=& \int\big| F_{P-\mathrm{N}} \big|\,\dd\leb
    \,\ =\,\ \int_0^\epsilon \big( \Phi(\epsilon)-\Phi(x)\big)\,\dd x \,\ \sim\,\ \frac{\epsilon^2}{2\sqrt{2\pi}}
\al
by (\ref{Eq:zeta_1=kappa_1_new},\ref{Eq:varkappa_1=L^1-enorm}) in the first step and, say,
de l'H\^opital in the last.
The commutative ring identity
\la                                       \label{Eq:ring_identity_P_N_n=2}
 P^{\ast 2} &=& \mathrm{N}^{\ast 2}+2\, (P-\mathrm{N})\ast \mathrm{N} + (P-\mathrm{N})^{\ast 2}, 
\al
which by the way is a special case of~\eqref{Eq:ring_identity} in the proof of 
Theorem~\ref{Thm:F_star_G_vs_H_star_H}, yields here in particular
\[
 \left\| P^{\ast 2} -\mathrm{N}^{\ast 2}\right\|_\mathrm{K}
 & \ge & F_{P^{\ast 2} - \mathrm{N}^{\ast 2}}(0) 
 \,\ =\,\ 2 \int  F_{P-\mathrm{N}}( 0-y)\phi(y)\,\dd y + (P-\mathrm{N})^{\ast 2}(\{0\})  \\
 &= &  2\int_{-\epsilon}^0 \big(\Phi(\epsilon)-\Phi(-y)\big)\phi(y)\,\dd y
    + \big(\Phi(\epsilon)-\tfrac{1}{2}\big)^2 
 \,\ \sim \,\ \frac{\epsilon^2}{2\pi}    + \frac{\epsilon^2}{2\pi} 
 \,\ =\,\ \frac{\epsilon^2}{\pi} \, ,
\]
say again by de l'H\^opital in the penultimate step.
Since, in the other direction, we have
\[
 \left\| P^{\ast 2} -\mathrm{N}^{\ast 2}\right\|_\mathrm{K}
 & \le & \left( 2 \,\sqrt{ \left\|\Phi'\right\|_\infty \zeta_1(P-\mathrm{N})   }\,\right)^2
  \,\ \sim \,\ \frac{\epsilon^2}{\pi} \, ,
\]
by~\eqref{Eq:F_star_G_vs_H_star_H} in the first step, 
and by $\left\|\Phi'\right\|_\infty=\frac{1}{\sqrt{2\pi}}$ 
and~\eqref{Eq:zeta_1_asymptotics_example_conv_ineq_sharp} in the second, the claim follows. 

\ref{part:Example_Convolution_inequality_sharp_b}
The first asymptotics claim is clear as in part~\ref{part:Example_Convolution_inequality_sharp_a}.
For the second using symmetry yields
$ \left\| P^{\ast 2} -\mathrm{N}^{\ast 2}\right\|_\mathrm{K}
 \ge F_{P^{\ast 2} - \mathrm{N}^{\ast 2}}(0)
  = \frac{1}{2}P^{\ast 2}(\{0\}) = 2\big(\Phi(\epsilon)-\frac{1}{2}\big)^2\sim\frac{\epsilon^2}{\pi}$,
and the corresponding upper asymptotic bounding follows using~\eqref{Eq:The_convolution_inequality_new}.
\end{proof}

This following instructive example is treated here in more detail than in the original sources
and in \citet[Examples 1.2 and 1.3]{Yaroslavtseva2008b}; of course one could go still further. 

\begin{Example}[\mbox{\citeposs{Zolotarev1972,Zolotarev1973}}
                        normal laws discretised near zero]              \label{Example:Zolotarev_1973}
For $\epsilon\in\mathopen]0,\infty\mathclose[$\,, let 
\[ 
 P &\coloneqq&   P_\epsilon \,\ \coloneqq\,\   \mathrm{N}(\,\cdot \setminus [-\epsilon,\epsilon])
  +  p \frac{\delta_{-a}+\delta_{a}}{2}
\]
with $p\coloneqq p_\epsilon \coloneqq \mathrm{N}([-\epsilon,\epsilon])$ and 
$a\coloneqq a_\epsilon\coloneqq (\frac{1}{p}\int_{-\epsilon}^\epsilon x^2\phi(x)\,\dd x)^\frac12$.
Let asymptotic comparisons in this example always refer to $\epsilon\downarrow0$,
with any other parameters $n$ or $r$ being fixed.
\smallskip 
\begin{parts}
\item{\rm \textbf{Simple properties.}} 
 Each $P$ is a symmetric law with all moments finite and with $\mu_2(P)=1$, 
 and hence in particular $P = \widetilde{P}\in\widetilde{\cP_3}$
 and $P-\mathrm{N} \in\cM_{4,3}$.  We have 
 \la                                         \label{Eq:p_sim_a_sim_Zolotarev_example}
  p &\sim& \frac{2\epsilon}{\sqrt{2\pi}}\,, \qquad 
   a \,\ \sim\,\ \frac{\epsilon}{\sqrt{3}} \,.
 \al   
\item{\rm \textbf{Asymptotics of CLT errors for $n$ small.}}
                            \label{part:Zolotarev_example_asymptotics}
 Let $n\in\{1,2,3,4\}$ and  $r\in\mathopen[0,\infty\mathclose[\,$,
\la        \label{Eq:Zolotarev_example_def_M}
  M &\coloneqq& \tfrac{1}{2}(\delta_{-1}+\delta_1) 
  - \tfrac{1}{2\sqrt{3}}  \leb(\,\cdot\,\cap\mathopen]-\sqrt{3},\sqrt{3}\mathclose[\,) \,,\\
 M_t &\coloneqq&  M(\tfrac{\cdot}{t}) \,\ =\,\ \tfrac12(\delta_{-t}+\delta_t)
        - \tfrac{1}{2\sqrt{3}t}  \leb(\,\cdot\,\cap\mathopen]-\sqrt{3}t,\sqrt{3}t\mathclose[\,)
           \quad \text{ for } t\in\mathopen]0,\infty\mathclose[\,.
\al
Let $\|\cdot\|$ be any of the following norms on, say,  $\cM_{4,3}\cap \cM_r$
subject to the indicated restrictions on $(n,r)$:
\[
 \|\cdot\|^{}_{\mathrm{K}} &:&  r =0\,, \\
 \nu_r       &:& n=1 \,\text{ or }\, r<5-n\,, \\
 \varkappa_r &:& r>0 \,\text{ and }\, ( n=1\,\text{ or }\, r < 5-n) \,, \\
 \zeta_r     &:& r\in\{0,1,2,3,4\} \,\text{ and }\, (n=1  \,\text{ or }\, r\le 4-n)\,.
\]
Then we have
\la                                     
          \label{Eq:widetile(P^n)-N_sim_n_le_4_Zolotarev-Example}
 \widetilde{P_{}^{\ast n}} -\mathrm{N}
  &\sim&  \big(\frac{2\epsilon}{\sqrt{2\pi}}\big)^n 
                \left(M_{\frac{a}{\sqrt{n}}}\right)^{\ast n} \quad \text{w.r.t. } \|\cdot\|\,, \\
  \left\| \widetilde{P_{}^{\ast n}}  \label{Eq:Zolotarev_example_norm_error_n}  
     -\mathrm{N} \right\| 
   &\sim& \big(\frac{2\epsilon}{\sqrt{2\pi}}\big)^n \big(\frac{\epsilon}{\sqrt{3n}}\big)^r
      \left\|  M^{\ast n} \right\|   \,. 
\al     
\item{\rm\textbf{Norms of $M$ and $M^{\ast 2}$.}}
 With $M$ from~\eqref{Eq:Zolotarev_example_def_M} and for $r\in [0,\infty[\,$, we have
 \la
  &&  \left\| M \right\|_{\mathrm{K}} = \frac{1}{2\sqrt{3}}\,, 
    \quad  \left\| M^{\ast 2} \right\|_{\mathrm{K}} = \frac{1}{4}\,, \\
  &&  \nu_r(M)= \frac{3^\frac{r}{2}}{r\!+\!1}+1   \,,   \quad    \label{Eq:Zolotarev_example_nu_r_M}
      \nu_0(M)=2\,, \quad \nu_0(M^{\ast 2}) = \frac{8-\sqrt{3}}{3} \,, \\
  && \varkappa_r(M)=                                            \label{Eq:Zolotarev_example_kappa_r_M}
  \frac{1}{r\!+\!1}\big( 3^\frac{r}{2}_{}+ \frac{2\sqrt{3}\!-\!3}{3}r -1 \big)
      \quad\text{ if }r>0 \,,    \\
  && \zeta_1(M) = \frac{5\sqrt{3}-6}{6}\, ,        \label{Eq:Zolotarev_example_zeta_r_M}
  \quad   \zeta_3(M) = \frac{3\sqrt{3}-4}{24}\,, \quad  \zeta_4(M) = \frac{1}{30}\,.
 \al 
\item{\rm \textbf{Specialisations.}} 
We have                          
\la
 \left\| P -\mathrm{N} \right\|_{\mathrm{K}}   \label{Eq:||_||_K_Zolotarev-example}
   &\sim& \frac{\epsilon}{\sqrt{3}\sqrt{2\pi}}\,,   \quad       
 \left\| \widetilde{P_{}^{\ast 2}}  -\mathrm{N} \right\|_{\mathrm{K}} \label{Eq:n=2_Kolmogorov_error_Zolotarev_example}
   \,\ \sim\,\ \frac{\epsilon^2}{2\pi}\,, \\
 \left\| \widetilde{P_{}^{\ast n}}  -\mathrm{N} \right\|_{\mathrm{K}}  
          \label{Eq:n=small_Kolmogorov_error_Zolotarev_example}
           \label{Eq:Zolotarev_example_norm_error_n_le_4_K_nu_0}  
   &\asymp& \nu_0(\widetilde{P_{}^{\ast n}} -\mathrm{N})
   \,\ \asymp \,\ \epsilon^n \quad \text{ for } n\in \{1,2,3,4\}\,, \\  
  \zeta_1(P-\mathrm{N})                       \label{Eq:zeta_1_Zolotarev-example}
   &=& \varkappa_1(P-\mathrm{N}) \,\ \sim \,\ \frac{5-2\sqrt{3}}{3\sqrt{2\pi}}\epsilon^2, \\
 \zeta_3(P -\mathrm{N})                 \label{Eq:zeta_3_Zolotarev-example}
    &=&       \frac{1}{6}\big(\nu_3(N)-\nu_3(P)\big) 
   \,\ \sim \,\  \frac{9-4\sqrt{3}}{108\sqrt{2\pi}}\epsilon^4, \\ 
  \qquad                                      \label{Eq:n=2_sharpness_BE_K_Z_Zolotarev_example}                                                                              
 \frac{\left\| \widetilde{P_{}^{\ast 2}}-\mathrm{N} \right\|_{\mathrm{K}}} 
      {\,\frac{1}{\sqrt{2}}\big(\zeta_1\!\vee\!\zeta_3\big)(P-\mathrm{N})\,}  
   &\sim&  \frac{\left\| \widetilde{P_{}^{\ast 2}}-\mathrm{N} \right\|_{\mathrm{K}}} 
      {\,\frac{1}{\sqrt{2}}\zeta_1(P-\mathrm{N})\,}  
   \,\ \rightarrow\,\  \sqrt{2}\, \frac{15+6\sqrt{3}}{13\sqrt{2\pi}}  \,\ =\,\ 1.1020\ldots \,,\\ 
  \nu_r(P-\mathrm{N})        \label{Eq:Zolotarev_example_r-norm_errors_n=1}
    &\asymp& \varkappa_r(P-\mathrm{N})
    \,\ \asymp\,\ \zeta_r(P-\mathrm{N}) \,\ \asymp\,\ \epsilon^{r+1} 
    \quad\text{ for }r\in\{1,2,3,4\} \,, \\    
 \varkappa_r(P-\mathrm{N})                    \label{Eq:varkappa_r_Zolotarev-example}
  &\sim&
    \frac{2}{\sqrt{2\pi}}\left(\frac{1}{r\!+\!1}
       +    \Big( \frac{2r}{r\!+\!1} -\sqrt{3}\Big) 3^{-\frac{r\!+\!1}{2}}\right)\epsilon^{r+1}                                    
  \quad\text{ for } r\in\mathopen]0,\infty\mathclose[   \,, \\
 \nu_r(P-\mathrm{N})                          \label{Eq:nu_r_Zolotarev-example} 
   &\sim& \frac{2}{\sqrt{2\pi}}\left(\frac{1}{r\!+\!1}+3^{-\frac{r}{2}}\right)\epsilon^{r+1}
    \quad\text{ for }r\in\mathopen[0,\infty\mathclose[\,, \\   
  \nu_0\!\left(\widetilde{P_{}^{\ast 2}}  -\mathrm{N} \right)  \label{Eq:n=2_nu_0_error_Zolotarev_example}
   &\sim& \frac{16-2\sqrt{3}}{3\pi}\epsilon^2 , \\ 
                            \label{Eq:n=2_sharpness_Shiganov_1_Zolotarev_example}
 \frac{\left\| \widetilde{P_{}^{\ast 2}}-\mathrm{N} \right\|_{\mathrm{K}}} 
      {\,\frac{1}{\sqrt{2}}\big(\nu_1\!\vee\!\nu_3\big)(P-\mathrm{N})\,}  
   &\sim&  \frac{\left\| \widetilde{P_{}^{\ast 2}}-\mathrm{N} \right\|_{\mathrm{K}}} 
      {\,\frac{1}{\sqrt{2}}\nu_1(P-\mathrm{N})\,}  
   \,\ \rightarrow\,\  \sqrt{2}\, \frac{2\sqrt{3}-3}{\sqrt{2\pi}}  \,\ =\,\ 0.26184  \ldots\,.
\al
\end{parts}
\end{Example} 
\begin{proof}
1. The claims up to $a\sim \frac{\epsilon}{\sqrt{3}}$ are obvious, 
and we have $a\in\mathopen]0,\epsilon\mathclose[$\,, using $x^2< \epsilon^2$ in the defining integral.
In what follows, we will 
omit the convolution symbol $\ast$\,, as explained at the beginning 
of section \ref{sec:zeta_distances}.

\smallskip 2. Let us first prove 
\la
  (P -\mathrm{N})^{n} &\sim&                      \label{Eq:P-N_sim_Zolotarev-Example}
  \big(\frac{2\epsilon}{\sqrt{2\pi}}\big)^n 
   M_{a}^{n} \quad\text{ w.r.t. }\nu_r \,, \text{ for every }n\in\N.
\al
We put $Q\coloneqq\frac12(\delta_{-1}+\delta_1)$ and $\mathrm{U}_t \coloneqq
 \frac{1}{2t}\leb(\,\cdot\,\cap\mathopen]-t,t\mathclose[\,)$ and get
\la
  t^{-r}\nu_r(M_t) &=& \nu_r(M) \,\ = \,\ \nu_r(Q) +\nu_r(\mathrm{U}_{\sqrt{3}})
   \,\ =\,\ 1+\tfrac{3^\frac{r}{2}}{r\!+\!1} \quad \text{ for }t\in\mathopen]0,\infty\mathclose[
\al
by~\eqref{Eq:homogeneity_on_cM_conclusion} in the first step, and by mutual singularity
due to discreteness and continuity in the second.
Hence, using also~\eqref{Eq:p_sim_a_sim_Zolotarev_example}, we get 
$\nu_r(pM_a) = p \nu_r(M_a) \asymp \epsilon^{r+1}$.
Together with 
\[
 \nu_r(P-\mathrm{N}-pM_a) 
  &=& \nu_r\big(p\mathrm{U}_{\sqrt{3}a} - \mathrm{N}(\,\cdot \cap [-\epsilon,\epsilon])\big) 
  \,\ = \,\ \int_{-\epsilon}^\epsilon |x|^r\left|\tfrac{p}{2\sqrt{3}a}-\phi(x)\right|\,\dd x
  \,\ \llcurly \,\ \epsilon^{r+1} \,,
\]
by~\eqref{Eq:p_sim_a_sim_Zolotarev_example} in the last step,
this yields~\eqref{Eq:P-N_sim_Zolotarev-Example} in case of $n=1$, 
by $\text{L.H.S.\eqref{Eq:P-N_sim_Zolotarev-Example}}\sim R-p\mathrm{U}_{\sqrt{3}a} 
\sim \text{R.H.S.\eqref{Eq:P-N_sim_Zolotarev-Example}}$,
using in the last step~\ref{Lem:sim_relations}\ref{part:sim_relations_multiplied} 
with $V_1\coloneqq\R $ and, say,  $V_2\coloneqq \cM_{r,0}$\,.

For a general $n$ we have 
\la                                                  \label{Eq:asymp_epsilon_(n+r)} 
  \text{R.H.S.\eqref{Eq:P-N_sim_Zolotarev-Example}} &\asymp& \epsilon^{n+r}
  \quad\text{ with respect to $\nu_r$,} 
\al
due to $\nu_r( M_a^{ n})  = a^r\nu_r(M^n)$ and
$\nu_r(M^n) = \nu_r\big( \sum\limits_{j=0}^n\binom{n}{j} Q_{}^{j}
               (-\mathrm{U}_{\sqrt{3}})^{n-j}\big)
  \ge \nu_r( Q_{}^{n}) >0$,
using discreteness of $Q^{n}$ and continuity of the other summands
in the penultimate step.

Hence, now for a general $n$ but restricting to the case of $r=0$, 
we get~\eqref{Eq:P-N_sim_Zolotarev-Example} inductively from the case of $n=1$, 
by using~\ref{Lem:sim_relations}\ref{part:sim_relations_multiplied} 
with $V_1\coloneqq V_2\coloneqq V \coloneqq\cM$, each norm being being $\nu_0$,
and with continuity of convolution due to~\eqref{Eq:nu_0_algebra_norm}.

Hence, using the scale invariance~\eqref{Eq:nu_0_and_Kolmogorov_scale_invariant}, 
and hence~\eqref{Eq:sim_on_cM_scale_invariant} for $\|\cdot\|\coloneqq\nu_0$
and $t\coloneqq\frac{1}{a}$, we get 
\la            \label{Eq:sim_a_cdot}
 (P -\mathrm{N})^{n}(a\cdot) \sim \big(\frac{2\epsilon}{\sqrt{2\pi}}\big)_{}^n  M_{a}^{n}(a\cdot) 
   \,\ =\,\ \big(\frac{2\epsilon}{\sqrt{2\pi}}\big)_{}^n  M_{}^{n}
\al
with respect to~$\nu_0$.
This is an asymptotic relation in the vector space 
$V \coloneqq \{M\in \cM_{0,0} : M(\R\!\setminus\![-c, c ])=0\}$
with $c\coloneqq \sup_{\epsilon\in \mathopen]0,\epsilon_0\mathclose]}\frac{\epsilon}{a}
\le n\left(\sqrt{3}+1\right)$ for $\epsilon_0$ small enough, 
by using~\eqref{Eq:p_sim_a_sim_Zolotarev_example}.
Since we  have $\nu_r \preccurlyeq \nu_0$ on $V$, 
and $  \nu_0(\text{R.H.S.\eqref{Eq:sim_a_cdot}}) \asymp \epsilon^n
\asymp \nu_r(\text{R.H.S.\eqref{Eq:sim_a_cdot}})$,
 we get~\eqref{Eq:P-N_sim_Zolotarev-Example}
as stated by applying~\ref{Lem:sim_relations}\ref{part:sim_relations_with_two_norms}.

\smallskip 3. 
Let $n\in\N$. Then, in generalisation of~\eqref{Eq:ring_identity_P_N_n=2}, 
we have the commutative ring identity
\[
 P^n -\mathrm{N}^n -(P-\mathrm{N})^n
  &=& (P-N)\Big(\sum_{j=0}^{n-1}P^{n-1-j}\mathrm{N}^j -  
  \sum_{j=0}^{n-1}\textstyle{\binom{n-1}{j}}P^{n-1-j}(-\mathrm{N})^j\Big) \\
  &=& (P-\mathrm{N})\mathrm{N}R_n
\]
with $R_n$ being some polynomial function, depending only on $n$, 
of the laws $P$ and $\mathrm{N}$,
for example $R_1=0$, $R_2= 2\delta_0$, $R_3 = 3P$, 
$R_4=4P^2-2P\mathrm{N}+ 2\mathrm{N}^2$.
Hence, for $r\in[0,\infty[$\,, $\epsilon\in\mathopen]0,1\mathclose]$, 
and with finite constants $c^{}_{r,n},c'^{}_{r,n},c''^{}_{r,n} $ 
depending only on $r$ and $n$, we get
\[
 \nu_r\big(  P^n -\mathrm{N}^n -(P-\mathrm{N}))^n \big)
  &\le& 2^{r\vee1} \big(\nu_0\!\vee\!\nu_r\big)((P-\mathrm{N})\mathrm{N})  
     \big(\nu_0\!\vee\!\nu_r\big)(R_n) \\
  &\le& c^{}_{r,n} \big(\nu_0\!\vee\!\nu_r\big)((P-\mathrm{N})\mathrm{N}) 
  \,\ \le \,\  c'^{}_{r,n}\,\zeta_4(P-\mathrm{N}) \\
  & \le & \tfrac{c'^{}_{r,n}}{4!}\, \nu_4 (P-\mathrm{N})  
  \,\ \le \,\  c''^{}_{r,n}\,\epsilon^5
\]
by~\eqref{Eq:nu_0_with_nu_r_submultiplicative} in the first step, 
boundedness of $\nu_r(P)$ and~\eqref{Eq:nu_0_with_nu_r_submultiplicative}
again in the second,
 \eqref{Eq:nu_r_smoothed_vs_zeta_k} and $(P-\mathrm{N})(\,\cdot\!\setminus\![-1,1])=0$
in the third, 
\eqref{Eq:zeta_vs_varkappa_etc} and $P-\mathrm{N}\in\cM_{4,3}$ in the fourth, 
and (\ref{Eq:P-N_sim_Zolotarev-Example},\ref{Eq:asymp_epsilon_(n+r)})
with~$1$ in place of $n$ 
in the fifth.
Hence, using now (\ref{Eq:P-N_sim_Zolotarev-Example},\ref{Eq:asymp_epsilon_(n+r)})
for the present $n$, we get 
\[
 P^n-\mathrm{N}^n &\sim& (P-\mathrm{N})^n \quad\text{w.r.t.~$\nu_r$ if $n=1$ or $r<5-n$}\,.
\]
This, combined with~\eqref{Eq:P-N_sim_Zolotarev-Example} and using 
the scaling behaviour~\eqref{Eq:nu_r_scaling_behaviour}, yields
the claim of~\eqref{Eq:widetile(P^n)-N_sim_n_le_4_Zolotarev-Example} 
in case of $\|\cdot\|=\nu_r$.

\smallskip 4. Let $\|\cdot\|$ be any of the norms as specified in~\ref{part:Zolotarev_example_asymptotics}.
Then
\[
 \|\text{R.H.S.\eqref{Eq:widetile(P^n)-N_sim_n_le_4_Zolotarev-Example}}\|
  &=&  \big(\frac{2\epsilon}{\sqrt{2\pi}}\big)^n \big(\frac{a}{\sqrt{n}}\big)^r\left\|M^n\right\|
  \,\ \asymp \,\ \nu_r( \text{R.H.S.\eqref{Eq:widetile(P^n)-N_sim_n_le_4_Zolotarev-Example}})\,,
\]
by the scaling properties (\ref{Eq:homogeneity_on_cM_conclusion},\ref{Eq:nu_r_scaling_behaviour}).
Hence, using~\ref{Lem:sim_relations}\ref{part:sim_relations_with_two_norms} and $\|\cdot\|\preccurlyeq\nu_r$,
we get~\eqref{Eq:widetile(P^n)-N_sim_n_le_4_Zolotarev-Example} from the above part~3 of this proof,
and then~\eqref{Eq:Zolotarev_example_norm_error_n} 
by~\ref{Lem:sim_relations}\ref{part:sim_relation_implies_same_for_norms} 
and~\eqref{Eq:p_sim_a_sim_Zolotarev_example}.
%

\smallskip 5.
We have $\nu_r(M)=\nu_r(\mathrm{U}_{\sqrt{3}})+\nu_r(Q)$, and hence~\eqref{Eq:Zolotarev_example_nu_r_M}
except for the value of $\nu_r(M^2)$. 
Let here $F\coloneqq F_M$, so $F(-x)=-F(x)$ for $x\in\R$, 
\[
 F(x) &=& -\frac{x}{2\sqrt{3}}(0\le x<1)+ \big(\,\frac{1}{2}-\frac{x}{2\sqrt{3}}\, \big)(1\le x\le \sqrt{3}\,)
   \quad \text{ for }x\in [0,\infty[\,,
\]
and hence $\|M\|^{}_{\mathrm{K}} = \sup_{x\in\R}|F(x)|
=\max\{\frac{1}{2\sqrt{3}},\frac{1}{2}- \frac{1}{2\sqrt{3}}\}=\frac{1}{2\sqrt{3}}$ and, if $r>0$,  
\[
 \varkappa_r(M)&=& 2\int_0^\infty\! rx^{r-1}|F(x)|\,\dd x
  = r\,\Big( \int_0^1\frac{x^r}{\sqrt{3}}\dd x + \int_1^{\sqrt{3}} x^{r-1}\big(1-\frac{x}{\sqrt{3}}\big)\dd x\Big)
  \,\ =\,\ \text{R.H.S.\eqref{Eq:Zolotarev_example_kappa_r_M}}, 
\]
and hence also $\zeta_1(M)=\varkappa_1(M)$ as claimed in~\eqref{Eq:Zolotarev_example_zeta_r_M}.

\smallskip 6. We have $ M^2 = \tfrac14\delta_{-2}+\tfrac12\delta_0+\tfrac14\delta_2 +f\leb$ with 
\[
 f(x) &\coloneqq& \tfrac{1}{2\sqrt{3}}\left( \big(1-\tfrac{|x|}{2\sqrt{3}}\big)_+ 
     -  \big( (\,|x+1|\le\sqrt{3}\,) +(\,|x-1|\le\sqrt{3}\,) \big)\right) \quad \text{ for } x\in\R\,.
\]
The function $f$ is even with 
\[
 f(x)&=& \frac{1}{12} \left\{\begin{array}{ll}
      -2\sqrt{3}-x   & \text{ if }0\le x\le\sqrt{3}-1,\\
      -x & \text{ if } \sqrt{3}-1 < x \le \sqrt{3}+1,   \\
      2\sqrt{3}-x& \text{ if } \sqrt{3}+ 1 < x \le 2\sqrt{3},  \\
        0 & \text{ if }x\ge2\sqrt{3} ,
   \end{array}
  \right.
\]
and we get
$ \nu_0(M^2) =  1+ \int\! |f|\,\dd\leb \,\ =\,\ \frac{8-\sqrt{3}}{3}$ and,
setting now $F\coloneqq F_{M^2}$, 
we have $M(\R)=0$ and by symmetry then $F(0)=\frac14$ and, 
using piecewise monotonicity in the second step,
\[
 \left\|M\right\|_{\mathrm{K}} &=& \sup\{ |F(x)| : x\in[0,\infty[\} 
  \,\ =\,\ \max\left\{F(0)\,, -F(2-)\,, F(2)\,, -F(\sqrt{3}+1) \right\}\\
  &=& \max\left\{ \frac{1}{4}\,, \frac{5-2\sqrt{3}}{12}\,, \frac{\sqrt{3}-1}{6}\,, \frac{\sqrt{3}-1}{12}\right\}
  \,\ = \,\ \frac{1}{4} \,.
\]

7. To compute $\zeta_3(M)$ und $\zeta_4(M)$
as claimed in~\eqref{Eq:Zolotarev_example_zeta_r_M}, 
we use Theorem~\ref{Thm:Cut_criteria}\ref{part:zeta_r=nu_r/r!}, applied to the present
$-M$ and with $r\in\{3,4\}$. 
We have $-M=f\mu$ with 
$\mu\coloneqq \leb + \delta_{-1} +\delta_1$ and 
$f\coloneqq \frac{1}{2\sqrt{3}}\1_{\mathopen]-\sqrt{3},\sqrt{3}\mathclose[\setminus\{-1,1\}} 
 -\frac{1}{2}\1_{\{-1,1\}}\,$,
and here  $S^-(f)=4$ and $f$ initially positive.

If now $r=3$, then, in the notation of 
Theorem~\ref{Thm:Cut_criteria}\ref{part:zeta_r=nu_r/r!},
we have $S^-(F_0)=S^-(f)=4=r-0+1$, 
so that condition $(C_0)$ is fulfilled, 
hence also $(C_3)$ by~\ref{Thm:Cut_criteria}\ref{part:zeta_r=translated_nu_r/r!}, 
and by symmetry and~\ref{Thm:Cut_criteria}\ref{part:zeta_r=nu_r/r!} then 
$\zeta_3(M)=\zeta_3(-M)= -\frac{1}{3!}\int|x|^3\,\dd M(x) = \frac{3\sqrt{3}-4}{24}$.

If instead $r=4$, then $S(F_0)=r-0$, hence condition $(B_0)$ is fulfilled, 
and then~\ref{Thm:Cut_criteria}\ref{part:zeta_r=mu_r/r!} yields 
$\zeta_4(M)=\zeta_4(-M)= -\frac{1}{4!}\int x^4\,\dd M(x) = \frac{1}{30}$.

Alternatively,  
$\zeta_4(M)$ is, essentially by  the very definition in~\eqref{Eq:Def_zeta_r}, 
the optimal error bound 
for the two-point Gauss quadrature on the interval $[-\sqrt{3},\sqrt{3}]$
for functions with their fourth derivative bounded in modulus by $1$, 
and hence, 
by \citet[p.~80, 3.5.19 and 3.5.21]{Olver2010NIST},
$\zeta_4(M)=(\,2\sqrt{3}\,)^{2n+1}\frac{\gamma_n}{(2n)!} = \frac{1}{30}$ 
with $\gamma_n=\frac{2^{2n+1}}{2n+1} \frac{(n!)^4}{((2n)!)^2}$
and $n=2$.

\smallskip 8.  Part (d) now follows easily.
\end{proof}

In the proof of Examples~\ref{Examples:zeta_3=mu_3}(a,b)
and~\ref{Example:beta_llcurly_zeta_1} below, we use:

\begin{Lem}             \label{Lem:Standard_normal_identies_and_asymptotics}
Below, the stated exact identities hold for $t\in\R$, 
and the $O(\ldots)$-relations hold for $t\in\mathopen]0,\infty\mathclose[$\,.
\la                       
   \Phi(-t) &=&                                     \label{Eq:Phi(-t)_asymptotic}
      \phi(t)\left(\frac{1}{t}-\frac{1}{t^3} +O\left(\frac{1}{t^5}\right)\right), \\  
   \int_{-t}^\infty x\phi(x)\,\dd x &=&  \phi(t)   \,,  \label{Eq:int_x_phi(x)_dx}    \\  
   \int_{-t}^\infty x^2\phi(x)\,\dd x &=& 1-t\,\phi(t) -\Phi(-t)   
      \,\ =\,\   1-t\,\phi(t) -\frac{\phi(t)}{t} + O\left(\frac{\phi(t)}{t^3}\right)\,, \\  
   \int_{-t}^\infty x^3\phi(x)\,\dd x &=& (t^2+2)\phi(t)    \,, \label{Eq:int_-t^infty_x^3_phi(x)_dx}  \\ 
   \int_{-\infty}^{-t} \Phi(x)\,\dd x                           \label{Eq:int_-infty^-t_Phi(x)_dx}
        &=& \phi(t)-t\,\Phi(-t)  
            \,\ =\,\  \frac{\phi(t)}{t^2}  +  O\left(\frac{\phi(t)}{t^4}\right)  \,, \\
   \int_{-t}^\infty (1-\Phi(x))\,\dd x &=& t\,\big(1-\Phi(-t)\big) +\phi(t) \label{Eq:int_(1-Phi(x))_dx}
     \,\ =\,\ t + O\left(\frac{\phi(t)}{t^2}\right) \,.
\al
\end{Lem}
\begin{proof}
 The exact identities, namely the respectively first identities in 
 (\ref{Eq:int_x_phi(x)_dx}--\ref{Eq:int_(1-Phi(x))_dx}), are obvious by differentiation. 
 The $O(\ldots)$-relation~\eqref{Eq:Phi(-t)_asymptotic} is well-known 
 to follow from writing $\Phi(-t)=\int_{-\infty}^{-t}\phi(x)\,\dd x
 = \frac{\phi(t)}{t} \int_{-\infty}^{0}\exp(-\frac{x^2}{2t^2})  \mathrm{e}^x\,\dd x$, 
 by the change of variables $x\mapsto \frac{x}{t}-t$ with $t>0$,
 and then using $1-z< \exp(-z)<1-z+\frac{z^2}{2}$ for $z\coloneqq\frac{x^2}{2t^2}>0$.
 The remaining $O(\ldots)$-relations follow easily. 
\end{proof}

\begin{proof}[Proof of Examples~\ref{Examples:zeta_3=mu_3}]
 We will use                            \label{page:Proof_of_Examples:zeta_3=mu_3}
 the general formulae $\sigma^2(P)=\mu_2(P)-\mu^2(P)$ for $P\in\Prob_2(\R)$ and 
 \la                            \label{Eq:mu_3_widetilde_P}
   \mu_3(\widetilde{P}) &=& \frac{\mu_3(\bdot{P})}{\sigma^3(P)}
    \,\ =\,\ \frac{1}{\sigma^3(P)}\big(\mu_3(P)-3\mu(P)\mu_2(P)+2\mu^3(P) \big)
    \quad\text { for }P\in\cP_3   \,,
 \al
 and also
 \la                              \label{Eq:nu_1_P_dot-N_le}
  \big| \nu_1(\bdot{P})- \nu_1(\mathrm{N})\big| 
   &\le &  
   \zeta_1(\bdot{P}-P\big) + \zeta_1(P-\mathrm{N}\big)  \quad\text{ for }P\in\Prob_1(\R)\,,
 \al
 which follows from~\eqref{Eq:zeta_vs_varkappa_etc} for $r\coloneqq1$ and  
 $M\coloneqq \bdot{P}-\mathrm{N}$, 
 together with the triangle inequality for~$\zeta_1$. 
 In each of the three parts,  obviously $P\in\cP_3\!\setminus\!\{\mathrm{N}\}$.

 In parts~\ref{part:Left-truncated_normal_laws} and~\ref{part:Left-winsorised_normal_laws},
 we have $P\stge \mathrm{N}$ in the sense of~\eqref{Eq:usual_stochastic_order},
 and hence
 (\ref{Eq:zeta_1=kappa_1},\ref{Eq:kappa(P-Q)_if_stle},\ref{Eq:lambda_r=mu_r_if_r_odd})
 and $\mu(\mathrm{N})=0$ yield
 \la                           \label{Eq:zeta_1(P-N)=mu(P)}
   \zeta_1(P-\mathrm{N}) &=&  \mu(P-\mathrm{N}) \,\ =\,\ \mu(\mathrm{P})
      \quad \text{ in parts ~\ref{part:Left-truncated_normal_laws} and~\ref{part:Left-winsorised_normal_laws}}.
 \al
 In parts \ref{part:Left-truncated_normal_laws} and
 \ref{part:Gamma_laws_and_power_transforms}, 
 obviously $\widetilde{f}-\phi$ is initially negative, 
 since $\widetilde{f}$ initially vanishes.
 In part~\ref{part:Left-winsorised_normal_laws}, obviously $\widetilde{F}-\Phi$ 
 is initially negative.

\smallskip\ref{part:Left-truncated_normal_laws} 
 Let $I
       \coloneqq \mathopen]-t,\infty\mathclose[$
 and $\widetilde{I} \coloneqq \{ x\in\R: \widetilde{f}(x) >0\}
  = \mathopen]-\frac{t+\mu(P)}{\sigma(P)} ,\infty\mathclose[\,$.
 We then have $S^-(\widetilde{f}-\phi)\le S^-\big((\widetilde{f}-\phi)|_{\widetilde{I}}\big) + 1$ and 
 \la                               \label{Eq:S^-_Left-truncated_normal_density_on_I} 
  S^-\big((\widetilde{f}-\phi)|^{}_{\widetilde{I}}\big)
  &= &  S^-\big(\log\left.\tfrac{\widetilde{f}}{\phi}\right|_{\widetilde{I}}\big) 
  \,\ = \,\ S^-(\text{a quadratic polynomial})
  \,\ \le \,\ 2\,,
 \al
 and hence Lemma~\ref{Lem:zeta_3_computation} here yields the claim up to~\eqref{Eq:zeta_3=mu_3_in_Examples}.
 
 %
%
%
 For $t\rightarrow\infty$, 
 using in several steps below Lemma~\ref{Lem:Standard_normal_identies_and_asymptotics}
 and also $t^k\phi(t)\rightarrow 0$ for each $k\in\Z$, we get 
 \la
  1-\mathrm{N}(I) &=& \Phi(-t) \,\ \sim\,\  \frac{\phi(t)}{t}\,,\nonumber \\
  \mu(P) &=& \frac{\int_{-t}^\infty x\,\phi(x)\,\dd x  } {\mathrm{N}(I)}
    \,\ =\,\ \frac{\phi(t)}{\Phi(t)} \,\ \sim\,\ \phi(t)\,,\nonumber \\
  1-\mu_2(P) &=& \frac{\mathrm{N}(I) - \int_{-t}^\infty  x^2\,\phi(x)\,\dd x } {\mathrm{N}(I)}
    \,\ =\,\ \frac{t\phi(t)}{\Phi(t)}
    \,\ \sim\,\ t\,\phi(t)\,, \nonumber\\  
  \mu_3(P) &=& \frac{\int_{-t}^\infty x^3\,\phi(x)\,\dd x}{\mathrm{N}(I)}
    \,\ =\,\ \frac{(t^2+2)\phi(t)}{\Phi(t)}
  \,\ \sim\,\ t^2\phi(t)\,,\nonumber \\
  1-\sigma^2(P) &=& 1-\mu_2(P) +\mu^2(P) \,\ \sim\,\ t\,\phi(t)\,,   \nonumber \\
  1-\sigma(P) &=& \frac{1-\sigma^2(P)}{1+\sigma(P)} \,\ \sim\,\ \tfrac{1}{2}t\,\phi(t)\,,\label{Eq:1-sigma(P)_asymp}\\
  \mu_3(\widetilde{P}) &=& \frac{\mu_3(\bdot{P})}{\sigma^3(P)}                    \nonumber
    \,\ =\,\ \frac{1}{\sigma^3(P)} \big(\mu_3(P)-3\mu(P)\mu_2(P)+2\mu^3(P) \big) \,\ \sim\,\ t^2\phi(t)\,,
\al
and hence
the first relation in~\eqref{Eq:zeta_3_1_truncated_normal_standardised}, using~\eqref{Eq:zeta_3=mu_3_in_Examples}.  

We further get, using~\eqref{Eq:zeta_1(P-N)=mu(P)} in the first step,
\la                                        \label{Eq:zeta_1(P-N)_left-truncated_normal}
 \zeta_1(P-\mathrm{N})  &=& \mu(P) \,\ \sim \,\  \phi(t)  \,.
\al
Next, \eqref{Eq:zeta_1_under_x_maspto_a+x} with $(M,a,b) \coloneqq  (P,0,-\mu(P))$ yields
\la                                  \label{Eq:zeta_1(P_dot-P)_left-truncated_normal}
 \zeta_1(\bdot{P}-P) &=& |\mu(P)| \,\ \sim \,\ \phi(t)\,.
\al

Starting from~\eqref{Eq:zeta_1_under_x_maspto_ax} with $(M,a,b)\coloneqq(\bdot{P},1,\frac{1}{\sigma(P)})$,
we get 
\[
 \zeta_1(\widetilde{P}-\bdot{P}) &=& \big|\frac{1}{\sigma(P)} - 1\big| \nu_1(\bdot{P})
    \,\ \sim \,\ \tfrac{1}{2}t\,\phi(t)\nu_1(\mathrm{N}) 
    \,\ =\,\ \tfrac{1}{\sqrt{2\pi}}\phi(t)
\]
by using (\ref{Eq:1-sigma(P)_asymp},\ref{Eq:nu_1_P_dot-N_le},\ref{Eq:zeta_1(P-N)_left-truncated_normal},%
\ref{Eq:zeta_1(P_dot-P)_left-truncated_normal}) in the second step,
and~\eqref{Eq:nu_1_and_nu_3_of_N} in the last.

Hence we have~\eqref{Eq:zeta_1_three_asymptotics_truncated_normal},
hence the second asymptotic equality in~\eqref{Eq:zeta_3_1_truncated_normal_standardised}
and the final claim of part~\ref{part:Left-truncated_normal_laws},  
and~\eqref{Eq:zeta_1_vee_zeta_3_over_Esseen1956_truncated_normal_etc}
follows from~(\ref{Eq:zeta_3=mu_3_in_Examples},\ref{Eq:zeta_3_1_truncated_normal_standardised}).

\smallskip\ref{part:Left-winsorised_normal_laws}  
 Let again $I\coloneqq I_t\coloneqq \mathopen]-t,\infty\mathclose[\,$,
 and 
 $\widetilde{I} \coloneqq \{ x\in\R: \widetilde{F}(x-) >0\}
  = \mathopen]-\frac{t+\mu(P)}{\sigma(P)} ,\infty\mathclose[\,$.
 We have 
 \la                              \label{Eq:Left-winsorised_normal_mu>0_sigma_<1}
   && \mu(P)\,>\,0\,, \quad \sigma(P) \, <\, 1\,, 
 \al
 with the second inequality following from 
 Corollary~\ref{Cor:Variance_under_winsorisation},
 and the first being even more obvious.
  
 We have
 $S^-( (\widetilde{F}|_{\widetilde{I}})'-\phi|_{\widetilde{I}})\le 2$
 as in~\eqref{Eq:S^-_Left-truncated_normal_density_on_I},
 and the function $\widetilde{I}\ni x \mapsto (\widetilde{F}-\Phi)'(x)
  = \sigma(P)\phi\big(\sigma(P)x+\mu(P)\big)- \phi(x)$ is finally positive 
  due to $\sigma(P)<1$.
 Hence, say by \citet[Lemma~2.8(b,a), with the present $\widetilde{I}$ 
  in the role of $I$ there]{MattnerShevtsova},
 $S^-((\widetilde{F}-\Phi)|_{\widetilde{I}})\le 2 $ and  $\widetilde{F}-\Phi$ is finally negative.
 Therefore $S^-(\widetilde{F}-\Phi)\le 2+1 = 3$
 and hence, $\widetilde{F}-\Phi$ being initially as well as finally negative,
 $S^-(\widetilde{F}-\Phi)\le 2$.
 
 We have $\mu_k(P)=(-t)^k\Phi(-t)+\int_{-t}^\infty x^k\phi(x)\,\dd x$ for $k\in\N$  
 and hence get, 
 using~(\ref{Eq:Phi(-t)_asymptotic}--\ref{Eq:int_-t^infty_x^3_phi(x)_dx},\ref{Eq:mu_3_widetilde_P}),
 \la
  \mu(P) &=& -t\,\Phi(-t) +\phi(t) \,\ \sim\,\ \frac{\phi(t)}{t^2}\,,\label{Eq:mu_leftwinsorised_normal} \\
  \mu_2(P) &=& t^2\Phi(-t) + 1-t\,\phi(t) - \Phi(-t) 
           \,\ = \,\ 1-2\frac{\phi(t)}{t} + O\left(\frac{\phi(t)}{t^3}\right)\,,\nonumber   \\
  \mu_3(P) &=& -t^3\Phi(-t) +(t^2+2)\phi(t) \,\ \sim\,\ 3\,\phi(t)\,, \nonumber   \\
  \sigma^2(P) &=& 1-2\frac{\phi(t)}{t} + O\left(\frac{\phi(t)}{t^3}\right)\,, \nonumber   \\
  1-\sigma(P) &=& \frac{1-\sigma^2(P)}{1+\sigma(P)} \,\ \sim\,\ \frac{\phi(t)}{t} \,,\label{Eq:sigma_leftwinsorised_normal}  \\
  \mu_3(\widetilde{P}) &\sim& \mu_3(P)\,\ \sim\,\  3\,\phi(t) \,,  \nonumber
 \al
 and hence, using~\eqref{Eq:zeta_3=mu_3_in_Examples}, the first relation 
 in~\eqref{Eq:zeta_3_1_truncated_normal_winsorised}.  
 
 We further get, using (\ref{Eq:zeta_1(P-N)=mu(P)},\ref{Eq:mu_leftwinsorised_normal},
\ref{Eq:zeta_1_under_x_maspto_a+x},%
\ref{Eq:zeta_1_under_x_maspto_ax},\ref{Eq:sigma_leftwinsorised_normal},%
\ref{Eq:nu_1_P_dot-N_le},\ref{Eq:nu_3_von_N}), 
 \[
  \zeta_1(P-\mathrm{N}) &=& \mu(P) 
    \,\ \sim\,\ \frac{\phi(t)}{t^2} \,, \qquad
  \zeta_1(\bdot{P}-P) \,\ =\,\ |\mu(P)| \,\ \sim \,\ \frac{\phi(t)}{t^2}\,, \\
  \zeta_1(\widetilde{P}-\bdot{P}) &=&  \big|\frac{1}{\sigma(P)} - 1\big| \nu_1(\bdot{P}) 
    \,\ \sim \,\ \tfrac{2}{\sqrt{2\pi}}\frac{\phi(t)}{t} 
 \]
 and hence the rest of~\eqref{Eq:zeta_3_1_truncated_normal_winsorised}. 
 
This concludes the proof of the present part, but we will 
continue here in the proof of 
Example~\ref{Example:Left-winsorised_normal_re_Nagaev} 
on page~\pageref{page:Proof_of_Example_Left-winsorised_normal_re_Nagaev}.
 
\smallskip\ref{part:Gamma_laws_and_power_transforms}  
 Here $\lambda^\frac{1}{\beta}$ is just a scale parameter, 
 and we may therefore assume $\lambda =1$ in what follows.
 The claim about the finiteness and then the value of $\nu_r(P)$
 in~\eqref{Eq:nu_r_power_transformed_gamma} is easily checked. 
 
 From now on we use the assumption $\alpha+\frac{2}{\beta}>0$ for
 the existence of $\mu\coloneqq\mu(P)$, $\sigma\coloneqq\sigma(P)$,
 and we let $\widetilde{f}$  denote the  $\leb$-density of $\widetilde{F}$ defined by 
$\widetilde{f}(x)\coloneqq \sigma f_{\Gamma_{\alpha,1,\beta}}(\sigma x+\mu)$ for $x\in\R$, 
so that $\{\widetilde{f}>0\} =\mathopen]-\frac{\mu}{\sigma},\infty\mathclose[ 
\eqqcolon \widetilde{I}$.
Let further $h(x) \coloneqq \log(\widetilde{f}(x)/\phi(x))$ for $x\in \widetilde{I}$, 
$\gamma\coloneqq\alpha\beta-1$,
and $g(t)\coloneqq \gamma -\frac{\mu}{\sigma^2}t+\frac{t^2}{\sigma^2} -\beta t^\beta$
for $t\in\mathopen]0,\infty[$, so that 
$  h'(x) = \frac{\sigma}{t}g(t)$  for $x\in \widetilde{I}$ and 
$t\coloneqq \sigma x+\mu$, and therefore $S^{-}(h')=S^{-}(g)$. 

We hence get 
\[
 S^{-}(\widetilde{f}-\phi)  
  &\le& S^{-}((\widetilde{f}-\phi)|^{}_{\widetilde{I}}) + (\beta>0)(\alpha\beta\le1)  \\
  &=& S^{-}(h) + (\beta>0)(\gamma\le 0)  \\
  &\le& S^{-}(g) + 1  +(\beta>0)(\gamma\le0) \\
  &\le&
   \left\{ 
    \begin{array}{ll}
      S^{-}(-\beta,\gamma,-\tfrac{\mu}{\sigma^2},\tfrac{1}{\sigma^2}) +1\ =\ 3
        &\text{ if }\beta<0,\\ 
      S^{-}(\gamma,-\tfrac{\mu}{\sigma^2},\tfrac{1}{\sigma^2}) +2 \ =\ 3
        &\text{ if }0<\beta\le2\text{ and }\gamma\le 0,\\ 
      S^{-}(\gamma,-\tfrac{\mu}{\sigma^2},\tfrac{1}{\sigma^2}) +1 \ =\ 3 
        &\text{ if }0<\beta\le2\text{ and }\gamma> 0,\\
      S^{-}(\gamma,-\tfrac{\mu}{\sigma^2},\tfrac{1}{\sigma^2},-\beta) +2 \ =\ 4 
        &\text{ if }\beta>2\text{ and }\gamma\le 0,\\  
      S^{-}(\gamma,-\tfrac{\mu}{\sigma^2},\tfrac{1}{\sigma^2},-\beta) +1 \ =\ 4 
        &\text{ if }\beta>2\text{ and }\gamma> 0  
    \end{array}
   \right.   \\
  &=& \left\{\begin{array}{ll} 3 &\text{ if }\beta\le 2,\\
                               4 &\text{ if }\beta> 2
             \end{array}
      \right.       
\]
by using in the third step Rolle's theorem,
as given for example in \citet[p.~510, Lemma 2.8(b)]{MattnerShevtsova},
to bound 
the number of sign changes of an absolutely continuous function by those
of its derivative,   
and using in the fourth step the theorem of Laguerre, as presented 
in \citet[pp.~46--47, number 77]{PolyaSzegoeII}, which bounds $S^{-}(g)$ by 
the number of sign changes of the at most four 
coefficients of $g$, ordered according to increasing exponents.

If $S^-(\widetilde{f}-\phi)\le 3$, and hence in particular if $\beta\le 2$,
the Lemma~\ref{Lem:zeta_3_computation}\ref{part:phi-f_tilde_3_sign_changes}
yields $S^-(\widetilde{f}-\phi)= 3$
and~\eqref{Eq:zeta_3=mu_3_in_Examples}.
If now $\beta>2$, then $\widetilde{f}-\phi$ is essentially not only initially
but also finally negative, and hence assuming $S^-(\widetilde{f}-\phi)\le 3$
would by Lemma~\ref{Lem:zeta_3_computation}\ref{part:phi-f_tilde_3_sign_changes}
yield $\mu_3(\widetilde{P})=0$ in contradiction to~\eqref{Eq:zeta_3=mu_3_in_Examples}.
Hence $S^-(\widetilde{f}-\phi)=4$ in case of $\beta>2$.

The (first) identity in~\eqref{Eq:mu_3_Gamma_power_transformed_standardised}
follows from~\eqref{Eq:mu_3_widetilde_P}, 
using $\mu_k(P)=\nu_k(P)$ for $k\in\{1,2,3\}$ and~\eqref{Eq:nu_r_power_transformed_gamma}, 
and specialises in case of $\beta=1$ to~\eqref{Eq:mu_3_Gamma_standardised}. 

For the $O(\ldots)$-claim in~\eqref{Eq:mu_3_Gamma_power_transformed_standardised}, 
we recall from \citet[pp.~135--136, (6) and (5$'$)]{TricomiErdelyi1951}
the asymptotic expansion
\[
  G(a,x) 
  & =& x^a\cdot\left(1+\binom{a}{2}x^{-1} + \frac{3a-1}{4}\binom{a}{3}x^{-2} 
                    + \binom{a}{2}\binom{a}{4}x^{-3}  + O(x^{-4})\right)
\]
for, say, $a\in\R$ fixed and variable real $x\ge 1\vee(-a+1)$, 
and conclude for $\beta\in\R\setminus\{0\}$ fixed and variable 
$\alpha \ge 1\vee (-\frac{2}{\beta}+1)$, setting $a\coloneqq\frac{1}{\beta}$ and 
$x\coloneqq \alpha$,
\[
 \sigma^2(P)
  &=& G(2a,x) - G_{}^2(a,x)  \\ 
  &=& x^{2a}\left( 
    1 + \binom{2a}{2}x^{-1} +  O(x^{-2})
    - \left(1+ \binom{a}{2}x^{-1} +  O(x^{-2}) \right)^2 \right) \\
  &=&  x^{2a}\left( \left(
    \binom{2a}{2}
   -2\binom{a}{2}
   \right)x^{-1} + O(x^{-2}) \right)  \\
 &= &  a^2 x^{2a-1} + O(x^{2a-2})
 \,\ = \,\  \frac{\alpha^{\frac{2}{\beta}-1}}{\beta^2} + O(\alpha^{\frac{2}{\beta}-2})\,.
\]
If $\alpha \ge 1\vee (-\frac{3}{\beta}+1)$, we similarly get 
\[  
 \mu_3(\bdot{P}) &=& G(3a,x)
      -3\,G(2a,x)G(a,x) 
      +2\,G^{3}_{}(a,x) \\
  &=&\textstyle x^{3a}
    \Big(
      \Big(
             \binom{3a}{2} - 3\left(\binom{2a}{2} + \binom{a}{2} \right) +6 \binom{a}{2}
      \Big) x^{-1}  \\
  &&\textstyle  
    +\Big( \frac{9a-1}{4}\binom{3a}{3}
           - 3\left(  \frac{3a-1}{4}\binom{a}{3} + \binom{2a}{2}\binom{a}{2}
                      +\frac{6a-1}{4}\binom{2a}{3}
              \right) 
  + 6\left( \frac{3a-1}{4}\binom{a}{3}+\binom{a}{2}^2 \right) 
     \Big)x^{-2}
    + O(x^{-3})
   \Big) \\ 
 &=& x^{3a}\left(
      (6a^4 -6a^3+\tfrac{9}{4}a^2-\tfrac{1}{4}a ) x^{-2}+O(x^{-3}) 
     \right)     
\]
and hence 
\[
 \zeta_3(\widetilde{P}-\mathrm{N}) 
  &=& \frac{1}{6} \frac{(6a^4 -6a^3+\tfrac{9}{4}a^2-\tfrac{1}{4}a ) x^{3a-2} +O(x^{3a-2}) }
                        { \big(a^2x^{2a-1} +O(x^{2a-2})\big)^\frac{3}{2} } \, ,
\]
and hence the final term in \eqref{Eq:mu_3_Gamma_power_transformed_standardised}.
The proof of \eqref{Eq:mu_3_Gamma_standardised} is obvious.

For~\eqref{Eq:zeta_1/zeta_3->...<1_for_gamma_laws}, 
let us write here $\Gamma_\alpha\coloneqq \Gamma_{\alpha,1}=\Gamma_{\alpha,1,1}$.
Then \citet[Theorem 4.2, (4.24)]{Esseen1958} yields 
$\sqrt{n}\zeta_1(\widetilde{\Gamma_n}-\mathrm{N}) =
 \sqrt{n}\zeta_1(\widetilde{\Gamma_1^{\ast n}}-\mathrm{N})
 \rightarrow \frac{2}{3\sqrt{2\pi\mathrm{e}}}\mu_3(\widetilde{\Gamma_1}) 
 =  \frac{4}{3\sqrt{2\pi\mathrm{e}}}$
for $\N\ni n \rightarrow\infty$. 
Now the theory of Edgeworth expansions, as used by Esseen,
extends easily from sequences $(P^{\ast n} : n\in\N)$ of convolution powers
with $P\in\cP_3$ to more general ``one-parameter semigroups''
$(P_\alpha : \alpha \in\Alpha)$,  with $\Alpha$ a subsemigroup of $(\,]0,\infty[\,,+)$
with $1\in\Alpha$, and  $P_\alpha\in\cP_3$ and $P_{\alpha+\beta}=P_\alpha\ast P_\beta$
for $\alpha,\beta\in\Alpha$.
Hence we get here analogously 
$\sqrt{\alpha}\zeta_1(\widetilde{\Gamma_\alpha}-\mathrm{N})
\rightarrow  \frac{4}{3\sqrt{2\pi\mathrm{e}}}$ even for $\mathopen]0,\infty\mathclose[ \ni \alpha\rightarrow\infty$,
and combined with~\eqref{Eq:mu_3_Gamma_standardised}
then~\eqref{Eq:zeta_1/zeta_3->...<1_for_gamma_laws}.
%
\end{proof}

\begin{proof}[Proof of Example~\ref{Example:Subbotin_laws} from page \pageref{page:Example_Subbotin_laws}]  
                                                  \label{page:Proof_of_Example_Subbotin}
We may assume 
$\alpha=1$. For $\beta\in\mathopen]0,\infty\mathclose[\,$,
we 
get $\nu_r(P_{\beta,1})=\Gamma(\frac{r+1}{\beta})/\Gamma(\frac{1}{\beta})$
for $r\in[0,\infty[\;$,
hence 
$P_{\beta,1}\in\Prob(\R)$, and, by symmetry,
$\nu_3(\widetilde{P_{\beta,1}}) = h(\frac{1}{\beta})$ with
\[
 h(x) &\coloneqq& \frac{\Gamma(4x) \Gamma(x)^\frac{1}{2}}{\Gamma(3x)^\frac{3}{2}} 
 \quad\text{ for }x\in \mathopen]0,\infty\mathclose[\,. 
\]
We get here  $S^{-}(\widetilde{f}-\phi) \le 2\cdot 2=4$ easily by Laguerre,  
$\mu_3(\widetilde{P})=0$ by symmetry and
hence in case of $\beta\ne2$ not $S^{-}(\widetilde{f}-\phi)\le3$
by Lemma~\ref{Lem:zeta_3_computation}, 
and hence
Theorem~\ref{Thm:Cut_criteria}(c,d)
yields~\eqref{Eq:zeta_3_Subbotin_-_N} up to the third expression, for $\beta<\infty$.
The case of $\beta=\infty$ follows easily,
using $\Gamma(x)\sim\frac{1}{x}$ for $x\rightarrow 0$.

For the monotonicity claim, we recall the digamma function expansion 
$\psi(x) = \frac{\Gamma'}{\Gamma}(x)=-\gamma+\sum_{k\in\N_0} \big(\frac{1}{k+1}-\frac{1}{k+x}\big)$
and get by a simple computation, using in particular $4+\frac{1}{2}-\frac{9}{2}=0$ to simplify,
\[
 (\log h)'(x) &=& 4\psi(4x)+\tfrac{1}{2}\psi(x) -\tfrac{9}{2}\psi(3x)
  \,\ =\,\  \sum_{k\in\N_0} \frac{3kx}{(k+4x)(k+3x)(k+x)} \,\ >\,\ 0\,.
\]
The remaining claims follow easily, the one concerning $\beta\downarrow 0$
say by Stirling's formula.
\end{proof}

\begin{proof}[Proof of Example~\ref{Example:beta_llcurly_zeta_1}]\label{page:Proof_of_Example_beta_llcurly_zeta_1}
For any $p,s\in\mathopen]0,\infty[\,$, the measure $P$   
defined by~\eqref{Eq:Def_example_beta_llcurly_zeta_1} is positive and symmetric, 
and for $k\in\N_0$ we get, recalling Lemma~\ref{Lem:Standard_normal_identies_and_asymptotics}
in what follows,
\[
 \nu_k(P) &=& \nu_k(\mathrm{N}) - 2\int_t^\infty x^k\phi(x)\,\dd x 
               -2\,\phi(t)\frac{t^{k+1}}{k+1}  + p\,s^k \\
  &=&\left\{\begin{array}{ll}
    1-2\,\Phi(-t)-2t\,\phi(t)+p  &\text{ if }k=0\,,\\
    \nu_1(\mathrm{N}) - 2\,\phi(t)-t^2\phi(t)+p\,s  &\text{ if }k=1\,,\\
    1-2\,\big(t\,\phi(t) +\Phi(-t)\big) -\tfrac{2}{3}t^3\phi(t)+p\,s^2  &\text{ if }k=2\,,\\
    \nu_3(\mathrm{N}) - 2\,(t^2+2)\phi(t) 
        -\tfrac{1}{2}t^4\phi(t) +p\, s^3  &\text{ if }k=3\,,
   \end{array}  \right.             
\]
with the 
case of $k=3$ 
included to prove below
~\eqref{Eq:zeta_3_for_Example_beta_llcurly_zeta_1}
for a later use in the proof of
Lemma~\ref{Lem:zeta_1_bounded_by_beta_and_zeta_3}.

The conditions $\nu_0(P)=\nu_2(P)=1$ are fulfilled exactly for 
\[
 p &=&  2\,\big(t\,\phi(t) +\Phi(-t) \big) \,\ \sim\,\  2t\,\phi(t)   \,, \\
 s &=& \Big(\tfrac{2}{p}\big(\tfrac{1}{3}t^3\phi(t)+t\,\phi(t)
     +\Phi(-t)\big)\Big)^\frac{1}{2} \,\ \sim\,\ \frac{t}{\sqrt{3}}\,,
\]
and hence we get 
\[
 \nu_0(P-\mathrm{N}) &=& \phi(t)\mathrm{N}(\,]-t,t[\,) +\mathrm{N}(\R\setminus\,]-t,t[\,)+ p \\
    &=& \phi(t)\big( 1-2\,\Phi(-t)\big) +2\,\Phi(-t) +p \,\ \sim\,\ 2t\,\phi(t)\,, \\
 \nu_1(P)- \nu_1(\mathrm{N})
    &=& -2\phi(t)-t^2\phi(t) +p\,s \,\ \sim\,\ \big( \tfrac{2}{\sqrt{3}}-1\big)t^2\phi(t) \,.    
\]
Hence, using the general inequalities $\beta \le \nu_0$ on $\cM$
and $\zeta_1(M)\ge \int |x| \,\dd M(x)$ for $M\in\cM$, the claim follows.

We  further observe for the proof of Lemma~\ref{Lem:zeta_1_bounded_by_beta_and_zeta_3}
that we get
\la                                       \label{Eq:zeta_3_for_Example_beta_llcurly_zeta_1} 
 \zeta_3(P-\mathrm{N}) &=&  \tfrac{1}{6}\big(\nu_3(\mathrm{N})-\nu_3(P)\big) \\ \nonumber
  & = & \tfrac{1}{6} \big( 2\,(t^2+2)\phi(t)  +\tfrac{1}{2}t^4\phi(t) -p\, s^3\big)
  \,\ \sim\,\ \big(\tfrac{1}{12}- \tfrac{1}{9\sqrt{3}} \big)   t^4\phi(t)  \,
\al
by using in the first step Theorem~\ref{Thm:Cut_criteria}(c,d) 
applied to $M\coloneqq P-\mathrm{N}$ and $r\coloneqq3$,
namely $(C_0) \Rightarrow (C_3)$ in~\ref{Thm:Cut_criteria}\ref{part:zeta_r=mu_r/r!}
with $\mu\coloneqq \leb +\delta_{-s}+\delta_s$ and
$f(x)\coloneqq f_M(x)\coloneqq \frac{p}{2}\1_{\{-s,s\}}(x)
 -\phi(|x|\!\vee\!t)\1_{\R\setminus\{-s,s\}}(x)$ for $x\in\R$,
and indeed $(-1)^0F_0=F_0=-f$ initially positive with $S^{-}(F_0)=4$.
\end{proof}

\begin{proof}[Proof of Example~\ref{Example:Left-winsorised_normal_re_Nagaev}]
We continue    \label{page:Proof_of_Example_Left-winsorised_normal_re_Nagaev} 
to use the notation and facts established in the proof
of Example~\ref{Examples:zeta_3=mu_3}\ref{part:Left-winsorised_normal_laws}.
 Using the commutative ring identity~\eqref{Eq:ring_identity_P_N_n=2}
 and here $\mathrm{N}^{\ast 2}=\mathrm{N}_{\sqrt{2}}\,$ and 
 $P-\mathrm{N}=\Phi(-t)\delta_{-t}-\mathrm{N}(\,\cdot\setminus I)$,
 and assuming from now on $x\ge -2t $ and $t>0$, 
 we get
  $\big(P-\mathrm{N}\big)^{\ast 2}(\,]-\infty,x]) =\big(P-\mathrm{N}\big)^{\ast2}(\,]-\infty,-2t])=0$
 and hence
 \[
  F^{\star 2}(x) &=& 
     \Phi\big(\frac{x}{\sqrt{2}} \big)  
      +2\, \Big( \Phi(-t)\Phi(x+t) - \int_{-\infty}^{-t} \Phi(x-y)\phi(y)\,\dd y \Big) \\
   &=& \Phi\big(\frac{x}{\sqrt{2}} \big) 
       + 2\int_{-\infty}^{-t}\big(\Phi(x+t)-\Phi(x-y)\phi(y)\big)\,\dd y  \\
   &=&  \Phi\big(\frac{x}{\sqrt{2}} \big)  
     + 2\frac{\phi(t)}{t} \int_{-\infty}^{0}\big(\Phi(x+t)-\Phi(x+t-\tfrac{y}{t}) \big)
       \mathrm{e}_{}^{-\frac{y^2}{2t^2}} \mathrm{e}^y \,\dd y    
 \]
 by the change of variables $y\mapsto \frac{y}{t}-t$ in the last step.
 Choosing now $x\coloneqq -t$ as to roughly maximise
 the modulus of the last integral for $t$ large, 
 and using $\big(\Phi(0)-\Phi(z)\big)\mathrm{e^{-z^2/2}}=-\frac{z}{\sqrt{2\pi}}+O(z^2)$ for
 $z\in\R$, we obtain
 \la                   \label{Eq:F^ast2_symptotics_Left-winsorised_normal}
   F^{\star 2}(-t) &=& \Phi\big(\frac{-t}{\sqrt{2}} \big)   
       - \frac{2}{\sqrt{2\pi}}\,\frac{\phi(t)}{t^2} 
       +   O\big( \frac{\phi(t)}{t^3} \big) \quad\text{ for } t\in\mathopen]0,\infty\mathclose[\,.
 \al
 
 Assuming still $t>0$ and setting now $x_t\coloneqq\frac{-t-2\mu}{\sqrt{2}\sigma}$,
 we obtain
 \[
  \widetilde{F^{\ast2}}(x_t) &=& F^{\ast 2}(\sqrt{2}\sigma x_t+2\mu)
   \,\ =\,\ F^{\ast 2}(-t)
 \]
 and, using~\eqref{Eq:Left-winsorised_normal_mu>0_sigma_<1}, 
 also $x_t < \frac{-t}{\sqrt{2}}$ and hence, 
 by (\ref{Eq:mu_leftwinsorised_normal},\ref{Eq:sigma_leftwinsorised_normal})
 in the last step,
 \[
  0&<& \Phi\big(\frac{-t}{\sqrt{2}}\big) - \Phi(x_t)
   \,\ <\,\ \phi\big(\frac{-t}{\sqrt{2}}\big)\,\Big(\frac{-t}{\sqrt{2}}-x_t \Big)
   \,\ =\,\ O\big( \mathrm{e}^{-\frac{t^2}{4}}_{} \phi(t) \big)\,,
 \]
 and therefore, using~\eqref{Eq:F^ast2_symptotics_Left-winsorised_normal}, 
 \[
  \widetilde{F^{\ast2}}(x_t) - \Phi(x_t) 
   &=& - \frac{2}{\sqrt{2\pi}}\,\frac{\phi(t)}{t^2} 
       +   O\big( \frac{\phi(t)}{t^3} \big) \quad\text{ for } t\in\mathopen]0,\infty\mathclose[\,
 \]
and hence, using $x_t \sim -\frac{t}{\sqrt{2}}$,
\[
 \frac{\sup\limits_{x\in\R}\big(1+|x|^3\big)\left| F_{\widetilde{P^{\ast 2}}}(x)-\Phi(x)\right|}
        {\big(\zeta_1\!\vee\!\zeta_3\big)(\widetilde{P}-\mathrm{N})}
  & \succcurlyeq & \frac{t^3\frac{\phi(t)}{t^2}}{\phi(t)} \,\, \rightarrow\,\ \infty
  \quad\text{ for } t\rightarrow\infty\,. 
\]  

\vspace{-1.5\baselineskip}
\end{proof}

\section*{Acknowledgements}
An invitation to the conference \textit{Esseen 100 Years} at Uppsala University, 
for which I thank the organisers and in particular Silvelyn Zwanzig, provided the impetus 
for proving in Summer~2018 and then presenting
Theorems~\ref{Thm:F_star_G_vs_H_star_H}, then for $H_1=H_2$ only, 
and \ref{Thm:Berry-Esseen_for_Z-close_to_normal}, 
then with a worse constant and still relying on \citeposs{Zolotarev1997} slightly erroneous computation;
see \citet{Mattner2018}.

I further thank Jochen Wengenroth for showing me the kind of example used in 
Remark~\ref{Rem:phi_discontinuous_on_K},
Gerd Christoph and Ludwig Paditz for helpful remarks,
Lea Willems for help with TikZ,
Patrick van Nerven for proofreading in particular the numerics in Example~\ref{Example:Discretised_normal_laws},
Alexander Zvonkin for sharing with us in \citet[p.~92]{BravermanEtAl2021}
Shubin's ``Three Principles for writing a mathematical paper'',
and Bero Roos and two referees for suggesting several corrections
and improvements based
upon the version \url{https://arXiv.org/abs/2210.04060v3}
of this paper.

{\footnotesize
    \newcommand{\forceindent}{\leavevmode{\parindent=1.2em\indent}}
    \def\bibpreamble{\forceindent
      Some journal title abbreviations are here avoided,
      following \citet{BondGreen2014},
      in particular if a title occurs only once or if
      \url{https://mathscinet.ams.org/msnhtml/serials.pdf}
      does not apply.
      %
      Links to full texts, if provided here,
       are {\color{OliveGreen} OliveGreen} if they are {\color{OliveGreen}free},
       {\color{YellowOrange} YellowOrange}  if merely {\color{YellowOrange}read-only free},
       and {\color{BrickRed} BrickRed}  if  {\color{BrickRed} paywalled},
       as experienced by us at various times 
       starting in March 2022.
     The asterisk *  marks 
     references 
     I have taken from secondary sources, without looking at the original.
      Links in {\color{Blue} blue} are for navigating within this paper.\\ }

} 

\end{document}